\documentclass[a4paper,11pt]{amsart}
\usepackage[latin1]{inputenc}

\usepackage{amsmath,amsfonts,amssymb}
\usepackage{dsfont}
\usepackage{latexsym}
\usepackage{hyperref}

\usepackage{amsthm}

\usepackage{fullpage}
\usepackage[all]{xy}
\usepackage{mathrsfs}
\usepackage{amsbsy}
\usepackage{stmaryrd}

\usepackage{graphicx}
\usepackage{stackrel}
\usepackage[english]{babel} 
\usepackage{xcolor}
\usepackage{relsize}
\usepackage[nosort]{cite}

\numberwithin{equation}{section}

\newtheorem{thm}{Theorem}[section]

\newtheorem{prop}[thm]{Proposition}
\newtheorem{lem}[thm]{Lemma}
\newtheorem{definition}[thm]{Definition}
\newtheorem{Rem}[thm]{Remark}

\title{Semiclassical limit of  cubic  Nonlinear \\ Schr\"odinger equations for mixed states}
\author{Daniel Han-Kwan}
\address{CNRS, Laboratoire de Math\'ematiques Jean Leray (UMR 6629), Nantes Universit\'e, 44322 Nantes Cedex 03, France.} \email{daniel.han-kwan@univ-nantes.fr}
\author{Fr\'ed\'eric Rousset}
\address{Universit\'e Paris-Saclay, CNRS, Laboratoire de Math\'ematiques d'Orsay,  91405 Orsay Cedex, France.}\email{frederic.rousset@universite-paris-saclay.fr}
\date{\today}

\newcommand{\norm}[1]{{\left\lVert#1\right\rVert}}

\newcommand{\eps}{\varepsilon}
\newcommand{\intr}[1]{\int_{\bf R^{n}}}

\newcommand{\intrd}[1]{\int_{\bf R^{2n}}}

\newcommand{\intt}[1]{\int_{\bf T^{n}}}

\newcommand{\R}{\mathbb{R}}
\renewcommand{\H}{\mathcal{H}}

\newcommand{\pa}{\partial}

\begin{document}

\begin{abstract}
In this work, we study the semiclassical limit of cubic Nonlinear Schr\"odinger equations for mixed states. We justify the limit to a singular Vlasov equation  (in which the force field is proportional to the gradient of the density),
for  data with finite Sobolev regularity whose velocity profiles satisfy a  quantum Penrose stability condition.
This latter  condition is always satisfied for small data (with a smallness condition independent of the semiclassical parameter)
both in the focusing and  the defocusing case, and  for small perturbations of a large class of  physically relevant examples in the defocusing case, such as
local Maxwellian-like profiles.

\end{abstract}

\maketitle



\section{Introduction}

\subsection{The semiclassical nonlinear Schr\"odinger and Hartree equations} 

We are interested in the semiclassical limit of the cubic nonlinear Schr\"odinger equation (NLS) 
modeling  the mean-field dynamics of  quantum particles.  We shall use the description of the system based on
 the evolution of  a self-adjoint nonnegative trace class operator $\gamma(t) \in \mathscr{L}(L^2(\R^d; \mathbb{C}))$ 
 which solves the following form of the cubic nonlinear Schr\"odinger equation: 
 \begin{equation}
\label{eq:hartree-cubic}
 \left\{
\begin{aligned}
&i \eps \partial_t \gamma = \left[ -\frac{\eps^2}{2} \Delta \pm \rho_\gamma, \gamma\right], \\
&\gamma|_{t=0}= \gamma^0.
\end{aligned}
\right.
\end{equation}
both in the defocusing ($+$) and focusing ($-$) case.
Here, $[\cdot, \cdot]$ denotes the commutator between two operators. The density $\rho_{\gamma}(t,x)$ is
 defined as  $\rho_\gamma(t,x) ={\gamma}(t,x,x)$, where ${\gamma}(t,\cdot,\cdot)$ is the Schwartz kernel of  $\gamma(t)$.
 The parameter $\eps \in (0,1]$ stands for a scaled Planck constant and the semiclassical limit $\eps \to 0$ corresponds to the transition from quantum to classical dynamics.
 
In the  special case of pure states where $\gamma(t)$ is a rank-one operator, we have
$\gamma(t,x,y)= u(t,x) \overline{u(t,y)}$ and up to a time dependent phase, it is equivalent for
 $\gamma$  to solve \eqref{eq:hartree-cubic} and for the complex wave function $u(t,x)$ to solve the one-particle cubic NLS  equation 
 \begin{equation}
\label{eq:NLS}
\left\{
\begin{aligned}
&i \eps \partial_t u + \frac{\eps^2}{2} \Delta u = \pm |u|^2 u, \qquad  x \in \R^d,  \\
&u|_{t=0}=u_0.
\end{aligned}
\right.
\end{equation}
Here, we focus on the  general  case of \emph{mixed states} described by  \eqref{eq:hartree-cubic}.


Our techniques actually allow to  study a natural generalization of~\eqref{eq:hartree-cubic}  involving a short-range pair potential. Let $V \in \mathscr{S}'(\R^d)$  be a real and even  potential, with Fourier transform   $\widehat{V} \in \mathscr{C}^\infty_b(\R^d)$ (meaning that $\widehat{V}$ and all its derivatives are uniformly bounded) and $\langle V,1\rangle \neq 0$. 
We shall consider a  scaled interaction potential  $V_\eps$   defined by
$\widehat{V_{\eps}}(\xi)= \widehat{V}(\eps \xi)$.
 When $V$ is an $L^1_{\text{loc}}$ function this yields $V_{\eps}= {1 \over \eps^d} V(\cdot/\eps)$.
We can then consider the  nonlinear Hartree (or Von Neumann) equation with short-range potential:
 \begin{equation}
\label{eq:hartree-intro}
i \eps \partial_t \gamma = \left[ -\frac{\eps^2}{2} \Delta + V_\eps  \ast \rho_\gamma, \gamma\right].
\end{equation}
Note that the Dirac mass $V= \pm \delta_{0}$  is covered by the assumptions and that it is invariant by our scaling so that 
 $V_{\eps}= V$  and  \eqref{eq:hartree-intro} reduces  to~\eqref{eq:hartree-cubic}. 
  Another physically relevant potential that is admissible is the screened Coulomb potential, corresponding to $\widehat{V} (\xi) = \frac{1}{1+|\xi|^2}$. Note though that the unscreened Coulomb potential, corresponding to $\widehat{V} (\xi) = \frac{1}{|\xi|^2}$ is not covered.
 With general potentials $V$, it is natural to refer to the case $\widehat{V}> 0$ as the defocusing case, and $\widehat{V}< 0$ as the focusing case.

The scaling  for the pair potential $V_\eps$ is natural and physically relevant. Let us mention at least
two motivations for this scaling.
\begin{itemize}
\item   Consider the unscaled Hartree equation with the    pair potential $V$:
 \begin{equation}
\label{eq:hartree-intro-unscaled}
i  \partial_t \Gamma = \left[ -\frac{1}{2} \Delta+ V \ast \rho_\Gamma, \Gamma\right].
\end{equation}
We consider  for $\eps>0$ an hyperbolic scaling, meaning that we set  $\gamma^\eps := \lambda_{1/\eps} \Gamma$, where $\lambda_{1/\eps}$ stands for the dilation of ratio $1/\eps$ both in time and space. This reads at
the kernel level
$$  \Gamma (t,x, y) = \gamma^\eps (\eps t, \eps x, \eps y).$$
Then $\gamma^\eps$ precisely solves ~\eqref{eq:hartree-intro}.  Roughly speaking this scaling means that we are trying to  describe a large scale, long time regime for the Hartree equation~\eqref{eq:hartree-intro-unscaled}.
\item Another motivation can be related to the understanding of the  mean-field limit for fermions,
in a scaling  which is the natural  counterpart   
 to the one used for bosons in order to derive
 the NLS equation \eqref{eq:NLS} as a mean-field model, see \cite{erdos-schlein0, erdos-schlein, lieb-seiringer} for example.
Starting with  the Hamiltonian operator associated with the evolution of $N$ fermions, $N \gg 1$, which reads
$$
H_N= \sum_{j=1}^N - \frac{1}{2} \Delta_{x_j} + \lambda(N) \sum_{i<j}^N V\left({x_i-x_j \over L}\right),
$$
where the parameter $\lambda(N)$ accounts for the strength of the potential energy, 
 $L$ is the typical length scale of interaction;
  $H_{N}$ is acting on the space of $L^2(\R^{dN}; \mathbb{C})$ functions  with anti-permutation
 symmetry.
Because of the antisymmetry, as a consequence of the Lieb-Thirring inequality, the typical kinetic energy of $N$ fermions confined in a volume of order one is at least of order $N^{1+2/d}$. Therefore, for the potential energy to play a significant role in the dynamics, one has to choose $\lambda(N)$ at least of order $N^{-1+2/d}$ (which is significantly larger than the usual mean field scaling for bosons, that is  $N^{-1}$). The choice $\lambda(N)=  N^{-1+2/d}$ (in dimension $d=3$) was specifically made in \cite{NS,Spohn,BPSS} (for other scalings, see \cite{BGGM}, \cite{PP}).  Since the  typical velocity of a particle is of order $N^{1/d}$, it is  natural  to rescale time so that to focus on short times of order $\eps:=N^{-1/d}$. After multiplication by $\eps^2$, the associated many-body Schr\"odinger equations then writes
\begin{equation}
\label{eq:Nbody-fermions}
i\varepsilon \partial_t \psi_{N,t} = \left( \sum_{j=1}^N - \frac{1}{2} \eps^2 \Delta_{x_j} +   \lambda(N) \eps^2  \sum_{i<j}^N V\left({x_i-x_j \over L}\right) \right) \psi_{N,t},
\end{equation}
with $\varepsilon = N^{-1/d}.$
Here, we specifically make the choice of the supercritical scaling $\lambda(N)= N^{   {2 \over d} }$ so that
$\lambda(N) \eps^2 = { 1 \over N} \eps^{-d}$, 
and $L= \eps$. 
Taking formally the limit $N\longrightarrow +\infty$  while now fixing $\eps$ and neglecting the exchange term, we end up with the 
Hartree equation in the scaling~\eqref{eq:hartree-intro} as an intermediate model.


\end{itemize}

 As $V_\eps$ converges in the sense of distributions to a Dirac mass, the  semiclassical limit of \eqref{eq:hartree-cubic} and~\eqref{eq:hartree-intro} are  similar, namely  we shall obtain the \emph{singular} Vlasov equation
\begin{equation}
\label{eq:vlasov-benney-intro}
\partial_t f + v\cdot \nabla_x f  - c_V \nabla_x\rho_f \cdot \nabla_v f=0, \quad (x,v) \in   \R^d \times \R^d.
\end{equation}
where $c_V:=\langle V,1\rangle$. In the case $c_V>0$, this  equation is known as as the Vlasov(--Dirac)--Benney equation \cite{Zak,Bardos}.  
The aim of this work is to justify the derivation of ~\eqref{eq:vlasov-benney-intro} from \eqref{eq:hartree-intro}, for a class of initial data with finite regularity.

\subsection{The semiclassical Wigner equation}
The Wigner formalism is particularly useful to uncover the link between the Hartree and the Vlasov equations. The (semiclassical) Wigner transform of an operator $\gamma$  is defined as
$$
W_{\eps}[\gamma](x,v) = \frac{1}{(2\pi)^d}\int_{\mathbb{R}^d}e^{-{i}v\cdot y} {\gamma}\left(x+\frac{\eps y}{2},x-\frac{\eps y}{2}\right) dy,  \quad (x,v)\in \mathbb{R}^d \times \mathbb{R}^d,
$$
where $ {\gamma}(\cdot,\cdot)$ denotes the Schwartz kernel of $\gamma$. Recall that the Wigner transform can be understood as the formal dual (or inverse) of the Weyl quantization, defined for a  symbol $a$ as
$$
\mathrm{Op}^{W,\eps}_a \varphi =  \frac{1}{(2\pi )^d} \int_{\R^d} \int_{\R^d} e^{{i}(x-y)\cdot \xi} a\left( \frac{x+y}{2}, \eps \xi\right) \varphi (y) \, dy  d\xi, \quad \varphi \in \mathscr{S}(\R^d),
$$
in the sense that
$$
 \langle W_{\eps}[\gamma], a\rangle=  \operatorname{Tr}(
\gamma \,  \mathrm{Op}^{W,\eps}_a), \quad  a \in \mathscr{D}(\R^d\times\R^d).
$$
Note that when $\gamma$ is self-adjoint, its Wigner transform is a real function.

When $\gamma_{\eps}$ solves  the Hartree equation~\eqref{eq:hartree-intro}, the Wigner transform of $\gamma_{\eps}$, denoted by $f_\eps$ solves the Wigner equation
\begin{equation}
\label{eq:wigner-intro}
\left\{
\begin{aligned}
&\partial_t f_\eps+ v\cdot\nabla_x f_\eps + B_\eps [ \rho_{f_\eps},f_\eps]=0,\\
&f_\eps|_{t=0}= f^0_\eps \left(:= W_{\gamma^0_\eps}\right),
\end{aligned}
\right.
\end{equation}
where $\rho_{f_\eps} =  \int_{\R^d} f_\eps \, dv$ and
\begin{equation}
\label{defB-intro}
\begin{aligned}
B_\eps [ \rho_{f_\eps},f_\eps](t,x,v)&=  \frac{i}{\eps} \left(V_\eps \ast{\rho}_{f_\eps} \left(x-\frac{\eps}{2i} \nabla_v \right) -V_\eps \ast{\rho}_{f_\eps} \left (x+\frac{\eps}{2i} \nabla_v\right) \right)  f_\eps \\
 & \hspace{-1cm} \Big(=:   \frac{i}{(2\pi)^d} \int_{\mathbb{R}^d}e^{i v\cdot \xi_v} \frac{1}{\eps}\left( V_\eps \ast{\rho}_{f_\eps} \left(x-\frac{\varepsilon \xi_v}{2}\right) -V_\eps \ast{\rho}_{f_\eps} \left(x+\frac{\varepsilon \xi_v}{2}\right) \right) \mathcal{F}_{v}f_{\eps}(t,x, \xi_{v}) \, d\xi_v.\Big)
\end{aligned}
\end{equation}
Formally, if $f_{\eps}$ converges to some $f$ sufficiently strongly, then, by Taylor expansion,  
we expect the convergence
$$
  B_\eps [ \rho_{f_\eps},f_\eps]  \xrightarrow{\eps \to 0} - c_V \nabla_x \rho_f \cdot \nabla_v f,
$$
so  that the formal limit of the semiclassical Wigner equation is indeed  the Vlasov equation~\eqref{eq:vlasov-benney-intro}.

\subsection{Previous justifications of semiclassical limits of the  Hartree equations}
We shall now review the literature on the analysis of the semiclassical limit of  the Hartree equation.
There are many available works that we can  roughly classify into three types:
results for the Hartree equation with unscaled pair potentials, results in the case of  pure states
(where \eqref{eq:NLS} is studied directly) focusing on WKB initial data,  and results in dimension one for pure states 
which rely on the integrable structure of \eqref{eq:NLS}.

The semiclassical Hartree equation with  unscaled pair potential $w$ reads
 \begin{equation}
\label{eq:hartree-intro0}
i \eps \partial_t \gamma = \left[ -\frac{\eps^2}{2} \Delta  + w \ast \rho_\gamma, \gamma\right], 
\end{equation}
in this case, the  formal limit is the Vlasov equation 
\begin{equation}
\label{eq:vlasov-intro}
\partial_t f + v\cdot \nabla_x f - \nabla_x w \ast \rho_f \cdot \nabla_v f=0
\end{equation}

\medskip

\noindent {\bf Smooth pair potentials.} For smooth  pair potentials $w$,  the Vlasov equation has been   derived directly from the N-body dynamics for fermions, in the pioneering works \cite{NS,Spohn}. The derivation from Hartree to Vlasov, with quantitative estimates,  in strong topologies, 
was subsequently obtained in \cite{Pul,AP,APPP,AKN1,AKN2,BPSS,GP}. Non trace-class data were treated in \cite{LS20}.

\medskip

\noindent {\bf Coulomb potential.} A physically important interaction kernel is the Coulomb potential (namely $w = \frac{1}{4\pi} \frac{1}{|x|}$ in dimension $d=3$), in which case~\eqref{eq:vlasov-intro} is referred to as the Vlasov-Poisson equation. A justification of the semiclassical limit for mixed states towards Vlasov-Poisson was obtained by \cite{LP} and \cite{MM}, using the Wigner transform. 
Their methods are based on weak compactness techniques and  the use  of the conservation laws of the equations. A general 1D result allowing pure states was subsequently proved in \cite{ZZM}.

Recently, new approaches providing quantitative estimates, with convergence rates, in the case of the Coulomb potential and even  more singular potentials (but  not as singular as the Dirac measure that we allow  here) were developed in  \cite{Saf-CMP,Saf,LS,CLS}. 
In the latter, the general idea is to consider the Weyl quantization of the solution to the Vlasov equation in view of applying stability estimates at the level of the Hartree equation~\eqref{eq:hartree-intro0}. Some regularity is required for the initial data.

We can also mention the recent \cite{IL} which uses a combination of the quantum Monge-Kantorovich distance of \cite{GP} with the kinetic Wasserstein distance of \cite{Iacobelli} to obtain stability estimates for solutions having bounded density.

\medskip

\noindent {\bf NLS.}  For the cubic NLS   (or for other power nonlinearities), all results  on the semiclassical limit  we are aware of deal  with \emph{pure states}, that is to say with (variants of)  the NLS equation~\eqref{eq:NLS}, and are most often restricted to the defocusing case. The WKB approximation for one-phase
 initial data, that is to say for initial data under the form
$$
u_0(x)  =\sqrt{\rho_0(x)} \exp \left( i \frac{S_0(x)}{\eps} \right),
$$
was justified in \cite{Ger93,Gre98}. 
Namely, \cite{Ger93} proved the semiclassical limit in the analytic class (a focusing nonlinearity is then allowed), while in \cite{Gre98}, the case of data with finite Sobolev regularity  in the defocusing case was treated. The justification of the  WKB approximation in the defocusing case  consists in proving that the solution to~\eqref{eq:NLS} can be written as
$$
u(t,x) = a^\eps(t,x) \exp \left( i \frac{S^\eps(t,x)}{\eps} \right),
$$
on a small but  uniform interval $[0,T]$, with $(|a^\eps|^2, \nabla S^\eps)$ converging in a Sobolev norm to
  $(\rho,u)$,  a smooth  solution to the following isentropic Euler equation,
\begin{equation}
\label{eq:euler}
 \left\{
\begin{aligned}
&\partial_t \rho + \mathrm{div} (\rho u) =0, \\
&\partial_t u + u \cdot \nabla u + \nabla_x \rho = 0, \\
&\rho|_{t=0}= \rho, \quad u|_{t=0} = \nabla S_0,
\end{aligned}
\right.
\end{equation}
often referred to as the shallow water equation. We emphasize that~\eqref{eq:euler} can be seen as a special case of the Vlasov-Benney equation~\eqref{eq:vlasov-benney-intro}, namely for monokinetic data of the form $f(t,x,v) = \rho(t,x) \otimes \delta_{v=u(t,x)}$.  
We refer to \cite{Zhang,LinZhang,AC,CR} for extensions based on the modulated energy method (or variants)  and  to the monographs \cite{Carles-book} and \cite{Zhang2} for a broader overview. 

In dimension one, relying on the integrability of the cubic NLS equation \eqref{eq:NLS}
 and the inverse scattering method, more results are available, in particular the description of the solution  after singularity
 formation in the Euler equation, we refer to  \cite{JLM99} and to the review \cite{Miller-dispersif} for example.

For the case of many phases, that is when considering a multiphase WKB initial data
$$
u_0(x)  =\int_M \sqrt{\rho_0^\alpha (x)} \exp \left( i \frac{S_0^\alpha(x)}{\eps} \right)  \, d \mu(\alpha),
$$
where $(M, \mu)$ is a given probability space,  instabilities, even in the defocusing case,  are expected and  the literature is much more scarce.
In \cite{Car24}, the WKB analysis of~\cite{Gre98} is extended to the case of a finite number of phases, as long as they do not interact. The work \cite{BB16} justified the semiclassical limit to the Vlasov-Benney equation for multiphase WKB data with uniform \emph{analytic} regularity, thus extending \cite{Ger93}. Although not explicitly stated, the result of  \cite{BB16} extends as well to the focusing case.


\subsection{The Vlasov--Benney equation and the Penrose stability condition}
In order to justify the semiclassical limit to \eqref{eq:vlasov-benney-intro} in finite regularity, an important
issue is related to  the well-posedness  theory in finite regularity of this class of equations.
The equations~\eqref{eq:vlasov-benney-intro} belong to the family of \emph{singular Vlasov equations} \cite{HKR} which display a loss of derivative at the level of the force,
in sharp contrast with the Vlasov-Poisson equation in which the force field rather gains one derivative.
Above all, owing to Cauchy-Kowalevskaya type theorems, they are locally well-posed in analytic category, see in particular \cite{Gre95,JN, MV
}. However, 
 they are in general \emph{ill-posed} in Sobolev spaces \cite{BN}, even in arbitrarily small time and with an arbitrary finite loss of derivatives and weights \cite{HKN,Bar} (note that these results are stated for the Vlasov-Benney equation
 i.e when $c_{V}=1$ but can be readily extended to all the equations~\eqref{eq:vlasov-benney-intro}).
Broadly speaking, ill-posedness in Sobolev spaces is related to a possible loss of hyperbolicity (akin to \cite{Met}) for \eqref{eq:vlasov-benney-intro} and is, from the physical point of view, due to instabilities that  occur at the linear level for some particular initial conditions; typical examples are the so-called \emph{two-stream} instabilities which appear around functions whose profile in velocity  displays two large bumps (or more). 
However, when these instabilities do not develop, one may expect the equation to be well-posed in finite regularity. A first result in this direction  is \cite{BB13}, where it was proved that  in dimension $d=1$, the Vlasov-Benney equation is indeed locally well-posed for Sobolev initial data which, for all $x$, display a one bump velocity profile (which is indeed linearly stable). As a matter of fact, most studies of the  Vlasov-Benney equation were motivated by  its relation  to the Vlasov-Poisson equation in   the  \emph{quasineutral limit} of plamas.  
Namely,  it appears as the formal limit  for
\begin{equation}\label{eq:VP-quasineutre}
\left\{
      \begin{aligned}
&\partial_t f_\eps +v \cdot \nabla_x f_\eps   - \nabla_{x} U_{\eps} \cdot \nabla_v f_\eps=0, \qquad (x,v) \in \mathbb{R}^d \times \R^d, \\
&(\mathrm{I}-\eps^2 \Delta_x)U_\eps= \int_{\R^d} f_\eps(t,x,v) \, dv-1,
\end{aligned}
    \right.
\end{equation}
a system modeling the dynamics of ions in a plasma, in which the small parameter $\eps \to 0 $ stands for the scaled Debye length, which is the typical length scale of electrostatic interaction. 
This scaling can also be  interpreted as a  hyperbolic scaling for the Vlasov-Poisson system,
 so that large time instabilities of the unscaled Vlasov-Poisson system may show up in times $\mathcal{O}(\eps)$
  in the scaled system \eqref{eq:VP-quasineutre}  and prevent the formal limit to hold in general in finite regularity \cite{HKH}. 
Nevertheless, in \cite{HKR}, we have justified the quasineutral limit to Vlasov-Benney for data with finite regularity satisfying a certain stability condition, which precisely allows to avoid these  instabilities.
Note that the analysis of \cite{HKR} is performed on the periodic torus, that is for $x \in \mathbb{T}^d$; nevertheless, it  can
be easily adapted to the whole space case $x \in \mathbb{R}^d$.
For other types of results regarding the quasineutral limit of plasmas from  various forms of the Vlasov-Poisson
system, we refer for example  to \cite{Gre95,GPIproc,GPI,Bre00,HK11} which deal either with analytic regularity
 or with monokinetic data (which are the counterpart  of the one phase WKB approximation for NLS that was previously mentioned) in order to avoid instabilities.  

 Let us explain the result of \cite{HKR}. Given a profile in velocity $v\mapsto \mathbf{f}(v)$,  consider  what we shall call  generically a  Penrose function
\begin{equation*}
\mathcal{P}_{\mathrm{VP}}(\gamma,\tau,\eta,\mathbf{f})=-\frac{1}{1+|\eta|^2}\int_{0}^{+\infty} e^{-(\gamma+i\tau)s} s |\eta|^2  (\mathcal{F}_v\mathbf{f}) (s \eta) ds,\quad \gamma>0, \tau \in \mathbb{R}, \eta \in \mathbb{R}^d,
\end{equation*}
where the convention for the Fourier transform will be specified in~\eqref{eq:fourier}. 
Here the subscript $\mathrm{VP}$ means that it is associated with the Vlasov-Poisson system~\eqref{eq:VP-quasineutre}. Given a function $f(x,v)$ we say that the Penrose stability condition is satisfied if
\begin{equation}
\label{eq:quasi-penrose}
\inf_{x\in \mathbb{R}^d}\underset{(\gamma,\tau,\eta)\in(0,+\infty)\times \mathbb{R}\times \mathbb{R}^d}{\inf} \vert 1- \mathcal{P}_{\mathrm{VP}}(\gamma,\tau,\eta,f(x,\cdot))\vert \geq c_0,
\end{equation}
for some $c_{0}>0$.
The Penrose stability condition for homogeneous profiles $\mathbf{f}(v)$ appeared in \cite{Pen}. It notably played a key role in asymptotic stability results, referred to as Landau Damping, for the Vlasov-Poisson equation posed in $\mathbb{T}^d\times \R^d$, see \cite{MV}.
The main result of \cite{HKR} is that the limit from~\eqref{eq:VP-quasineutre} to Vlasov-Benney holds for a sequence of uniformly smooth (but of finite regularity) initial data, satisfying the Penrose stability condition~\eqref{eq:quasi-penrose}, also uniformly $\eps$.
As a corollary of the analysis of \cite{HKR}, the Vlasov-Benney system appears to be locally well-posed  in any dimension  for finite regularity initial conditions satisfying the  Penrose stability condition~\eqref{eq:quasi-penrose}.

The stability condition~\eqref{eq:quasi-penrose}, though necessary for the justification of the quasineutral limit (as the formal limit is wrong when it is violated \cite{HKH}), is however non-optimal for what concerns the well-posedness of the Vlasov-Benney equation. 
In the work \cite{CHKR}   in collaboration with K. Carrapatoso, we   prove that the Vlasov-Benney equation is indeed locally well-posed in finite regularity under the optimal  condition
\begin{equation}
\label{eq:penrose}
\inf_{x\in \mathbb{R}^d}\underset{(\gamma,\tau,\eta)\in(0,+\infty)\times \mathbb{R}\times \mathbb{R}^d}{\inf} \vert 1- \mathcal{P}_{\mathrm{VB}}(\gamma,\tau,\eta,f(x,\cdot))\vert \geq c_0
\end{equation}
for some $c_{0}>0$, where 
\begin{equation*}
\mathcal{P}_{\mathrm{VB}}(\gamma,\tau,\eta,\mathbf{f})=- \int_{0}^{+\infty} e^{-(\gamma+i\tau)s} s |\eta|^2  (\mathcal{F}_v \mathbf{f}) (s\eta) ds,\quad \gamma>0, \tau \in \mathbb{R}, \eta \in \mathbb{R}^d.
\end{equation*}
This condition is more natural since it can be derived through a direct stability analysis of the Vlasov-Benney equation,
whereas the condition \eqref{eq:quasi-penrose} is dependent of the approximation process used
 to construct the solution (namely the quasineutral limit process considered in \eqref{eq:VP-quasineutre}).
 By a continuity argument, the condition~\eqref{eq:quasi-penrose} implies~\eqref{eq:penrose} but  one can find examples with two bumps where $1- \mathcal{P}_{\mathrm{VP}}$ vanishes
 whereas ~\eqref{eq:penrose} holds.
 Note that we shall not use the existence results of   \cite{CHKR} or \cite{HKR}  here. As a byproduct of  our  main result, 
 we  obtain an existence result  for  \eqref{eq:vlasov-benney-intro} in finite regularity under another Penrose type
 stability condition adapted to the semiclassical  Wigner equation \eqref{eq:wigner-intro}.

\subsection{Main result}
We shall now present the main result of this paper. In order to state it, 
we need to introduce  appropriate functional spaces to measure regularity and localization (which are adapted to  the semiclassical Wigner equation
 \eqref{eq:wigner-intro}), and to introduce our stability condition.
 
Let us first  define  the   vector fields
\begin{equation}
\label{defvecVX}
V_{\pm}= \eps \nabla_{x} \pm 2  i v, \quad X_{\pm}= \eps \nabla_{v} \pm 2  i x.
\end{equation}
Note that they depend on $\eps$ but that we omit this dependence for notational convenience.
We shall use  that these vector fields have good commutation properties with  the linear part of the Wigner equation.
They correspond to natural differentiation and multiplication by weights at the level of operators,  that is to say when acting on 
$\gamma = \mathrm{Op}_{f}^{W, \eps} $. Indeed, by  definition of the Wigner transform, we  observe that 
 \begin{equation}
  \label{linkgamma}   V_{+} f= W^\eps [2 \eps \nabla \gamma], \quad V_{-} f= W^\eps[ 2 \eps \gamma \nabla  ], \quad X_{+} f= W^\eps [2x \gamma], \quad X_{-} f= W^\eps[  2\gamma x ].
 \end{equation}
We shall work with the following weighted Sobolev spaces based on $V_\pm, X_\pm$. 
\begin{definition}
Let $m, r \in \mathbb{N}$. 
\begin{itemize}
\item For a  function $f(x,v)$ on $\mathbb{R}^{2d}$, we define the $\H^0_r$ norm as
\begin{equation}
\label{defH0r} \|f \|_{\H^0_{r}} = \sum_{\substack{ |\beta | + | \beta' |  \leq r \\ |\gamma| + | \gamma'| \leq r} }\|V_{+}^{\beta}X_{-}^{\beta'}
 V_{-}^{\gamma} X_{+}^{\gamma'} f \|_{L^2(\R^{2d})},
 \end{equation}
 where  $\beta,\beta', \gamma,\gamma'  \in \mathbb{N}^d$, and the $\H^m_r$ norm as
\begin{equation}
\label{defHmr} \|f \|_{\H^m_{r}}= \sum_{| \alpha | \leq m} \| \partial^\alpha_{x,v} f \|_{\H^0_{r}}, 
\end{equation}
 where  $\alpha \in \mathbb{N}^{2d}$.
 \item For a  function $\rho(x)$ on $\mathbb{R}^d$, we define the $H^0_r$ norm as 
  $$\| \rho \|_{H^0_{r} }= \sum_{|\beta | \leq r } \| (\eps \partial_x)^\beta  \rho \|_{L^2(\mathbb{R}^d)},$$
  where  $\beta\in \mathbb{N}^d$,  and the $H^m_r$ norm as
 $$ \|\rho\|_{H^m_{r}} = \sum_{ | \alpha | \leq m} \| \partial^\alpha_x \rho \|_{H^0_{r}},$$
 where $\alpha \in \mathbb{N}^d$.
  \end{itemize}
     Note  that all these norms depend on $\eps$, but this dependence is never specified.
  \end{definition}

 For  the convergence result, it  will be convenient to rely on   standard weighted Sobolev spaces
  which do not depend on $\eps$.
\begin{definition}
\label{def:Hmr-droit}
The weighted Sobolev space $\mathrm{H}^m_{r}$ of functions $f(x,v)$ on $\R^{2d}$, for $m,r\in \R$ is associated with the norm
$$
\|f \|_{\mathrm{H}^m_{r}}=\| \langle v \rangle^r (I-\Delta_{x,v})^{m/2} f \|_{L^2(\R^{2d})},
$$
where $\langle \cdot \rangle = \sqrt{1+|\cdot|^2}$. We shall also denote $\mathrm{H}^m_{r,v}$ for the analogous space pertaining to functions of $v$ only.
\end{definition}

Note that the spaces $\mathrm{H}^m_{r}$ will be used only in Section \ref{sec:end-conv}. Let us observe
(see Lemma \ref{LemHmrHmr} for details) that we have  for some $C>0$ independent of $\eps\in (0, 1]$ the relation
$$ \| \cdot \|_{\mathrm{H}^m_{r}} \leq C \|\cdot \|_{\mathcal{H}^m_{r}}$$
for $m$, $r$ nonnegative integers.

\bigskip

We finally introduce the relevant stability condition for the semiclassical limit.
Throughout this paper, the Fourier transform on $\mathbb{R}^n$ for all $n \in \mathbb{N}\setminus\{0\}$, that will be denoted indifferently by  $\mathcal{F}(u)$ or $\widehat{u}$, will be normalized as
\begin{equation}
\label{eq:fourier}
\mathcal{F}(u)(\xi)=\widehat{u}(\xi)= \int_{\mathbb{R}^n} u(y) e^{-i \xi\cdot y} dy.
\end{equation}

\begin{definition} Given a profile
$\mathbf{f}(v)$, we define its \emph{quantum Penrose} function by
\begin{equation}\label{Penrose-quantique}
\mathcal{P}_{\mathrm{quant}}(\gamma,\tau,\eta,\mathbf{f})=- 2 \widehat{V}( \eta)  \int_{0}^{+\infty}  e^{-(\gamma+i\tau)s} {\sin\left(\frac{s |\eta|^2}{2}\right)} (\mathcal{F}_v  \mathbf{f}) (s \eta ) ds,\quad \gamma>0, \tau \in \mathbb{R}, \eta \in \mathbb{R}^d.
\end{equation}
We say that a  function $f(x,v)$ satisfies  for a given  $c_{0}>0$ the  $c_0$ \emph{quantum Penrose stability}   if  the following inequality holds 
\begin{equation}\label{QCrit}
\underset{x \in \mathbb{R}^d}{\inf}\underset{(\gamma,\tau,\eta)\in(0,+\infty)\times \mathbb{R}\times \mathbb{R}^d}{\inf} \vert1- \mathcal{P}_{\mathrm{quant}}(\gamma,\tau,\eta,{f}(x,\cdot))\vert \geq c_0
\end{equation}
and  that  $f$ satisfies the \emph{quantum Penrose stability} condition if it satisfies the $c_0$ quantum Penrose stability condition for some $c_0>0$.
\end{definition}

The main result of this work is a derivation of the singular Vlasov equation~\eqref{eq:vlasov-benney-intro} from the Wigner  equation~\eqref{eq:wigner-intro} in the semiclassical limit $\eps\to0$. 
 The result is achieved for a family of initial data with uniform bound in the weighted Sobolev space $\H^m_r$  (with $m,r$ large enough) and that satisfy a uniform quantum Penrose stability condition. 

\begin{thm}
\label{thm:main}
 Let  $r\geq  2d +  2 \lfloor d/2  \rfloor  + 8$ and $m \geq \min \left(10d +  d/2+ 14 + r , 3d + 6 + 2r\right)$. 
Let 
 $(f^0_{\eps})_{\eps \in (0,1]} $ a real-valued family of initial data for \eqref{eq:wigner-intro} that satisfies
 the following assumptions.


\medskip

\noindent {\bf A1. Uniform weighted Sobolev regularity.} There is $M_0>0$ such that
\begin{equation}
\label{hypintrouniforme}
\sup_{\eps \in (0,1]} \| f^0_\eps \|_{\H^{m}_{r}}  \leq M_0.
\end{equation}

\noindent  {\bf A2. Uniform quantum Penrose stability.}  The family  $(f^0_\eps)_{\eps \in (0,1]}$ satisfies the $c_0$ quantum Penrose stability condition \eqref{QCrit} for some $c_0>0$ independent of $\eps$.

\medskip

Then there exist $T>0$ and $\eps_0>0$ such that, for all  $\eps \in (0,\eps_0)$, there is  a unique solution $f_\eps \in C([0,T]; \H^{m}_{r})$ to~\eqref{eq:wigner-intro} such that the following properties hold.

\medskip

\noindent $\bullet$ {\bf Uniform bounds.}  There exists $M>0$  such that
for all $\eps \in (0,\eps_0)$,
\begin{equation}
\label{estuniformeintro}
\norm{f_\eps}_{L^{\infty}\left(0,T;\mathcal{H}^{m-1}_{r}\right)}+\norm{\rho_{f_{\eps}}}_{L^{2}\left(0,T;H^{m}_{r}\right)} \leq M.
\end{equation}

\noindent $\bullet$ {\bf Convergence to singular Vlasov.} Assume in addition  that $f^0_\eps \rightarrow f^0 $ in $
L^2(\mathbb{R}^{2d)}$.
Then, 
there exists  $f \in C([0,T]; \mathrm{H}^{m-1}_{r})$ with  $\rho_{f} \in L^2(0, T; H^m)$,
 solution   to~\eqref{eq:vlasov-benney-intro} with initial datum  $f^0$ such that  the following convergences hold:
\begin{equation}
\label{convergenceintro}
\lim_{\eps \to 0} \left( \sup_{[0,T]} \| f_{\eps} - f\|_{\mathrm{H}^{m-1-\delta}_{r-\delta}} + \| \rho_{f_\eps} -\rho_f\|_{L^2(0,T;H^{m-\delta})}\right)=0,
\end{equation}
for any $\delta>0$.
\end{thm}

We shall explain the general strategy for the proof of Theorem \ref{thm:main} in section \ref{sec:strategy}. Let us first provide a few comments.

\begin{itemize}

\item Note that $f_{\eps}$ stays  in $\H^m_{r}$ for all $t \in [0, T]$. Nevertheless, we are only 
able to obtain uniform in $\eps$ estimates for the $\H^{m-1}_{r}$ norm of $f_{\eps}$, and only $\rho_{f_{\eps}}$
 can be controlled with the maximal regularity when we measure it in the $L^2$ norm in time.
 Note that we  have kept track of regularity  and localization in the above result but did not try to optimize it.
 We obtain a strong convergence result.  By Sobolev embedding, \eqref{convergenceintro}
  implies in particular  convergence in $L^\infty_{x,v}$.  
We could  also obtain weights in $x$ in the convergence result (thanks to the vector fields $X_\pm$) but we have chosen not to dwell on this aspect.

\item 
We have chosen to state everything   in terms of  $f_{\eps}$, since we shall perform the proof at the level of
the Wigner equation \eqref{eq:wigner-intro},  nevertheless, by using the properties
of the Wigner transform and  \eqref{linkgamma}, the assumptions and results can be translated
at the operator level.  Note that we assumed that $f_{0}^\eps$ is real but that we did not assume
that $f_{0}^\eps$ is non-negative to take into account a well-known flaw of the Wigner transform: 
a non-negative self-adjoint operator yields a real Wigner transform but not necessarily a non-negative one.

If one starts from a given family of self-adjoint non-negative  trace class operators $\gamma^{0}_{\eps}$
as initial conditions for \eqref{eq:hartree-intro}, by setting  $ f_{0}^\eps:=  W_\eps[\gamma^{0}_{\eps}]$
  the uniform regularity assumptions \eqref{hypintrouniforme} bearing on $f^0_\eps$  follows from a uniform control in weighted Hilbert-Schmidt norm of  commutators of $\gamma^0_{\eps}$ with $\frac{x}{\eps}$ and $\frac{\nabla}{\eps}$, namely, for some $C>0$ 
\begin{equation*}
\label{eq:hyp-operator}
\sup_{\eps \in (0,1]}\eps^{-d} \| \langle x \rangle^r \langle\nabla \rangle^r [a_1, [\ldots,[a_\ell, \gamma^0_{\eps}]\cdots]]\langle x \rangle^r \langle\nabla \rangle^r \|_{\mathrm{HS}}^2\leq C,
\end{equation*}
for all $\ell = 0, \cdots, m$, for all choices of $a_i = \frac{x}{\eps}$ or $\nabla$. Here,  $\| \cdot \|_{\mathrm{HS}}$ stands for the Hilbert-Schmidt norm.
According to  \cite[Remark 3) after Theorem 2.1]{BPSS}, smooth superpositions of fermionic coherent states naturally satisfy this assumption. 
Note that pure states do not satisfy it.

Asking for analytic regularity would mean to ensure that
$$
\sup_{\eps \in (0,1]}\eps^{-d} \Big\| \langle x \rangle^r \langle\nabla \rangle^r \underbrace{\left[\frac{\nabla}{\eps}, \left[\ldots,\left[\frac{\nabla}{\eps}, \gamma^0_{\eps}\right]\cdots\right]\right]}_{\ell \text{ times }}\langle x \rangle^r \langle\nabla \rangle^r \Big\|_{\mathrm{HS}}^2\leq C,
$$
 holds for all $\ell \in \mathbb{N}$.

The convergence result  \eqref{convergenceintro} also implies the following. Denoting by $\gamma_{f_\eps}$ (resp. $\gamma_f$) the Weyl quantization of $f_\eps$ for all $\eps \in (0,1]$ (resp. of $f$), $\gamma_{f_\eps}$ satisfies the Hartree equation~\eqref{eq:hartree-intro} associated with the initial condition $\gamma^0_{\eps}$ on $[0,T]$ and 
$$
\lim_{\eps \to 0} \sup_{[0,T]} \eps^{-d} \|   [a_1, [\ldots,[a_{\ell},\gamma_{f_\eps}- \gamma_{f}]\cdots]]  \|^2_{\mathrm{HS}} =0,
$$
for all $\ell =0, \cdots, m-2$ and all choices of $a_i = \frac{x}{\eps}$ or $\frac{\nabla}{\eps}$.


\item{\bf Other nonlinearities.}
We can handle other smooth nonlinearities for NLS, with essentially the same analysis. We only need to introduce the appropriate quantum Penrose stability condition. Namely, consider the nonlinear Hartree equation
$$
i \eps \partial_t \gamma = \left[ -\frac{\eps^2}{2} \Delta + \Psi(\rho_\gamma), \gamma\right], 
$$
where $\Psi \in \mathscr{C}^\infty(\R)$, which corresponds to the mixed state version of the NLS equation 
$$i \eps \partial_t u + \frac{\eps^2}{2} \Delta u = \Psi( |u|^2) u.$$ The quintic case for instance corresponds to $\Psi(x)= \pm x^2$. (Note that a convolution with short-range pair potential as done in the cubic case may also be considered.)
The Penrose function for a function $\mathbf{f}(v)$ (with density $\rho_{\mathbf{f}}= \int_{\R^d} \mathbf{f}  \,dv$) reads in this case
\begin{equation}\label{Penrose-quantique-other}
\mathcal{P}_{\mathrm{quant}}(\gamma,\tau,\eta,\mathbf{f})=- 2 \Psi' (\rho_{\mathbf{f}})\int_{0}^{+\infty}  e^{-(\gamma+i\tau)s} {\sin\left(\frac{s |\eta|^2}{2}\right)} (\mathcal{F}_v  \mathbf{f}) (\eta s) ds,\quad \gamma>0, \tau \in \mathbb{R}, \eta \in \mathbb{R}^d,
\end{equation}
while the limit is the singular Vlasov equation
$$
\partial_t f + v\cdot \nabla_x f  -  \nabla_x \Psi(\rho_f) \cdot \nabla_v f=0. 
$$

\item 
An analogue of Theorem~\ref{thm:main}  restricted to the case of smooth and fastly decaying  pair potentials,  namely  $V$
 in the Schwarz class  $\mathscr{S}(\mathbb{R}^d)$  has previously been obtained in collaboration with T. Chaub in \cite{Chaub}.  The  cubic NLS  case \eqref{eq:hartree-cubic} is therefore not covered by this previous result.  The fact that high  frequencies $\eps | \xi_{x}|\gtrsim
 1$ can  be  controlled
  by the fast decay of  $\widehat V$ and  of its derivatives is crucially used in many steps of the proof in  \cite{Chaub}.
  We   follow the same  general strategy (itself inspired by \cite{HKR}), but
  in order to handle general pair potentials, the method has to be significantly improved and sharpened.
   In particular, we  perform  the analysis  in the weighted spaces $\H^m_r$ defined above  which are really tailored for the cubic Wigner equation~\eqref{eq:wigner-intro}, instead of the standard weighted Sobolev spaces $\mathrm{H}^m_{r}$. 
\end{itemize}

\subsection{The quantum Penrose stability condition}

To the best of our knowledge, the (homogeneous) quantum Penrose  stability condition was first introduced in the mathematical literature by \cite{LS14} in the context of asymptotic stability of space invariant equilibria of the Hartree equation.
 We refer to \cite{CHP,CdS,Mal,Had,You,BHS} for recent developments on this topic. Such results can be understood as the quantum analogue of asymptotic stability results (usually referred to as Landau damping) for nonlinear Vlasov equations in the whole space, see for example \cite{BMM, HKNR, HNX, BMMlin,HKNRlin, IPWW}. 
The recent work \cite{Smith} established a connection between these  quantum and classical results (see also \cite{HadamaHong}).
In the physical literature, the quantum Penrose function is often referred to as the Lindhard function \cite{Lin} and was already identified to play a key role in the stability of space invariant quantum gases, see e.g. \cite[Chapter 4]{GV}.
For the same reason as  in the study of  the quasineutral limit, the fact that  the quantum Penrose stability condition plays a prominent role in the study of~\eqref{eq:hartree-intro} is natural.
Indeed, because of  the hyperbolic scaling related to the semiclassical  Hartree equation \eqref{eq:hartree-intro}, the linear instabilities
which may show up in the large time behavior of the unscaled Hartree equation are now expected to 
occur in times $\mathcal{O}(\eps)$. An adaptation of the analysis of  \cite{HKH}  yields that the Penrose condition is necessary  to justify the semiclassical limit on  times $\mathcal{O}(1)$ in finite regularity.

It is important to note that the quantum Penrose condition~\eqref{QCrit} implies the Penrose condition~\eqref{eq:penrose} (see the upcoming Lemma~\ref{lem:Penrose}), so that the last part of Theorem~\ref{thm:main} is in  agreement with the sharp  local well-posedness result of \cite{CHKR}.

The quantum Penrose stability condition is open with respect to strong enough topologies, in the sense that if it is satisfied for a function $f$, it is also satisfied for all functions in the vicinity (in a strong enough topology) of $f$. This comes from a stability inequality:  for  $m>3$ and $r>d/2$ and for any two profiles $f(v),g(v)$,  we have
\begin{equation}
\label{perturbpenrose}
\sup_{(\gamma,\tau, \eta) \in (0,+\infty)\times \R \times \R^d} |\mathcal{P}_{\mathrm{quant}}(\gamma, \tau, \eta, f) - \mathcal{P}_{\mathrm{quant}}(\gamma, \tau, \eta, g)| \lesssim \|\widehat{V}\|_\infty \| f-g \|_{\mathrm{H}^m_{r,v}} .
\end{equation}
where for functions $h(v)$ the  $H^m_{r, v}$ norm is defined as 
$$
\|h\|_{\mathrm{H}^m_{r, v}}=\| \langle v \rangle^r (I-\Delta_{v})^{m/2} h \|_{L^2(\R^d)}.
$$

In particular,  the quantum Penrose condition always  holds for  data in $L^\infty_x \mathrm{H}^m_{r,v}$ ($m>3, r>d/2$), satisfying a smallness condition involving  the pair potential $V$. Namely, there exists a constant $c_d>0$ such that, if
$$
c_d \|\widehat{V}\|_\infty  \| f(x,v) \|_{L^\infty_x \mathrm{H}^m_{r,v}}< 1, \quad m>3, r>d/2,
$$
then $f(x,v)$ satisfies the quantum Penrose stability condition.
Therefore in the focusing case, Theorem~\ref{thm:main} consequently holds for initial data of this kind. This is, as far as we know, the first  class of examples for which   the semiclassical limit  for NLS  in finite regularity
 can be justified.

Note that another direct consequence of \eqref{perturbpenrose}  is that for $f \in \mathscr{C}^{0}_0(\mathbb{R}^d; \mathrm{H}^m_{r,v})$ (the space of continuous functions,  converging to zero at infinity
 in $x$, with values in $\mathrm{H}^m_{r,v})$,   $m>3, r>d/2$,  we have by a finite covering argument  that $f$ satisfies  the quantum Penrose condition
  if and only if for every $x \in \mathbb{R}^d$, the profile $f(x,\cdot)$ satisfies the Penrose condition  
$$ \inf_{(\gamma,\tau,\eta)\in(0,+\infty)\times \mathbb{R}\times \mathbb{R}^d}  \vert1- \mathcal{P}_{\mathrm{quant}}(\gamma,\tau,\eta,{f}(x,\cdot))\vert  >0.$$
There are thus  interesting cases where   large data are allowed, at least in the {\it defocusing} case, that is when assuming that $\widehat{V} \geq 0$: indeed, in this case, the quantum Penrose condition is satisfied  for non-negative initial data that are radial decreasing in $v$ in dimension $d= 1, \, 2$, and only radial in $v$ in dimensions $d\geq 3$, see \cite{MV,BMM,LS14,Mal}.
For instance, in the defocusing case, the quantum Penrose stability condition together with \eqref{hypintrouniforme} holds  
for the following inhomogeneous distribution of :
\begin{itemize}
\item  Boltzmann gases 
$$
\varphi (x,v)=  \rho(x) e^{\frac{-|v-u(x)|^2 - \mu(x)}{T(x)}},
$$
\item Fermi gases
$$ \varphi (x,v)=  {\rho(x) \over  e^{\frac{|v-u(x)|^2 - \mu(x)}{T(x)}} + 1},$$
\item Bose gases 
$$  {\rho(x) \over  e^{\frac{|v-u(x)|^2 - \mu(x)}{T(x)}} - 1},$$
\end{itemize}
 where $\rho$,  $u$, $\mu$ and $T $ are  bounded,  smooth enough,   $\rho$   positive,    decaying to zero at infinity quickly enough, 
  $\inf_{\R^d} T>0$ and $\mu$ is such that $\sup_{\mathbb{R}^d} \mu <0$ in the third case.
  More generally for a given function $F: \mathbb{R}^d \times [0, + \infty)\rightarrow \mathbb{R}_+$
   smooth, sufficiently decaying at infinity and such that $F(x, \cdot)$ is decreasing in dimension $d=1, \, 2$
   for every $x$, the distribution 
   $$ \varphi(x,v)= F\left(x, {|v-u(x)|^2 \over T(x) } \right)$$
   matches the regularity assumption \eqref{hypintrouniforme} and the quantum Penrose stability condition.
   
   Owing to \eqref{perturbpenrose}, one can also add to  these examples an arbitrary  small enough perturbation.

%
%

\subsection{Notations}

We first provide a convenient notation that will be systematically used in the paper. 
We will often write the variable $y$ or $z=(x,v) \in \mathbb{R}^{d} \times \mathbb{R}^{d}$ to handle both variables $x$ and $v$ at the same time; in some specific cases, we use $x$ and $v$ to highlight their specific role.
Likewise, we denote the dual variable $\xi = (\xi_x,\xi_v) \in \mathbb{R}^{d} \times \mathbb{R}^{d}$, writing $\xi_x$ or $\xi_v$ only when required.

Given a function $u_\eps$, the subscript $\eps$ refers to a dependence with respect to $\eps$ of the function $u_\eps$. Most of the time, to simplify the expressions, when this dependence is not singular, it will be dismissed, while keeping in mind that the main focus will be to obtain estimates which are uniform with respect to $\eps$. 

Given a function $u(t,z,\xi)$, with $z \in \R^n$ to be seen as the physical variable (in practice, $z= x,v$ or $(x,v)$) and $\xi \in \R^n$ its dual Fourier variable, the notation $u^\eps$ means that we evaluate $u$ at the point $(t,z,\eps\xi)$:
\begin{equation}
\label{def:epsenhaut}
u^\eps(t,z,\xi)=u(t,z,\eps \xi).
\end{equation} 
In the case of multiple variables, for example for a function $u(t,z,y,\xi,\eta)$, all dual variables are rescaled, meaning that $u^\eps(t,z,y,\xi,\eta)=u (t,z,y,\eps\xi,\eps\eta)$.

We use in this work different types of pseudodifferential calculus. 
\begin{itemize}
\item We consider standard pseudodifferential operators with the following notation. Let $y= x, v$ or $(x,v) \in \R^n$, and $\xi_y \in \R^n$ be its dual Fourier variable. Given $a(y, \xi_y)$ a scalar or vectorial symbol, we denote by $a(y, D_y)$ the associated pseudodifferential operator (where $D_y$ can be understood as $\frac{1}{i} \nabla_y$), defined by the formula
\begin{equation}\label{PseudoDef-nonsemi}
a(y,D_y)u := \frac{1}{(2\pi)^n}\int_{\xi_y} e^{i y\cdot \xi_y} a(y,\xi_y)  \widehat{u}(\xi_y) d\xi_y, \quad u \in \mathscr{S}(\R^n).
\end{equation}
This notation allows to explicitly indicate the variables with respect to which the pseudodifferential calculus is performed.
With the notation~\eqref{def:epsenhaut}, the operator $a^\eps(y,D_y)$ denotes the associated semiclassical pseudodifferential operator. In particular, observe that the operator $B_\eps$ appearing in the Wigner equation~\eqref{eq:wigner-intro} can be recast as
\begin{equation}
\label{defB-re}
B_\eps [ \rho_{f_\eps},f_\eps]=  \frac{i}{\eps} \left( [V_\eps \ast{\rho}_{f_\eps}]^\eps \left(x-\frac{1}{2} D_v \right) -[V_\eps \ast{\rho}_{f_\eps}]^\eps \left (x+\frac{1}{2}D_v\right) \right)  f_\eps.
\end{equation}

\item We will furthermore use a pseudodifferential calculus for operator-valued symbols, meaning that given a separable Hilbert space $\mathbf{H}$ and a symbol $\mathrm{L}(y,\xi_y) \in \mathscr{L}(\mathbf{H})$, we consider the pseudodifferential operator $\operatorname{Op}_{\mathrm{L}}$ defined by the formula
\begin{equation}\label{PseudoDef-operator}
\operatorname{Op}_{\mathrm{L}} u := \frac{1}{(2\pi)^n}\int_{\xi_y} e^{i y\cdot \xi_y} \mathrm{L}(y,\xi_y)  \widehat{u}(\xi_y) d\xi_y, \quad u \in \mathscr{S}(\R^n; \mathbf{H}).
\end{equation}
This will be used in the case   $\mathbf{H}= L^2(0, T)$.
\item We will finally use a pseudodifferential calculus with parameter $\gamma>0$ (for functions of time and space). To avoid any confusion, the associated pseudodifferential operators will be referred to with bold letters, with the symbol as subscript. For a symbol $a(x,\gamma,\tau,\eta)$ on $\mathbb{R}^d\times (0,+\infty)\times \mathbb{R}\times\mathbb{R}^d $ and $u$, we denote by $\mathbf{Op}^{\gamma}_a$ (resp. $\mathbf{Op}^{\eps,\gamma}_a$) the operator
\begin{equation}\label{PseudoParamDef}
\begin{aligned}
\mathbf{Op}^{\gamma}_a u &:=    \frac{1}{(2\pi)^{d+1}} \int_{\tau} \int_{\xi}  e^{{i}\left( x\cdot \xi +\tau t \right)}  a(x, \gamma, \tau, \xi)  \mathcal{F}_{t,x} u(\tau,\xi) d\xi d\tau,  \quad u \in \mathscr{S}(\R \times \R^d)\\
\mathbf{Op}^{\eps,\gamma}_a u &:=   \frac{1}{(2\pi)^{d+1}}  \int_{\tau} \int_{\xi}  e^{{i}\left( x\cdot \xi +\tau t \right)}  a(x, \eps\gamma, \eps\tau, \eps\xi)  \mathcal{F}_{t,x} u(\tau,\xi) d\xi d\tau.
\end{aligned}
\end{equation}
\end{itemize}

The integer
$$
k_d :=  \lfloor d/2 \rfloor +2
$$
will appear many times in the analysis, in particular in the pseudodifferential estimates. We will use very often this notation, sometimes recalling its definition to ease readability.

Finally, throughout this work, $\Lambda$ will stand for a generic continuous nondecreasing function
with respect to all its arguments, that may change from line to line but that stays independent of $\varepsilon$. 
We also use the notation $\cdot \lesssim \cdot$ for $\cdot \leq C \cdot$ where $C$ is an harmless number which does not
depend on $\eps \in (0, 1]$.

\section{Strategy and organization of the paper}
\label{sec:strategy}

As already mentioned, the well-posedness of singular Vlasov equations such as \eqref{eq:vlasov-benney-intro}
is a subtle question and this class of equations is  not known to admit a weak-strong stability principle. Consequently, the strategy of the recent works \cite{LS} or \cite{CLS} in the Coulomb case, which consists in lifting  a weak-strong stability estimate for Vlasov to the level of Hartree (or Wigner)
 does not seem to be possible. 
  However,  as seen from \cite{HKR}, a stability estimate turns out to hold for smooth enough solutions to \eqref{eq:vlasov-benney-intro}, as long as one of the two satisfies the Penrose stability condition~\eqref{eq:penrose}. 

\subsection{Strategy} The proof of Theorem~\ref{thm:main} relies on a generalization of this principle to the Wigner equation~\eqref{eq:wigner-intro}. To this aim, we need to be able to propagate uniform regularity at the level of the semiclassical Wigner equation~\eqref{eq:wigner-intro}. The main goal is to  get the uniform estimates \eqref{estuniformeintro}.


\medskip 

\noindent $\bullet$ {\bf Propagation of uniform regularity. Bootstrap.} 
The first step is to establish a suitable local well-posedness theory for  the Wigner equation \eqref{eq:wigner-intro}
which is adapted to our purpose.
We shall obtain in Lemma~\ref{lembilinB} several properties of the bilinear  operator $B$ defined in~\eqref{defB-re}; in particular,  continuity estimates in the  weighted Sobolev spaces $\mathcal{H}^m_{r}$, which rely on  commutation properties with the vector fields $V_\pm$. As a consequence, we obtain  the local well-posedness of the Wigner equation in $\H^m_r$ spaces, for $m$ and $r$ larger than $d/2$ (see Proposition~\ref{prop-localexistence}).
The motivation for the use of the  weighted Sobolev spaces $\mathcal{H}^m_{r}$  is the following.
Note that at first sight, we need to propagate a sufficient amount of  weights in $v$ as in the classical case  in order for the density to be defined. Nevertheless, 
without assuming  decay of the Fourier transform of the pair potential $V$,  as previously done  in  \cite{Chaub}, the usual weight
 $v$  is  not  convenient for the analysis  of the Wigner equation, since it does not commute well with $B$: it produces  a loss of an additional
 $\eps$ derivative on $\rho_{f}$. It turns out that the density can be  also controlled from the control of enough powers
 of the vector fields  
  $V_{\pm}$ acting on $f$, see \eqref{embedH02}.  They are better choices since they  have  better  commutations properties with $B$.
  Some control of  the localization in $x$ is  also needed later in the analysis. The weight $x$
   is then  also not well suited for the Wigner equation since it produces a weight $v$ when commuted with the
   free transport $v \cdot \nabla_{x}$.  The natural objects are  instead the vector fields
    $X_{\pm}$ which produce the vector fields $V_{\mp}$ when commuted with the free transport operator.

Then the proof of the main result of the paper, namely Theorem~\ref{thm:main}, relies on a bootstrap argument that we  set up in  Section~\ref{bootstrap}.   For some $m,r \in \mathbb{N}$  and $M>0$ large enough,  we define 
$$
\mathcal{N}_{m,r}(t,f):=\norm{f}_{L^{\infty}\left(0,t;\mathcal{H}^{m-1}_{r}\right)}+\norm{\rho}_{L^{2}\left(0,t;H^{m}_r\right)},
$$
where $\rho(t,x)$ stands for $\rho_{f}(t,x)$ 
and
$$
T_\eps := \sup\left\{T \geq 0,  \, \mathcal{N}_{m,r}(T,f) \leq M\right\}.
$$
The goal is to show that there exist $T_\ast>0$ and $\eps_0>0$, such that $\forall \eps \in (0,\eps_0]$, $T_\eps \geq T_\ast$.  This corresponds to the first part of  Theorem~\ref{thm:main}.
The control of $\mathcal{N}_{m,r}(T_\ast,f)$ will  eventually lead, by a compactness argument, to a derivation of the singular Vlasov equation~\eqref{eq:vlasov-benney-intro}.
The bootstrap argument is formalized in Theorem~\ref{thm:main2}. 
By an energy estimate in the weighted Sobolev spaces $\H^m_r$, 
a control of $\|f \|_{L^\infty(0,T;\H^{m-1}_r)}$ for $T\leq T_\eps$ directly follows from the one  of  $\|\rho \|_{L^2(0,T;H^m_r)}$, which means that the latter  is the key quantity to control.
This main difficulty can thus be seen as similar 
 as the one in  \cite{HKR} for the quasineutral limit\eqref{eq:VP-quasineutre} (see also the introduction of \cite{Chaub} for a presentation on a toy model) and thus we will follow a related strategy.
 The analysis in \cite{HKR} for the estimate of the density without loss of derivatives relies strongly
 on the properties of the  average in time and velocity of the solutions of  the  transport operator
 \begin{equation}
 \label{transport:intro} \partial_{t} + v \cdot \nabla_{x}  - \nabla_{x} \rho(t,x) \cdot \nabla_{v}
 \end{equation} 
 for a given $\rho(t,x)$ smooth enough. One of the main difficulty here will be to develop an appropriate
 quantum analogue where the transport operator is basically replaced,  
 again for a given  $\rho$, 
 by the Wigner operator
$$
\mathscr{T}= \partial_t  + v\cdot \nabla_x + B[\rho,\cdot]
$$
 where we recall that $B$ is defined in~\eqref{defB-re}. In particular, in the case of the cubic NLS, 
 this operator is under the form
 $$ \partial_{t} +  v \cdot \nabla_{x}  +  {i \over \eps}\left( \rho(t, x-{\eps \over 2} D_{v}) -  \rho(t, x+{\eps \over 2} D_{v})\right).$$


\medskip

\noindent $\bullet$ {\bf The extended Wigner system.} 
We thus aim at estimating $\partial_x^\alpha \rho$, for $|\alpha|=m$, in $H^0_r$. As observed in~\cite{HKR,Chaub}, it is not sufficient to apply $\partial^\alpha_x$ for all $|\alpha|=m$ to the Wigner equation~\eqref{eq:wigner-intro}, as this procedure involves terms of type $B[\partial^{\alpha'}_x \rho, \partial^{\alpha''}_x f]$, with $|\alpha'|+ |\alpha''|=m$, $|\alpha'|=1$, which we do not control uniformly on $[0,T_\eps)$. The idea is to consider the full vector of higher derivatives $\mathrm{F}=(\partial^\alpha_x \partial^\beta_v f)_{|\alpha|+|\beta|=m}$, which is shown to satisfy a  pseudodifferential system of the form
\begin{equation}
\label{eq:wigner-higher-intro}
\mathscr{T} \mathrm{F}  +\mathcal{M} \mathrm{F} + { 1 \over \eps} b^\eps_{f}(x,v,D_x)\mathrm{V}_{\rho_\mathrm{F}}={R},
\end{equation}
where $\mathcal{M} $ is a certain matrix-valued pseudodifferential operator, $b_{f}$ a certain symbol related to  $B$ (see \eqref{SinForm}), and $\mathrm{V}_{\rho_\mathrm{F}}$ stands for the vector  $(\partial^\alpha\mathrm{V}_\eps \ast\rho_{\mathrm{F}})_{|\alpha|=m}$. In the right-hand side, $R$ is a well-controlled remainder on $[0,T_\eps)$. The system~\eqref{eq:wigner-higher-intro} is referred to as the \emph{extended Wigner system}. Section~\ref{sec:higher} is precisely dedicated to this second preliminary step.


\medskip
\noindent $\bullet$ {\bf Parametrix for the extended Wigner operator.} By fairly standard arguments, the operator $\mathscr{T}+ \mathcal{M}$ generates a strongly continuous propagator  ${U}_{t,s}$ on $\H^0_{r,0}$
(which is the variant of $\H^0_{r}$ which involves only powers of $V_{\pm}$). 
 The key point to control the regularity of the density will be to prove a quantum analogue
 of the averaging Lemma with gain of one derivative proven in \cite{HKR} in the Vlasov-Benney case.
 However,   contrary to its analogue for the Vlasov case for which the method of characteristics can be naturally  used
to provide  an explicit representation, and eventually to justify that in small time the effect of the free transport
is dominant (see \cite{HKR}), 
we do not have  here at our disposal an explicit tractable  representation formula.
A systematic idea consists in building a parametrix for the extended Wigner operator. 
To simplify, let us neglect the zero order term $\mathcal{M}$ and focus 
on  the scalar operator $\mathscr{T}$. We thus study the linear   semiclassical pseudodifferential equation
$$
\mathscr{T}=\partial_t f + \frac{i}{\eps} a^\eps (t,x,v, D_{x},  D_v)
$$
with symbol $a$ defined by
$$
a(t,z,\xi)=v\cdot\xi_x +{\left(V_\eps \ast \rho \left(t,x-\frac{\xi_v}{2}\right)-V_\eps \ast \rho \left(t,x+\frac{\xi_v}{2}\right)\right)}, 
$$
and the parametrix we look for naturally takes the form of a Fourier Integral Operator (see e.g. \cite{Robert,Zwo}):
\begin{equation}\label{eq:FIO-intro}
{U}^{\mathrm{FIO}}_{t,s} u (z)=\frac{1}{(2\pi)^{2d}} \int_{\xi}\int_{y} e^{\frac{i}{\varepsilon} \left( \varphi_{t,s}^{\eps} (z,\xi)-\langle y,\eps\xi\rangle\right)}b_{t,s}^{\eps}(z,\xi)  u(y)  \, dy d\xi, \quad u \in \mathscr{S}(\R^{2d}),
\end{equation}
where $\varphi$ is the phase and $b$ the amplitude of the FIO. More specifically, we ask that ${U}^{\mathrm{FIO}}_{t,s}$ is such that we have the expansion
\begin{equation}\label{eq:approx-intro}
{U}_{t,s} = {U}^{\mathrm{FIO}}_{t,s} + \eps {U}^{\mathrm{rem}}_{t,s},
\end{equation}
where both ${U}^{\mathrm{FIO}}_{t,s}$ and ${U}^{\mathrm{rem}}_{t,s}$ must be linear continuous operators on $\H^0_{r,0}$, with uniform bound with respect to $\eps$. The term~$\eps {U}^{\mathrm{rem}}_{t,s}$ can be considered as a remainder in the analysis. We are hence led to develop $V_\pm$-weighted $L^2$ continuity results for FIO operators of the form~\eqref{eq:FIO-intro}, with phases satisfying appropriate properties; this is achieved in Appendix~\ref{SFIO}.

 For \eqref{eq:approx-intro} to hold, the phase $\varphi$
must satisfy the eikonal equation
\begin{equation}\label{eq:HJ-intro}
 \left\{
    \begin{aligned}
        &\partial_t{\varphi}_{t,s}+a(t,z, \nabla_z \varphi_{t,s}) =0, \quad z= (x,v) \in \R^{2d}, \, \xi = (\xi_x,\xi_v) \in \R^{2d},\\
	&\varphi_{s,s}(z,\xi)=z\cdot \xi,
    \end{aligned}
\right.
\end{equation}
while $b$ must solve a first order linear equation with coefficients depending on $\nabla_v \varphi$. 
When $a(t,z,\xi)=v\cdot\xi_x$, the eikonal equation~\eqref{eq:HJ-intro} reduces to the free transport equation; the solution is  then explicit, given by $\varphi_{t,s}^{\mathrm{free}}(x,\xi) = (x-(t-s) v)\cdot \xi_x + v \cdot \xi_v$.
One cornerstone of the  proof is the fact that the phase $\varphi_{t,s}$ is close enough (in a precise sense to be specified) to the free phase~$\varphi_{t,s}^{\mathrm{free}}$, see Proposition~\ref{prop-HJ1}. 

Actually, as we need to study the extended Wigner system~\eqref{eq:wigner-higher-intro}, we are enforced to build an approximation to the propagator  associated with $\mathscr{T}+ \mathcal{M}$, 
which leads to the study of a matrix-valued Fourier Integral Operator 
\begin{equation}\label{eq:FIO-matrix-intro}
{U}^{\mathrm{FIO}}_{t,s} u(z)=\frac{1}{(2\pi)^{2d}} \int_{\xi} \int_{y} e^{\frac{i}{\varepsilon} \left( \varphi_{t,s}^{\eps} (z,\xi)-\langle y,\eps\xi\rangle\right)}\mathrm{B}_{t,s}^{\eps}(z,\xi)  u(y)  \, dy d\xi,
\end{equation}
where the amplitude $\mathrm{B}_{t,s}$ is here a matrix. 
We provide complete details to this procedure in Section~\ref{sec:param}.

\medskip

{

\noindent $\bullet$ {\bf  Quantum averaging lemma.}
We are  led to study the following averaging operator
\begin{equation}
\label{def:averagingquantum}
\mathcal{U}_{[\Phi,b,G]}: \quad \varrho(t,x) \mapsto \int_{0}^t \int_{v}\int_{\xi} e^{{i}  \Phi_{t,s} (z,\xi)}b_{t,s}(z,\xi) \widehat{ B[\varrho,G_{t,s}] }(z,\xi)  d\xi  dv ds,
\end{equation}
where the phase $\Phi$ satisfies certain model properties that are verified by the free phase  $\varphi_{t,s}^{\mathrm{free}}$ 
and by $\frac{1}{\eps}\varphi_{t,s}^\eps$ where $\varphi$ is the phase of the FIO from the previous step. Direct estimates for the operator $B$ seem to indicate that the operator $\mathcal{U}_{[\Phi,b,G]}$ is not uniformly bounded with respect to $\eps$ as an operator on $L^2(0,T;L^2(\R^d))$.

In \cite{HKR}, 
we have considered the averaging operator with kernel $H$
\begin{equation}
\label{def:averagingclassical}
\mathcal{U}^{\mathrm{free}}_H: \quad \varrho(t,x) \mapsto \int_0^t \int_{v} \nabla_x \varrho(s,x-(t-s)v)  \cdot H(s,t,x,v) \, dv ds,
\end{equation}
 that is related to the resolution of \eqref{transport:intro} with a special type of source terms
  adapted to the obtention of  a priori estimates for ~\eqref{eq:vlasov-benney-intro}.
We proved that despite the apparent loss of derivative in $x$, $\mathcal{U}^{\mathrm{free}}_H$ is bounded on $L^2(0,T ; L^2(\R^d))$ for all $T>0$, as soon as the  kernel $H(s,t,x,v)$ is sufficiently regular. 
This can be seen as a kinetic averaging lemma, in the spirit of \cite{GLPS}, but tailored for singular Vlasov equations such as Vlasov-Benney.  As a matter of fact, the operator~\eqref{def:averagingclassical} is related to  the operator~\eqref{def:averagingquantum}, when considering the case of the  free phase $\Phi_{t,s}=\varphi_{t,s}^{\mathrm{free}}$ and an amplitude $b_{t,s}$ which does not depend on $\xi$, though a quantum effect remains through the operator $B$.

We shall provide a quantum counterpart of the result of \cite{HKR}, pertaining to the operator $\mathcal{U}_{[\Phi,b,G]}$.   Namely, we shall prove that thanks to fine properties of the phase $\Phi$, if $b,G$ are sufficiently smooth and decaying, then $\| \mathcal{U}_{[\Phi,b,G]} \|_{\mathscr{L}(L^2(0,T;L^2(\R^d)))} \leq C$, uniformly in $\eps$.

The proof of the averaging lemma for \eqref{def:averagingclassical} in the classical case in  \cite{HKR} is based on writing
\begin{align*}
\mathcal{U}^{\mathrm{free}}_H(\varrho)(t,x)&= \frac{1}{(2\pi)^d} \int_0^t \int_{v} \int_\xi \int_{y}   e^{i \left[(x-(t-s)v)- y \right]\cdot \xi } \nabla_x \varrho(s,y)  \cdot H(s,t,x,v) \,  dy d\xi dv ds \\
&= \int_0^t \int_\xi \int_{y}  e^{i (x- y )\cdot \xi } \nabla_x \varrho(s,y)  \cdot \mathcal{F}_v H(s,t,x,(t-s)\xi) \,  dy d\xi ds 
\end{align*}
and then using Bessel-Parseval's formula together with a variant of Schur's test. In the quantum case, as the phase may not be  linear, we cannot proceed similarly. Our approach  for studying $\mathcal{U}_{[\Phi,b,G]}$ consists first in noticing that it can be recast as a pseudodifferential operator in space, associated with an operator-valued symbol in $\mathscr{L}(L^2(0,T))$, that is to say
$$
\mathcal{U}_{[\Phi,b,G]}(\varrho)(\cdot,x)= \int_{\eta} e^{i x\cdot \eta} \mathrm{L}(x,\eta)  \mathcal{F}_x(\varrho) (\cdot,\eta) \, d\eta,
$$
where for all $x,\eta \in \R^d$, $\mathrm{L}(x,\eta) \in \mathscr{L}(L^2(0,T))$. Explicitly we have
\begin{align*}
&[ \mathrm{L}(x,\eta) u] (t)\\
&=2 \int_0^t \int_{v}\int_{\xi=(\xi_x,\xi_v)} e^{-ix\cdot \eta} e^{{i}   \Phi_{t,s} (z,\xi)}  b_{t,s}(z,\xi) \mathcal{F}_{x,v}G(\xi_x-\eta,\xi_v) \frac{1}{\eps} \sin\left(\frac{\eps \xi_v\cdot \eta} {2}\right)   \widehat V(\eps \eta) u(s,\eta)  \,  d\xi dv ds.
\end{align*}
Then the boundedness of $\mathcal{U}_{[\Phi,b,G]}$ on $L^2(0,T;L^2(\R^d))$ follows from showing that $\mathrm{L}$ is a symbol in a class such that a Calder\'on-Vaillancourt theorem for operator-valued symbols can be applied. That $\mathrm{L}$ satisfies suitable properties follows from non-stationary phase estimates, crucially relying on the fine properties of the phase.  
Section~\ref{sec:quantum} is dedicated to this development.

}

\medskip

\noindent $\bullet$ {\bf Quantum Penrose stability.} 
Using the parametrix and~\eqref{eq:approx-intro}, after several reductions, some of them crucially involving applications of our  quantum averaging lemma, we show that $ \partial^\alpha_x \rho $ for $|\alpha|=m$ satisfies an equation of the form
$$
\left(\mathrm{I} -\mathbf{Op}^{\eps,\gamma}_{\mathcal{P}_{\mathrm{quant}}}\right)[e^{-\gamma t} \partial^\alpha_x \rho] = R,
$$
where $\gamma>0$ is a parameter, the symbol $\mathcal{P}_{\mathrm{quant}}(x,\gamma,\tau,\xi)$, defined as
$$
\mathcal{P}_{\mathrm{quant}}(x,\gamma,\tau,\xi) =- 2  \widehat{V}(\xi) \int_0^{+\infty} e^{-(\gamma+i\tau)s} {\sin\left(\frac{s |\xi|^2}{2}\right)} \mathcal{F}_v{f^{0}_\eps}(x, \xi s)\, ds
$$
is nothing but the quantum Penrose function introduced in~\eqref{Penrose-quantique} for $f^0_\eps(x,\cdot)$, 
and $R$ is a controlled remainder.
The quantum Penrose stability condition~\eqref{QCrit} precisely means that the symbol $1-\mathcal{P}_{\mathrm{quant}}$ is uniformly away from $0$, which leads to an uniform estimate of $\|\rho \|_{L^2(0,T;H^m_r)}$, owing to pseudodifferential calculus with (large enough) parameter $\gamma>0$.
This ultimate step is led in Section~\ref{sec:penrose}, and the proofs of Theorem~\ref{thm:main2} and finally Theorem~\ref{thm:main} are completed in Section~\ref{sec:end}.

\subsection{Organization of the paper} This paper is structured as follows.
Section~\ref{sec:preliminary} is mostly dedicated to the local well-posedness theory for the Wigner equation in the $\H^m_r$ spaces (for $m,r>d/2$), and to the setup of the bootstrap argument.  In Section~\ref{sec:higher}, in view of obtaining higher order estimates for the density, we derive the so-called extended Wigner system that is satisfied by derivatives of the solution to the Wigner equation. 
In Section~\ref{sec:param}, we obtain and study a parametrix for the extended Wigner propagator, that takes the form of a Fourier Integral Operator.  This FIO is shown to be bounded in the weighted $\H^0_{r,0}$ spaces. Several fine properties of its phase of are also provided.
In Section~\ref{sec:quantum}, we establish quantum averaging lemmas for a class of  operators related to the latter FIO.  Sections~\ref{sec:penrose} and~\ref{sec:end} correspond to the final stages of the proof of Theorem~\ref{thm:main}. In Section~\ref{sec:penrose}, we apply the parametrix and quantum averaging lemmas to reduce the problem of deriving higher estimates for the density to the study of a semiclassical pseudodifferential equation. 
Finally, the bootstrap is concluded in Section~\ref{sec:end} and the convergence statement is also justified.

The paper ends with the Appendix ~\ref{sec:pseudo} where  several useful results of continuity for pseudodifferential operators and Fourier Integral Operators are collected and proved.
 Section~\ref{Pseudo} is dedicated to a Calder\'on-Vaillancourt result for operator-valued symbols. Section~\ref{SFIO} provides continuity results for FIO, especially in the weighted $\H^0_{r,0}$ spaces, for phases satisfying appropriate properties.
 Eventually, in Section~\ref{sec:PseuDiff}, we present some elements of pseudodifferential calculus with (large) parameter.



\section{Preliminaries for the Wigner equation}
\label{sec:preliminary}

\subsection{Functional inequalities in weighted Sobolev spaces}

 Recall the definition of the vector fields $V_\pm$, $X_\pm$ in \eqref{defvecVX}. As shorthand, we shall sometimes write
 \begin{equation}
 \label{defvecZ} Z_{+}= ( V_{+}, X_{-}), \quad Z_{-}= ( V_{-}, X_{+}), 
 \end{equation} 
 and for $\gamma=(\gamma_x, \gamma_v) \in \mathbb{N}^{d} \times \mathbb{N}^d$, we set
 $$
 Z_+^\gamma =  X_-^{\gamma_x} V_+^{\gamma_v}, \quad  Z_-^\gamma =  X_+^{\gamma_x} V_-^{\gamma_v},
 $$
 so that the ${\H^0_{r}}$ norm as defined in~\eqref{defH0r}  can be recast as
 \begin{equation}
 \label{defH0rbis}  \|f \|_{\H^0_{r}} = \sum_{\substack{ |\beta|  \leq r, |\gamma|  \leq r \\ \beta,  \gamma \in \mathbb{N}^{2d}}} \|Z_{+}^{\beta} Z_{-}^\gamma f \|_{L^2(\mathbb{R}^{2d})}.
 \end{equation}
 In our analysis, we shall also sometimes use another version of weighted spaces where only the vector fields $V_{\pm}$ are involved.
 
 \begin{definition}
 For $r\in \mathbb{N}$, we define the $\H^{0}_{r, 0}$ norm as
 \begin{equation}
 \label{defH0rbisbis}
  \|f\|_{\H^{0}_{r, 0}}
  =  \sum_{\substack{ |\beta|  \leq r, |\gamma|  \leq r \\ \beta,  \gamma \in \mathbb{N}^{d}}} \| V_{+}^\beta V_{-}^{\gamma} f
   \|_{L^2(\mathbb{R}^{2d})}.
 \end{equation}
 \end{definition}
 Note that clearly,  we have the relation $ \| \cdot \|_{\H^0_{r, 0}} \leq  \| \cdot \|_{\H^0_{r}}.$

\bigskip

In the next lemma, we state some properties of the norms $\H^0_r, \H^0_{r,0}$ and $H^0_r$ that will be crucial for the study of the Wigner equation.
  \begin{lem}
  \label{LemHmrHmr}
      The exists $C>0$ such that for every $\eps \in (0, 1]$ and every $ f \in \H^0_{r}$, we have that:
    \begin{align}
   &  \label{embedH01}\| \langle x\rangle^r f \|_{L^2(\mathbb{R}^{2d})} + \| \langle v \rangle^r f \|_{L^2(\mathbb{R}^{2d})}
    + \eps^r \|f \|_{H^r(\mathbb{R}^{2d})}  \leq C \|f \|_{\H^0_{r}}, \quad \forall r \in \mathbb{N},  \\
   & \label{embedH02}   \|  \rho \|_{H^0_{r}} \leq C \|  f \|_{\H^{0}_{r, 0}}, \quad \rho(x)= \int_{\mathbb{R}^d} f(x,v) \, dv, \quad   \forall  \,   r>d/2.
  \end{align}
  \end{lem}
  \begin{Rem}
  Note that an immediate consequence of \eqref{embedH02} is that for all integers $m \geq 0$ and $r>d/2$, it holds
  \begin{equation}
  \label{embedH02m}
   \| \rho \|_{H^m_{r}} \leq C \|f\|_{\H^m_{r}}, \quad \rho(x)= \int_{\mathbb{R}^d} f(x,v) \, dv.
  \end{equation}
  \end{Rem}
  \begin{proof}
  For \eqref{embedH01}, we may just observe that 
  $$ x = \frac{1}{4i}( X_{+}- X_{-}), \quad v = \frac{1}{4i}(V_{+}- V_{-}),  \quad \eps \nabla_{x}= \frac{1}{2}( V_{+}+ V_{-}), \quad
   \eps \nabla_{v}= \frac{1}{2}( X_{+}+ X_{-}).$$
   For \eqref{embedH02}, note first that \eqref{embedH01} combined with the Cauchy-Schwarz inequality only yields
   $$    \| \rho \|_{H^0_{r}} \lesssim  \|f \|_{\H^0_{r+s}},$$ 
    for $s >d/2$, which therefore displays a loss  in terms 
    of the parameter for the weight, in comparison with the claimed \eqref{embedH02}. We thus need something more subtle. 
    
    Let $|\alpha|\leq r$. We first write
    $$
    \| (\eps \partial_x)^\alpha \rho \|_{L^2}^2 = \int_{\R^d}  \left| (\eps \partial_x)^\alpha  \mathcal{F}_v f (x,0) \right|^2 \, dx.
    $$
    Then introduce the function $g$ such that for all $x,\eta \in \R^d$, setting $y_\pm = \frac{1}{2\eps} x \pm   \frac{1}{2} \eta$,
       $$
    \mathcal{F}_v f (x,\eta) = g(y_+,y_-).  
    $$
    (Note that the transform $(x,\eta) \mapsto (y_+, y_-)$ is indeed invertible, with determinant equal to $(-2\eps)^{-d}$.)
    This choice is made so that, denoting $    \widetilde{V}_\pm  = \eps \nabla_x \pm 2 \nabla_\eta$ (which corresponds to the action of $V_\pm$ in the Fourier space in $v$),
\begin{equation}
\label{eq:identityVpm}
    \widetilde{V}_\pm  f (x,\eta)  = \nabla_{y_\pm} g(y_+,y_-).
    \end{equation}
    It follows that
     \begin{align*}
\int_{\R^d}  \left| (\eps \partial_x)^\alpha  \mathcal{F}_v f (x,0) \right|^2 \, dx
 &=  \int_{\R^d}  \left| (\eps \partial_x)^\alpha \left[ g\left(\frac{x}{2\eps},\frac{x}{2\eps}\right)\right]  \right|^2 \, dx \\
  &\lesssim  (2\eps)^d \sum_{\beta\leq \alpha}  \int_{\R^d}  \left| \partial_{y_+}^\beta  \partial_{y_-}^{\alpha-\beta}  g \left(x,x\right) \right|^2 \, dx.
    \end{align*}
     Let us write $g(y_+,y_-)= g_{1}(y_+,y_-)+ g_{2}(y_+,y_-)$ where $\widehat{g_{1}}(\xi_{y_+}, \xi_{y_-})= \widehat{g}(\xi_{y_+},\xi_{y_-}) \mathbf{1}_{|\xi_{y_-}|\leq \xi_{y_+}|}$ and $\widehat{g_{2}} = \widehat{g}- \widehat{g_{1}}$.
     By definition, their Fourier transform in $y_+, y_-$ are such that $\widehat{g_{1}}(\xi_{x}, \xi_{y_-})$ is supported in $| \xi_{y_-} | \leq |\xi_{x}|$
      and $\widehat{g_{2}}(\xi_{x},  \xi_{y_-})$ in  $| \xi_{y_-} | \geq |\xi_{x}|$. By Sobolev embedding with respect to the second variable, we have that for every $y_+,y_- \in \R^d$,
      $$  |(\partial_{y_+}^\beta  \partial_{y_-}^{ \alpha - \beta} g_{1})(y_+, y_-) | \lesssim   \| (\partial_{y_+}^\beta  \partial_{y_-}^{ \alpha - \beta} g_{1})(y_+, \cdot),
      \|_{H^s_{y_-}(\mathbb{R}^d)}$$
      for any $s>d/2$. This yields in particular 
      $$ \| (\partial_{y_+}^\beta  \partial_{y_-}^{ \alpha - \beta} g_{1})(x,x) \|_{L^2(\mathbb{R}^{d})} 
       \lesssim  \| (\partial_{y_+}^\beta  \partial_{y_-}^{ \alpha - \beta} g_{1})
      \|_{L^2_{y_+}(\mathbb{R}^d, H^s_{y_-}(\mathbb{R}^d))} \lesssim   \|  \xi_{y_+}^\beta \xi_{y_-}^{\alpha - \beta}( 1 + |\xi_{y_-}|^s) \widehat{g_{1}}\|_{L^2(\mathbb{R}^{2d})}.$$
      Since $\widehat{g}_{1}$ is supported in  $| \xi_{y_-} | \leq |\xi_{y_+}|$, we therefore get the inequality
      $$ \| (\partial_{y_+}^\beta  \partial_{y_-}^{ \alpha - \beta} g_{1})(x,x) \|_{L^2(\mathbb{R}^{d})}  \lesssim \| | \xi_{y_+}|^{| \alpha  |} ( 1 +  |\xi_{y_-}|^s) \widehat{g}\|_{L^2(\mathbb{R}^{2d})}$$
      and hence since  $r>d/2$, we can  choose $s=r$ 
      and we get that
      $$ \| (\partial_{y_+}^\beta  \partial_{y_-}^{ \alpha - \beta} g_{1})(x,x) \|_{L^2(\mathbb{R}^{d})} 
       \lesssim\sum_{| \alpha'| \leq r, |\alpha''| \leq r} \|  \partial_{y_+}^{\alpha'}   \partial_{y_-}^{\alpha''} g \|_{L^2(\mathbb{R}^{2d})}.$$
       We can use a symmetric argument to estimate  $ \| (\partial_{y_+}^\beta  \partial_{y_-}^{ \alpha - \beta} g_{2})(x,x) \|_{L^2(\mathbb{R}^{d})}$
        by the same quantity, and consequently we obtain that
$$    \| (\eps \partial_x)^\alpha \rho \|_{L^2(\mathbb{R}^d)}  \lesssim   \eps^{d/2}  \sum_{| \alpha'| \leq r, |\alpha''| \leq r} \|\partial_{y_+}^{\alpha'}  \partial_{y_-}^{\alpha''} g \|_{L^2(\mathbb{R}^{2d})}.$$ 
By a final reverse change of variable and using \eqref{eq:identityVpm}, we have
$$
\int_{\R^{2d}} | \partial_{y_+}^{\alpha'}  \partial_{y_-}^{\alpha''} g (y_+,y_-) |^2  \, d{y_+} d{y_-} 
 = \int_{\R^{2d}} | \widetilde{V}_+^{\alpha'}  \widetilde{V}_-^{\alpha''} \mathcal{F}_v f  (x,\eta) |^2  \, \frac{d {x} d{\eta} }{(2\eps)^d},
$$
and we deduce by Bessel-Parseval that 
$$
 \| (\eps \partial_x)^\alpha \rho \|_{L^2(\mathbb{R}^d)}  \lesssim \sum_{|\alpha'|, |\alpha''|\leq r} \| {V}_+^{\alpha'}  {V}_-^{\alpha''}  f \|_{L^2(\mathbb{R}^{2d})} = \| f\|_{\H^0_{r,0}},
$$ 
hence \eqref{embedH02} holds.

  \end{proof}

  We finally state commutator properties of the vector fields $X_\pm$ and $V_\pm$ with the free transport operator $ \mathcal{T}_{0}$:
   $$ \mathcal{T}_{0} := \partial_{t}+ v \cdot \nabla_{x}.$$
    \begin{lem}
  \label{lemcomfree}
  We have the identities:
  \begin{align*}
  &\nabla_{x}\mathcal{T}_{0}= \mathcal{T}_{0} \nabla_{x}, \quad  \nabla_{v} \mathcal{T}_{0}= \mathcal{T}_{0} \nabla_{v} + \nabla_{x},\\
 &  V_{\pm} \mathcal{T}_{0}= \mathcal{T}_{0} V_{\pm}, \quad X_{\pm} \mathcal{T}_{0}= \mathcal{T}_{0} X_{\pm} + V_{\mp}.
  \end{align*}
  \end{lem}

\subsection{Local well-posedness of the Wigner equation}
In this subsection, we discuss the local well-posedness in $\H^m_r$  of the Wigner equation 
\begin{equation}\label{eq:wigner}
 \left\{
    \begin{aligned}
        &\partial_t{f}+v\cdot \nabla_x f +B[\rho, f]  = 0,\\
	&f (0,x,v) = f^0 (x,v),
   \end{aligned}
\right.
\end{equation} 
where $\rho(t,x) = \int_{\R^d} f(t,x,v) \, dv$.
Recalling~\eqref{defB-re}, we have
\begin{equation}
\label{eq:K1}
B[\rho, f] =\frac{i}{\eps} a_{\rho}(t,x,D_v)f , \qquad a_\rho(x,\xi_v) := \mathrm{V}_\rho \left(t,x-\frac{\xi_v}{2}\right)- \mathrm{V}_\rho \left(t, x+\frac{\xi_v}{2}\right),
\end{equation}
 with the notation
\begin{equation}
\label{def:Vrho}
\mathrm{V}_\rho =  V_\eps \ast \rho. 
\end{equation}
In the following, we will often use the notation
$$
\mathscr{T}= \partial_t  + v\cdot \nabla_x + B[\rho,\cdot],
$$
so that~\eqref{eq:wigner} recasts as  $\mathscr{T} f =0$.

It is well-known that $B[\rho,f]$ can be recast in equivalent ways,  this will allow us to choose the most convenient form 
 depending on the type of estimates we perform. 
\begin{lem}
\label{lem:Keps}

The following identities holds for all $\rho \in \mathscr{S}(\R^d), f \in \mathscr{S}(\R^d\times \R^d)$:
\begin{align}
\label{SinForm} B[\rho,f]&= \frac{1}{\eps} b^\eps_{f} (x,v, D_{x}) (\mathrm{V}_\rho),   \qquad b_{f}(x,v,  \xi_x) :=2  \int_{\xi_v}  e^{i v\cdot \xi_v} \sin\left(\frac{ \xi_v\cdot \xi_x} {2}\right) \mathcal{F}_vf(x,\xi_v)d\xi_v,
\end{align}
\begin{equation}
\label{fouriersinform}
\mathcal{F}_{x,v}\left( B[\rho,f]\right) (\xi)=  {(2\pi)^d}   \int_{\eta}   \frac{2}{\eps}  \sin\left(\frac{\eps (\xi_{x}- \eta) \cdot \xi_{v}} {2 }\right)
\widehat{\mathrm{V}_\rho} (\xi_{x} - \eta) \mathcal{F}_{x,v} f( \eta, \xi_{v}) \, d\eta, \quad \xi= (\xi_{x}, \xi_{v}) \in \R^{2d}.
\end{equation}

\end{lem}

\begin{Rem} The symbol $b_{f}$ can be recast as
\begin{equation}
 \label{eq:K3}
b_{f}(x,v,  \xi_x)=\int_{-1/2}^{1/2} (\xi_x\cdot \nabla_v) f(x,v+ \lambda \xi_x) \, d\lambda.
\end{equation}
\end{Rem}
\begin{proof}
For the first identity, we note that by inverse Fourier transform,
$$
\mathrm{V}_\rho \left(x\pm\frac{\varepsilon \xi_v}{2}\right)  = \frac{1}{(2\pi)^d} \int_{\xi_x} e^{i\xi_x \cdot(x\pm \frac{\varepsilon \xi_v}{2})} \widehat{\mathrm{V}_\rho} (\xi_x)  \, d \xi_x,
$$
and consequently, we get
\begin{align*}
B[\rho,f] (x,v)
&= \frac{2}{(2\pi)^d}  \int_{\xi_v}\int_{\xi_x}  e^{-i (v\cdot \xi_v-x\cdot \xi_x)} \frac{1}{\eps}\sin\left(\frac{-\eps \xi_v\cdot \xi_x} {2}\right) \mathcal{F}_vf(x,-\xi_v)  \widehat{\mathrm{V}_\rho}(\xi_x) d\xi_x d\xi_v 
\\&= \frac{2}{(2\pi)^d} \int_{\xi_v}\int_{\xi_x}  e^{i (v\cdot \xi_v+x\cdot \xi_x)} \sin\left(\frac{\eps \xi_v\cdot \xi_x} {2}\right) \mathcal{F}_vf(x,\xi_v) \widehat{\mathrm{V}_\rho} (\xi_x) d\xi_x d\xi_v  ,
\end{align*}
which yields~\eqref{SinForm}. The second identity \eqref{fouriersinform}  follows by  taking the Fourier transform in $(x,v)$ of~\eqref{SinForm} and using again the Fourier inverse formula.  
\end{proof}

It turns out convenient to see $B$ as a bilinear operator as defined below.
\begin{definition}
Let  $B: \mathscr{S}(\R^d)\times  \mathscr{S}(\R^d\times \R^d) \rightarrow  \mathscr{S}(\R^d\times \R^d)$ be the  bilinear operator  defined by its Fourier transform
\begin{equation}
\label{expressiondeB}
\left(\widehat{ B[F,f] }\right) (\xi)= {(2\pi)^d}  \int_{\eta}  \frac{2}{\eps} \sin\left(\frac{\eps (\xi_{x}- \eta) \cdot \xi_{v}} {2}\right) \widehat{V_{\eps}}(\xi_{x}- \eta)
\widehat{F} (\xi_{x} - \eta)  \widehat f( \eta, \xi_{v}) \, d\eta, 
\end{equation}
for $\xi= (\xi_{x}, \xi_{v}) \in \mathbb{R}^{2d}$.
\end{definition}
In the following, it will be sometimes useful to use the decomposition
\begin{equation}
\label{defBpm}
 B[F, f]=  B_{+}[F, f] - B_{-}[F, f]
\end{equation}
where
\begin{align}
\label{defB+}
 \widehat{B_{+}[F, f] }(\xi) & =  {(2\pi)^d}  \int_{\eta}  \frac{1}{i  \eps} e^{i\frac{\eps (\xi_{x}- \eta) \cdot \xi_{v}} {2}} \widehat {V_{\eps}}(\xi_{x}- \eta)
\widehat F (\xi_{x} - \eta)  \widehat f( \eta, \xi_{v}) \, d\eta, \\
\label{defB-} \widehat{B_{-}[F, f] }(\xi) &= {(2\pi)^d}    \int_{\eta}  \frac{1}{i \eps} e^{- i\frac{\eps (\xi_{x}- \eta) \cdot \xi_{v}} {2}} \widehat {V_{\eps}}(\xi_{x}- \eta)
\widehat F (\xi_{x} - \eta)  \widehat f( \eta, \xi_{v}) \, d\eta.
\end{align}
The energy estimates  and  local well-posedness theory in $\H^m_r$  thus rely on continuity properties of the bilinear operator $B$ in $\H^m_r$.
They are established in the next lemma.   
The estimates we wish to prove  rely on improved commutator properties satisfied by $V_+$ with respect to $B_+$ (respectively $V_-$ with respect to $B_-$); these key properties, given in~\eqref{comvectB2}--\eqref{comvectB3} below, further justify the use of these vector fields in the weighted spaces.
 
\begin{lem}
\label{lembilinB}
The operator $B$ satisfies the following properties.

\noindent $\bullet$ {\bf Identities.} It holds
 \begin{equation}
 \label{Breal} \overline{B[F,f]}= B(\overline{F}, \overline{f}),
 \end{equation}
 so that $B[F,f]$ is real-valued if $F, f$ are. Moreover, for $F, f$ real-valued, it holds
 \begin{equation}
 \label{Bskew} \langle B[F, f],  f\rangle= 0,
 \end{equation}
 where $\langle\cdot,\cdot \rangle$ stands for the $L^2$ scalar product.

\noindent $\bullet$ {\bf Weighted Sobolev estimates.} 
For any integer $s>d/2$ and any nonnegative integer $r$, we have
\begin{equation}
 \label{estB1}
 \|B[\partial^\alpha F,  \partial^{\alpha'}f] \|_{\mathcal{H}^0_{r}} \lesssim { 1 \over \eps} \|F \|_{H^s_{r}} \|f\|_{\mathcal{H}^s_{r}}, \quad
 \forall  \alpha, \alpha' \in \mathbb{N}^{2d}, \, | \alpha |+ |\alpha'| \leq s; \\
\end{equation}
\noindent $\bullet$ {\bf Commutation estimates.} For any integer $s>d/2$ and any nonnegative integer $r$, we have
\begin{equation}
\label{estB4}
  \| Z_{+}^\beta Z_{-}^\gamma   B[F,f] - B[F, Z_{+}^\beta  Z_{-}^\gamma f]\|_{L^2} \lesssim   \| \nabla_{x} F \|_{H^s_{r-1}} \|f \|_{\mathcal{H}^0_{r}},  \quad \forall \beta, \gamma \in \mathbb{N}^{2d}, \,  
  |\beta |\leq r, \,  | \gamma | \leq r.
 \end{equation}
Moreover, if $s >1+ d/2$, we have  
\begin{equation}\label{estB3} 
 \| \partial^\alpha (B[F,f]) - B[F, \partial^\alpha f] \|_{\mathcal{H}^0_{r}} \lesssim  \| \nabla_{x} F \|_{H^s_{r}},
  \|f \|_{\mathcal{H}^s_{r}}, \quad   \forall  \alpha\in \mathbb{N}^{2d}, \,  |\alpha  | \leq s,
  \end{equation}
  
\noindent $\bullet$ {\bf Uniform weighted Sobolev estimates.} 
For any integer $s>3+ d/2$, there holds
  \begin{equation}
  \label{estBreste}\| B [\partial^{\alpha'} F, \partial^{\alpha''} f] \|_{\mathcal{H}^0_{r}} \lesssim \|  F \|_{H^s_{r}}
   \|f\|_{\mathcal{H}^{s-1}_{r}}, \quad  \forall \alpha' \in \mathbb{N}^d, \alpha'' \in  \mathbb{N}^{2d}, \,  |\alpha'| + |\alpha''|= s, \, 2 \leq  |\alpha'| \leq s-1.
   \end{equation}
 \end{lem}
 \begin{Rem}
 \label{remHmr}
 Note that from Leibnitz formula, we have the expansion
 $$ \partial^\alpha( B[\rho, f])= \sum_{\alpha'+\alpha''= \alpha} C_{\alpha', \alpha"} B[\partial^{\alpha'}\rho, \partial^{\alpha''}f], \quad \forall \alpha \in \mathbb{N}^{2d}, $$
 where the $C_{\alpha', \alpha"}$ are numerical coefficients and therefore,  \eqref{estB1} yields
 \begin{equation}
 \label{estB1bis}
 \|B[\rho, f]\|_{\mathcal{H}^s_{r}} \lesssim { 1 \over \eps} \|\rho\|_{H^s_{r}} \|f\|_{\mathcal{H}^s_{r}}, \quad
 \forall  s>d/2, r \geq 0,
 \end{equation}
 which means that $B$ is a bounded operator from $H^s_{r} \times \mathcal{H}^s_{r}$ to $\mathcal{H}^s_{r}$, but with a norm that is non-uniform  in $\eps$.
 \end{Rem}
 
 
\begin{proof}[Proof of Lemma \ref{lembilinB}]

We start with the proof of the identities.

\noindent$\bullet$ {\bf Proof of \eqref{Breal} and \eqref{Bskew}.} 
We recall that we have assumed that the pair interaction potential is real and even so that
 its Fourier transform is also real and even.
 To prove \eqref{Breal}, we write
\begin{align*}
\widehat{\overline{ B[F,f] } } (\xi) & = \overline{\widehat{B[F, f]}(-\xi)}=  {(2\pi)^d}  \int_{\eta}  \frac{2}{\eps} \sin\left(\frac{\eps (\xi_{x}+\eta) \cdot \xi_{v}} {2}\right) \widehat{V_{\eps}}(\xi_{x}+ \eta)
 \overline{\widehat {F} (- \xi_{x} - \eta)} \, \overline{  \widehat{f}( \eta, -\xi_{v}) } \, d\eta \\
 &=   {(2\pi)^d}  \int_{\eta}  \frac{2}{\eps} \sin\left(\frac{\eps (\xi_{x}+ \eta) \cdot \xi_{v}} {2}\right)  \widehat{V_{\eps}}(\xi_{x}+ \eta)
\widehat{\overline{ F}} ( \xi_{x}  + \eta)   \widehat{ \overline{ f}}(- \eta, \xi_{v}) \, d\eta =  \widehat{ B[\overline{F}, \overline{f}] }(\xi),
\end{align*}
where for the final identity we have just used the change of variable $\eta \mapsto -\eta$ in the integral.

Finally, for \eqref{Bskew}, using Bessel-Parseval, for $F, f$ real-valued, we write  that
\begin{align*} \langle B[F,f], f\rangle&= (2\pi)^{2d} \int_{\xi}
 \widehat{ B[F,f]}(\xi) \overline{\widehat{f}(\xi)} \, d\xi \\
 &= (2\pi)^{3d} \int_{\xi} \int_{\eta}  \frac{2}{\eps} \sin\left(\frac{\eps (\xi_{x}- \eta) \cdot \xi_{v}} {2}\right)  \widehat{V_{\eps}}(\xi_{x} -  \eta)
\widehat F (\xi_{x} - \eta)  \widehat f( \eta, \xi_{v}) \widehat f( - \xi) \, d\eta d\xi.
\end{align*}
By exchanging the roles of $\eta$ and $\xi_{x}$  and by oddness of the sin function, we infer that
$$
\langle B[F,f], f\rangle = -\langle f,  B[F,f]\rangle,
$$
hence the proof of the second identity.

\bigskip

\noindent$\bullet$ {\bf Convolution inequalities.} For the estimates, we shall first  establish a useful  elementary  convolution estimate.
Define the bilinear operator  $K: \mathscr{S}(\R^d)\times  \mathscr{S}(\R^d\times \R^d) \rightarrow  \mathscr{S}(\R^d\times \R^d)$ as
$$ K[h, g](\xi) =\int_{\eta} h(\xi_{x}- \eta) g(\eta, \xi_{v}) \, d\eta, \quad \xi=(\xi_x, \xi_{v}) \in \mathbb{R}^{2d}.$$
 Then, we have that
 \begin{equation}
 \label{usefulK}
  \|K[h,g]\|_{L^2} \lesssim \|h\|_{L^1} \|g\|_{L^2}, \quad   \|K[h,g]\|_{L^2} \lesssim\|h\|_{L^2} \|g\|_{L^2_{\xi_{v}} L^1_{\xi_{x}}}.
  \end{equation}
  Indeed, from the Young inequality for the convolution in $\xi_{x}$, we first get that for any $\xi_{v}$
  $$  \|K[h,g](\cdot, \xi_{v})\|_{L^2} \lesssim \|h\|_{L^1} \|g(\cdot, \xi_{v})\|_{L^2_{\xi_x}}.$$ Alternatively, we can also get
  $$ \|K[h,g](\cdot, \xi_{v})\|_{L^2} \lesssim   \|h\|_{L^2} \|g(\cdot, \xi_{v}) \|_{L^1_{\xi_{x}}}$$
  and we conclude in both cases by taking the $L^2$ norm in $\xi_{v}$.

  \bigskip
  
\noindent$\bullet$ {\bf Proof of \eqref{estB1}.} We start with \eqref{estB1}.
  Let us recall that by assumption, the pair potential $V$ satisfies $\widehat{V} \in \mathscr{C}^\infty_b(\R^d)$, so that any occurence of $\widehat{V}$ or of its derivatives in the estimates can be directly bounded.
For $r=0$, thanks to the Bessel-Parseval identity, we observe that
we just need to study $ \|  \widehat{B[\partial^\alpha F, \partial^{\alpha'} f]}(\xi) \|_{L^2}$. From the definition \eqref{expressiondeB}, we first  have the rough estimate 
 $$| \widehat{B[\partial^\alpha F, \partial^{\alpha'} f]}(\xi)|
  \lesssim { 1 \over \eps}  K[ | \xi_{x}|^{|\alpha|}|\widehat F|,  |\xi|^{|\alpha'|}|\widehat f|](\xi)
   \lesssim  { 1 \over \eps}   K[ \langle  \xi_{x}\rangle^{s}|\widehat F|, |\widehat f|](\xi) + { 1 \over \eps }K[ |\widehat F|, \langle  \xi \rangle^{s} |\widehat f|](\xi).$$
   Consequently, by using \eqref{usefulK}, we obtain that 
  $$ \|   \widehat{B[\partial^\alpha F, \partial^{\alpha'} f]} \|_{L^2} \lesssim {1 \over \eps} \left(\| \langle \xi \rangle^s \widehat F \|_{L^2} \|\widehat f\|_{L^2_{\xi_{v}}L^1_{\xi_{x}}} 
   +  \|\widehat F\|_{L^1} \| \langle \xi \rangle^s\widehat f\|_{L^2}\right).$$
   By using  Bessel-Parseval and  observing that 
   \begin{equation}
   \label{bizarre}
     \|\widehat F\|_{L^1} \lesssim \|F\|_{H^s(\mathbb{R}^d)}, \quad \| \widehat f \|_{L^2_{\xi_{v}}L^1_{\xi_{x}}} \lesssim \|f\|_{H^s(\mathbb{R}^{2d})}, 
    \quad s>d/2
    \end{equation}
   we finally get
   \begin{equation}
   \label{OKpourK}
  \|  B [\partial^\alpha F, \partial^{\alpha'} f] \|_{L^2}   \lesssim  { 1 \over \eps }\|F\|_{H^s} \|f\|_{\mathcal{H}^s_{0}}.
   \end{equation}
   This yields \eqref{estB1} for $r=0$.
   
   \begin{Rem}
   \label{remBpm0}
   Note that by repeating the above arguments, we also have that 
   \begin{equation}
   \label{OKpourKpm}
   \|  B_{\pm}[\partial^\alpha F, \partial^{\alpha'} f] \|_{L^2}   \lesssim  { 1 \over \eps }\|F\|_{H^s} \|f\|_{\mathcal{H}^s_{0}}, \, | \alpha |+ | \alpha '|=s >d/2
   \end{equation}
   for the operators $B_\pm$ defined in~\eqref{defB+}--\eqref{defB-}. 
     \end{Rem}

 To get \eqref{estB1}  for integers $r> 0$, recalling the definitions \eqref{defH0r} and \eqref{defvecZ}, we need to estimate
 $$ \|Z_{+}^\beta Z_{-}^\gamma B[ \partial^\alpha F, \partial^{\alpha'} f]\|_{L^2}, \quad |\beta|+ | \gamma | \leq r.$$
 We first observe that we have the following commutator formula for every $F, \, f$
 \begin{align}
 \label{comvectB1}
   & X_{\pm} B_{+}[F, f] = B_{+}[F, X_{\pm}f],  \quad   X_{\pm} B_{-}[F, f] = B_{-}[F, X_{\pm}f], \\
  \label{comvectB2}   & V_{-} B_{+}[F, f]=  B_{+}[ 2 \eps \nabla_{x} F, f] + B_{+}[ F, V_{-}f], \quad V_{-} B_{-}[F, f]= B_{-}[F, V_{-}f], \\
  \label{comvectB3}  & V_{+} B_{-}[F,f]= B_{-} [ 2 \eps \nabla_{x} F, f] + B_{-}[ F, V_{+}f], \quad V_{+} B_{+}[F, f]= B_+[F, V_{+}f].
 \end{align}
 In other words, both $X_{+}$ and $X_{-}$ commute with $B_{+}[F, \cdot]$ and $B_{-}[F, \cdot]$ so that they also commute with $B$,
 whereas $V_{+}$  (resp. $V_{-}$) only commutes with $B_{+}[F, \cdot]$ (resp. with $B_{-}[F, \cdot]$).
 
 To get \eqref{estB1}, it suffices to estimate separately the terms involving $B_{+}$ and those involving $B_{-}$. We shall perform  the estimate for
 $$  \|Z_{+}^\beta Z_{-}^\gamma B_{+}[ \partial^\alpha F, \partial^{\alpha'} f]\|_{L^2}$$ 
 the other one being similar. By using the above commutator formulas, we observe that the term $Z_{+}^\beta Z_{-}^\gamma B[ \partial^\alpha F, \partial^{\alpha'} f]$ can be expanded as
 \begin{equation}
 \label{expandZB1} Z_{+}^\beta Z_{-}^\gamma B_{+}[ \partial^\alpha F, \partial^{\alpha'} f]= 
  \sum_{\substack{\gamma_{1}' + \gamma_{1}^{''}= \gamma_1 \\ |\beta_{1}|+ |\beta_{2}|\leq r,  |\gamma_{1}|+ |\gamma_{2}| \leq r}}  C_{\gamma_{1}, \gamma_{1}',\gamma_2, \beta_1,\beta_2}
   B_{+} \left[ (\eps \partial)^{\gamma_{1}'} \partial^\alpha F, V_{+}^{\beta_{1}} X_{-}^{\beta_{2}} V_{-}^{\gamma_{1}''} X_{+}^{\gamma_{2}} \partial^{\alpha'}f,
   \right]\end{equation}
   where the $C_{\gamma_{1}, \gamma_{1}',\gamma_2, \beta_1,\beta_2}$ are numerical coefficients.
  Since  the commutators $[\nabla_{x}, X_{\pm}]$ and $[\nabla_{v}, V_{\pm}]$ are constant (with constant independent of $\eps$), 
    we are reduced to estimating in $L^2$  terms under the form
    $$ B_{+} \left[ \partial^\alpha (\eps \partial)^{\beta_{1}}  F,  \partial^{\alpha''} V_{+}^{\beta_{1}'} X_{-}^{\beta_{2}'} V_{-}^{\gamma_{1}'''} X_{+}^{\gamma_{2}'} f
   \right], $$
   where  $ | \alpha''| \leq  |\alpha'|,$  $|\gamma_{1}'''| \leq |\gamma_{1}''|$, $|\gamma_{2}'| \leq | \gamma_{2}|$,  $| \beta_{i}'| \leq | \beta_{i}|$, $i=1, \, 2$.
   By using Remark \ref{remBpm0}, we thus obtain
   \begin{align*} \left\| B_{+} \left[ \partial^\alpha (\eps \partial)^{\beta_{1}}  F,  \partial^{\alpha''} V_{+}^{\beta_{1}''} X_{-}^{\beta_{2}'} V_{-}^{\gamma_{1}'} X_{+}^{\gamma_{2}'} f
   \right] \right\|_{L^2} &\lesssim {1 \over \eps} \| (\eps \partial)^{\beta_{1}}  F\|_{H^s} \| V_{+}^{\beta_{1}''} X_{-}^{\beta_{2}'} V_{-}^{\gamma_{1}'} X_{+}^{\gamma_{2}'} f \|_{\mathcal{H}^s_{0}} \\ &\lesssim { 1 \over \eps } \|F\|_{H^s_{r}} \| f \|_{\mathcal{H}^s_{r}}.
   \end{align*}
   This ends the proof of \eqref{estB1}.
   
   \bigskip
   
\noindent$\bullet$ {\bf Proof of \eqref{estB4}.}  To obtain \eqref{estB4}, we note that \eqref{expandZB1} yields
   \begin{equation}
   \label{expandZB+}\left[ Z_{+}^\beta Z_{-}^\gamma,  B_{+}[F,  \cdot ] f \right]=  \sum_{\substack{\gamma_{1}' + \gamma_{1}^{''} = \gamma_{1}, \gamma_{1}' \neq 0 \\ \cdots}} C_{\gamma_{1}, \gamma_{1}',\ldots}
   B_{+} \left[ (\eps \partial)^{\gamma_{1}'} F, V_{+}^{\beta_{1}} X_{-}^{\beta_{2}} V_{-}^{\gamma_{1}''} X_{+}^{\gamma_{2}} f\right].
   \end{equation}
Similarly, by using again the commutator relations \eqref{comvectB1}, \eqref{comvectB2}, \eqref{comvectB3}, we have
   \begin{equation}
   \label{expandZB-}  \left[ Z_{+}^\beta Z_{-}^\gamma,  B_{-}[F,  \cdot ] f \right]=  \sum_{\substack{\beta_{1}' + \beta_{1}^{''} = \beta_{1}, \beta_{1}' \neq 0 \\ \cdots }} C_{\beta_{1}, \beta_{1}', \ldots }
   B_{-} \left[ (\eps \partial)^{\beta_{1}'} F, V_{+}^{\beta_{1}''} X_{-}^{\beta_{2}} V_{-}^{\gamma_{1}} X_{+}^{\gamma_{2}} f\right].
   \end{equation}
   We can estimate separately the contributions of the two sums in $L^2$.  We shall only give the details for the estimate of
   $$  \left\|B_{+} \left[ (\eps \partial)^{\gamma_{1}'} F, V_{+}^{\beta_{1}} X_{-}^{\beta_{2}} V_{-}^{\gamma_{1}''} X_{+}^{\gamma_{2}} f\right]\right\|_{L^2}$$
   where $ \gamma_{1}' + \gamma_{1}^{''} = \gamma_{1}, \, \gamma_{1}' \neq 0, $ $|\beta_{1}|+ |\beta_{2}|\leq r$, 
    $|\gamma_{1}| + | \gamma_{2}| \leq r.$ Let us set $g = V_{+}^{\beta_{1}} X_{-}^{\beta_{2}} V_{-}^{\gamma_{1}''} X_{+}^{\gamma_{2}} f$.
    Since $\gamma_{1}'\neq 0$, we have by using the Fourier transform that
    $$ \left| \widehat{ B_{+} \left[ (\eps \partial)^{\gamma_{1}'}F, g\right]} (\xi) \right|
    \leq K\left[ |\xi|  | \eps \xi|^{ |\gamma_{1}'|-1} |\widehat{F}|,  | \widehat{g}| \right](\xi)$$
    and we deduce by using the first estimate of \eqref{usefulK} and \eqref{bizarre} that for $s>d/2$:
    $$ \left\| B_{+} \left[ (\eps \partial)^{\gamma_{1}'}F, g\right]\right\|_{L^2}
     \lesssim \| \langle \eps \partial \rangle^{r-1} \nabla F \|_{H^s} \|g \|_{L^2} \lesssim \| \nabla F\|_{H^s_{r-1}} \|f\|_{\mathcal{H}^0_{r}}.$$
     This concludes the proof of \eqref{estB4}.
   \begin{Rem}
   \label{remarkr00}
   Note that by using similar estimates,  we can  also get  the following variants which will be also useful.
   A version of \eqref{estB4}  where  only  the vector fields $V_{\pm}$ are involved, holds:
   \begin{equation}
   \label{comutqueV} \|V_{+}^\beta V_{-}^\gamma B [F, f] - B[F, V_{+}^\beta V_{-}^\gamma f] \|_{L^2} \lesssim \| \nabla F \|_{H^s_{r-1}} \| f\|_{\H^0_{r,0}}, \quad
   \forall (\beta, \gamma) \neq (0, 0), \, | \beta| \leq r, \, | \gamma| \leq r.
   \end{equation}
   Moreover, we also have the following  variants of \eqref{estB1} which are  useful when either $F$ or $f$  is smoother:
   \begin{equation}
   \label{utilepourFIO}
    \|B[F, f]\|_{\H^0_{r, 0}} \lesssim  {1 \over \eps }\| F\|_{H^{s}_{r}} \|f \|_{\H^0_{r, 0}}, \quad s>d/2, \, r\in \mathbb{N}, 
   \end{equation}
    \begin{equation}
   \label{utilepourFIO2}
    \|B[F, f]\|_{\H^0_{r, 0}} \lesssim  {1 \over \eps }\| F\|_{H^{0}_{r}} \|f \|_{\H^s_{r, 0}}, \quad s>d/2, \, r\in \mathbb{N}.
   \end{equation}
Similar estimates also hold for $B^\pm$.
   \end{Rem}  
   
 \noindent$\bullet$ {\bf Proof of \eqref{estB3}.}   To prove \eqref{estB3}, we need to estimate
   $$ \left\| B[\partial^{\alpha'} F, \partial^{\alpha''}f] \right\|_{\mathcal{H}^0_{r}}, \quad |\alpha'| + | \alpha''| \leq s, \, \alpha' \neq 0.$$
   By using again the expansion \eqref{expandZB1}, we have to estimate  three types of terms:
    \begin{align*}
    I& =  \left\| B[\partial^{\alpha'} F,  Z_{+}^\beta Z_{-}^\gamma \partial^{\alpha''}f] \right\|_{L^2}, \quad | \beta | \leq r, \, | \gamma | \leq r, \\
    II & = \left\|  B_{+} \left[ (\eps \partial)^{\gamma_{1}'} \partial^{\alpha'}F, V_{+}^{\beta_{1}} X_{-}^{\beta_{2}} V_{-}^{\gamma_{1}''} X_{+}^{\gamma_{2}} \partial^{\alpha''} f\right]\right\|_{L^2}, \,  | \gamma_{1}'|+ | \gamma_{1}''| + | \gamma_{2}| \leq r, \, | \beta_{1}| + | \beta_{2}| \leq r, \, \gamma_{1}' \neq 0, \\
    III & = \left\|  B_{-} \left[ (\eps \partial)^{\beta_{1}'} \partial^{\alpha'}F, V_{+}^{\beta_{1}''} X_{-}^{\beta_{2}} V_{-}^{\gamma_{1}} X_{+}^{\gamma_{2}} \partial^{\alpha''} f\right]\right\|_{L^2}, \,  | \beta_{1}'|+ | \beta_{1}''| + | \beta_{2}| \leq r, \, | \gamma_{1}| + | \gamma_{2}| \leq r, \, \beta_{1}' \neq 0.
    \end{align*}
  To estimate $I$, by commuting $\partial$ and $X_{\pm}$, $V_{\pm}$, we observe that it suffices to estimate
   $$ \widetilde{I}=   \left\| B[\partial^{\alpha'} F,  \partial^{\alpha'''} (Z_{+}^\beta Z_{-}^\gamma )f] \right\|_{L^2}, \quad | \beta | \leq r, \, | \gamma | \leq r, \,
    |\alpha'|+ | \alpha'''| \leq s, \, \alpha' \neq 0.$$
    Let us set $g= Z_{+}^\beta Z_{-}^\gamma f$. By using the expression \eqref{expressiondeB}
     and the inequality $|\sin u| \leq |u|$,  we get that
     \begin{multline*}
      \left| \widehat{ B[ \partial^{\alpha'} F, \partial^{\alpha''} g ]}(\xi)\right|
     \lesssim K\left[ | \xi_{x}|^{ 1 + |\alpha'|} | \widehat F|,  |\xi|^{  1 + |\alpha''|} | \widehat{g}|\right](\xi) \\
      \lesssim       K\left[ | \xi_{x}| \langle\xi \rangle^{s} | \widehat F|,  | \xi | | \widehat{g}|\right](\xi) + K\left[ | \xi_{x}|^{ 2} | \widehat F|,  \langle \xi\rangle^{ s } | \widehat{g}|\right](\xi).
       \end{multline*}
       Consequently, by using \eqref{usefulK} and \eqref{bizarre} with $s-1$  (since we are assuming that $s>1+d/2$), we obtain that
       $$  \widetilde{I} \lesssim \| \nabla_{x} F \|_{H^{s}} \|g\|_{\mathcal{H}^s_{0}} \lesssim  \| \nabla_{x} F \|_{H^{s}} \|f\|_{\mathcal{H}^s_{r}}.$$ 
       This yields 
       \begin{equation}
       \label{IcomB}
        I \lesssim\| \nabla_{x} F \|_{H^{s}} \|f\|_{\mathcal{H}^s_{r}}.
        \end{equation}
  To estimate $II$,  as before, we can commute the derivatives with the vector fields so that we have to estimate
  $$  \widetilde{II} = \left\|  B_{+} \left[  \partial^{\alpha'} (\eps \partial)^{\gamma_{1}'} F,  \partial^{\alpha'''}V_{+}^{\beta_{1}} X_{-}^{\beta_{2}} V_{-}^{\gamma_{1}''} X_{+}^{\gamma_{2}}  f\right]\right\|_{L^2},  $$
  with $  | \gamma_{1}'|+ | \gamma_{1}''| + | \gamma_{2}| \leq r, \, | \beta_{1}| + | \beta_{2}| \leq r, \, \gamma_{1}' \neq 0$ and $|\alpha'|+ | \alpha'''| \leq s.$ 
 Consequently, by using again  \eqref{OKpourKpm}
  we obtain
  $$ \widetilde{II} \lesssim { 1 \over \eps} \|  (\eps \partial)^{\gamma_{1}'} F \|_{H^s} \|f\|_{\mathcal{H}^s_{r}}.$$
  Since $| \gamma_{1}'| >0$, this yields
  \begin{equation}
  \label{IIcomB}
   II \lesssim \|  \nabla_{x}F \|_{H^s_{r-1}} \|f\|_{\mathcal{H}^s_{r}}.
   \end{equation}
   In a similar way, we obtain that
    \begin{equation}
  \label{IIIcomB}
   III \lesssim \|  \nabla_{x}F \|_{H^s_{r-1}} \|f\|_{\mathcal{H}^s_{r}}.
   \end{equation}
   The estimate \eqref{estB3} then follows from  a combination of \eqref{IcomB}, \eqref{IIcomB}, \eqref{IIIcomB}.
  
  \bigskip
  
\noindent$\bullet$ {\bf Proof of \eqref{estBreste}.}   To get \eqref{estBreste}, we can again  after commutator with the vector fields reduce the estimate to the one of  three types of terms
   \begin{align*}
    I_{\mathrm{rem}}& =  \left\| B[\partial^{\alpha'} F,  \partial^{\alpha''} Z_{+}^\beta Z_{-}^\gamma f] \right\|_{L^2}, \quad | \beta | \leq r, \, | \gamma | \leq r, \\
    II_{\mathrm{rem}} & = \left\|  B_{+} \left[  \partial^{\alpha'} (\eps \partial)^{\gamma_{1}'}F,  \partial^{\alpha''} V_{+}^{\beta_{1}} X_{-}^{\beta_{2}} V_{-}^{\gamma_{1}''} X_{+}^{\gamma_{2}} f\right]\right\|_{L^2}, \,  | \gamma_{1}'|+ | \gamma_{1}''| + | \gamma_{2}| \leq r, \, | \beta_{1}| + | \beta_{2}| \leq r, \, \gamma_{1}' \neq 0, \\
    III_{\mathrm{rem}} & = \left\|  B_{-} \left[ \partial^{\alpha'} (\eps \partial)^{\beta_{1}'} F, \partial^{\alpha''}  V_{+}^{\beta_{1}''} X_{-}^{\beta_{2}} V_{-}^{\gamma_{1}'} X_{+}^{\gamma_{2}} f \right]\right\|_{L^2}, \,  | \beta_{1}'|+ | \beta_{1}''| + | \beta_{2}| \leq r, \, | \gamma_{1}| + | \gamma_{2}| \leq r, \, \beta_{1}' \neq 0
    \end{align*} 
    with  $2\leq |\alpha'| \leq s-1$ and $\alpha''$ is now such that $| \alpha '' | \leq s - |\alpha'|.$ 
    
    For $I_{\mathrm{rem}}$,  by setting $g= Z_{+}^\beta Z_{-}^\gamma f$, we obtain similarly to above  that 
     \begin{align*}
      \left| \widehat{ B[ \partial^{\alpha'} F, \partial^{\alpha''} g ]}(\xi)\right|
    &\lesssim K\left[ | \xi_{x}|^{ 1 + |\alpha'|} | \widehat F|,  |\xi|^{  1 + |\alpha''|} | \widehat{g}|\right](\xi) \\
      &\lesssim  K\left[ | \xi_{x}|^{ 3} | \widehat F|,  \langle \xi\rangle^{ s-1} | \widehat{g}|\right](\xi) + 
       K\left[ | \xi_{x}| \langle\xi \rangle^{s-1} | \widehat F|,  | \xi | \langle \xi \rangle | \widehat{g}|\right](\xi).
       \end{align*}
      Note that  we have used that $|\alpha'| \geq 2$ and that $|\alpha'| \leq s-1$ to handle the case $\alpha''=0$.
       Consequently, by using \eqref{usefulK} and \eqref{bizarre} with $s-3$ (since we are assuming here  that $s>3+d/2$), we obtain that
       $$  I_{\mathrm{rem}} \lesssim \|  F \|_{H^{s}} \|g\|_{\mathcal{H}^{s-1}_{0}} \lesssim  \|  F \|_{H^{s}} \|f\|_{\mathcal{H}^{s-1}_{r}}.$$ 
       
      For $II_{\mathrm{rem}}$,   we can set again $g=  V_{+}^{\beta_{1}} X_{-}^{\beta_{2}} V_{-}^{\gamma_{1}''} X_{+}^{\gamma_{2}} f$
       and $G= (\eps \partial)^{\gamma_{1}'}F$ with $\gamma_{1}'  \neq 0$ so that we have to estimate
       $$  \left\|  B_{+} \left[  \partial^{\alpha'} G,  \partial^{\alpha''} g\right]\right\|_{L^2}.$$
       We have assumed that $|\alpha'| \geq 2$. When additionally $\alpha'' \neq 0$, we can write
       $$  \left\|  B_{+} \left[  \partial^{\alpha'} G,  \partial^{\alpha''} g\right]\right\|_{L^2}
       =   \left\|  B_{+} \left[  \partial^{\widetilde{\alpha}'} \partial^{e'} G,  \partial^{ \widetilde{\alpha}''} \partial^{e''} g\right]\right\|_{L^2}$$
       where $|e'|=2, \, |e''|=1$ and thus $|\widetilde{\alpha}'|+ | \widetilde{\alpha}''|\leq s-3$. Since $s-3>d/2$, we can use \eqref{OKpourKpm} with $s-3$
        instead of $s$, this yields
      $$  \left\|  B_{+} \left[  \partial^{\alpha'} G,  \partial^{\alpha''} g\right]\right\|_{L^2} \lesssim { 1 \over \eps}  \| \nabla^2 G \|_{H^{s-3}}
      \|\nabla_{x,v} g\|_{\mathcal{H}^{s-3}_{0}}.$$
      When $\alpha''=0$, we can rely on the assumption that $|\alpha'| \leq s-1$  and just use that
      $$ | \widehat{B_{+} \left[  \partial^{\alpha'} G,  g\right]}(\xi)| \lesssim { 1 \over \eps} K\left[|\xi|^{|\alpha'} |\widehat G|, |\widehat{g}| \right](\xi)$$
      and the second inequality in \eqref{usefulK} and \eqref{bizarre}.
      Since $s-3>d/2$ and $| \alpha' | \leq s-1$, this yields 
      $$  \|B_{+} \left[  \partial^{\alpha'} G,  g\right]\|_{L^2} \lesssim { 1 \over \eps} \| G \|_{H^{s-1}} \|g\|_{\mathcal{H}^{s-3}_{0}}.$$
       We thus obtain in all cases  the estimate 
     \begin{align*} II_{\mathrm{rem}} \lesssim   \|B_{+} \left[  \partial^{\alpha'} G,  \partial^{\alpha''}g\right]\|_{L^2} \lesssim { 1 \over \eps} \| G \|_{H^{s-1}} \|g\|_{\mathcal{H}^{s-2}_{0}}
        \lesssim  \|  F \|_{H^{s}_{r-1}}  \|f\|_{\mathcal{H}^{s-2}_{r}}
       \end{align*}
        since $\gamma_{1}' \neq 0.$ A similar estimate can be obtained for $III_{\mathrm{rem}}$ and consequently, \eqref{estBreste} follows,
and this ends the proof of the proposition.
  \end{proof}

We conclude this subsection with the local well-posedness result in $\mathcal{H}^{m}_{r}$ for $m,r>d/2$.

\begin{prop} 
\label{prop-localexistence}
The Wigner equation \eqref{eq:wigner} is locally well-posed in $\mathcal{H}^{m}_{r}$ for all integers $m,r>d/2$:  if $f^{0} \in \mathcal{H}^{m}_{r}$, there exists $T>0$ (which may depend on $\varepsilon$) such that there is a unique  solution $f \in \mathscr{C}([0,T];\mathcal{H}^{m}_{r})$ of \eqref{eq:wigner}.
Moreover, if $f_{0}$ is real-valued, $f$ also is.
\end{prop}
\begin{proof}
For the existence part,  using the characteristics of the free transport, it is equivalent to solve  the integral  equation
\begin{equation}
\label{pff0} f(t,x,v)= f^0(x-vt, v) - \int_{0}^t B\left[\rho_{f}(s), f(s)\right] (x-(t-s)v, v) \, ds.
\end{equation}
Defining the bilinear operator
$$ \mathcal{B}[ g  ,f ] (t,x,v)= - \int_{0}^t B\left[ \rho_{g} (s), f(s)\right] (x-(t-s)v, v) \, ds, $$
a solution is therefore given by a fixed point of the map $f \mapsto f^{0}(x-vt, v) +  \mathcal{B}[ f, f]$.
Note that it holds
\begin{equation*}
  \|\mathcal{B}[ g,f ](t) \|_{\mathcal{H}^m_{r}} \lesssim ( 1 + t^{m}) \int_{0}^t \| B\left[ \rho_{g} (s), f(s)\right]\|_{\mathcal{H}^m_{r}}.
\end{equation*}
By using \eqref{estB1bis}, we have
\begin{equation}
\label{pff1}   \|\mathcal{B}[ g,f ](t) \|_{\mathcal{H}^m_{r}} \lesssim {1 \over \eps} ( 1 + t^{m}) \int_{0}^t \| \rho_{g} (s) \|_{H^m_{r}} \|f(s) \|_{\H^m_{r}}\, ds.\end{equation}
Consequently, thanks to the estimate \eqref{embedH02m},   we get 
the bilinear estimate
 $$  \|\mathcal{B}[g,f ]\|_{L^\infty(0, T; \mathcal{H}^m_{r})} \lesssim { 1 \over \eps}( 1 + T^{m}) T  \|g \|_{L^\infty(0, T;\mathcal{H}^m_{r})}
    \| f \|_{L^\infty(0,T; \mathcal{H}^m_{r})}.$$  
 This allows to get existence for small times thanks to the Banach fixed point Theorem  and also uniqueness of the solution.
 Note finally that if $f^{0}$ is real, then thanks to \eqref{Breal}, $\overline{f}$ is solution of the same equation with the same initial data,
 so $f = \overline f$  by uniqueness and $f$ is real.
\end{proof}

\subsection{The bootstrap argument}\label{bootstrap}

The proof of Theorem~\ref{thm:main} relies on a bootstrap argument that we initiate in this final subsection. For $m, r > d/2$ (to be fixed large enough), thanks to Proposition~\ref{prop-localexistence}, there exists a maximal lifespan $T^*>0$ and a maximal solution $f \in \mathscr{C}\left([0,T^{*});\mathcal{H}^{m}_{r}\right)$ to the Wigner equation
 \eqref{eq:wigner}. Though $f$ depends on $\eps$, we de not specify it explicitely for the sake of readability.
 In the same way, $\rho$ will now stand for  $\rho_{f}$.
 
   For $t \in [0,T^{*})$, consider the functional
$$
\mathcal{N}_{m,r}(t,f):=\norm{f}_{L^{\infty}\left(0,t;\mathcal{H}^{m-1}_{r}\right)}+\norm{\rho}_{L^{2}\left(0,t;H^{m}_{r}\right)}.
$$
The functional $ \mathcal{N}_{m,r}(t,f )$ is well-defined and is continuous with respect to $t$ on $[0,T^*)$. This allows  to consider for some parameter $M>0$ to be chosen  appropriately later,  
\begin{equation*}
T_\eps=\sup \left\{T\in[0,T^{*}),\mathcal{N}_{m,r}(T,f)\leq M\right\}.
\end{equation*}
By taking $M$ large enough (at the very least $M > \| f^0\|_{\mathcal{H}^{m}_{r}}$), we have by continuity that $T_\eps>0$. The goal is to show that,  up to choosing  the value of  $M$  large enough (but  independent of $\varepsilon$), $T_\eps$ is uniformly bounded from below by some time $T^{\#}>0$.  This is formalized in the following statement.

\begin{thm}
\label{thm:main2}
With the same assumptions as in Theorem~\ref{thm:main},
there exist $M>0$, $\eps_0>0$ and $T^{\#}>0$, such that, for all $\eps \in (0,\eps_0]$, there is a unique solution $f \in C([0,T^{\#}];\H^{m}_{r})$ of the Wigner equation \eqref{eq:wigner}. Furthermore the following estimate holds: 
$$
\sup_{\eps\in(0,\eps_0]} \mathcal{N}_{m,r}(T^{\#},f)\leq M.
$$
\end{thm}
 This corresponds to the first part of Theorem~\ref{thm:main}; 
 the convergence statement is a consequence which will be obtained in Section~\ref{sec:end}. 
 
 \bigskip
 
 Note that from the definitions of $T_\eps$ and $T^*$, the following alternative holds: either $T_\eps=T^{*}$, or $T_\eps < T^{*}$ and $\mathcal{N}_{m,r}(T_\eps,f )= M$. Let us analyze the first case. If $T_\eps=T^{*}=+\infty$, then $\mathcal{N}_{m,r}(T ,f )\leq M$ for every $T>0$ and therefore Theorem~\ref{thm:main2} holds automatically; we thus only need to study the subcase $T_\eps=T^{*}<+\infty$. As a matter of fact,  this subcase is impossible. Indeed, the following estimate holds.

\begin{lem}\label{PE}
Assume that $T_\eps< + \infty$, then the solution $f $  of \eqref{eq:wigner},satisfies, for some $C>0$ independent of $\varepsilon$, the estimate 
\begin{equation*}
\underset{[0,T_\eps)}{\sup}\norm{f (t)}_{\mathcal{H}^{m}_{r}}\lesssim  (1 + T_\eps^{m} ) \norm{f ^{0}}_{\mathcal{H}^{m}_{r}}\exp\left[\frac{C(1 + T_\eps^{m} )}{\varepsilon} T_\eps^{1/2} M\right].
\end{equation*}
\end{lem}
\begin{proof}
By \eqref{pff0}, \eqref{pff1}, we have that
 for $t\in[0,T_\eps)$,
\begin{equation*}
\norm{f(t)}_{\mathcal{H}^{m}_{r}}\lesssim   (1 + T_\eps^{m} ) \left( \norm{f^{0}}_{\mathcal{H}^{m}_{r}} +\int_{0}^{t}  \frac{1}{\varepsilon} \norm{\rho(s)}_{H_r^{m}} \norm{f(s)}_{\mathcal{H}^{m}_{r}} ds \right).
\end{equation*} 
From the Gronwall inequality, we deduce that
$$  \norm{f(t)}_{\mathcal{H}^{m}_{r}}\lesssim   (1 + T_\eps^{m} )  \norm{f^{0}}_{\mathcal{H}^{m}_{r}}  \exp\left( { C (1 + T_\eps^{m} ) \over \eps }\int_{0}^t \| \rho(s) \|_{H^m_{r}}\, ds\right)$$ 
for some $C>0$ independent of $\eps$ 
and the lemma follows from an application of the Cauchy-Schwarz inequality,
since we have $\mathcal{N}_{m,r}(T_\eps,f)\leq M$.
\end{proof}
Applying Lemma~\ref{PE} in the subcase $T_\eps=T^{*}<+\infty$, we obtain that 
\begin{equation*}
\underset{[0,T^{*})}{\sup}\norm{f(t)}_{\mathcal{H}^{m}_{r}}\lesssim(1 + (T^*)^{m} ) \norm{f ^{0}}_{\mathcal{H}^{m}_{r}}\exp\left[\frac{C(1 + (T^*)^{m} )}{\varepsilon} (T^*)^{1/2} M\right].
\end{equation*}
This means that the solution could be continued beyond $T^{*}$, which contradicts its definition. As a result, this case cannot occur and we can therefore focus on the second case.

\bigskip

\fbox{Namely, until the end of the paper, we assume that $T_\eps < T^{*}$ and $\mathcal{N}_{m,r}(T_\eps,f )=M$. 
 }
 \bigskip

 The goal  will be  to 
 find some time $T^{\#}>0$ independent of $\varepsilon$, such that 
\begin{equation*}
\mathcal{N}_{m,r}(T^{\#},f )<M.
\end{equation*}
This  will prove that $T_\eps>T^{\#}>0$.
To this end, we need to  uncover an improved  estimate of $\mathcal{N}_{m,r}(T,f )$ for sufficiently small $T<T_\eps$. 

We first provide a control of the term $\norm{f}_{L^{\infty}\left(0,t;\mathcal{H}^{m-1}_{r}\right)}$ which can be obtained by an energy estimate  and the bilinear estimates of Lemma~\ref{lembilinB}.
\begin{lem}
\label{lem-1ertermeN}
Assume that $r>d/2$ and that $m >2 + d/2$. The solution $f$  of \eqref{eq:wigner} satisfies for all $T\in[0,T_\eps]$  the estimate
\begin{equation}
\label{eq:estim-m-1}
\underset{[0,T]}{\sup}\norm{f}_{\mathcal{H}^{m-1}_{r}}\lesssim \norm{f^{0}}_{\mathcal{H}^{m-1}_{r}} \:+ \sqrt{T}\Lambda (T,M).
\end{equation}
\end{lem}
To prove this estimate, we   first need to commute derivatives and  the vector fields $Z_\pm$  with the  Wigner equation  \eqref{eq:wigner}. 
 
 \begin{definition} 
 \label{eq:defI1}
 For $\alpha=(\alpha_{x}, \alpha_{v})  \in  \mathbb{N}^{2d}$, $i=1,\ldots,d$, if $\alpha_{v, i} \neq 0$, we define  $  {\alpha}^{i,+,-} =( {\alpha}^{i,+,-}_{x}, {\alpha}^{i,+,-}_{v})\in \mathbb{N}^{2d}$ by
\begin{equation}
{\alpha}^{i,+,-}_{x,j}= \alpha_{x,j} + \delta_{i,j}, \quad  {\alpha}^{i,+,-}_{v,j}= \alpha_{v,j} - \delta_{i,j} \quad j=1,\ldots,d.
\end{equation}
\end{definition}
Note that we have $| {\alpha}^{i,+,-}|= | \alpha|.$
\begin{proof}[Proof of Lemma~\ref{lem-1ertermeN}] 
 Since  $\mathscr{T}f=0$, we get 
 for $| \alpha | \leq  m-1$, $ \alpha= (\alpha_{x}, \alpha_{v}) \in \mathbb{N}^{2d}$,  that
 \begin{equation}
 \label{comTdalpha} \mathscr{T} \partial^\alpha f = - \sum_{j=1}^d \alpha_{v, j} \partial^{{\alpha}^{j,+,-}} f - \left[ \partial^\alpha, B[\rho, \cdot] \right] f.
 \end{equation}
 Next, for $\beta, \gamma \in \mathbb{N}^{2d}$, $|\beta| \leq r$, $|\gamma| \leq r$, we obtain that 
\begin{equation}
\label{energypoids}
\mathscr{T}  Z_{+}^\beta Z_{-}^\gamma  \partial^\alpha f = - \sum_{i=1}^4 \mathscr{S}_{i}, 
\end{equation}
where
\begin{align*}
& \mathscr{S}_{1}=  \sum_{j=1}^d \alpha_{v, j} Z_{+}^\beta Z_{-}^\gamma \partial^{{\alpha}^{j,+,-}} f, \qquad \mathscr{S}_{2}=  Z_{+}^\beta Z_{-}^\gamma  \left[ \partial^\alpha, B[\rho, \cdot] \right] f,  \\
&  \mathscr{S}_{3}= \left[ Z_{+}^\beta Z_{-}^\gamma, B[\rho, \cdot] \right] \partial^\alpha f, \qquad \mathscr{S}_{4}= - \left[  Z_{+}^\beta Z_{-}^\gamma, v \cdot \nabla_{x} \right]\partial^\alpha f.
\end{align*}
Thanks to the identities in Lemma \ref{lemcomfree}, $\mathscr{S}_4$ can be expanded as 
$$
\mathscr{S}_{4} =\sum_{|\widetilde \beta| \leq r, \, | \widetilde \gamma|\leq r} C_{\beta, \gamma, \widetilde \beta, \widetilde \gamma}  Z_{+}^{\widetilde\beta} Z_{-}^{\widetilde\gamma} \partial^\alpha f,
$$
where  $C_{\beta, \gamma, \widetilde \beta, \widetilde \gamma}$ are numerical coefficients. 
We first clearly get that
$$ \| \mathscr{S}_{1}\|_{L^2}  \lesssim \| f \|_{\mathcal{H}^{m-1}_{r}}.$$
For $\mathscr{S}_{2}, $ we can use \eqref{estB3} with $s=m-1$ since $m>2+d/2$, this yields
$$  \| \mathscr{S}_{2}\|_{L^2} \lesssim  \| \nabla_{x} \rho \|_{H^{m-1}_{r}} \|f\|_{\H^{m-1}_{r}}.$$
For $\mathscr{S}_{3}$, we use \eqref{estB4}, again with $s=m-1$, we also  obtain
$$  \| \mathscr{S}_{3}\|_{L^2} \lesssim  \| \nabla_{x} \rho \|_{H^{m-1}_{r-1}} \|\partial^\alpha f\|_{\H^{0}_{r}}
 \lesssim   \| \nabla_{x} \rho \|_{H^{m-1}_{r}} \| f\|_{\H^{m-1}_{r}}.$$
 Finally, from the expansion of $\mathscr{S}_{4}$, we  have the estimate
 $$  \| \mathscr{S}_{4}\|_{L^2} \lesssim \|f\|_{\H^{m-1}_{r}}.$$
By taking the $L^2$ scalar product of \eqref{energypoids} with $ Z_{+}^\beta Z_{-}^\gamma \partial^\alpha f$, 
 we   get from standard integration by parts and the above estimates  that
 $$ \frac{1}{2} { d \over dt}  \| Z_{+}^\beta Z_{-}^\gamma \partial^\alpha f \|_{L^2}^2 + \langle B[  \rho, Z_{+}^\beta Z_{-}^\gamma \partial^\alpha f], Z_{+}^\beta Z_{-}^\gamma \partial^\alpha  f \rangle
  \lesssim (1 + \| \rho \|_{H^m_r} )\| f \|_{\mathcal{H}^{m-1}_{r}}^2.$$  
  Since $f$ and thus $\rho_{f}$ are real-valued, the second term in the left hand side vanishes thanks to \eqref{Bskew}. 
  By integrating in time and  summing on the multi-indices, we get
  $$ \|f(t) \|_{\mathcal{H}^{m-1}_{r}}^2 \lesssim \|f^0 \|_{\mathcal{H}^{m-1}_{r}}^2 + \int_{0}^t ( 1 + \| \rho(s) \|_{H^m_r})
  \,  \|f(s) \|_{\mathcal{H}^{m-1}_{r}}^2\, ds$$
  and therefore, we infer from the Gronwall inequality that
  $$  \|f(t) \|_{\mathcal{H}^{m-1}_{r}}^2 \lesssim \|f^0 \|_{\mathcal{H}^{m-1}_{r}}^2  \exp\left(C \int_{0}^t (1 + \|\rho(s) \|_{H^m_r} ) \, ds\right),$$
    for some $C>0$ independent of $\eps$.
Since  $\mathcal{N}_{m,r}(T,f)\leq M$, by the Cauchy-Schwarz inequality, this yields
  $$  \|f(t) \|_{\mathcal{H}^{m-1}_{r}} \lesssim \|f^0 \|_{\mathcal{H}^{m-1}_{r}} \exp\left( {C \over 2 }(T + M T^{1 \over 2})\right).$$ 
  The result follows since 
  $ e^x \leq 1 + xe^x$, for  $x \geq 0$.

\end{proof}


\section{The extended Wigner system}
\label{sec:higher}

As set up in the previous section, we work on the interval $[0,T_\eps]$, where
$$
T_\eps=\sup \left\{T\in[0,T^{*}), \, \mathcal{N}_{m,r}(T,f)\leq M\right\}.
$$
in which we recall $\mathcal{N}_{m,r}(T,f )= \norm{f}_{L^{\infty}\left(0,T;\mathcal{H}^{m-1}_{r}\right)}+\norm{\rho}_{L^{2}\left(0,T;H_r^{m}\right)}$.
With the aim to estimate $\norm{\rho}_{L^{2}\left(0,T;H^{m}_r\right)}$, we look for an equation satisfied by $\pa^\alpha_x \rho$, for $|\alpha|=m$. To this end, it seems natural to apply the operator $\pa^\alpha_x$ to the Wigner equation and integrate with respect to $v$. However this approach is not directly conclusive since commutators with $B$  involve the control of terms of the form $\pa^{\widetilde \alpha}_x \pa^{\widetilde\beta}_v f$ with $|\widetilde\alpha|+|\widetilde\beta|=m$ and $|\widetilde\beta|=1$, which are not controlled by $\mathcal{N}_{m,r}(T,f)$ and thus cannot be estimated uniformly with respect to $\eps$. To bypass this issue, as in \cite{HKR} for the case of the Vlasov equation, we look for an equation for the full vector of  higher derivatives $(\pa^\alpha_x \pa^\beta_v f)_{|\alpha|+|\beta|=m}$. 


\subsection{Applying derivatives to the Wigner equation} 
The aim is now to uncover the structure of the system satisfied by the partial derivatives
 $ \partial^\alpha_{x,v} f $ for $ \alpha = (\alpha_{x}, \alpha_{v}) \in \mathbb{N}^{2d }, \, | \alpha | = m$.
 Let us define $\mathcal{E}_{m}= \{  \alpha= (\alpha_{x}, \alpha_{v}) \in \mathbb{N}^{2d}, \, | \alpha | = m\}, $ and $N_{m}= \mathrm{card} (\mathcal{E}_{m}).$
 It turns out  convenient to fix a parametrization of $\mathcal{E}_{m}$  by $ \llbracket 1,  N_{m}\rrbracket $, denoted by 
  $ \alpha : \llbracket 1,  N_{m}\rrbracket \rightarrow   \mathcal{E}_{m}$
  and to define a vector
   $\mathrm{F} \in \mathbb{R}^{N_{m}}$ such that
   $\mathrm{F}_{i}= \partial^{\alpha(i)} f$.
   We choose the parametrization with the additional property that
   $\alpha_{v,i}= 0$ for all  $i \in \llbracket 1, n_{m}\rrbracket$ where
   $$ n_{m}= \mathrm{card} \{ \alpha \in \mathbb{N}^d, \, | \alpha | =m\},$$
   so that the first $n_{m}$ components of $\mathrm{F}$ correspond to partial derivatives in $x$ only.
  
\begin{lem}\label{RReste} 
There exist $(c_{p,k}), (d_{p,k}) \in \mathcal{M}(\R^{N_m})$, and a function $\beta: \llbracket 1,N_m\rrbracket \times \llbracket1,N_m\rrbracket \rightarrow \mathbb{N}^d$, with $|\beta|=2$, such that the following holds. 
Define the matrix-valued pseudodifferential operator (with symbol in $ \mathcal{M}(\R^{N_{m}})$) by
\begin{equation}
\label{def:M}
\mathcal{M}\mathrm{F}= \mathfrak{m}_\rho^\eps(t,x, D_{v}) \mathrm{F}, \quad
( \mathfrak{m}_\rho)_ {p,k}(t, x, \xi_{v})= c_{p,k} + d_{p,k}  \int_{-{ 1 \over 2}}^{ 1 \over 2} \partial_x^{\beta(p,k)}  \mathrm{V}_{\rho}(t, x + \lambda \xi_{v})\, d\lambda.
 \end{equation}
Let  $f$ be the solution of the Wigner equation \eqref{eq:wigner}. The vector $\mathrm{F} = (\partial^{\alpha(i)} f)_{i \in \llbracket1, N_m\rrbracket}$  satisfies the system
\begin{equation}\label{IJ}
\mathscr{T} \mathrm{F}  +\mathcal{M}\mathrm{F} + { 1 \over \eps} b^\eps_{f}(x,v,D_x)\mathrm{V}_{\rho_\mathrm{F}}=\mathcal{R},
\end{equation} 
where $\mathrm{V}_{\rho_\mathrm{F}}=  (\partial^{\alpha(i)}  \mathrm{V}_{\rho})_{i \in \llbracket1, N_m\rrbracket} $.
Moreover, the remainder $\mathcal{R}$ satisfies 
\begin{equation}\label{EstV}
\norm{\mathcal{R}}_{L^{2}(0,T;\mathcal{H}^{0}_{r})}\leq  \Lambda(T,M).
\end{equation}

\end{lem}

In the following, it will be also convenient to  use the notation
$$ (B[ \mathrm{V}_{\rho_\mathrm{F}}, f ])_{i} = B[ (\mathrm{V}_{\rho_\mathrm{F}})_{i}, f], \quad i \in \llbracket 1, N_m\rrbracket,$$
so that we can write
\begin{equation}
\label{Betendu}
{ 1 \over \eps}  b^\eps_{f} (x,v,D_x)\mathrm{V}_{\rho_\mathrm{F}} = B[ \mathrm{V}_{\rho_\mathrm{F}}, f],
\end{equation}
thanks to Lemma \ref{lem:Keps}. To summarize, we can recast~\eqref{IJ} as
\begin{equation}\label{eq:extended-wigner}
\mathscr{T} \mathrm{F}  +\mathfrak{m}_{\rho}^\eps(t,x, D_{v}) \mathrm{F} + B[ \mathrm{V}_{\rho_\mathrm{F}}, f]=\mathcal{R}.
\end{equation}

\begin{proof}[Proof of Lemma~\ref{RReste}]
By  further expanding  \eqref{comTdalpha} (in the case $\beta= \gamma=0$), we obtain that for all $\alpha= (\alpha_x,\alpha_v) \in \mathbb{N}^{2d}$,
\begin{equation*}
\mathscr{T}(\partial^\alpha f )+\mathbf{1}_{\vert \alpha_{x}\vert =m}B[\partial^{\alpha} \rho, f]+ \mathcal{P} _{\alpha}=\mathcal{R}_{\alpha},
\end{equation*} 
in which, recalling Definition~\ref{eq:defI1} for ${\alpha}^{j,+,-}$,
\begin{align*}
\mathcal{P}_{\alpha}&=  \sum_{j=1}^d  \alpha_{v,j} \partial^{{\alpha}^{j,+,-}}f + \sum_{\substack{\gamma\leq \alpha_x \\ |\gamma |= 1}}  c_{\alpha, \gamma} B[\partial^\gamma_x \rho ,\partial^{\alpha_x - \gamma }_x \partial^{\alpha_v}_v f],\\
\mathcal{R}_{\alpha}&= - \sum_{k=2}^{\min(\vert \alpha_{x}\vert ,m-1)}\sum_{\substack{\sigma \leq \alpha_x \\ \vert\sigma\vert=k}}c_{\alpha ,\sigma}\: B[\partial^{\sigma}_x \rho, \partial^{\alpha_{x}-\sigma}_x \partial^{\alpha_{v}}_v f], \end{align*}
where $c_{\alpha,\gamma}$, $c_{\alpha,\sigma}$ are numerical coefficients. In the case $\vert \alpha_{x}\vert \leq 1$, we note that $\mathcal{R}_{\alpha}=0$. For $\mathcal{P} _{\alpha}$, according to \eqref{eq:K1}, and using a Taylor expansion, we have
\begin{multline*}
B[\partial^\gamma \rho ,\partial^{\alpha_x - \gamma }_x \partial^{\alpha_v}_v f]
\\
=  - i 
\sum_{|\gamma'|=1} \int_{-1/2}^{1/2} \int_{\xi_{v}} \int_{w}  e^{i (v-w) \cdot \xi_{v}} \partial_x^{\gamma+\gamma'} \mathrm{V}_{\rho}(x + \eps \lambda \xi_{v})  \cdot \partial^{\gamma'} _w \partial^{\alpha_x - \gamma }_x \partial^{\alpha_v}_v f(x, w) \, dw d\xi_{v}   d\lambda.
\end{multline*}
and we can therefore use the indexing explained in the beginning of the subsection to write the contribution of such terms  of $\mathcal{P}_{\alpha}$ as in~\eqref{def:M}.

To conclude, it remains to estimate $\norm{\mathcal{R}_{\alpha}}_{L^{2}(0,T;\mathcal{H}^{0}_{r})}$ in order to show that $\mathcal{R}_{\alpha}$ is a controlled remainder.
We only need observe that all terms in the sum are under the form
$ B[ \partial_{x}^\sigma \rho, \partial_{x}^{\alpha_x- \sigma} \partial_{v}^{\alpha_v} f ]$ with $2\leq |\sigma|\leq m-1$, so that the estimate follows from \eqref{estBreste}  (since $m>3+d/2$). We get 
\begin{equation*}
\norm{\mathcal{R}_{\alpha}}_{\mathcal{H}^{0}_{r}}\lesssim   \norm{f}_{\mathcal{H}^{m-1}_{r}}\norm{\rho}_{H^{m}_r}.
\end{equation*}
Consequently, by definition of $T_\eps$, we have for every $0\leq T< T_\eps$ that
\begin{equation*}\label{EstRLD2}
\norm{\mathcal{R}}_{L^{2}(0,T;\mathcal{H}^{0}_{r})}\leq \Lambda(T,M),
\end{equation*}
hence concluding the proof.
\end{proof}

\begin{definition}
We shall refer to the system \eqref{IJ} in the sequel as the \emph{extended Wigner system}.  The (matrix-valued) operator $ \mathscr{T} +\mathcal{M}$ will be called the \emph{extended Wigner operator}.
\end{definition}

\subsection{The propagator associated with the extended Wigner operator}

In this subsection, we provide some properties of the propagator associated with the extended  Wigner operator operator $\mathscr{T} +\mathcal{M}$, which will allow to use the Duhamel formula to express the solution to   the extended Wigner system~\eqref{IJ}. 

\begin{lem} 
\label{lem-semigroup}
For all matrix-valued map $\mathrm{G}^0(x,v)  \in \H^0_{r, 0}$ and for every $s \in [0, T_\eps]$,   there exists a unique solution on $[0,T_\eps]$ to the problem 
\begin{equation}\label{TransportG}
(\mathscr{T} +\mathcal{M}) \mathrm{G} = 0, \quad \mathrm{G} |_{t=s}= \mathrm{G}^0(x,v).
\end{equation}
This solution is denoted by $U_{t,s} \mathrm{G}^0$ and $U_{t,s}$ is referred to as the propagator associated with the extended  Wigner operator $\mathscr{T} +\mathcal{M}$. It satisfies the uniform estimate, for all $T\in[0,T_\eps]$,
\begin{equation}\label{MajOpU}
\sup_{0\leq t,s\leq T} \norm{U_{t,s}}_{\mathscr{L}(\mathcal{H}^0_{r, 0})}\leq \Lambda(T,M).
\end{equation}
\end{lem}

\begin{proof} 
The equation \eqref{TransportG} can be at first  seen as a forced  free  transport equation, so that 
it is equivalent to 
\begin{equation}
\label{pff0-re} 
\begin{aligned}
G(t,x,v)&= G^0(x-v(t-s), v) - \int_{s}^t B\left[\rho(\tau), G(\tau)\right] (\tau, x-v(t-\tau), v) \, d\tau \\
&\qquad - \int_s^t \mathcal{M} G (\tau,x-v(t-\tau), v)  \, d\tau \\
\end{aligned}
\end{equation}
A local  solution can thus be obtained in short time  by a fixed argument as in the proof of  Proposition \ref{prop-localexistence}. Note that
 since $V_\pm$ commute with the free transport operator, by  using~\eqref{utilepourFIO} in  Remark \ref{remarkr00} (as $m> d/2 +2$), we have
 $$\| \mathcal{M}G (\tau,x-v(t-\tau), v)\|_{\H^0_{r,0}}  \leq \Lambda(T, M) \| G(\tau)\|_{\H^0_{r,0}} $$
 and by an estimate similar   to~\eqref{pff1}, it holds
 $$\left\| B\left[\rho(\tau), G(\tau)\right] (\tau, x-v(t-\tau), v)  \right\|_{\H^0_{r,0}}
 \leq  { 1 \over \eps} \Lambda (T, M)  \| G(\tau)\|_{\H^0_{r,0}}.
 $$
 We can then justify that the unique  local solution can  be continued on the whole  $[0, T_\eps]$ by using the Gronwall Lemma.


To obtain the uniform estimate~\eqref{MajOpU},  we proceed by energy estimates
 as in Lemma~\ref{lem-1ertermeN}. We  once again rely on the fact  that $V_{\pm}$ commute with the free transport operator,  on~\eqref{Bskew} and  on the bilinear estimate \eqref{comutqueV} to treat the contribution of $B[\rho,G]$.
We thus get as in the proof of  Lemma~\ref{lem-1ertermeN} that for $|\beta|, |\gamma|\leq r$, and all $p,k \in \llbracket1,N_m\rrbracket$,  and all $T \in [0,T_\eps]$,
 $$ \frac{1}{2} { d \over dt}  \| V_{+}^\beta V_{-}^\gamma \partial^\alpha G_{p,k} \|_{L^2}^2 
  \lesssim \Lambda(T,M) \| G \|_{\mathcal{H}^{0}_{r}}^2.$$  
Summing on all $\beta,\gamma$ and all $p,k$ yields the claimed result.

 \end{proof}

 \medskip
 
Applying Lemma~\ref{lem-semigroup}, we get that the solution to the exended Wigner system \eqref{IJ} can be rewritten as 
\begin{equation}\label{LinT}
\mathrm{F}=U_{t,0}  \mathrm{F}^0 - { 1 \over \eps }\int_0^t U_{t,s} b_{f}^\eps(s, x, v,  D_{x}) \mathrm{V}_{\rho_\mathrm{F}}\, ds  +  \int_{0}^t U_{t,s}\mathcal{R}(s)\,ds.
\end{equation}
Integrating with respect to $v$, we obtain a system of equations for $\rho_{\mathrm{F}}= (\pa^{\alpha(i)} \rho)_{i=1,\cdots,N_m}$, which is the starting point for obtaining $H^m_{r}$ estimates for $\rho$. The next goal of the analysis is to recast in a more tractable form the  Duhamel
term
\begin{equation}
\label{eq:need}
{ 1 \over \eps}\int_0^t \int_v U_{t,s} b_{f}^\eps(s,x, v,  D_{x})  \mathrm{V}_{\rho_\mathrm{F}}  \, dv ds
\end{equation}
and in particular prove that it is uniformly bounded in $L^2(0,T; H^0_r)$, for $T$ small enough. 


\section{Parametrix for the extended Wigner system}\label{sec:param}

To handle~\eqref{eq:need}, we need to put forward a smoothing effect due to the integration in time and velocity, in the style of kinetic averaging lemmas \cite{GLPS,DPLM,AL,JLT}; specifically we shall provide a quantum analogue of the averaging lemma of \cite{HKR}, that we briefly alluded to in the sketch of proof of Section~\ref{sec:strategy}.
The propagator  $U_{t,s}$ (introduced in Lemma~\ref{lem-semigroup}) associated with the extended  Wigner system~\eqref{IJ} is formally related to the propagator of the transport equation associated  with the Vlasov equation in the semiclassical limit $\eps\to 0$, which  suggests that a quantum analogue of \cite{HKR} may hold. However, a direct perturbative analysis is  not possible and,  in the way we have obtained it,    $U_{t,s}$  is  a too abstract object to be  useful   to perform a precise analysis. 

In this section, our goal is to build an explicit approximation of $U_{t,s}$, i.e. a \emph{parametrix} for the extended Wigner operator, under the form of a Fourier Integral Operator (FIO). We specifically look for a matrix-valued operator ${U}^{\mathrm{FIO}}_{t,s} \in \mathscr{L}(\H^0_{r,0})$ satisfying
\begin{equation}\label{eq:approx}
{U}_{t,s} = {U}^{\mathrm{FIO}}_{t,s} + \eps {U}^{\mathrm{rem}}_{t,s},
\end{equation}
where ${U}^{\mathrm{rem}}_{t,s} \in \mathscr{L}(\H^0_{r,0})$, so that terms due to $ \eps {U}^{\mathrm{rem}}_{t,s}$ will be considered as remainders  thanks to the gain of the  factor $\eps$.
According to~\eqref{eq:approx}, the study of~\eqref{eq:need} will then be reduced to that of 
\begin{equation}
\label{eq:need2}
{1 \over \eps }\int_0^t \int_v U_{t,s}^{\mathrm{FIO}}(b^\eps_{f}(s, x, D_{x}) \mathrm{V}_{\rho_\mathrm{F}})  \, dv ds,
\end{equation}
which will be the focus of the forthcoming Section~\ref{sec:quantum}. 
%
%
%

\subsection{General scheme of the construction}

From~\eqref{def:M}, 
the extended Wigner operator $\mathscr{T} + \mathcal{M}$ is a
pseudodifferential operator under the form 
\begin{equation}
\label{eq:higher-wigner-pseudo}
 \mathscr{T} + \mathcal{M} = \partial_t   + v\cdot \nabla_{x} +  \frac{i}{\eps}a_{\rho}^\eps(t,x, D_{v}) + \mathfrak{m}_\rho^\eps(t,x, D_{v}),
\end{equation}
where we recall
$a_{\rho}(t,x,\xi_v)=\mathrm{V}_\rho \left(t,x-\frac{\xi_v}{2}\right)-\mathrm{V}_\rho  \left(t,x+\frac{\xi_v}{2}\right)$
and the (matrix-valued) symbol  $\mathfrak{m}_\rho(t,x, \xi_{v})$ is defined in \eqref{def:M}.

It is thus natural to look for  a parametrix ${U}^{\mathrm{FIO}}_{t,s}$ under the form of a FIO, that is to say
\begin{equation}\label{FIOGen}
U_{t,s}^{\mathrm{FIO}}u(z)=\frac{1}{(2\pi)^{2d}}\int_{\xi}  \int_{y} e^{\frac{i}{\varepsilon} \left( \varphi_{t,s}^{\eps} (z,\xi)-\langle y,\eps\xi\rangle\right)}\mathrm{B}_{t,s}^{\eps}(z,\xi)  u(y) dyd\xi,
\end{equation}
where $\varphi$ is a phase and $\mathrm{B}$ a (matrix-valued) amplitude. 

Before getting into the details of the construction of the FIO, we state  a general lemma which will allow
to get the decomposition \eqref{eq:approx}.

\begin{lem} 
\label{lem:reduction-propagator}
Let $T\in [0,T_\eps]$. Assume that there exist two operators  $U_{t,s}^{\mathrm{FIO}}$ and $V_{t,s}^{\mathrm{rem}}$ such that, for  some
 $r \in \mathbb{N}$  and some $C>0$,
\begin{equation}\label{eq:MajFIO}
\sup_{0\leq t,s\leq T} \norm{U_{t,s}^{\mathrm{FIO}}}_{\mathscr{L}(\mathcal{H}^0_{r, 0})}+ \sup_{0\leq t,s\leq T} \norm{V_{t,s}^{\mathrm{rem}}}_{\mathscr{L}(\mathcal{H}^0_{r, 0})} \leq  C,
\end{equation}
and which satisfy for all $0\leq s,t\leq T$, the equation
\begin{equation}\label{FIOeq}
 \left\{
    \begin{array}{ll}
    (\mathscr{T}+\mathcal{M})U_{t,s}^{\mathrm{FIO}}=\eps V_{t,s}^{\mathrm{rem}},\\
        U_{s,s}^{\mathrm{FIO}} = \mathrm{I}.
 \end{array}
\right.
\end{equation}
Then defining for all $0\leq s,t\leq T$
\begin{equation}
\label{eq:lem-defU}
U_{t,s}^{\mathrm{rem}}   := -\int_s^t U_{t,\tau} V_{\tau,s}^{\mathrm{rem}} \,d\tau ,
\end{equation}
we have
\begin{equation}\label{MajUr}
 \sup_{0\leq t,s\leq T} \norm{U_{t,s}^{\mathrm{rem}}}_{\mathscr{L}(\mathcal{H}^0_{r, 0})} \leq  C^2 T,
\end{equation}
and it holds
\begin{equation}
\label{eq:lem-appr}
{U}_{t,s} = {U}^{\mathrm{FIO}}_{t,s} + \eps {U}^{\mathrm{rem}}_{t,s}.
\end{equation}
\end{lem}


\begin{proof}
Let $u^0 \in \H^0_{r, 0}$ and  $t,s\in [0,T_\eps]$. 
Introducing
$$u^{\mathrm{rem}} (t,s,z)= U_{t,s}u^0(z)-U_{t,s}^{\mathrm{FIO}}(t,0)u^0(z),$$
we infer that $u^{\mathrm{rem}}$ satisfies 
$$
 (\mathscr{T}+\mathcal{M})u^{\mathrm{rem}}= -\eps V_{t,s}^{\mathrm{rem}}u^0(z), \quad u^{\mathrm{rem}}(s,s,\cdot)=0, 
$$
that we can solve using Lemma~\ref{lem-semigroup} as 
$$
u^{\mathrm{rem}} (t,s,z)=-\eps \int_s^t U_{t,\tau} V_{\tau,s}^{\mathrm{rem}}u^0(z) ds.
$$
We can thus define the operator $U^{\mathrm{rem}}_{t,s}$ by the formula~\eqref{eq:lem-defU}
and by construction, the equality \eqref{eq:lem-appr} holds, while  the bound~\eqref{MajUr} directly follows from \eqref{eq:MajFIO}.
\end{proof}

In the following subsections, the goal  in summary will be
\begin{itemize}
\item to construct a  Fourier Integral Operator $U^{\mathrm{FIO}}_{t,s}$ such that \eqref{FIOeq} holds;
\item  to show that the properties of the phase and of the amplitude  ensure \eqref{eq:MajFIO};
\item to derive sharp  properties of the phase  which will allow to prove a quantum averaging lemma in the next section.
\end{itemize}

Note that the construction of a parametrix for an operator such as \eqref{eq:higher-wigner-pseudo} will follow  fairly standard steps.
Nevertheless, compared to the general theory, see for example  \cite{Zwo}, here we want to construct a parametrix which is valid globally on the phase space  (see also \cite{Cha,HR} in the elliptic case), and above all, to obtain precise continuity estimates in the weighted space $\mathcal{H}^0_{r,0}$ .
This will be possible thanks to the specific form of the  symbol  of the extended Wigner system.
 Note that we also want to perform this analysis in finite regularity and to quantify the required regularity for $\rho$ though we shall not try to optimize it.

\subsection{Eikonal equation,  transport equation and properties of the phase}\label{HJSec}

 First recall that $\mathrm{V}_{\rho}= V_{\eps}*_{x} \rho$ and $\widehat{V} \in \mathscr{C}^\infty_b(\R^d)$,  so that by definition of $T_\eps$, it holds
\begin{equation}\label{InegPot}
\sup_{\eps\in (0,1]}\norm{\mathrm{V}_{\rho}}_{L^2(0,T_\eps;H^m_{r})}\leq \Lambda(T,M).
\end{equation}
As expected  (see again \cite{Zwo}), in order to construct an appropriate FIO parametrix associated with $\mathscr{T}$,  the phase has to solve
  the following eikonal equation, which is an Hamilton-Jacobi equation:
\begin{equation}\label{H-J}
 \left\{
    \begin{array}{ll}
        \partial_t{\varphi}_{t,s}+ v \cdot \nabla_{x} \varphi_{t,s} + a_{\rho}(t,x, \nabla_v \varphi_{t,s}) =0, \quad z=(x,v),\xi \in  \R^{2d}, \\
        \varphi_{s,s}(z,\xi)=z\cdot \xi,
    \end{array}
\right.
\end{equation}
 where $a_{\rho}(t,x,\xi_v)=\mathrm{V}_\rho \left(t,x-\frac{\xi_v}{2}\right)-\mathrm{V}_\rho  \left(t,x+\frac{\xi_v}{2}\right)$.
 We first  gather, in the following proposition, the existence, uniqueness and  regularity properties for~\eqref{H-J}.

\begin{prop}
\label{prop-HJ1}
Let $p\geq 2$ be an integer such than   $m\geq \lfloor d/2\rfloor + p+2$. 
There exists a positive time $T(M) >0 $ such that for all $s \in [0,\min(T(M),T_\eps)]$, there is a unique solution $\varphi_{t,s} \in \mathscr{C}^2([0,\min(T(M),T_\eps)]^2 \times \mathbb{R}^{2d} \times \mathbb{R}^{2d})$ to  \eqref{H-J}. Morevover, $\varphi_{t,s}$ satisfies for all  $z, \xi \in \R^{2d}$ and all $0\leq t,s \leq \min(T_\eps,T(M))$ the  estimates
\begin{align}
\label{estphi}&\underset{{\substack{\vert \alpha\vert + \vert \beta\vert \leq p
 }}}{\sup} \left|\partial^{\alpha}_z\partial^{\beta}_{\xi}  \Big[ \varphi_{t,s}(z,\xi)- (x-(t-s)v) \cdot  \xi_x - v \cdot  \xi_v\Big]\right|\leq  1, \\
\label{estdzphi}&\underset{{\substack{\vert \alpha\vert  \leq p 
 }}}{\sup} \left|\partial^{\alpha}_z  \Big[\varphi_{t,s} (z,\xi)- (x-(t-s)v) \cdot  \xi_x - v \cdot  \xi_v\Big]\right|\leq  \vert  \xi_v\vert +\frac{1}{2} |t-s| \vert  \xi_x\vert,
\end{align}
and
\begin{equation}
\label{esthessphi}
\norm{ (\partial_z\partial_{\xi} \varphi_{t,s} -\mathrm{I} ) }_{L^{\infty}_{z,\xi}}\leq \frac{1}{2}.
\end{equation}
\end{prop}

Note that we obtain  a local existence result for the Hamilton-Jacobi equation \eqref{H-J} with estimates which are uniform with respect to $(z, \xi) \in \mathbb{R}^{2d} \times \mathbb{R}^{2d}$; this is due to the specific structure of the Hamiltonian. 

The estimates \eqref{estphi} and \eqref{esthessphi} will be crucial to ensure the boundedness of the Fourier integral operator $U^{\mathrm{FIO}}$ on $\mathcal{H}^0_{r}$
by using  continuity results proved in  Appendix~\ref{SFIO}.
The estimates \eqref{estphi}--\eqref{estdzphi} will be instrumental in the proof of the quantum averaging Lemma.

\medskip

Finally, the  amplitude $\mathrm{B}_{t,s}(z,\xi)$ will solve the  following first order linear equation:
\begin{equation}\label{Eqb}
 \left\{
    \begin{array}{ll}
        \partial_t \mathrm{B}_{t,s}+ v\cdot \nabla_x \mathrm{B}_{t,s}+\nabla_{\xi_v}a_{\rho} \left(t,x,\nabla_v\varphi_{t,s} \right)\cdot\nabla_v \mathrm{B}_{t,s}
        + \mathcal{N}_{t,s}\mathrm{B}_{t,s}=0,\\
        	\mathrm{B}_{s,s}(z,\xi)=\mathrm{I},
    \end{array}
\right.
\end{equation}
where $\mathcal{N}_{t,s} =\mathcal{N}_{1,t,s}   + \mathcal{N}_{2,t,s} $ with
\begin{align*}
 \mathcal{N}_{1,t,s} &:=  { 1 \over 2} \nabla_{v} \cdot \left[\left( \nabla_{\xi_v} a_\rho \right) \left(t,x,\nabla_v\varphi_{t,s}  \right) \right] ,\\
\mathcal{N}_{2, t,s}  &:= \mathfrak{m}_\rho(t,x, \nabla_v\varphi_{t,s}),
\end{align*}
where we recall the matrix $\mathfrak{m}_\rho$ is defined in \eqref{def:M}.
The existence, uniqueness and  regularity properties for~\eqref{Eqb} are gathered in the following proposition.
\begin{prop}\label{EstB}
Let $p\geq 2$ be an integer such than   $m\geq \lfloor d/2\rfloor + p+3$. 
Let  $T(M)>0$ be given by Proposition~\ref{prop-HJ1}. For all $s \in [0,\min(T(M),T_\eps)]$, there exists a unique solution $\mathrm{B}_{t,s}\in \mathscr{C}^1([0,\min(T(M),T_\eps)]^2 \times \R^{2d} \times \R^{2d})$ to~\eqref{Eqb}. Moreover    $B$ satisfies the following estimates:
\begin{align}\label{regStrongB}
\sup_{0\leq t,s\leq T} \sup_{\vert \alpha\vert + \vert \beta\vert \leq p} \norm{\partial^{\alpha}_z\partial^{\beta}_{\xi} \mathrm{B}_{t,s}}_{L^{\infty}_{z ,\xi}}&\leq \Lambda(T,M), \qquad T \in [0,\min(T(M),T_\eps)] \\
\label{SmallTimeB}
 \sup_{\vert \alpha\vert+\vert \beta\vert \leq p-1}   \norm{ \partial^{\alpha}_z\partial^{\beta}_{\xi}\left(\mathrm{B}_{t,s}-\mathrm{I}\right) }_{L^{\infty}_{z,\xi}} &\leq |t-s|\Lambda(T,M), \qquad t,s \in [0,\min(T(M),T_\eps)].
\end{align}

\end{prop}

\medskip

The proof of Propositions~\ref{prop-HJ1} and \eqref{EstB} are postponed to Subsections~\ref{HJSecProof1}--\ref{HJSecProof2}--\ref{sec-proofEstB}.
We point out that we shall obtain in Lemma~\ref{lemestimtildephi} a sharp version of the estimates of Proposition \ref {prop-HJ1}, that will be important in the final stage of the proof.

%

\subsection{Construction of the parametrix}
Thanks to Propositions  \ref{prop-HJ1} and \ref{EstB}, we can build the required FIO.

\begin{prop}\label{ConsFIO}
Let $T(M)>0$ be given by Proposition~\ref{prop-HJ1}.
Let $\varphi_{t,s}(z,\xi)$ be given by  Proposition  \ref{prop-HJ1} and $\mathrm{B}_{t,s}$ be given by Proposition \eqref{EstB}.  
Then, the (matrix-valued) Fourier Integral Operator $U^{\mathrm{FIO}}_{t,s}$  defined by  
\begin{equation}
\label{def-FIO-prop}
U_{t,s}^{\mathrm{FIO}}u=\frac{1}{(2\pi)^{2d}}\int_{\xi}   \int_{y} e^{\frac{i}{\varepsilon} \left( \varphi^{\eps} _{t,s}(z,\xi)-\langle y,\eps\xi \rangle \right)}\mathrm{B}^{\eps}_{t,s}(z,\xi)  u(y) dyd\xi,
\end{equation}
satisfies  for all $s,t \in [0,\min(T(M),T_\eps)]$ 
 the equation \begin{equation*}
 \left\{
    \begin{array}{ll}
        (\mathscr{T}+\mathcal{M})U_{t,s}^{\mathrm{FIO}}=\eps V_{t,s}^{\mathrm{rem}},\\
        	U_{s,s}^{\mathrm{FIO}} = \mathrm{I},
 \end{array}
\right.
\end{equation*}
where $V^{\mathrm{rem}}_{t,s}$ is an operator that satisfies the bound 
\begin{equation}
\label{MajVr}
\sup_{0\leq t,s\leq \min(T(M),T_\eps)} \norm{V_{t,s}^{\mathrm{rem}}}_{\mathscr{L}(\mathcal{H}^0_{r,0})} \leq  C,
\end{equation}
for $C>0$ independent of $\eps$.
\end{prop}

\begin{proof}[Proof of Proposition~\ref{ConsFIO}]
Let $T \in [0, \min(T(M),T_\eps)]$.
Let   $U^{\mathrm{FIO}}$  be a FIO under the form~\eqref{def-FIO-prop}. 
It will be convenient to use a more precise notation:  for a phase $\varphi$ and an amplitude $A$, we denote by 
$ I_{\varphi}[A]$ the semiclassical FIO defined by
$$   I_{\varphi}[A] u(z)= { 1 \over (2 \pi)^{d}} \int_{\xi} e^{\frac{i}{\eps} \varphi^\eps (z, \xi)} A^\eps(z, \xi) \, \widehat{u}(\xi) \, d\xi$$
so that 
$$ U^{\mathrm{FIO}}_{t,s} u= I_{\varphi_{t,s}}[\mathrm{B}_{t,s}] u.$$
Let us study the action of $\mathscr{T}+ \mathcal{M}$ on $U^{\mathrm{FIO}}$. 
By using \eqref{eq:higher-wigner-pseudo}, we get that
\begin{multline}
\label{expandparameps}
( \mathscr{T} + \mathcal{M}) U^{\mathrm{FIO}}_{t,s} u= {i \over \eps} I_{\varphi_{t,s}^\eps}\left[( \partial_{t} \varphi_{t,s} + v \cdot \nabla_{x} \varphi_{t,s}) \mathrm{B}_{t,s}
 + A_{\eps} \right]u \\+ I_{\varphi_{t,s}}\left[ \partial_{t} \mathrm{B}_{t,s}^\eps + v \cdot \nabla_{x} \mathrm{B}_{t,s}^\eps + M_{\eps} \right]u.
\end{multline}
where
\begin{align}
\label{defAeps}
A_{\eps}(t,s, z, \xi)&= e^{\frac{-i}{\varepsilon} \varphi_{t,s} (z,\xi)} a_{\rho}^\eps(t, x, D_{v})\left(e^{\frac{i}{\varepsilon}  \varphi_{t,s}  (z,\xi)}\mathrm{B}_{t,s}(z,\xi)\right), \\
\label{defMeps}
 M_{\eps}(t,s, z, \xi) &=  e^{\frac{-i}{\varepsilon} \varphi_{t,s}(z,\xi)} \mathfrak{m}_\rho^\eps(t, x, D_{v})\left(e^{\frac{i}{\varepsilon}  \varphi_{t,s}  (z,\xi)}\mathrm{B}_{t,s}(z,\xi)\right).
\end{align}
We shall next look for an expansion of $A_{\eps}$ and $M_{\eps}$ under the form
\begin{align}
\label{expandAeps}
{ 1 \over \eps} A_{\eps}(t,s,z, \xi ) &=  { 1 \over \eps} A_{-1}(t,s,z, \xi) + A_{0}(t,s,z, \xi)+ \eps A_{\mathrm{rem}}(t,s,z,\xi),  \\
\label{expandMeps}  M_{\eps}(t,s,z,\xi) &=  M_{0}(t,s,z,\xi)+ \eps M_{\mathrm{rem}}(t,s,z,\xi).
\end{align}
Note that $A_{-1}, A_0, A_{\mathrm{rem}}, M_0, M_{\mathrm{rem}}$ may all depend on $\eps$ but we shall not write explicitly this dependence for the sake of readability.

\medskip

\noindent $\bullet${\bf Expansions of $A_{\eps}$ and $M_{\eps}$.}
We have 
\begin{align}
\nonumber  A_{\eps}(t,s,z, \xi)&=\frac{1}{(2\pi )^{d}}\int_{\eta_v}\int_{w}e^{i \left(v-w\right)\cdot\eta_v}  e^{\frac{-i}{\varepsilon}\left(  \varphi _{t,s} (z,\xi)-\varphi_{t,s}  (x,w,\xi)\right)}a_{\rho}^\eps(t, x,\eta_v) \mathrm{B}_{t,s}(x,w,\xi) \, dw d\eta_v, \\
\label{defMeps2} M_{\eps}(t,s,z,\xi)&= \frac{1}{(2\pi )^{d}}\int_{\eta_v}\int_{w}e^{i \left(v-w\right)\cdot\eta_v}  e^{\frac{-i}{\varepsilon}\left(  \varphi _{t,s} (z,\xi)-\varphi_{t,s}  (x,w,\xi)\right)}\mathfrak{m}_\rho^\eps(t, x,\eta_v) \mathrm{B}_{t,s}(x,w,\xi) \, dw d\eta_v.
\end{align}
By a Taylor expansion with respect to the middle point $(v+w)/2$, we can write that
\begin{multline*}
\varphi _{t,s} (x,v,\xi)-\varphi _{t,s} (x,w,\xi) \\=\nabla_v\varphi _{t,s} (x, {v+ w \over 2},\xi)\cdot(v-w)
+ R^0_{t,s}(z,w,\xi) [v-w, v-w]\cdot (v-w),
\end{multline*}
where 
\begin{equation}
\label{def-R0}
 R^0_{t,s}(z,w,\xi)={1 \over 8}\int_{-1}^1 \int_{0}^1 \sigma_{1}^2 (1- \sigma_{2}) D_{v}^3 \varphi_{t,s}(x, { v+w \over 2}  + \sigma_{1} \sigma_{2} {v-w \over 2}) d\sigma_{1} d\sigma_{2},
\end{equation}
and we have denoted $R^0_{t,s}(z,w,\xi) [v-w, v-w]\cdot (v-w)= R^0_{t,s}(z,w,\xi) [v-w, v-w,v-w]$.
Let us first study the expansion of $A_{\eps}$.
By using  the change of variable 
$$\eta_v' :=\eta_v-\frac{1}{\eps}\left(\nabla_v\varphi _{t,s} (z,\xi)+R^0_{t,s}(z,w,\xi)[v-w, v-w]\right), $$ 
we obtain 
\begin{multline*}
A_{\eps}=\frac{1}{(2\pi )^{d}}\int_{\eta_v}\int_{w}e^{i \left(v-w\right)\cdot\eta_v}\\a_{\mathcal{\rho}}\left(x,\nabla_v\varphi _{t,s} (x, {v+ w \over 2 },\xi)  +\eps\eta_v+R^0_{t,s}(z,w,\xi)[v-w, v-w]\right)  \mathrm{B} _{t,s}(x,w,\xi) \, dw d\eta_v .
\end{multline*}
We can then  use  again a Taylor expansion  to write 
\begin{align*}
&a_{\mathcal{\rho}}\left(x,\nabla_v\varphi _{t,s} (x, {v+w \over 2},\xi) +  \eps \eta_v +R^0_{t,s}(z,w,\xi)[v-w, v-w]\right)\\ 
&\qquad\qquad =a_{\rho}\left(x,\nabla_v\varphi  _{t,s} ({v+w \over 2} ,\xi)\right)+\eps\eta_v \cdot \nabla_{\xi_{v}} a_{\mathcal{\rho}}\left(x,\nabla_v\varphi  _{t,s} ({v+w \over 2} ,\xi)\right) \\
&\qquad\qquad +R^1_{t,s}(z,w,\xi, \eps \eta_v), 
\end{align*}
where, recalling that $R^0$ is defined in~\eqref{def-R0}, 
\begin{equation}
\label{def:R1}
\begin{aligned}
R^1_{t,s}(z,w,\xi,&\eta_v) \\
&= \int_{0}^1 \nabla_{\xi_{v}}a_{\rho}\left(x,\nabla_v\varphi _{t,s} (x, {v+w \over 2},\xi) +   \eta_v + \sigma \{R_0\}_{t,s}(z,w,\xi)[v-w, v-w]\right) \, d \sigma \\
\qquad & \mbox{\hspace{2cm}}  \cdot  R^0_{t,s}(z,w,\xi)[v-w, v-w] \\
&+ \int_{0}^1 (1- \sigma) D_{\xi_{v}}^2 a_\rho\left(x,\nabla_v\varphi _{t,s} (x, {v+w \over 2},\xi) +  \sigma  \eta_v \right) [\eta_{v}, \eta_{v}] \, d\sigma.
\end{aligned}
\end{equation}
The key point in this expression is that $R^1_{t,s}(z, w, \xi, \eta_{v})$ can be seen as a bilinear  form, either in $\eta$ or  $v-w$, with bounded coefficients: more precisely, it  can be put  under the form 
\begin{equation}\label{expanR} 
R^{1}_{t,s}(z,w,\xi,\eta_v)=
\sum_{|\alpha |= 2}  c_{\alpha,t,s}(z,w,\xi,\eta_v)  \eta^\alpha + d_{\alpha,t,s}(z,w,\xi,\eta_v) (v-w)^\alpha,
\end{equation}
in which the coefficients $ c_{\alpha,t,s}$, $ d_{\alpha,t,s}$ satisfy the estimate 
\begin{align}\label{estCttDer}  
\sup_{t,s \in [0,T]} | \partial_{z, w, \xi, \eta_{v} }^\gamma  c_{\alpha,t,s}(z,w,\xi, \eta_{v}) |
 &\leq  \sup_{[0, T] } \Lambda (\|\rho \|_{W^{|\gamma|+2, \infty}},  \| \nabla^2 \varphi \|_{W^{|\gamma|, \infty}} ),   \\ 
 \label{estDttDer}
  \sup_{t,s \in [0,T]}  | \partial_{z, w, \xi }^\gamma \partial_{\eta_{v} }^{\gamma'}  d_{\alpha,t,s}(z,w,\xi, \eta_{v})|
 &\leq \sup_{[0, T]} \Lambda (\|\rho \|_{W^{|\gamma|+|\gamma'|+1, \infty}},  \| \nabla^2 \varphi \|_{W^{|\gamma|+1, \infty}}) \langle v- w \rangle^{ 2 | \gamma|}. 
\end{align}
%
Going back to $A_{\eps}$, we write 
$$
A_{\eps}=A_{-1}+\eps A_{0} + \eps^2 A_{\mathrm{rem}},
$$
with
\begin{align}
\label{def-I21}
A_{-1}(t,s,z,\xi)&:=a_{\rho}\left(x,\nabla_v\varphi _{t,s} (z,\xi)\right)\mathrm{B}_{t,s}(z,\xi), \\
\label{definutile} A_{0}(t,s,z,\xi)&:= \\
 \nonumber   \frac{1}{\eps}\frac{1}{(2\pi)^{d}}&\int_{\eta_v}\int_{w}e^{i \left(v-w\right)\cdot\eta_v}\eps\eta_v \cdot\nabla_{\xi_v}a_{\rho}\left(x,\nabla_v\varphi  _{t,s} ({v+ w \over 2 },\xi)\right) \mathrm{B} _{t,s}(x,w,\xi) \, dw d\eta_v,  \\
\label{def-I24}
 A_{\mathrm{rem}}(t,s,z,\xi)&:=\frac{1}{\eps^2}\frac{1}{(2\pi )^{d}}\int_{\eta_v}\int_{w}e^{i \left(v-w\right)\cdot\eta_v} R^{1}_{t,s}(z,w,\xi, \eps \eta_v)\mathrm{B} _{t,s}(x,w,\xi) \, dw d\eta_v.
\end{align}
We can further simplify the expression of   $A_{0}$ in \eqref{definutile} by resorting to integrations by parts in $w$ (one may also directly recognize  the   Weyl quantization  in the variables $(v, \eta_{v})$ of the symbol 
$ \eta_{v} \cdot \nabla_{\xi_{v}} a_\rho(x, \nabla_{v} \varphi_{t,s}(x, v, \xi)),$
$\xi$ and $x$ being parameters, acting on $B_{t,s}$  seen as a function of $v$).
This yields
\begin{equation}
\label{def-I22}
A_{0}(t,s,z, \xi)=\frac{1}{i}\nabla_{\xi_v}a_{\rho}\left(x,\nabla_v\varphi_{t,s} (z,\xi)\right)\cdot \nabla_v\mathrm{B}_{t,s}(z,\xi) +\frac{1}{2 i}  \nabla_{v} \cdot \left(  \nabla_{\xi_v}a_{\rho}\left(x,\nabla_v\varphi_{t,s}  (z,\xi)\right) \right)\mathrm{B}_{t,s}(z,\xi).
\end{equation}

\medskip

Similarly, we obtain an  expansion in powers of $\eps$ of $M_{\eps}$ defined in \eqref{defMeps}, by using \eqref{defMeps2}. 
 This is slightly easier since we only need to expand
 at first order. For example, for the phase, we can write
 \begin{align*}
\varphi_{t,s}  (x,v,\xi)-\varphi_{t,s}  (x, w  ,\xi)&= \nabla_v\varphi_{t,s}  (x,w,\xi)\cdot  (v-w )+R^2_{t,s}(z,\xi, \eps \eta)(v-w)\cdot (v-w)
\end{align*}
with 
\begin{equation}\label{ExpR}
R^2_{t,s}(z,w, \xi)= \int_0^1 D_{v}^2 \varphi _{t,s} (x,w + \sigma (v-w) ,\xi) \,  \sigma   d\sigma,
\end{equation}
and  we have denoted $R^2_{t,s}(z,\xi, \eps \eta)(v-w)\cdot (v-w)=R^2_{t,s}(z,\xi, \eps \eta)[v-w,v-w]$.
This yields 
$$
M_{\eps}(t,s,z,\xi)= M_{0}(t,s,z,\xi) + \eps M_{\mathrm{rem}}(t,s,z,\xi),
$$
with
\begin{align}
\label{def-I21bis}
M_{0}(t,s,z,\xi)&:=\mathfrak{m}_\rho \left(t,x,\nabla_v\varphi _{t,s}  (z,\xi)\right)\mathrm{B} _{t,s}(z,\xi), \\
\label{def-I21bisbis} M_{\mathrm{rem}}(t,s,z, \xi)&:=\frac{1}{\eps}\frac{1}{(2\pi)^{d}}\int_{\eta_v}\int_{w}e^{i \left(v-w\right) \cdot \eta_{v}} R^{3}_{t,s} (z, w, \xi, \eps \eta_{v})\mathrm{B} _{t,s}(x,w,\xi) \, dw d\eta_v,
\end{align}
where
\begin{multline*} R^{3}_{t,s}(z, w, \xi, \eta_{v})= \\ \int_{0}^1
 D_{\xi_{v}} \mathfrak{m}_\rho \left(t,x, \nabla_{v} \varphi_{t,s}(z,\xi) + \sigma( \eta_{v} + R^{2}_{t,s}(z, w, \xi) (v-w))\right) \cdot(\eta_{v}
  + R^{2}_{t,s}(z, w, \xi) (v-w))\,  d\sigma. 
  \end{multline*}
This time, we can expand  $ R^{3}_{t,s}$
  as  a linear form in $\eta$ and $v-w$:
 \begin{equation}\label{expanR3} 
R^{3}_{t,s}(z,w,\xi,\eta_v)=
\sum_{|\alpha |= 1}  c_{\mathfrak{m}, \alpha,t,s}(z,w,\xi,\eta_v)  \eta^\alpha + d_{\mathfrak{m}, \alpha,t,s}(z,w,\xi,\eta_v) (v-w)^\alpha,
\end{equation}
in which the coefficients $ c_{\mathfrak{m}, \alpha,t,s}$, $ d_{\mathfrak{m}, \alpha,t,s}$ satisfy the estimate
\begin{multline*}
\sup_{t,s \in [0,T]} | \partial_{z, w, \xi}^\gamma \partial_{\eta_{v} }^{\gamma'} c_{\mathfrak{m},\alpha,t,s}(z,w, \eta, \xi) |
+ | \partial_{z, w, \xi}^\gamma \partial_{\eta_{v} }^{\gamma'}d_{\mathfrak{m},\alpha,t,s}(z,w, \eta, \xi)| \\
 \leq  \sup_{[0, T] } \Lambda (\|\rho \|_{W^{|\gamma|+ |\gamma'|+3, \infty}},  \| \nabla^2 \varphi \|_{W^{|\gamma|, \infty}} ) \langle v- w \rangle^{  | \gamma|}. 
\end{multline*}

\medskip 

\noindent $\bullet$ {\bf Expression of the remainder.} 
By 
choosing $\varphi$ as the solution to the eikonal equation~\eqref{H-J}, we obtain by  using \eqref{def-I21} that
$$
\left(\partial_{t}   \varphi _{t,s}   + v \cdot \nabla_{x} \varphi_{t,s}\right) B_{t,s}   +A_{-1}=0,
$$
which cancels the terms of order $-1$ in $\eps$ in \eqref{expandparameps}, while choosing $ \mathrm{B}$ as the solution to~\eqref{Eqb} precisely yields that
$$
\partial_t  \mathrm{B} _{t,s}+  v \cdot \nabla_{x}B_{t,s} + i A_{0} + M_{0}=0,
$$
by using \eqref{def-I22}--\eqref{def-I21bis}, 
which cancels the terms of order $0$ in $\eps$.
Consequently, we have obtained  that 
\begin{align*}
(\mathscr{T}+\mathcal{M})U_{t,s}^{\mathrm{FIO}}=  \eps V_{t,s}^{\mathrm{rem}}
\end{align*}
where $ V_{t,s}^{\mathrm{rem}}$ is the semiclassical Fourier Integral Operator defined by 
$$  V_{t,s}^{\mathrm{rem}} = -  I_{\varphi}[ i A_{\mathrm{rem}} + M_{\mathrm{rem}}],$$
that is to say
\begin{equation}
\label{expressionVr}  V_{t,s}^{\mathrm{rem}} u(z) = { 1 \over (2 \pi)^{d}} \int_{\mathbb{R}^{2d}} e^{\frac{i}{\eps} \varphi^\eps_{t,s}(z, \xi)}\left( i A_{\mathrm{rem}}^\eps(t,s,z, \xi) + M_{\mathrm{rem}}^\eps (t,s, z,  \xi) \right) \widehat{u}(\xi) \, d \xi,
\end{equation}
where $A_{\mathrm{rem}}$ and $M_{\mathrm{rem}}$ are defined by \eqref{def-I24}, \eqref{expanR} and \eqref{def-I21bisbis}, \eqref{expanR3}, respectively.

%
%

\medskip

\noindent $\bullet${\bf Study of the remainder operator ${V}^{\mathrm{rem}}_{t,s}$.} 
To conclude the proof,  we need to prove  that ${V}^{\mathrm{rem}}$ is acting as  a bounded operator on $\mathcal{H}^0_{r, 0}$. Appendix \ref{SFIO} contains continuity results for FIO that are tailored for the present problem. Specifically, we shall apply Proposition~\ref{MajFIO}. 
Note that the  required estimates for the phase, namely \eqref{Hyppsi} and \eqref{hypFIO-weighted1},  clearly follow from Proposition \ref{prop-HJ1} with $p=2d + 2r +1$.
It remains to prove that the amplitude $iA_{\mathrm{rem}} + M_{\mathrm{rem}} $ matches  the required estimates~\eqref{hypFIO-weighted1}.

\begin{lem}\label{Maj-I_{2,4}}
The following estimates hold for $A_{\mathrm{rem}}$ and $M_{\mathrm{rem}}$:
\begin{multline}
\label{estimAr}
\sup_{t,s \in [0,T]}\sup_{\vert \alpha\vert \leq p_0}  \norm{\langle \eps \nabla_x \rangle^r  \langle \eps \nabla_{\xi_v} \rangle^r  \partial^{\alpha}_{z,\xi}A_{\mathrm{rem}} (t,s,z,\xi)}_{L^{\infty}_{z,\xi} }\\ \leq 
\sup_{[0, T]} \Lambda\left( \| \rho \|_{W^{3p_{0}+ 4d+7, \infty}_{2r}},\| B_{t,s}\|_{W^{p_{0}+ d+1, \infty}_r }, \| \nabla^2 \varphi \| _{W^{p_{0}+ d+4, \infty}_r } 
  \right), 
\end{multline}
\begin{multline}
\label{estimMr}
\sup_{t,s\in [0,T]}\sup_{\vert \alpha\vert \leq p_0}  \norm{\langle \eps \nabla_x \rangle^r  \langle \eps \nabla_{\xi_v} \rangle^r  \partial^{\alpha}_{z,\xi}M_{\mathrm{rem}}(t,s,z,\xi)}_{L^{\infty}_{z,\xi}} \\ \leq \sup_{[0, T]} \Lambda \left( \| \rho \|_{W^{  2p_{0}+
3d+7, \infty}_{2r}},\|B_{t,s} \|_{W^{p_{0}+ d+1, \infty}_r}, \| \nabla^2 \varphi \|_{W^{p_{0}+d+2, \infty}_r }\right),
\end{multline}
     where we have denoted for $k \in \mathbb{N}$,  $\| \cdot \|_{W^{k,\infty}_r} := \|  \langle \eps \nabla_x \rangle^r  \langle \eps \nabla_{\xi_v} \rangle^r  \cdot \|_{W^{k,\infty}}$.
     
\end{lem}

\begin{proof}[Proof of Lemma~\ref{Maj-I_{2,4}}]
By using \eqref{def-I24} and \eqref{expanR}, we can write the  decomposition
$$i A_{\mathrm{rem}}(t,s,z, \xi) =   \mathfrak{I}_{{\mathrm{rem}}, 1} (t,s,z,\xi)+  \mathfrak{I}_{{\mathrm{rem}}, 2}(t, s, z, \xi)$$
 where 
\begin{align*}
 \mathfrak{I}_{{\mathrm{rem}}, 1} &=  { 1 \over (2 \pi)^d } \sum_{|\alpha | = 2} \int_{\eta_{v}} \int_{w} e^{i \left(v-w\right)\cdot\eta_v}  c_{\alpha, t, s}(z,w,\xi,  \eps \eta_{v}) \eta_{v}^\alpha  \mathrm{B} _{t,s}(z,w,\xi) \, dw d\eta_v, \\
 \mathfrak{I}_{{\mathrm{rem}}, 2} &= { 1 \over (2 \pi)^d } \sum_{|\alpha | = 2} \int_{\eta_{v}} \int_{w} e^{i \left(v-w\right)\cdot\eta_v}  d_{\alpha, t, s}(z,w,\xi,  \eps \eta_{v}) (v-w)^\alpha  \mathrm{B} _{t,s}(z,w,\xi) \, dw d\eta_v.
 \end{align*}
 By integrating by parts in the integrals, 
 we can rewrite
 \begin{align*}
 \mathfrak{I}_{{\mathrm{rem}}, 1}(t,s,z, \xi) & = {- 1 \over (2 \pi)^d } \sum_{|\alpha | = 2} \int_{\eta_{v}} \int_{w} e^{i \left(v-w\right)\cdot\eta_v}  \partial_{w}^\alpha (c_{\alpha, t, s}  \mathrm{B} _{t,s})(z,w,\xi, \eps  \eta_{v}) \, dw d\eta_v,\\
  \mathfrak{I}_{{\mathrm{rem}}, 2}(t,s,z,\xi) & ={- 1 \over (2 \pi)^d } \sum_{|\alpha | = 2}\int_{\eta_{v}} \int_{w} e^{i \left(v-w\right)\cdot\eta_v}  \partial_{\eta_{v}}^\alpha \left[ (d_{\alpha, t, s}  \mathrm{B} _{t,s})(z,w,\xi, \eps  \eta_{v}\right] \, dw d\eta_v.
\end{align*}
 We shall focus on the  estimate of   $\mathfrak{I}_{{\mathrm{rem}}, 2}$, the estimate of $\mathfrak{I}_{{\mathrm{rem}}, 1}$ is slightly easier to obtain since
 the derivatives of $c_{\alpha}$ with respect to the $ (z,w,\xi)$ variables do not produce powers of $v-w$.
In order to take advantage of the oscillatory nature of the integrals,  we define   the operators
$$ \mathcal{L}_{w}= { 1 \over 1 + | \eta_{v}|^2 }( 1   + i  \eta_{v} \cdot  \nabla_{w})$$
 and
 $$  \mathcal{L}_{\eta_{v}}= { 1 \over 1 + |v-w|^2 }( 1   - i  (v-w) \cdot \nabla_{\eta_{v}})$$
which are such that
 $$  \mathcal{L}_{w}e^{i \left(v-w\right)\cdot\eta_v} = e^{i \left(v-w\right)\cdot\eta_v} , \quad  \mathcal{L}_{\eta_v}e^{i \left(v-w\right)\cdot\eta_v} = e^{i \left(v-w\right)\cdot\eta_v}.$$
 We thus get the identity 
 $$  \mathfrak{I}_{{\mathrm{rem}}, 2}(t,s,z,\xi)  ={ -1 \over (2 \pi)^d } \sum_{|\alpha | = 2}  \int_{\eta_{v}} \int_{w} e^{i \left(v-w\right)\cdot\eta_v}
  ( \mathcal{L}_{\eta_{v}}^T)^{N_{\eta_{v}}}
  ( \mathcal{L}_{w}^T)^{N_{w}}\partial_{\eta_{v}}^\alpha [(d_{\alpha, t, s}  \mathrm{B} _{t,s})(z,w,\xi, \eps  \eta_{v})] \, dw d\eta_v,$$
  where $ \mathcal{L}_{w}^T, \mathcal{L}_{\eta_{v}}^T$ stand for the formal $L^2$ transpose of $\mathcal{L}_{w}, \mathcal{L}_{\eta_v}$, and  for $N_{w}$ and $N_{\eta_{v}} $  integers to be chosen large enough.
  Note that $( \mathcal{L}_{w}^T)^{N_{w}}$ (resp.  $(\mathcal{L}_{\eta_{v}}^T)^{N_{\eta_{v}}}$) is a differential operator with  coefficients that
   decay like $1/\langle \eta_{v} \rangle^{N_{w}}$ (resp.    $1/\langle v-w \rangle^{N_{\eta_{v}}})$.
     By using the estimate \eqref{estDttDer}, we therefore get that
   \begin{multline*} | \mathfrak{I}_{{\mathrm{rem}}, 2}(t,s,z,\xi)|
    \leq  \int_{\eta_{v}}\int_{w} { 1 \over \langle \eta_{v} \rangle^{ N_{w}}} { 1 \over \langle v-w \rangle^{N_{\eta_{v}}- 2 N_{w}}} d\eta_{v} dw\, \\ 
    \times \Lambda \left( \|B \|_{W^{N_{w} , \infty} }, \| \nabla^2 \varphi \|_{W^{ N_{w}+ 3, \infty}}, \| \rho \|_{W^{ N_{\eta_{v}}+ N_{w}+ 3, \infty}}\right).
     \end{multline*}
More generally, for $| \alpha | \leq p_{0}$, we obtain
\begin{multline*} |  \langle \eps \nabla_x \rangle^r  \langle \eps \nabla_{\xi_v} \rangle^r \partial_{z, \xi}^\alpha \mathfrak{I}_{{\mathrm{rem}}, 2}|
    \leq  \int_{\eta_{v}}\int_{w} { 1 \over \langle \eta_{v} \rangle^{ N_{w}}} { 1 \over \langle v-w \rangle^{N_{\eta_{v}}- 2 N_{w}- 2 p_{0}}} d\eta_{v} dw\, \\ 
    \times \Lambda \left( \| B \|_{W^{N_{w}  +p_{0}, \infty}_r}, \| \nabla^2 \varphi \|_{W^{ N_{w}+ p_{0}+3, \infty}_r}, \| \rho \|_{W^{ N_{\eta_{v}}+ N_{w}+p_{0}+ 3, \infty}_{2r}}\right).
     \end{multline*}
 We may thus choose $N_{w}= d+1$, $N_{\eta_{v}}=  3 (d+1) + 2 p_{0}$ to get the  claimed estimate \eqref{estimAr}. 
 
 The estimates \eqref{estimMr} for $M_{\mathrm{rem}}$ can be deduced from  \eqref{def-I21bisbis}--\eqref{expanR3} by using similar arguments.
 This concludes the proof of Lemma \ref{Maj-I_{2,4}}.
\end{proof}

We are finally in position to end the proof of Proposition~\ref{ConsFIO}.  Taking $p_0= 2(1+d)$ in Lemma \ref{Maj-I_{2,4}},  using $m> 10d +  d/2+ 13 + r $ and $m\geq 3d + 6 + 2r$,  indeed shows that the condition \eqref{hypFIO-weighted1} in the assumptions of Proposition \ref{MajFIO} is satisfied. Hence, we can apply Proposition \ref{MajFIO} to infer that 
$$
\sup_{t,s\in [0,\min(T(M), T_\eps)]} \norm{V_{t,s}^{\mathrm{rem}}}_{\mathscr{L}(\mathcal{H}^0_{r, 0})}\lesssim 1,
$$
and the proof is complete.

\end{proof}

It remains to prove Proposition~\ref{prop-HJ1} and Proposition \ref{EstB}.

\subsection{Proof of  Proposition~\ref{prop-HJ1}, Part I: existence and uniqueness of a smooth solution}\label{HJSecProof1}

In this subsection, we show the existence and uniqueness of a smooth solution to the Hamilton-Jacobi equation~\eqref{H-J}.
Note that since we assume that $ m \geq 1 + \lfloor{d \over 2} \rfloor  + p+ 1$, we have  by  \eqref{InegPot}, Sobolev embedding 
 and Cauchy-Schwarz that 
 \begin{equation}
 \label{regularitehamilton} \|a_{\rho} \|_{L^1 (0, T; W^{p+1, \infty}_{x, \xi_{v}})} \leq T^{1 \over 2} \Lambda (T, M), \quad
   \|a_{\rho} \|_{L^\infty (0, T; W^{p, \infty}_{x, \xi_{v}})} \leq  \Lambda (T, M).
 \end{equation}
 Thanks to Lemma \ref{PE}, for $T <T_{\eps}$,  we know that $ f \in \mathscr{C}([0, T]; \H^m_{r})$ (though the estimate in this space
 depends on $\eps$) and also by using the equation \eqref{eq:wigner} that $f\in \mathscr{C}^1([0, T]; \H^{m-1}_{r-1})$. This yields
 by Sobolev embedding (using the notation $\mathscr{C}^{k}_b$ for $k-$differentiable bounded functions) that 
 \begin{equation}
 \label{regularho}  \nabla_{(x, \xi_{v})} a_{\rho} \in  \mathscr{C}^1([0, T];\mathscr{C}^1_{b}( \mathbb{R}^d \times \mathbb{R}^d)) \cap
  \mathscr{C}^0([0, T];  \mathscr{C}^2_{b}( \mathbb{R}^d \times \mathbb{R}^d)),
 \end{equation}
  assuming  $p \geq 2$.

The proof will be  based on the method of characteristics (see e.g. \cite{Zwo} for a closely related, more geometric approach). Here the properties of the Hamiltonian $a$ 
which is defined by
\begin{equation}
\label{defhamiltonian} a(t, z, \xi)=  v\cdot \xi_{x} + a_{\rho}(t, x, \xi_{v}),
\end{equation}
where we recall the notation~\eqref{eq:K1} for $a_\rho$, will allow to get a global in $z, \xi$ result. 
To motivate the use of the bicharacteristics, 
let us consider a curve, parametrized by time $t$, denoted by $(Z_{t,s}(z,\xi))_t$ in $\mathbb{R}^{2d}$ with
$Z_{s,s}(z,\xi)=z$. Then, given a solution $\varphi_{t,s}$ to the Hamilton-Jacobi equation $\partial_t \varphi_{t,s} + a(t,z,\nabla_z \varphi_{t,s})=0$ on some interval $[0,T]$,  let us  set  $ \Xi_{t,s}(z,\xi):= \nabla_z \varphi_{t,s}(Z_{t,s}(z,\xi), \xi)$.
Differentiating this relation with respect to time $t$, we thus  have
$$
\partial_t \Xi_{t,s} = (\partial_t  \nabla_z \varphi_{t,s})( Z_{t,s})  + \partial_t Z_{t,s} \cdot \nabla_z \Xi_{t,s}.
$$
On the other hand, differentiating the Hamilton-Jacobi equation with respect to $z$ and evaluating at the point $Z_{t,s}$, we obtain that
$$(\partial_t  \nabla_z  \varphi_{t,s})(Z_{t,s})  +  \nabla_{\xi} a(t,Z_{t,s},\Xi_{t,s})\cdot \nabla_z \Xi_{t,s} + \nabla_{z} a(t,Z_{t,s} ,\Xi_{t,s})= 0.$$ 
We therefore see that imposing $Z_{t,s}$ to solve
$$
\partial_t Z_{t,s}  =  \nabla_{\xi} a(t,Z_{t,s},\Xi_{t,s}),
$$
the vector field $\Xi_{t,s}$ must satisfy
$$
\partial_t \Xi_{t,s} = -  \nabla_{z} a(t,Z_{t,s} ,\Xi_{t,s}).
$$
Finally, as we require $\varphi_{s,s}(z,\xi) =z \cdot \xi$, we get that $ \Xi_{s,s}(z,\xi)=\xi$. 

\begin{Rem}{This argument shows that if a solution $\varphi_{t,s}$ exists on $[0,T]$ and is at least $\mathscr{C}^2$, then it must be unique.} Indeed, if we have two solutions $\varphi_1$ and  $\varphi_2$, we can associate with them the vector fields $(Z^1_{t,s},\Xi^1_{t,s})$ and  $(Z^2_{t,s},\Xi_{t,s}^2)$ and show that they satisfy the same (regular enough) differential equation with the same initial condition. Therefore,  they must be equal, leading to $\nabla_z(\varphi_1-\varphi_2)_{t,s}(Z_{t,s})=0$. Consequently,
\begin{equation*}
 \left\{
    \begin{array}{ll}
	\partial_t (\varphi_1-\varphi_2)_{t,s}(Z_{t,s}, \xi)=0,\\
	 (\varphi_1)_{s,s}(z, \xi)=(\varphi_2)_{s,s}(z, \xi)=z\cdot\xi.
    \end{array}
\right.
\end{equation*}
Finally, provided that $z \mapsto Z_{t,s}(z,\xi)$ is a  diffeomorphism (which will be proved  in the upcoming Lemma~\ref{invcar}),  we infer that $\varphi_1 \equiv \varphi_2$.
\end{Rem}
We are therefore naturally led to consider the bicharacteristics curves associated with the Hamiltonian $a$.
\begin{definition}
\label{def:bichar}
The bicharacteristics 
$$Z_{t,s}(z,\xi)=(Z_{t,s}^x(z,\xi),Z_{t,s}^v(z,\xi)), \qquad \Xi_{t,s}(z,\xi)=(\Xi^x_{t,s}(z,\xi),\Xi^v_{t,s}(z,\xi)),$$
are the curves in $\R^{2d}$ solving the system
\begin{equation}\label{carop}
 \left\{
    \begin{aligned}
        &\partial_t Z_{t,s}= \nabla_\xi a(t,Z_{t,s},\Xi_{t,s}), \qquad Z_{s,s}=z,\\
	&\partial_t \Xi_{t,s}= -\nabla_z a(t,Z_{t,s},\Xi_{t,s}), \qquad \Xi_{s,s}=\xi. 
    \end{aligned}
\right.
\end{equation}
\end{definition}
The bicharacteristics exist and  are uniquely defined on the interval of time $[0,T]$  thanks to \eqref{regularho} and  the Cauchy-Lipschitz Theorem. Indeed, the fact that they exist on the whole time interval comes from the structure of the vector field
in \eqref{carop}:  it is made of a linear part and a nonlinear bounded part. We also get from the Cauchy-Lipschitz Theorem with parameter that 
$(Z, \Xi) \in \mathscr{C}^2([0, T ]\times  \mathbb{R}^{2d}_{z} \times \mathbb{R}^{2d}_{\xi})$ (note that with this notation for regularity
we do not claim boundedness). 

To show  the existence of a solution to~\eqref{H-J}, we first introduce the following function $\psi_{t,s}(z,\xi)$:
\begin{equation}
\label{eq-varphialongX}
\begin{aligned}
\psi_{t,s}(z,\xi)&=z\cdot\xi +\int_s^t -a(\tau,Z_{\tau,s},\Xi_{\tau,s})+ \Xi_{\tau,s}\cdot \nabla_{\xi} a(\tau, Z_{\tau,s},\Xi_{\tau,s}) d\tau\\&=z\cdot \xi +\int_s^t -a_{\rho}(\tau,Z_{\tau,s},\Xi_{\tau,s})+ \Xi^v_{\tau,s}\cdot \nabla_{\xi_v} a_{\rho}(\tau, Z_{\tau,s},\Xi_{\tau,s}) d\tau.
\end{aligned}
\end{equation}
From the regularity of the bicharacteristics and $a_{\rho}$, we also get that
\begin{equation}
\label{psiregularity} \psi \in    \mathscr{C}^2([0, T ]^2 \times \mathbb{R}^{2d}_{z} \times \mathbb{R}^{2d}_{\xi}).
\end{equation} 
We then want to define a function $\varphi_{t,s}$ such that $\varphi_{t,s}(Z_{t,s}(z,\xi),\xi)=\psi_{t,s}(z,\xi)$. Before proving that such a function is indeed a solution of~\eqref{H-J}, we start by showing that we can inverse the space characteristics $z\mapsto Z_{t,s}(z,\xi)$. This is the purpose of the next lemma.

\begin{lem}\label{invcar}
There exists a positive time $T(M)>0$ such that the function $z\mapsto Z_{t,s}(z ,\xi)$ is a global  diffeomorphism for all $s,t\in[0,\min (T(M),T_\eps)]$ and all $\xi \in \R^{2d}$. 
\end{lem}

\begin{proof}[Proof of Lemma~\ref{invcar}]

Applying $\nabla_z$ to the bicharacteristics equations \eqref{carop} yields
\begin{equation} \label{car2}
 \left\{
    \begin{array}{ll}
      \nabla_z  Z_{t,s}(z,\xi)=\mathrm{I} + \int_s^t \nabla_z \left(\nabla_\xi a(\tau,Z_{\tau,s},\Xi_{\tau,s})\right) d\tau, \\
	\nabla_z \Xi_{t,s}(z,\xi)= -\int_s^t \nabla_z \left( \nabla_z a(\tau,Z_{\tau,s},\Xi_{\tau,s}) \right) d\tau.
    \end{array}
\right.
\end{equation}
For $T\in (0,T_\eps]$, we deduce from \eqref{regularitehamilton} that 
$$
 \vert \left(\nabla_z Z_{t,s}(z,\xi),  \nabla_z \Xi_{t,s}(z,\xi \right))\vert  \leq 1+  \Lambda(T, M)  \int_{s}^t \vert \left(  \nabla_z Z_{\tau s}(z,\xi), \nabla_z  \Xi_{\tau,s}(z,\xi) \right) \vert\, d\tau$$
 and hence, we get from the Gronwall Lemma that
\begin{equation}
\label{nablazcara} \sup_{t,s \in [0,T]} \vert \left(\nabla_z Z_{t,s}(z,\xi),  \nabla_z \Xi_{t,s}(z,\xi \right))\vert  \leq e^{T \Lambda(T,M)}.
\end{equation}
 Going back to \eqref{car2}, we then  deduce that for all $t,s \in [0,T]$,
\begin{align*}
\sup_{t,s \in [0,T]}  \vert \nabla_z Z_{t,s}(z,\xi) -\mathrm{I}\vert & \leq  T \Lambda(T, M).
\end{align*}
Therefore, we can find a time $T(M)>0$ small enough such that for all $s,t\in[0,\min (T(M),T_\eps)]$,
\begin{equation}\label{Diffeo}
\norm{\nabla_z Z_{t,s}(z,\xi)-\mathrm{I} }_{L^\infty_{z, \xi}}\leq \frac{1}{2}.
\end{equation}
As a result, for all $s,t \in [0,\min (T(M),T_\eps)]$ and all $\xi \in \R^{2d}$, the map $z \mapsto Z_{t,s}(z,\xi)$ is a  small $\mathscr{C}^1$  perturbation of the identity and hence a global diffeomorphism.
\end{proof}
 We  define $Y_{t,s}(z,\xi)$ as the inverse of $Z_{t,s}(z,\xi)$, i.e. $Y_{t,s}(z,\xi)$ is the vector field satisfying for all $z, \xi \in \mathbb{R}^{2d}$,
\begin{equation}\label{ConsY}
 Z_{t,s}( Y_{t,s}(z, \xi), \xi)= z.
\end{equation}
Note that  we  get from \eqref{Diffeo},  the regularity of $Z$ and  the Implicit Function Theorem that  $Y \in \mathscr{C}^2([0, T ]^2 \times  \mathbb{R}^{2d}_{z} \times \mathbb{R}^{2d}_{\xi})$ for  $T \leq \min (T(M),T_\eps).$
As a consequence we can properly define $\varphi$ from the formula~\eqref{eq-varphialongX} as:
\begin{equation}\label{psi}
{\varphi}_{t,s}(z,\xi)={\psi}_{t,s}(   Y_{t,s}(z, \xi), \xi)
\end{equation}

We are in position to show that $\varphi_{t,s}$ as defined in \eqref{psi} satisfies~\eqref{H-J}. 
On the one hand, by using the chain rule and  the definition of the bicharacteristics \eqref{carop}, 
we have
\begin{equation}
\label{unefacon}
\frac{d}{dt}\left(\varphi_{t,s}(Z_{t,s},\xi)\right)= \partial_t\varphi_{t,s}(Z_{t,s},\xi)+ \nabla_\xi a(t,Z_{t,s},\Xi_{t,s})\cdot\nabla_z \varphi_{t,s}(Z_{t,s},\xi),
\end{equation}
while on the other hand, by  differentiating ~\eqref{eq-varphialongX} with respect to time, we have
\begin{equation}
\label{lautrefacon}
 \frac{d}{dt}\left(\varphi_{t,s}(Z_{t,s},\xi)\right)= -a(t,Z_{t,s},\Xi_{t,s})+ \Xi_{t,s}\cdot \nabla_{\xi} a(t, Z_{t,s}, \Xi_{t,s}).
\end{equation}
To conclude, it  only remains to check that $\nabla_z \varphi_{t,s}(Z_{t,s},\xi) =\Xi_{t,s}$ for all $s,t \in [0,\min (T(M),T_\eps)]$.
By injecting this property into \eqref{unefacon}, \eqref{lautrefacon},  we shall obtain
\begin{equation}\label{finp}
 \partial_t\varphi_{t,s}(Z_{t,s},\xi)+a(t,Z_{t,s},\nabla_z\varphi_{t,s}(Z_{t,s},\xi))=0.
\end{equation}
Differentiating~\eqref{eq-varphialongX} with respect to $z$, we get
\begin{equation}
\label{lemmepsi2}
\begin{aligned}
\nabla_z\left(\varphi_{t,s}(Z_{t,s},\xi)\right)&= \xi +\int_s^t \Big( -  \nabla_z Z_{\tau,s}  \cdot \nabla_z a(\tau,Z_{\tau,s},\Xi_{\tau,s})  -\nabla_\xi a(\tau,Z_{\tau,s},\Xi_{\tau,s})\cdot \nabla_z \Xi_{\tau,s}  
\\
&+ \nabla_\xi a(\tau,Z_{\tau,s},\Xi_{\tau,s})\cdot   \nabla_z \Xi_{\tau,s}  + \Xi_{\tau,s} \cdot \partial_{\tau}  \nabla_z Z_{\tau,s}(z,\xi)   \Big) \, d\tau \\ 
&=\xi + \int_s^t \left( \partial_\tau {\Xi}_{\tau,s} \cdot \nabla_z{Z}_{\tau,s}  +{\Xi}_{\tau,s}\cdot    \partial_\tau  \nabla_z   {Z}_{\tau,s} \right) \, d\tau \\
&=  {\Xi}_{t,s}\cdot  \nabla_z Z_{t,s},
\end{aligned}
\end{equation}
by definition of the bicharacteristics~\eqref{carop}. 
We therefore infer that
$$
\partial_z Z_{t,s}    \, \left( {\Xi}_{t,s} - \nabla_z \varphi_{t,s}(Z_{t,s},\xi)\right) =0
$$
where here $\partial_{z}Z$ stands for  the jacobian matrix  with respect to the $z$ variable.
By using Lemma \ref{invcar}, this implies that  ${\Xi}_{t,s} =\nabla_z \varphi_{t,s}(Z_{t,s},\xi)$ for all $s,t \in [0,\min (T(M),T_\eps)]$ and all $z, \xi \in \R^{2d}$. 
This ends the proof of the first part of Proposition~\ref{prop-HJ1}: we have proven for every  $s \in [0, T]$ the existence of a unique classical
$ \mathscr{C}^2([0 , T]^2 \times \mathbb{R}^{2d}_z\times\mathbb{R}^{2d}_{\xi})$ solution of the Hamilton-Jacobi equation. 

Note that  we can easily deduce a  first  quantitative estimate for the derivatives of order two of the phase.
\begin{lem}
\label{derivee2phi}
For every $T \leq \min (T(M), T_{\eps})$, we have  the estimate
$$  \sup_{t,s \in [0, T]} \|  \nabla_{(z, \xi)}^2 \varphi_{t,s} \|_{L^\infty_{z, \xi}} \leq \Lambda(T, M).$$
\end{lem}

\begin{proof}
By using again \eqref{carop}, as in the previous Lemma, we have
\begin{equation*}
 \left\{
    \begin{array}{ll}
      \nabla_\xi Z_{t,s}(z,\xi)= \int_s^t \nabla_\xi \left(\nabla_\xi a(\tau,Z_{\tau,s},\Xi_{\tau,s})\right) d\tau, \\
	\nabla_\xi \Xi_{t,s}(z,\xi)= I-\int_s^t \nabla_\xi \left( \nabla_z a(\tau,Z_{\tau,s},\Xi_{\tau,s}) \right) d\tau.
    \end{array}
\right.
\end{equation*}
and hence, we obtain again from \eqref{regularitehamilton} and the Gronwall lemma that
\begin{equation}
\label{nablaxicara} \sup_{t,s \in [0,T]} \vert \left(\nabla_\xi Z_{t,s}(z,\xi),  \nabla_\xi \Xi_{t,s}(z,\xi \right))\vert  \leq e^{T \Lambda(T,M)}.
\end{equation}
  By using \eqref{ConsY} and  \eqref{Diffeo}, we deduce that we also have
  $$   \sup_{t,s \in [0,T]} \left\vert \nabla_{(z, \xi)} Y_{t,s}(z,\xi) \right\vert  \leq 
 \Lambda(T, M)$$
  if $T \leq \min(T(M), T_{\eps}).$
Since we have proven beforehand that
$$   \nabla_{z} \varphi_{t,s}(z, \xi)= \Xi_{t,s}( Y_{t, s}(z,\xi), \xi), $$
we then deduce that
\begin{equation}
\label{manquequedxixi}  \sup_{t,s \in [0, T]} \| \nabla_{(z, \xi)} \nabla_{z} \varphi_{t,s} \|_{L^\infty_{z, \xi}} \leq \Lambda(T, M).
\end{equation}
To get the estimate for $\nabla_{\xi}^2 \varphi_{t,s}$, we use directly that $\varphi_{t,s}$ satisfies the Hamilton-Jacobi equation \eqref{H-J}.
 For $| \alpha |= 2$, we have that
 $$ \partial_{t} \partial_{\xi}^\alpha \varphi_{t,s} + v \cdot \nabla_{x} \partial^\alpha_{\xi} \varphi_{t,s}
  + \nabla_{\xi_{v}} a_{\rho}(t, x, \nabla_{v}\varphi_{t,s}) \cdot \nabla_{v} \partial^\alpha_{\xi} \varphi_{t,s} =  R$$
  where  by using \eqref{manquequedxixi}, we have
  $$  \sup_{t \in [0, T]}\|R(t) \|_{L^\infty_{z, \xi}} \lesssim  \sup_{t \in [0, T]}  \|a_{\rho} \|_{W^{2, \infty}_{z, \xi}} \| \nabla_{v} \nabla_{\xi}\varphi_{t,s} \|_{L^\infty}^2
   \leq   \Lambda(T, M).$$
   By $L^\infty$ estimates for transport equations, we thus deduce 
   $$  \sup_{t,s \in [0, T]}  \|\nabla_{\xi}^2 \varphi_{t,s} \|_{L^\infty_{z, \xi}}  \leq T \Lambda(T, M),$$
   which concludes the proof of the lemma.

\end{proof}

\subsection{Proof of Proposition~\ref{prop-HJ1} Part II: estimates  of the phase}
\label{HJSecProof2}

We shall now prove the estimates \eqref{estphi} and  \eqref{estdzphi}. It is convenient to set
\begin{equation}
\label{phiperturb}
\varphi_{t,s}(z, \xi)= (x- (t-s) v) \xi_{x}+ v \cdot \xi_{v} + \widetilde{\varphi}_{t, s}(z, \xi).
\end{equation}
Note that we still have the regularity $ \widetilde \varphi_{t,s} \in   \mathscr{C}^2([0, T ] \times  \mathbb{R}^{2d}_{z} \times \mathbb{R}^{2d}_{\xi})$ and that $\widetilde \varphi_{t,s}$ solves the perturbed equation
\begin{equation}
\label{celledephitilde}
        \partial_t \widetilde {\varphi}_{t,s}+ v \cdot \nabla_{x}\widetilde  \varphi_{t,s} + a_{\rho}\left(t,x,  \xi_{v} - (t-s) \xi_{x} +   \nabla_v \widetilde  \varphi_{t,s}\right) =0, \quad 
       \widetilde \varphi_{s,s}(z,\xi)=0.
 \end{equation}
 We shall prove that:
 \begin{lem}
 \label{lemestimtildephi}
 For every $T \leq \min (T(M), T_{\eps})$, we have  the estimates
\begin{equation}
\label{tildephiunif}
  \| \widetilde   \varphi_{t,s} \|_{W^{p, \infty}_{z, \xi}} \leq T^{1 \over 2 } \Lambda(T, M),  \quad \forall t, \, s \in [0, T], 
  \end{equation}
  and
  \begin{equation}
\label{tildephidegreun}
   \| \widetilde   \varphi_{t,s}(\cdot, \xi)  \|_{W^{p, \infty}_{z}} \leq T^{1 \over 2 } \Lambda(T, M) (|\xi_{v}| + (t-s) |\xi_{x}|),\quad \forall  t, \, s \in [0, T],
    \, \forall \xi \in \mathbb{R}^d. 
  \end{equation}
 \end{lem}
 Once \eqref{tildephiunif} and \eqref{tildephidegreun} are established, \eqref{estphi}, \eqref{estdzphi} and \eqref{esthessphi}
 directly follow from the definition of $\widetilde \varphi$  by choosing $T(M)$ small enough.
 
 \begin{proof}[Proof of Lemma \ref{lemestimtildephi}]
 We first prove \eqref{tildephiunif}.
 Integrating \eqref{celledephitilde} along the  characteristics of free  transport, we first get that 
 \begin{equation}
 \label{pasnablatildephi} \| \widetilde \varphi_{t,s} \|_{L^\infty_{z,\xi}} \leq  \int_{s}^t  \|a_{\rho}(\tau)\|_{L^\infty} \, d\tau  \leq T \Lambda(T,M).
 \end{equation}
Taking the gradient in \eqref{celledephitilde}, and 
using $L^\infty$ estimates for the transport equation \eqref{eqdalphatildephi}, we then also get
$$  \| \nabla_{(z, \xi)} \widetilde \varphi_{t,s} \|_{L^\infty_{x,\xi}} \leq  \Lambda(T, M)  \int_{s}^t  \| \nabla_{(z, \xi)} \widetilde \varphi_{\tau,s}  \|_{L^\infty_{x,\xi}}\, d\tau  + T \Lambda(T, M)$$
and hence
\begin{equation}
\label{nablatildephi}  \| \nabla_{(z, \xi)} \widetilde \varphi_{t,s} \|_{L^\infty_{x,\xi}} \leq  T\Lambda(T, M)
\end{equation}
from the Gronwall inequality.  We then write  that for $| \alpha|  \geq  2$, 
\begin{equation}
\label{eqdalphatildephi} \partial_{t} \partial ^\alpha  \widetilde \varphi_{t,s} + v \cdot \nabla_{x} \partial^\alpha \widetilde  \varphi_{t,s}
  + \nabla_{\xi_{v}} a_{\rho}(t, x,  \xi_{v} - (t-s) \xi_{x} + \nabla_{v} \widetilde \varphi_{t,s}) \cdot \nabla_{v} \partial^\alpha \widetilde  \varphi_{t,s} =  R_{\alpha},
  \end{equation}
  where $R_{\alpha}$ is a commutator term.
   For $| \alpha | = 2$, we have 
  $$  \|R_{\alpha}(t) \|_{L^\infty_{z, \xi}} \lesssim ( \|\widetilde \varphi_{t,s} \|_{W^{2, \infty}_{z, \xi}} + \|\widetilde \varphi_{t,s} \|_{W^{2, \infty}_{z, \xi}}^2) ( 1   +  \|a_{\rho}(t) \|_{W^{2, \infty}_{z, \xi}})   
   +  \|a_{\rho}(t) \|_{W^{2, \infty}_{z, \xi}}.$$
    Note that  from \eqref{phiperturb} and  Lemma \ref{derivee2phi} we already have
 that
 \begin{equation}
 \label{lordre2phi}   \| \nabla_{(z, \xi)}^2 \widetilde \varphi_{t,s} \|_{L^\infty_{z,\xi}} \leq   \Lambda(T, M),
 \end{equation}
and consequently,  by also using  \eqref{nablatildephi} and \eqref{pasnablatildephi}, we obtain the estimate
  $$  \|R_{\alpha}(t) \|_{L^\infty_{z, \xi}} \leq \Lambda(T, M).$$
By $L^\infty$ estimates for the transport equation \eqref{eqdalphatildephi}, this yields 
  $$ \| \partial^\alpha \widetilde \varphi_{t,s} \|_{L^\infty_{z, \xi}} \leq T \Lambda (T, M), \quad \forall | \alpha | = 2.$$
  The estimates for $| \alpha|= k $, $ 3 \leq k \leq p$ then follow by induction. Indeed, for $k  \geq 3$, we can again write the equation \eqref{eqdalphatildephi}. Assuming that the wanted estimates hold for all $|\alpha|=k-1$, we have
  \begin{align*}
    \|R_{\alpha}(t) \|_{L^\infty_{z, \xi}} &\leq   \|\widetilde \varphi_{t,s} \|_{W^{k, \infty}_{z, \xi}}  ( 1   +  \|a_{\rho}(t) \|_{W^{2, \infty}_{z, \xi}})
   \Lambda(  \|\widetilde \varphi_{t,s} \|_{W^{k-1, \infty}_{z, \xi}})  +   \Lambda(  \|\widetilde \varphi_{t,s} \|_{W^{k-1, \infty}_{z, \xi}}) \|a_{\rho}(t) \|_{W^{k, \infty}_{z, \xi}} \\ 
   &\leq \Lambda(T, M)  \|\widetilde \varphi_{t,s} \|_{W^{k, \infty}_{z, \xi}} + \Lambda(T, M) \|a_{\rho}(t) \|_{W^{k, \infty}_{z, \xi}},
   \end{align*}
  where we have used  in  the last inequality  the estimates for $|\alpha|=k-1$. 
  It follows from $L^\infty$ estimates for the transport equation \eqref{eqdalphatildephi} and Gronwall's inequality that
  $$ \|\widetilde \varphi_{t,s} \|_{W^{k, \infty}_{z, \xi}}  \leq   T^{1 \over 2}\Lambda(T, M)  \|a_{\rho} \|_{L^2(0, T; W^{k+1, \infty}_{z, \xi})} 
  \leq T^{1 \over 2} \Lambda(T, M),$$
  thanks to \eqref{regularitehamilton}. This concludes the proof of \eqref{tildephiunif}.

We can now prove \eqref{tildephidegreun}.
For this estimate, we shall use more precisely the structure of $a_{\rho}$ in \eqref{celledephitilde}, which we recall is given by 
\begin{equation}
\label{arhostructure}
 a_{\rho}(t, x, \xi_{v})= \mathrm{V}_{\rho}(t, x- {\xi_{v} \over 2}) - \mathrm{V}_{\rho}(t, x+ {\xi_{v} \over 2}).
 \end{equation}
By a Taylor expansion, we get that 
\begin{multline*}    \partial_t \widetilde {\varphi}_{t,s}+ v \cdot \nabla_{x}\widetilde  \varphi_{t,s} + b_{\rho}\left(t,x,  \xi_{v} - (t-s) \xi_{x} +    \nabla_v \widetilde  \varphi_{t,s}\right) \cdot \nabla_{v} \widetilde \varphi_{t,s} \\ =-  b_{\rho}\left(t,x,  \xi_{v} - (t-s) \xi_{x} +  \nabla_v \widetilde  \varphi_{t,s}\right) \cdot (\xi_{v}-(t-s) \xi_{x}),
\end{multline*}
where 
\begin{equation}
\label{defbrhotildephi}
 b_{\rho}(t, x, \xi_{v})= - \int_{-{1 \over 2}}^{1 \over 2} \nabla_{x} V_{\rho}(t, x + \sigma \xi_{v})\, d\sigma.
 \end{equation}
Consequently, integrating along the characteristics of the vector field
$$  v \cdot \nabla_{x} +  b_{\rho}\left(t,x,  \xi_{v} - (t-s) \xi_{x} +    \nabla_v \widetilde  \varphi_{t,s}(z, \xi)\right) \cdot \nabla_{v}, $$
we get that for all $\xi \in \R^{2d}$,
$$ \|\widetilde \varphi_{t,s}( \cdot, \xi) \|_{L^\infty_{z}} \lesssim\int_{s}^t  \| \nabla_{x} \mathrm{V}_{\rho}(\tau)\|_{L^\infty}
 (|\xi_{v}| + |\tau-s| \xi_{x}) \, d\tau 
  \leq T \Lambda(T, M)  (|\xi_{v}| + |t-s| \xi_{x}).$$
 For higher order derivatives, we write for all $1\leq | \alpha | \leq p$  that
 \begin{multline}
 \label{eqtildephihighorder}    \partial_t \partial^\alpha_{z} \widetilde {\varphi}_{t,s}+ v \cdot \nabla_{x}\widetilde \partial^\alpha_{z}  \varphi_{t,s} +
 \nabla_{\xi_{v}} a_{\rho}\left(t, x,  \xi_{v} - (t-s) \xi_{x} +    \nabla_v \widetilde  \varphi_{t,s}\right)\cdot \nabla_{v} \partial^\alpha_z \widetilde \varphi
 \\ = R_{\alpha}  -  \partial_{z}^\alpha b_{\rho}\left(t,x,  \xi_{v} - (t-s) \xi_{x} +  \nabla_v \widetilde  \varphi_{t,s}\right) \cdot (\xi_{v}-(t-s) \xi_{x}),
 \end{multline} 
 where $R_{\alpha}$ is again a commutator term which can be estimated, for all $\xi \in \R^{2d}$,  by 
 $$ \|R_{\alpha}(t, \xi)\|_{L^\infty_{z}} \leq \Lambda\left( \| \nabla_{z} \widetilde  \varphi_{t,s}\|_{W^{p-1, \infty}_{z, \xi}} \right)
 (1+\| a_{\rho}(t)\|_{W^{p, \infty}_{z,\xi}}) \| \nabla_{z} \widetilde  \varphi_{t,s}(\cdot, \xi)\|_{W^{p-1, \infty}_{z}}.$$
 This implies by using \eqref{tildephiunif} which is already established and \eqref{regularitehamilton} that for all $\xi \in \R^{2d}$
 $$ \|R_{\alpha}(t, \xi)\|_{L^\infty_{z}} \leq \Lambda(T, M) \| \nabla_{z} \widetilde  \varphi_{t,s}(\cdot, \xi)\|_{W^{p-1, \infty}_{z}}.$$
  By integrating \eqref{eqtildephihighorder} along the characteristics of  $v \cdot \nabla_{x} + \nabla_{\xi_{v}} a_{\rho}(t,x, \nabla_v \varphi_{t,s}) \cdot \nabla_{v}$, 
  we obtain that for all $\xi \in \R^{2d}$,
  \begin{align*}
  \| \nabla_{z} \widetilde \varphi_{t,s} (\cdot, \xi) \|_{W^{p-1, \infty}_{z}} &\leq  \Lambda(T, M) \int_{s}^t   \| \nabla_{z} \widetilde  \varphi_{\tau,s}(\cdot, \xi)\|_{W^{p-1, \infty}_{z}}\, d\tau   + \| \rho   \|_{L^1(0,T;W^{p+1, \infty})}  (| \xi_{v}| + |t-s| | \xi_{x}|) \\
   &\leq \Lambda(T,M) \int_{s}^t   \| \nabla_{z} \widetilde  \varphi_{\tau,s}(\cdot, \xi)\|_{W^{p-1, \infty}_{z}}\, d\tau 
    + T^{1 \over 2} \Lambda(T, M)( | \xi_{v}| + |t-s| | \xi_{x}|),
  \end{align*}
  where we have used \eqref{defbrhotildephi}. 
  We finally obtain from the Gronwall inequality that 
  $$ \| \nabla_{z} \widetilde \varphi_{t,s} (\cdot, \xi) \|_{W^{p-1, \infty}_{z}} \leq   T^{1 \over 2} \Lambda(T, M)( | \xi_{v}| + |t-s| | \xi_{x}|).$$
  This ends the proof of \eqref{tildephidegreun}.

\end{proof}

\subsection{Proof of Proposition~\ref{EstB}}\label{sec-proofEstB}
Let $T(M)>0$ be the positive time provided by Proposition~\ref{prop-HJ1}.
Thanks to the regularity estimates \eqref{InegPot} and \eqref{estphi}, the equation~\eqref{Eqb} can be seen as a transport equation with  coefficients in $L^\infty(0,T;W^{p,\infty}_{z,\xi})$, plus an operator of order $0$ which is just a multiplication by a matrix also bounded in $L^\infty(0,T;W^{p,\infty}_{z,\xi})$. Therefore, the existence and uniqueness of the solution $\mathrm{B}_{t,s}$ on  $[0,\min(T(M),T_\eps)]$ follows by standard arguments.
For all $z,\xi \in \R^{2d}$, let $\widetilde{Z}_{t,s}(z,\xi)=(\widetilde{X}_{t,s}(z,\xi), \widetilde{V}_{t,s}(z,\xi))_t$ be the characteristics associated with the vector field $z=(x,v) \mapsto (  v, \nabla_{\xi_v}a_{\rho} \left(t,x,\nabla_v\varphi_{t,s} \right))$, with $(\widetilde{X}_{s,s}(z,\xi), \widetilde{V}_{s,s}(z,\xi))=z$. By the Duhamel formula, it holds for $0\leq s\leq t \leq \min(T(M),T_\eps)$,
\begin{equation}
\label{eq:BB}
 \mathrm{B}_{t,s}(z,\xi) = \mathrm{I} - \int_s^t [\mathcal{N}(\tau)  \mathrm{B}_{\tau,s} ](\widetilde{Z}_{\tau,t}(z,\xi), \xi) \, d \tau.
\end{equation}
The estimate \eqref{regStrongB} thus follows from this equation, arguing by induction (similarly to the proof of Lemma~\ref{lemestimtildephi}). 
The estimate~\eqref{SmallTimeB} then rely also on~\eqref{eq:BB}, using~\eqref{regStrongB}. 

\section{Quantum averaging lemmas}\label{sec:quantum}

In this section, we develop one of the key aspects of the proof, which is a quantum version of the averaging lemma with gain of one derivative from~\cite{HKR}. 
We recall that we intend to study the term 
\begin{equation}
\label{eq:quantum-FIO}
\int_v  \int_0^t U_{t,s}^{\mathrm{FIO}} B[\partial_{x}^{\alpha(i)}V_\eps \ast \rho, f] \, dsdv, \qquad |\alpha(i)|=m,
\end{equation}
with $B$ defined in \eqref{eq:K1} and that a naive uniform estimate relying on Lemma~\ref{lembilinB} would require a control $m+1$ derivatives of $\rho$, which we do not have; this apparent loss of derivative reflects the singularity of the Vlasov-Benney equation~\eqref{eq:vlasov-benney-intro}.

\begin{definition}\label{defAL}
Let $T>0$.
Let  $\Phi_{t,s}(z,\xi)$  be a real-valued phase, we shall say that it matches the assumption
 $(\boldsymbol{A}_{p})$ for some $p \geq 0$ if  for all $t,s \in [0,T]$, $z=(x,v),\xi=(\xi_x,\xi_v) \in \R^{2d}$, we have the estimates
\begin{equation}
\label{eq:boundPsi}
\begin{aligned}
&\underset{{\substack{0\leq\vert \alpha\vert + \vert \beta\vert \leq p
 }}}{\sup} \left|\partial^{\alpha}_z\partial^{\beta}_{\xi}   \nabla_\xi \Big[ \Phi_{t,s}(z,\xi)- (x-(t-s)v) \cdot  \xi_x - v \cdot \xi_v\Big]\right|\leq  1, \\
&\underset{{\substack{0\leq\vert \alpha\vert + \vert \beta\vert \leq p
 }}}{\sup} \left|\partial^{\alpha}_z\partial^{\beta}_{\xi} \nabla_z\Big[\Phi_{t,s} (z,\xi)- (x-(t-s)v) \cdot  \xi_x - v \cdot \xi_v\Big]\right|\leq   \langle  \xi_v \rangle  +\frac{1}{2} \langle (t-s)   \xi_x \rangle.
\end{aligned}
\end{equation}
Let $b_{t,s}(z,\xi)$ and $G_{t,s}(\xi)$ be  given smooth  amplitudes and kernels.
 We denote by $\mathcal{U}_{[\Phi, b, G]}$ the operator defined by
 $$
\mathcal{U}_{[ \Phi, b, G]} (\varrho) (t,x)=\frac{1}{(2\pi )^{2d}} \int_{v}\int_{0}^t\int_{\xi} \int_{y} e^{i \Phi_{t,s} (z,\xi) }b_{t,s}(z,\xi)  \widehat{B[\varrho,G_{t,s}]}(\xi) d\xi ds dv.
$$
\end{definition}
Let us recall that  thanks to \eqref{expressiondeB}, we have that
\begin{equation}
\label{defBrecalled} \left(\widehat{ B[\rho,G_{t, s}] }\right) (\xi)= (2\pi)^d \int_{\eta}  \frac{2}{\eps} \sin\left(\frac{\eps (\xi_{x}- \eta) \cdot \xi_{v}} {2}\right) \widehat{V_{\eps}}(\xi_{x}- \eta)
\widehat{\rho} (s, \xi_{x} - \eta)  \widehat{G_{t, s}}( \eta, \xi_{v}) \, d\eta,
\end{equation}
for  $\xi= (\xi_{x}, \xi_{v}) \in \mathbb{R}^{2d}$.
Note that the operator $\mathcal{U}_{[ \Phi, b, G]}$ thus depends on $\eps$ through the definition of $B$.

\begin{definition} For $l, \, p \in \mathbb{N}$  we consider the norm $\| \cdot \|_{T, l,p}$ defined as 
$$
\norm{G}_{T, l,p}:= \sup_{t, \, s \in [0, T]}\sum_{0\leq\vert\alpha\vert  \leq p} \norm{\langle  \xi \rangle^{l}  \partial^\alpha_{\xi} \widehat{G_{ t, s}}(\xi)}_{L^{\infty}_{\xi}}
$$
and we set
$$ 
 \norm{ b}_{L^\infty_{T} W^{k,\infty}_{z,\xi}} = \sup_{s, \, t \in [0, T]} \| b_{ t,s} \|_{W^{k, \infty}_{z,\xi}}.$$
 \end{definition}
\begin{Rem}
\label{rem:norme-av}
Note that we can use the  norms $\H^m_{r}$  to control these norms by using  a Sobolev embedding
 in $\xi$ and  \eqref{embedH01}. We have:
$$
\norm{G}_{T, l,p} \lesssim \sup_{s,\,  t \in[0, T]} \| G_{t,s}\|_{\H^{l}_{p+k}}
$$
for all $k>d$.
\end{Rem}
Let us recall the notation  $k_{d}=\lfloor d/2\rfloor+2$ that will be systematically used throughout this section. The main quantum averaging lemma is stated in the following result. 
\begin{thm}\label{M-Pro}
  For every $T_{0}>0$, 
there exists $C_{0}>0$ such that for every $T \in [0, T_{0}]$,  if the assumption  ($\boldsymbol{A}_{4k_{d} + d + 4}$) holds, 
 we have  for every  $\eps \in (0,1)$ that
$$
\norm{\mathcal{U}_{[ \Phi, b, G]}}_{\mathscr{L}(L^2(0,T;L^2(\R^d)))}\leq C_{0}   \norm{b}_{L^\infty_{T} W^{d+4k_{d}+4,\infty}_{z,\xi}}  
 \norm{\langle \eps \nabla_x \rangle^{k_{d}} \langle \eps\nabla_v \rangle^{k_{d}} G}_{T, 3k_{d}+ 2 d + 6, k_{d}+ d+2}.
 $$
 \end{thm}

This result will notably be used in the following situations:
\begin{itemize}
\item When $\Phi$ is the phase associated with the free transport operator, that is when
$$
\Phi_{t,s}(z,\xi) = (x-(t-s)v) \cdot   \xi_x + v \cdot    \xi_v.
$$
The estimates~\eqref{eq:boundPsi} are then clear, the right hand side vanishes.  Note that even  in  this case,  in $\mathcal{U}_{[ \Phi, b, G]}$  there is  still a quantum contribution  through the $\sin$ term in the definition of $B$
and the dependence of $b_{t,s}$ in the $\xi$ variable.
\item  When $\Phi$ is the phase associated with the FIO constructed in the previous section, that is when
$$
\Phi_{t,s}(z,\xi) = { 1 \over \eps } \varphi_{t,s}^\eps(z,\xi) \quad  (=  { 1 \over \eps}\varphi_{t,s}(z,\eps\xi)),
$$
where $\varphi$ satisfies the eikonal equation~\eqref{H-J}.
The estimates~\eqref{eq:boundPsi} are then a consequence of Proposition~\ref{prop-HJ1}, hold for $T= \min (T(M), T_\eps)$ 
 and are uniform in $\eps$.  Indeed, the first set of estimates directly follow from~\eqref{estphi}. For the second set of estimates, when there is at least one derivative in $\xi$, we can simply use  \eqref{estphi}  and the fact that $\varphi$ is evaluated at $\eps \xi$, so that we gain a factor $\eps$. When there is no derivative in $\xi$ at all, we can use  \eqref{estdzphi}.

\end{itemize}
In view of applications to~\eqref{eq:quantum-FIO}, this result can be used  for the amplitude $(\mathrm{B}^\eps_{t,s})_{i,j}$ of the FIO constructed in the previous section. The required regularity assumptions  come from Lemma~\ref{EstB}. The kernel $G$ will typically be the solution $f_\eps$ to the Wigner equation itself. Theorem~\ref{M-Pro} thus shows that the loss of derivative in~\eqref{eq:quantum-FIO} is only apparent.


\subsection{Proof of Theorem \ref{M-Pro}}
In the  proof, we shall only denote $\mathcal{U}$ instead of $\mathcal{U}_{[\Phi, b, G]}$.
We shall rewrite $\mathcal{U}$ as a pseudodifferential operator with an operator-valued symbol. 
 By using the expression \eqref{defBrecalled}, we have 
\begin{equation}
\label{expan}
\begin{aligned}
\mathcal{U}(\varrho) (t,x)
&=2 \int_{\eta}e^{i x\cdot \eta} \Bigg[\int_{0}^t \Bigg( \int_{v}\int_{\xi} e^{-ix\cdot \eta} e^{i \Phi_{t,s} (z,\xi) } b_{t,s}(z,\xi)\\
& \quad\quad\quad \quad\quad\quad\quad\quad\quad   \mathcal{F}_{x,v}G_{t,s}(\xi_x-\eta,\xi_v)\left(  \frac{1}{\eps} \sin\left(\frac{\eps \xi_v\cdot \eta} {2}\right)  \widehat{V}(\eps \eta)\right) d\xi  dv  \Bigg) \widehat \varrho (s,\eta) ds \Bigg] d\eta,
\end{aligned}
\end{equation}
where $z=(x,v)$ in the above expression.
This can be seen  as a pseudodifferential operator acting on $\varrho$, that is to say, 
\begin{equation}
\mathcal{U}(\varrho) =\mathrm{Op}_{\mathrm{L}} (\varrho),
\end{equation}
for the  quantization~\eqref{PseudoDef-operator} and where, for $x, \eta \in \R^d \times \R^d$, $\mathrm{L}(x,\eta)$ is an operator-valued symbol acting on $L^2(0, T)$.  Namely $\mathrm{L}(x,\eta) : L^2(0,T) \to L^2(0,T)$ is the operator defined by
\begin{equation}
\label{def:L}
(\mathrm{L}(x,\eta)  \Upsilon)(t) = 2 \int_0^t H_{t,s}(x,\eta) \Upsilon(s) \, ds, \quad \Upsilon \in L^2(0,T),
\end{equation}
with
\begin{equation}
\label{def:H}
H_{t,s}(x,\eta) = \int_{v}\int_{\xi} e^{-ix\cdot \eta} e^{i    \Phi_{t,s} (z,\xi)}  b_{t,s}(z,\xi) \mathcal{F}_{x,v}G_{t,s}(\xi_x-\eta,\xi_v) \frac{1}{\eps} \sin\left(\frac{\eps \xi_v\cdot \eta} {2}\right)   \widehat{V}(\eps \eta)  \, d\xi dv  .
\end{equation}
We will prove that $H_{t,s}$ is a well-defined oscillating integral, and that it  enjoys the following key estimate.
\begin{prop} 
\label{prop-finH}
With the same notations as in Theorem~\ref{M-Pro}, for every $\ell \in \mathbb{N}$,  $0 \leq \vert \alpha\vert,\vert \beta\vert \leq k_{d}$,
if  the Assumption ($\boldsymbol{A}_{4k_{d} + d + 4}$) holds,  then
 for all $x, \eta \in \R^d$ and all $s,t \in [0,T]$, we have 
\begin{multline}\label{MajorationH}
|\partial^\alpha_x \partial^\beta_\eta H_{t,s}(x,\eta)| \lesssim  \\
   \norm{b}_{L^\infty_{T}W^{d+4k_{d}+4,\infty}_{z,\xi}}  
 \norm{\langle \eps\nabla_x \rangle^{k_{d}} \langle \eps\nabla_v \rangle^{k_{d}} G}_{2\ell + 3k_{d}+ 2 d + 4, k_{d} + d +2}\left( \int_{\xi_x} \frac{\langle  \xi_x\rangle }{  \langle \xi_x-\eta \rangle^{\ell}\langle (t-s)\xi_x\rangle^2}  \, d \xi_x  \right).
\end{multline}
 
\end{prop}
The proof of this proposition is technical and is left to the following subsection. 
Let us explain how it leads to a proof of Theorem~\ref{M-Pro}.
We have from \eqref{def:L} that
$$ (\partial_{x}^\alpha \partial_{\eta}^\beta L(x, \eta)   \Upsilon)(t) = 2 \int_0^t  \partial_{x}^\alpha \partial_{\eta}^\beta H_{t,s}(x,\eta) \Upsilon(s) \, ds, $$
therefore, by using the Schur test, we deduce that
$$   \left\|    \partial^\alpha_x \partial^\beta_\eta \mathrm{L}(x, \eta) \right\|_{\mathscr{L}(L^2(0,T))}^2 
 \leq  \sup_{0\leq t \leq T}\left(\int_0^t  | \partial^\alpha_x \partial^\beta_\eta H_{t,s}(x,\eta)| \, ds \right) \sup_{0\leq s \leq T}\left( \int_s^T  | \partial^\alpha_x \partial^\beta_\eta H_{t,s}(x,\eta)|  \, dt\right).$$
Thanks to Proposition~\ref{prop-finH}, taking $\ell=d+1$ it holds
\begin{multline*}
\sup_{0\leq t \leq T}\left(\int_0^t  | \partial^\alpha_x \partial^\beta_\eta H_{t,s}(x,\eta)| \, ds \right)  
\lesssim   \left( \sup_{0\leq t \leq T} \int_{\xi_x} \frac{1}{\langle \xi_x - \eta\rangle^r} \int_0^t  \frac{\langle\xi_x\rangle}{\langle (t-s)\xi_x\rangle^2}  \, ds  d \xi_x  \right)  \\
\times \norm{b}_{L^\infty_{T}W^{d+4k_{d}+4,\infty}_{z,\xi}}\norm{\langle \eps\nabla_x \rangle^{k_{d}} \langle \eps\nabla_v \rangle^{k_{d}} G}_{3k_{d}+ 2 d + 6, k_{d} + d +2} .
\end{multline*} 
Since we have 
$$  \int_{\xi_x} \frac{1}{\langle \xi_x - \eta\rangle^{d+1}} \int_0^t  \frac{\langle\xi_x\rangle}{\langle (t-s)\xi_x\rangle^2}  \, ds  d \xi_x  \lesssim    \int_{\xi_x} \frac{1}{\langle \xi_x - \eta\rangle^{d+1}} d\xi_x \int_0^{+\infty}  \frac{1}{\langle\tau\rangle^2}  \, d \tau \lesssim   1, $$
we get that
$$ \sup_{0\leq t \leq T}\left(\int_0^t  | \partial^\alpha_x \partial^\beta_\eta H_{t,s}(x,\eta)| \, ds \right)  \\
\lesssim\norm{b}_{L^\infty_{T}W^{d+4k_{d}+4,\infty}_{z,\xi}}\norm{\langle \eps\nabla_x \rangle^{k_{d}} \langle \eps\nabla_v \rangle^{k_{d}} G}_{3k_{d}+ 2 d + 6, k_{d} + d +2}.$$
With similar arguments, again by Proposition~\ref{prop-finH}, we also  infer
\begin{align*}
\sup_{0\leq s \leq T}\left( \int_s^T  | \partial^\alpha_x \partial^\beta_\eta H_{t,s}(x,\eta)|  \, dt\right)  \lesssim  \norm{b}_{L^\infty_{T}W^{d+4k_{d}+4,\infty}_{z,\xi}}\norm{\langle \eps\nabla_x \rangle^{k_{d}} \langle \eps\nabla_v \rangle^{k_{d}} G}_{3k_{d}+ 2 d + 6, k_{d} + d +2}.
\end{align*}
We conclude that
\begin{equation}
\max_{|\alpha|,|\beta| \leq k_{d}}\sup_{x,\eta \in \R^d} \left\|    \partial^\alpha_x \partial^\beta_\eta \mathrm{L} \right\|_{\mathscr{L}(L^2(0,T))} \lesssim  \norm{b}_{L^\infty_{T}W^{d+4k_{d}+4,\infty}_{z,\xi}}\norm{\langle \eps\nabla_x \rangle^{k_{d}} \langle \eps\nabla_v \rangle^{k_{d}} G}_{3k_{d}+ 2 d + 6, k_{d} + d +2} .
\end{equation}
Therefore, by the Cald\'eron-Vaillancourt theorem for operator-valued symbols (see  Proposition~\ref{SymTim} with  $\mathbf{H}= L^2(0,T)$), we obtain the desired result, namely
$$
\norm{\mathcal{U}}_{\mathscr{L}(L^2(0,T;L^2(\R^d)))}\lesssim  \norm{b}_{L^\infty_{T}W^{d+4k_{d}+4,\infty}_{z,\xi}}\norm{\langle \eps\nabla_x \rangle^{k_{d}} \langle \eps\nabla_v \rangle^{k_{d}} G}_{3k_{d}+ 2 d + 6, k_{d} + d +2} .
$$

\subsection{Proof of Proposition \ref{prop-finH} }\label{SubHTec}


We shall  prove  that $H$, as defined in~\eqref{def:H} is a well-defined oscillatory integral  in $v$ and $\xi$, thanks to a non-stationary phase argument. 
This is where various bounds from below  for certain derivatives of the phase are crucial. As a matter of fact, the absolute convergence in $\xi$ can be easily  ensured thanks to the decay of $\mathcal{F}_{x,v}G(\xi_x - \eta, \xi_v)$; however to obtain appropriate uniform estimates with respect to $\eta$, a special treatment is required for the decay in $\xi_x$.

To this end, it is convenient to distinguish between a low and a high frequency regime (in $\eta$): introducing a cut-off function $\chi \in \mathscr{C}^\infty(\R_+; \R_+)$ with $\chi \equiv 1$ on $[0,1]$ and $\chi \equiv 0$ on $[2,+\infty)$, we write
\begin{align*}
H_{t,s}(x,\eta) &= H_{t,s}(x,\eta) \chi\left(\eps |\eta|\right) + H_{t,s}(x,\eta) \left[1- \chi\left(\eps |\eta|\right)\right] \\
&=: H_{t,s}^-(x,\eta) + H_{t,s}^+(x,\eta),
\end{align*}
and shall argue differently according to the regime in $\eta$.
We shall focus on the case $d \geq 2$ in the following, the case $d=1$ being a simple adaptation (we just have to notice that
the direction orthogonal to $\eta$ considered below is empty).

\subsubsection{The low $\eta$ regime}

We start by studying the term $H^-_{t,s}(x,\eta) = H_{t,s}(x,\eta) \chi\left(\eps |\eta|\right) $, which roughly corresponds to $\{\eps |\eta| \leq 2\}$. In this regime, we can consider the  operator
\begin{equation}
\label{eq:Lxiv}
\mathcal{L}_{\xi_v} = \frac{\lambda   - i  \nabla_{\xi_v} \Phi_{t,s}(z,\xi)\cdot \nabla_{\xi_v}}{ \lambda\   +|\nabla_{\xi_v} \Phi_{t,s}(z,\xi)|^2},
\end{equation}
where, according to~\eqref{eq:boundPsi}, choosing $\lambda>0$ large enough, the following bound from below holds:
\begin{equation}
\label{eq:below-xiv}
\lambda   +|\nabla_{\xi_v} \Phi_{t,s}(z,\xi)|^2 \geq C  \langle v \rangle^2.
\end{equation}
By construction $\mathcal{L}_{\xi_v} e^{i  \Phi_{t,s} (z,\xi)}=e^{i    \Phi_{t,s} (z,\xi)}$. Let $p_1$ be an integer to be chosen large enough. Thanks to this identity, we have
\begin{equation}
\label{def:H2}
\begin{aligned}
&H^-_{t,s}(x,\eta) \\&= \int_{v}\int_{\xi} e^{-ix\cdot \eta} \mathcal{L}_{\xi_v}^{p_1} \left(e^{i  \Phi_{t,s} (z,\xi)} \right) b_{t,s}(z,\xi)  \mathcal{F}_{x,v}G_s(\xi_x-\eta,\xi_v) \frac{1}{\eps} \sin\left(\frac{\eps \xi_v\cdot \eta} {2}\right)   W^-(\eps \eta)  \, d\xi dv   
\\& = \int_{v}\int_{\xi} e^{-ix\cdot \eta} e^{i    \Phi_{t,s} (z,\xi)}\left(\mathcal{L}_{\xi_v}^T\right)^{p_1} \left( b_{t,s}(z,\xi)  \mathcal{F}_{x,v}G_s(\xi_x-\eta,\xi_v) \frac{1}{\eps} \sin\left(\frac{\eps \xi_v\cdot \eta} {2}\right) W^-(\eps \eta)  \right)  \, d\xi dv ,
\end{aligned}
\end{equation}
where we have set $W^-(\eps \eta)= \widehat{V}(\eps \eta) \chi\left(\eps| \eta|\right)$ and $\mathcal{L}_{\xi_v}^T$ is the formal adjoint of $\mathcal{L}_{\xi_v}$. 

\begin{lem}\label{MajConv}For every integer $p_1\geq 1$ and all $l_x, l_v>0$,  if the
assumption  ($\boldsymbol{A}_{p_{1}}$)  holds, we have the estimate
\begin{multline*}
\bigg\vert \left(\mathcal{L}_{\xi_v}^T\right)^{p_1} \left({b}_{t,s}(z,\xi)\mathcal{F}_{x,v}G_s(\xi_x-\eta,\xi_v) \frac{1}{\eps} \sin\left(\frac{\eps \xi_v\cdot \eta} {2}\right) W^-(\eps \eta) \right)   \bigg \vert   \\ 
\lesssim {1 \over \eps}  \frac{1}{\langle  \xi_x -\eta \rangle^{l_x} \langle  \xi_v \rangle^{l_v}\langle v \rangle^{p_1}} \norm{b}_{L^\infty_{T}W^{p_1,\infty}_{z,\xi}}
\norm{G}_{T, l_x+ l_v,p_1}.
\end{multline*}

\end{lem}

\begin{proof}[Proof of Lemma~\ref{MajConv}]
The adjoint of $\mathcal{L}_{\xi_v}$ reads, for a smooth function $u$, as
\begin{align*}
\mathcal{L}_{\xi_v}^T u &=\frac{\lambda  + i  \Delta_{\xi_v}  \Phi_{t,s} + i   \nabla_{\xi_v} \Phi_{t,s} \cdot\nabla_{\xi_v}}{\lambda    +\vert\nabla_{\xi_v} \Phi_{t,s} \vert ^2 \ } u  +  i  \frac{\nabla_{\xi_v}  \Phi_{t,s}\cdot  \nabla_{\xi_v}  \vert\nabla_{\xi_v}  \Phi_{t,s}\vert ^2 }{\left(\lambda    +\vert\nabla_{\xi_v}  \Phi_{t,s}\vert ^2 \right)^2}  u.
\end{align*}
By induction, we obtain the expansion
\begin{equation}\label{expanL-xiv} 
\left(\mathcal{L}_{\xi_v}^T\right)^{p_1} = \sum_{\vert \alpha\vert \leq p_1} c_{t,s}^{\alpha}(z,\xi) \partial_{\xi_v}^{\alpha},
\end{equation}
where $c_{t,s}^{\alpha}$ involves at most $p_{1}$ derivatives of $\nabla_{\xi_{v}} \Phi_{t,s}$.
Moreover, since we have  the lower bound  \eqref{eq:below-xiv} and thanks to the assumption ($\boldsymbol{A}_{p_{1}}$) also the upper bound
\begin{equation}
\label{eq:upper1-low-xiv}
\sup_{0\leq |\alpha|\leq p_1} | \partial^\alpha_{\xi_v} \nabla_{\xi_{v}} \Phi_{t,s} (z,\xi)| \lesssim  {  \langle v \rangle},
\end{equation}
we obtain for 
the functions $c_{t,s}^{\alpha,\beta}$  the  estimate 
\begin{equation}\label{EstC-xiv}
\begin{aligned}
\vert c_{t,s}^{\alpha}(z,\xi)\vert &\lesssim       \frac{1}{ \langle v  \rangle^{p_1}},
\end{aligned}
\end{equation}
as long as the  ($\boldsymbol{A}_{p_{1}}$) assumption is matched.

\begin{Rem} Let us record for later use that since  we also have 
$$
\sup_{0\leq |\alpha_z| + |\alpha_\xi|\leq p_1+  2 k_{d} +p_{v}} | \partial^{\alpha_z}_{z} \partial^{\alpha_\xi}_{\xi} \partial_{\xi_v} \Phi_{t,s} (z,\xi)| \lesssim    \langle v\rangle,
$$
when the assumption   ($\boldsymbol{A}_{p_{1}+ 2 k_{d}+p_{v}}$) is matched for some $p_{v} \in \mathbb{N}$, then 
the derivatives $ \partial_{z}^{\alpha_{z}} \partial_{\xi_v}^{\alpha_{\xi_{v}}} c_{t,s}^{\alpha}$ also  satisfy an estimate similar to~\eqref{EstC-xiv}, namely
\begin{equation}\label{EstC-xiv-deriv}
\begin{aligned}
\sup_{0\leq |\alpha_z| + |\alpha_\xi|\leq  2 k_{d} +p_{v} }  \vert \partial^{\alpha_z}_{z}  \partial^{\alpha_\xi}_{\xi}  c_{t,s}^{\alpha}(z,\xi)\vert  \lesssim    \frac{1}{ \langle v  \rangle^{p_1}}.
\end{aligned}
\end{equation}
\end{Rem}
Introducing
\begin{equation}\label{def:A-}
\mathcal{A}_{t,s}^-(z, \xi, \eta):=\left(   {b}_{t,s}(z,\xi)  \mathcal{F}_{x,v}G_{t,s}(\xi_x-\eta,\xi_v) \frac{1}{\eps} \sin\left(\frac{\eps \xi_v\cdot \eta} {2}\right)  W^-(\eps \eta)\right),
\end{equation}
  we have
  $$
H^-_{t,s}(x,\eta) = \int_{v}\int_{\xi} e^{-ix\cdot \eta} \left(\mathcal{L}_{\xi_v}^T\right)^{p_1} \mathcal{A}^-_{t,s}\,  \, d\xi dv.
$$
We can  control the action of $\partial_{\xi_v}^{\alpha}$ on $\mathcal{A}_{t,s}^-$
   by using the Leibniz formula. When 
 $\partial_{\xi_v}$ acts on $\sin\left(\frac{\eps \xi_v\cdot \eta}  {2} \right)$,  a power of  $\eps \eta$ appear, which  can  be absorbed thanks to the  
  potential  $W^-$  (a  more involved procedure will be required in the high $\eta$ regime). This observation leads to the estimate 
\begin{equation}\label{MajHH}
\bigg\vert \partial_{\xi_v}^{\alpha}
 \mathcal{A}_{t,s}^- \bigg\vert
\lesssim   {1 \over \eps}
		  \sum_{0\leq \alpha'\leq \alpha}  |\pa^{\alpha'}_{\xi_v}  \mathcal{F}_{x,v}G_{t,s}(\xi_x-\eta,\xi_v)|   
		   \norm{b}_{L^\infty_{T}W^{p_1,\infty}_{z,\xi}} 
\end{equation}
and by~\eqref{EstC-xiv}, we eventually obtain
\begin{equation}
\label{conclulemH}
\bigg\vert \left(\mathcal{L}_{\xi_v}^T\right)^{p_1}  \mathcal{A}_{t,s}^-
 \bigg \vert  \\ \lesssim {1 \over \eps }\frac{1}{\langle  \xi_x-\eta \rangle^{l_x}\langle  \xi_v \rangle^{l_v}\langle v\rangle^{p_1} }\norm{b}_{L^\infty_{T}W^{p_1,\infty}_{z,\xi}} \norm{G}_{T, l_x + l_v ,p_1},
\end{equation}
which concludes the proof of the lemma.
\end{proof}

Therefore, by the identity~\eqref{def:H2} of Lemma~\ref{MajConv}, choosing $p_1, l_x, l_v$ all strictly larger than $d$, $H_{t,s}^-$ can be turned into an absolutely converging integral which justifies the definition of $H_{t,s}^-$ as an oscillating integral. 

%

In the next lemma, we study the action of derivatives with respect to $x$ and $\eta$ applied to $H^-$.
\begin{lem}
\label{lem-paH}
For every  $\ell>d$, $p_1>d$,  $p_{v}>k_{d}$,  if the assumption  $(\boldsymbol{A}_{p_{1}+ p_{v}+ 2k_{d}})$ holds, then  for every $0 \leq \vert \alpha\vert,\vert \beta\vert \leq k_{d}$, 
 $\partial^{\alpha}_{\eta}\partial^{\beta}_{x}H_{t,s}^-$ can be rewritten  under the  form 
\begin{equation}\label{Hder2}
\partial^{\alpha}_{\eta}\partial^{\beta}_{x}H^-_{t,s}(x,\eta)=  \int_{v} \int_{\xi}e^{-ix\cdot \eta}  e^{i \Phi_{t,s}(z,\xi)} d_{t,s}^{\alpha,\beta}(z,\xi,\eta) \, d\xi dv, 
\end{equation}
where $d^{\alpha,\beta}$ satisfies
\begin{multline}\label{EstD}
| d_{t,s}^{\alpha,\beta}(z,\xi,\eta)| \\
\lesssim \frac{\langle \xi_x \rangle  }{  \langle v \rangle^{p_1}\langle \xi_v\rangle^\ell \langle \xi_x-\eta \rangle^\ell \langle (t-s)\xi_x\rangle^{p_{v}-k_{d}}} 
\norm{b}_{L^\infty_{T}W^{p_1+ p_{v}+  2k_{d},\infty}_{z,\xi}} \norm{\langle \eps\nabla_v\rangle^{k_{d}-1}  G}_{T,  2\ell +k_{d} + p_{v}+ 2, p_1+k_{d}}.
\end{multline}
\end{lem}

\begin{proof}[Proof of Lemma~\ref{lem-paH}]
For a given function $ a(z, \xi, \eta), $
 we observe that
 \begin{align*}
\partial_{\eta_j}  \int_{v} \int_{\xi} e^{-ix\cdot \eta}  e^{ i \Phi_{t,s}(z,\xi)} a(z,\xi,\eta)\, d\xi d \eta &=\int_{v} \int_{\xi} \partial_{\eta_j}\left( e^{-ix\cdot \eta}  e^{ i \Phi_{t,s}(z,\xi)} \right) a(z,\xi,\eta) \, d\xi dv  
\\&+\int_{v} \int_{\xi}e^{-ix\cdot \eta}  e^{ i \Phi_{t,s}(z,\xi)} \partial_{\eta_j} a(z,\xi,\eta) \, d\xi dv .
\end{align*}
Since we can write
\begin{multline*}
\partial_{\eta_j}\left(e^{-ix\cdot \eta}  e^{ i \Phi_{t,s}(z,\xi)}\right)=-(\partial_{{\xi_x}_j} + (t-s) \partial_{\xi_{v_{j}}})\left(e^{-ix\cdot \eta}   e^{ i  \Phi_{t,s}(z,\xi)}\right)  \\+i \left(\partial_{{\xi_x}_j} \Phi_{t,s}(z,\xi)-(x_j - (t-s) v_{j})\right)\left(e^{-ix\cdot \eta}  e^{ i  \Phi_{t,s}(z,\xi)}\right)
 \\+ i(t-s)\left( \partial_{\xi_{v_{j}}} \Phi_{t,s}- v_{j}\right)\left(e^{-ix\cdot \eta}e^{ i  \Phi_{t,s}(z,\xi)}\right) ,
\end{multline*}
we obtain by integration by parts in $\xi_x$ and $\xi_{v}$ that 
\begin{align}
\label{identiteXetaj}
\partial_{\eta_j} \int_{v} \int_{\xi} e^{-ix\cdot \eta}  e^{ i \Phi_{t,s}(z,\xi)} a(z,\xi,\eta)\, d\xi d \eta 
&=  \int_{v} \int_{\xi} e^{-ix\cdot \eta}  e^{ i \Phi_{t,s}(z,\xi)} X_{\eta_{j}}a(z,\xi,\eta)d\xi d \eta
\end{align}
where the vector field $X_{\eta_{j}}$ is defined as 
\begin{equation}
\label{defXetaj}
X_{\eta_{j}}= \partial_{{\xi_x}_j} + (t-s) \partial_{\xi_{v_{j}}} + \partial_{\eta_{j}} + i \left(\partial_{{\xi_x}_j} \Phi_{t,s}(z,\xi)-(x_j - (t-s) v_{j})\right)
+  i(t-s)\left( \partial_{\xi_{v_{j}}} \Phi_{t,s}- v_{j}\right).
\end{equation}
In a similar way, we can write
\begin{equation}
\label{identiteXxj}
\partial_{x_j}  \int_{v} \int_{\xi} e^{-ix\cdot \eta}  e^{ i \Phi_{t,s}(z,\xi)} a(z,\xi,\eta)\, d\xi d \eta 
=  \int_{v} \int_{\xi} e^{-ix\cdot \eta}  e^{ i \Phi_{t,s}(z,\xi)} X_{x_{j}}a(z,\xi,\eta)d\xi d \eta,
\end{equation}
where
\begin{equation}
\label{defXxj}
X_{x_{j}}= \partial_{x_j} + i \left( \xi_{x_{j}}- \eta_{j}\right) + i\left( \partial_{x_{j}} \Phi_{t,s}(z, \xi)- \xi_{x_{j}}\right).
\end{equation}
From these observations, we thus get that
$$ \partial^{\alpha}_{\eta}\partial^{\beta}_{x}H^-_{t,s}(x,\eta)=  \int_{v} \int_{\xi}e^{-ix\cdot \eta}  e^{i \Phi_{t,s}(z,\xi)}
  X_{\eta}^\alpha  X_{x}^\beta (\mathcal{L}_{\xi_{v}}^T)^{p_{1}} \mathcal{A}_{t,s}^-\, d \xi dv$$
where we have set
$$ X_{\eta}^\alpha = X_{\eta_{d}}^{\alpha_{d}} \cdots X_{\eta_{1}}^{\alpha_{1}}, \quad X_{x}^\beta = X_{x_{d}}^{\beta_{d}} \cdots X_{x_{1}}^{\beta_{1}}.$$
Since we have  the upper bounds
\begin{equation}
\label{eq:upper2-low}
\begin{aligned}
& \sup_{0\leq |\alpha_z| +  |\alpha_{\xi}|\leq  2  k_{d}+p_{v}} \left|\partial^{\alpha_z}_{z}  \partial^{\alpha_\xi}_{\xi} \left(\partial_{\xi_{x_j}} \Phi_{t,s}(z,\xi)-(x_j- (t-s)v_{j})\right)\right| \\
&\qquad \qquad \qquad \qquad \qquad \qquad \qquad \qquad \qquad  +   \left|\partial^{\alpha_z}_{z}  \partial^{\alpha_\xi}_{\xi} \left(\partial_{\xi_{v_j}} \Phi_{t,s}(z,\xi)-v_j \right)\right| \lesssim  1, \\
 &  \sup_{0\leq |\alpha_z| +  |\alpha_{\xi}|\leq   2 k_{d}+p_{v}}  \left|\partial^{\alpha_z}_{z}  \partial^{\alpha_\xi}_{\xi} \left(\partial_{x_j}  \Phi_{t,s}(z,\xi)-\xi_{x_j}\right)\right| \lesssim \langle \xi_v\rangle +  \langle (t-s)  \xi_{x} \rangle
\end{aligned}
\end{equation}
when  the $(\boldsymbol{A}_{2k_{d}+p_{v}})$ assumption is matched,  we can
expand  $X_{x}^\beta X_{\eta}^\alpha$  by using the definitions \eqref{defXetaj}, \eqref{defXxj} and the  Leibniz formula under the form
\begin{equation}
\label{XetaXxexpansion}
 X_{\eta}^\alpha  X_{x}^\beta=  \sum_{\substack{0\leq\vert\gamma\vert\leq k_{d} \\ 0\leq \vert\sigma\vert+\vert\rho\vert + \vert \mu \vert\leq k_{d}}} e_{t,s}^{\alpha,\beta,\gamma,\sigma,\rho}(z,\xi,\eta)\partial^{\gamma}_x\,\partial^{\sigma}_{\xi_x}\,\partial^{\rho}_\eta \left( (t-s) \partial_{\xi_{v}}\right)^\mu
 \end{equation}
where we have for the coefficients the estimate
 \begin{equation}
\label{estim-eij}
\sup_{0\leq |\alpha_v| \leq p_{v}} \vert  \partial^{\alpha_v}_{v} e_{t,s}^{\alpha,\beta,\gamma,\sigma,\rho}(z,\xi,\eta)\vert \lesssim   \left(\langle\xi_v\rangle +   \langle (t-s) \xi_{x}  \rangle\right)^{k_{d}} \langle \xi_x-\eta\rangle^{k_{d}}. 
\end{equation}
We finally introduce  the  operator
\begin{equation}
\label{eq:Lv}
\mathcal{L}_v = \frac{\lambda\langle \xi_v \rangle^2  - i  \nabla_v \Phi_{t,s}(z,\xi)\cdot \nabla_v}{ \lambda \langle \xi_v \rangle^2   +|\nabla_v \Phi_{t,s}(z,\xi)|^2},
\end{equation}
where, according to~\eqref{eq:boundPsi}, for  $\lambda>0$ large enough, the following bound from below holds
\begin{equation}
\label{eq:below-v}
\lambda  \langle \xi_v \rangle^2  +|\nabla_v \Psi_{t,s}(z,\xi)|^2 \geq  {1 \over 2}\left( \langle (t-s)\xi_x  \rangle^2 + \langle  \xi_v \rangle^2 \right).
\end{equation}
By also using  the upper bound provided by~\eqref{eq:boundPsi}, 
 $\mathcal{L}_v$ can be seen as a first order differential operator in $v$ whose coefficients and their derivatives  are bounded by 
 $${\langle \xi_{v}\rangle  \over  \langle \xi_{v}  \rangle+\langle (t-s) \xi_{x} \rangle}.$$
By construction it holds $ \mathcal{L}_ve^{i  \Phi_{t,s} (z,\xi)}=e^{ i    \Phi_{t,s} (z,\xi)}$. 
We therefore have
$$
\partial^{\alpha}_{\eta}\partial^{\beta}_{x} H^-_{t,s}(x,\eta) = 
  \int_{v}\int_{\xi} e^{-ix\cdot \eta} e^{i    \Phi_{t,s} (z,\xi)} 
  (\mathcal{L}_v^T)^{p_{v}}  X_{x}^\beta X_{\eta}^\alpha  (\mathcal{L}_{\xi_{v}}^T)^{p_{1}} \mathcal{A}_{t,s}^-\, d \xi dv.  
$$
and we  set 
\begin{equation*}
d^{\alpha,\beta}_{t,s}(z,\xi,\eta) =   (\mathcal{L}_v^T)^{p_{v}}  X_{x}^\beta X_{\eta}^\alpha  (\mathcal{L}_{\xi_{v}}^T)^{p_{1}} \mathcal{A}_{t,s}^-,
\end{equation*}
so that \eqref{Hder2} holds.

The adjoint of $\mathcal{L}_{v}$ reads when acting on a  smooth function $u$, as 
\begin{align*}
\mathcal{L}_{v}^T u &=\frac{\lambda \langle \xi_v\rangle^2 + i  \Delta_{v}  \Phi_{t,s} + i  \nabla_{v} \Phi_{t,s}  \cdot\nabla_{v}}{\lambda \langle\xi_v\rangle^2 +\vert\nabla_{v} \Phi_{t,s} \vert ^2 \ } u  +  i  \frac{\nabla_{v}  \Phi_{t,s}\cdot  \nabla_{v}    \vert\nabla_{v}  \Phi_{t,s}\vert ^2}{\left(\lambda \langle \xi_v\rangle^2+\vert\nabla_{v}  \Phi_{t,s}\vert ^2 \right)^2}  u.
\end{align*}
We therefore  obtain an  expansion
\begin{equation}\label{expanL} 
\left(\mathcal{L}_{v}^T\right)^{p_{v}} = \sum_{ \vert \beta\vert \leq p_{v}} c_{t,s}^{\beta}(z,\xi) \partial_v^{\beta},
\end{equation}
where 
the functions $c_{t,s}^{\beta}$ satisfy the  estimate 
\begin{equation}\label{EstC}
\begin{aligned}
\vert c_{t,s}^{\beta}(z,\xi)\vert 
\lesssim    \frac{\langle \xi_v \rangle^{p_{v}}}{ \left( \langle \xi_{v} \rangle+  \langle (t-s)\xi_x  \rangle\right)^{p_{v}}}.
\end{aligned}
\end{equation}
As a result, combining~\eqref{EstC-xiv-deriv}, ~\eqref{estim-eij} and~\eqref{EstC}, we  can write
\begin{multline}\label{MajH1}
d_{t,s}^{\alpha,\beta}(z,\xi,\eta)=  \\ \sum_{\substack{0\leq\vert\gamma\vert\leq k_{d} \\ 0\leq \vert\sigma\vert+\vert\rho\vert + \vert \mu \vert\leq k_{d}  }} \sum_{0\leq \vert \alpha'\vert \leq p_1} \sum_{0\leq \vert \beta'\vert\leq p_{v}} f_{t,s}^{\alpha,\beta, \gamma,\sigma,\rho,\lambda, \alpha',\beta'}(z,\xi,\eta)  \partial^{\gamma}_x\,\partial^{\sigma}_{\xi_x}\partial^{\rho}_\eta \left( (t-s) \partial_{\xi_{v}}\right)^\mu  \partial_{\xi_v}^{\alpha'} \partial_v^{\beta'}  \mathcal{A}_{t,s}^-,
\end{multline}
where the $f_{t,s}^{\alpha,\beta, \gamma,\sigma,\rho,\mu,\alpha',\beta'}$ satisfy the estimate
\begin{equation}\label{MajH2}
\bigg\vert f_{t,s}^{\alpha,\beta, \gamma,\sigma,\rho,\mu, \alpha',\beta'}(z,\xi,\eta)  \bigg  \vert  \lesssim \frac{ \langle \xi_v\rangle^{p_{v}}
\langle \xi_x-\eta\rangle^{k_{d}}}{ \langle  v  \rangle^{p_1}\left( \langle \xi_{v} \rangle + \langle  (t-s) \xi_x  \rangle \right)^{p_{v}-k_{d}}}
 \lesssim  \frac{ \langle \xi_v\rangle^{p_{v}}
\langle \xi_x-\eta\rangle^{k_{d}}}{ \langle  v  \rangle^{p_1} \langle  (t-s) \xi_x  \rangle^{p_{v}-k_{d}}} .
\end{equation}
Finally, there only remains to study the action of $\partial^{\gamma}_x\,\partial^{\sigma}_{\xi_x}\partial^{\rho}_\eta  \partial_{\xi_v}^{\alpha'+\mu} \partial_v^{\beta'} $ on 
$\mathcal{A}^-_{t,s}$, recalling~\eqref{def:A-}. 
Once again, we can use the Leibniz formula.  
There is no issue for the derivatives in 
 $\xi_x$,  $x$ and $v$. For the derivatives in $\eta$, when they fall on the potential $W^-$, we actually gain a power of $\eps$, 
 and for the derivatives in $\xi_{v}$ and $\eta$ when they fall on the  $\sin $ term, we use that
 $$ \left| {1 \over \eps} \partial_{\xi_{v}}^{\widetilde{\alpha}}  \partial_{\eta}^{\widetilde{\beta}} \left(\sin\left(\frac{\eps \xi_v\cdot \eta} {2}\right)\right) \right|
 \lesssim  \mathrm{1}_{| \widetilde \beta | \geq 1} \langle \xi_{v}\rangle \langle \eps \xi_{v}\rangle^{| \widetilde{\beta}| -1} \langle \eps \eta\rangle^{\widetilde{ \alpha}}
  + \mathrm{1}_{| \widetilde{\beta}|=0,\,  | \widetilde{\alpha}| \geq 1}
   | \eta| \, | \eps \eta|^{| \widetilde{\alpha}|-1} + \mathrm{1}_{| \widetilde{\beta}|=0, \, | \widetilde{\alpha}| = 0}
   | \xi_{v}|  | \eta|, $$
   where in the latter, we have used that $|\sin x|\leq |x|$ to absorb the prefactor $\eps^{-1}$.

 We can then rely on  the potential $W^-(\eps \eta)$ to absorb the powers of $\eps| \eta|$. Note that because of this property,
 in this regime, we do not need to use that in \eqref{MajH1} we have $((t-s) \partial_{\xi_{v}})^\mu$ instead of  $ \partial_{\xi_{v}}^\mu$.
  Since   $|\eta|\leq \langle \xi_x \rangle \langle \xi_x -\eta\rangle$, this
   yields the  estimate 
\begin{multline}\label{MajH3}
\bigg\vert \partial_{\xi_v}^{\alpha'}\,\partial_{v}^{\beta'}\,\partial^{\gamma}_x\,\partial^{\sigma}_{\xi_x}\,\partial^{\rho}_\eta\left(   {b}_{t,s}(z,\xi)  \mathcal{F}_{x,v}G_{t,s}(\xi_x-\eta,\xi_v) \frac{1}{\eps} \sin\left(\frac{\eps \xi_v\cdot \eta} {2}\right)  W^-(\eps \eta)\right)\bigg \vert \\
\lesssim   \langle \xi_x \rangle   \langle \xi_x -\eta\rangle  \langle \xi_v\rangle
		  \sum_{0\leq |\alpha''| +|\beta''|\leq p_1 + k_{d}}  \langle \eps\xi_v\rangle^{k_{d}-1} |\pa^{\alpha''}_{\xi_x} \pa^{\beta''}_{\xi_v}  \mathcal{F}_{x,v}G_{t,s}(\xi_x-\eta,\xi_v)|   
		   \norm{b}_{L^\infty_{T}W^{p_1+ p_{v}+2k_{d},\infty}_{z,\xi}}. 
\end{multline}
Combining \eqref{MajH1}, \eqref{MajH2} and \eqref{MajH3}, we obtain that
\begin{multline}
| d_{t,s}^{\alpha,\beta}(z,\xi,\eta)|  \\ \lesssim  
\frac{
\langle \xi_x \rangle }{\langle v \rangle^{p_1} \langle \xi_v\rangle^\ell \langle \xi_x-\eta\rangle^\ell \langle (t-s) \xi_x  \rangle^{p_{v}-k_{d}}}  \norm{b}_{L^\infty_{T}W^{p_1+ p_{v}+  2k_{d},\infty}_{z,\xi}} \norm{\langle \eps\nabla_v\rangle^{k_{d}-1}  G}_{T, 2\ell+ k_{d}+p_{v}+2, p_{1}+k_{d}},
\end{multline}
hence the lemma.
\end{proof}

We can conclude the argument for the low $\eta$ regime. By Lemma~\ref{lem-paH}, choosing $p_1=d+1$, $p_{v}= k_{d}+2$, we have
\begin{equation*}
\partial^{\alpha}_{\eta}\partial^{\beta}_{x}H^-_{t,s}(x,\eta)=  \int_{v} \int_{\xi}e^{-ix\cdot \eta}  e^{i \Phi_{t,s}(z,\xi)} d_{t,s}^{\alpha,\beta}(z,\xi,\eta)\, d\xi dv  
\end{equation*}
and we apply~\eqref{EstD} to directly integrate with respect to $v$ and $\xi_v$ and get
\begin{multline}
\label{eq:Hconclusion-low}
\vert \partial^{\alpha}_{\eta}\partial^{\beta}_{x}H^-_{t,s}(x,\eta)\vert
\\
 \lesssim \left(\int_{\xi_x} \frac{\langle  \xi_x\rangle }{  \langle \xi_x-\eta \rangle^{\ell}\langle (t-s)\xi_x\rangle^2} d\xi_x \right) \norm{b}_{L^\infty_{T}W^{d+3k_{d}+3,\infty}_{z,\xi}}  
 \norm{\langle \eps\nabla_v\rangle^{k_{d}-1}  G}_{T,  2 \ell + 2 k_{d} + 4,d+k_{d}+1}, 
 \end{multline}
hence the claimed estimate.

\subsubsection{The high $\eta$ regime}

We now study the term $H^+_{t,s}(x,\eta) = H_{t,s}(x,\eta) \left[1- \chi\left(\eps |\eta|\right)\right]$, which corresponds to the region $\{\eps|\eta|\geq 1\}$. The treatment of this regime is more technically involved and we additionally  need to distinguish between a low and high velocity regime. 
As in this regime $\eta \neq 0$, we can define coordinates adapted to $\eta$ by setting for all $y \in \R^d$,
\begin{align*}
y_\parallel = \left( \frac{\eta}{|\eta|} \cdot y\right)  \frac{\eta}{|\eta|} , \quad y_\perp = y- y_\parallel.
\end{align*}
Let us also denote
$$
\nabla_\parallel =  \frac{\eta}{|\eta|} \left(\frac{\eta}{|\eta|} \cdot \nabla_{\xi_v}\right) , \quad \nabla_\perp = \nabla_{\xi_v} - \nabla_\parallel.
$$
Setting $W^+(\eps\eta)= \widehat{V}(\eps \eta) \left[1- \chi\left(\eps |\eta|\right)\right] $ and 
$$
{b}^-_{t,s}(z,\xi) = {b}_{t,s}(z,\xi)  \chi\left( \frac{|v_\parallel|}{\eps |\eta|}\right), \quad {b}^+_{t,s}(z,\xi) = {b}_{t,s}(z,\xi) \left[1-  \chi\left( \frac{|v_\parallel|}{\eps |\eta|}\right)\right], 
$$
we write
\begin{equation}
\label{def:H-high}
\begin{aligned}
H^+_{t,s}(x,\eta)  
 &= \int_{v}\int_{\xi} e^{-ix\cdot \eta} e^{i    \Phi_{t,s} (z,\xi)} {b}^-_{t,s}(z,\xi)\mathcal{F}_{x,v}G_s(\xi_x-\eta,\xi_v) \frac{1}{\eps} \sin\left(\frac{\eps \xi_v\cdot \eta} {2}\right)   W^+(\eps \eta)  \, d\xi dv   \\
 &+ \int_{v}\int_{\xi} e^{-ix\cdot \eta} e^{i    \Phi_{t,s} (z,\xi)} {b}^+_{t,s}(z,\xi)\mathcal{F}_{x,v}G_s(\xi_x-\eta,\xi_v) \frac{1}{\eps} \sin\left(\frac{\eps \xi_v\cdot \eta} {2}\right)   W^+(\eps \eta) \, d\xi dv   \\
 &=: H^{+,-}_{t,s}(x,\eta)+ H^{+,+}_{t,s}(x,\eta).
\end{aligned}
\end{equation}
In the following, we will systematically use that since $\chi'(z)=0$ for $z \in [0,1]\cup (2,+\infty)$, for all multi-indices $\alpha, \beta, \gamma$,
$$
\left|\pa^\alpha_v \pa^\beta_\eta  \chi\left( \frac{|v_\parallel|}{\eps |\eta|} \right) \pa^\gamma_\eta W^+(\eps \eta)\right| \lesssim 1,
$$
and that derivatives of $ \eta/| \eta|$ are uniformly bounded  on the support of $W^+$.

We can then define the two vector fields that we shall use instead of $\mathcal{L}_{\xi_{v}}$,  
$$
\mathcal{L}_\parallel = \frac{\lambda    - i  \nabla_\parallel \Phi_{t,s}(z,\xi)\cdot \nabla_\parallel}{ \lambda   +|\nabla_\parallel \Phi_{t,s}(z,\xi)|^2}, \quad \mathcal{L}_\perp = \frac{\lambda    - i  \nabla_\perp \Phi_{t,s}(z,\xi)\cdot \nabla_\perp}{ \lambda    +|\nabla_\perp \Phi_{t,s}(z,\xi)|^2},
$$
where $\lambda>0$ is a large enough constant, independent of $\eps$ such that, according to~\eqref{eq:boundPsi},
\begin{align}
\label{eq:below-parallel}
 \lambda   +|\nabla_\parallel \Phi_{t,s}(z,\xi)|^2 \geq C  \langle v_\parallel\rangle^2, \\
\label{eq:below-perp}
\lambda   +|\nabla_\perp \Phi_{t,s}(z,\xi)|^2 \geq C  \langle v_\perp \rangle^2.
\end{align}
By construction, we have $\mathcal{L}_\parallel   e^{i   \Phi_{t,s} (z,\xi)}  = \mathcal{L}_\perp   e^{ i    \Phi_{t,s} (z,\xi)} =   e^{i   \Phi_{t,s} (z,\xi)}$.

We shall also use the vector field $\mathcal{L}_v$ as defined in~\eqref{eq:Lv} and the vector fields $X_{\eta}$, $X_{x}$ defined
in \eqref{defXetaj}, \eqref{defXxj}.

\bigskip

\noindent $\bullet$ {\bf Study of $H^{+,-}$: the high $\eta$, low $v$ regime.}
In this regime we only need the vector field $\mathcal{L}_\perp$ and not $\mathcal{L}_\parallel$. Since $\nabla_\perp \sin\left(\frac{\eps \xi_v\cdot \eta} {2}\right) =0$,  we will not  get powers of $\eps \eta$ to absorb.
 Let $p_\perp>0$ be an integer to be fixed later.
Arguing as in~\eqref{expanL-xiv}, it follows  that
\begin{equation}\label{expanL-perp} 
\left(\mathcal{L}_\perp^T\right)^{p_\perp} = \sum_{\vert \alpha\vert  \leq p_\perp} c_{\perp, t,s}^{\alpha}(z,\xi)\partial_{{\xi_v}_\perp}^{\alpha} ,
\end{equation}
where we set $(\partial_{{\xi_v}_\perp})_{j}= (\nabla_{\perp})_{j}, $ $ j \in \llbracket 1,  d \rrbracket$.
By using the lower bound~\eqref{eq:below-perp} and the upper bound
\begin{equation}
\label{eq:upper-phi-perp}
\sup_{0\leq |\alpha_z| +  |\alpha_\xi|\leq  2k_{d}+p_{v}} \sup_{0\leq |\alpha| \leq p_\perp} \left| \pa^{\alpha_z}_z \pa^{\alpha_\xi}_\xi \pa^{\alpha}_{{\xi_v}} \nabla_{{\xi_{v}}_\perp} \Phi_{t,s} \right| \lesssim  \langle v_\perp\rangle,
\end{equation}
if the assumption ($\boldsymbol{A}_{2k_{d} +p_{\perp}+ p_{v}}$) is matched, 
 we then have that the functions $c_{\perp,t,s}^{\alpha}$  satisfy the estimate
\begin{equation}\label{EstC-perp}
\sup_{0\leq |\alpha_z| +  |\alpha_\xi|\leq 2 k_{d}+p_{v}} \vert  \pa^{\alpha_z}_z \pa^{\alpha_\xi}_\xi  c_{\perp_,t,s}^{\alpha}(z,\xi)\vert \lesssim 
  \frac{1} {  \langle v_\perp \rangle^{p_\perp}}.
\end{equation}
We can then  write that
\begin{align*}
 H^{+,-}_{t,s}(x,\eta) & = 
 \int_{v}\int_{\xi} e^{-ix\cdot \eta} e^{\frac{i}{\varepsilon}   \Psi_{t,s} (z,\xi)} \left(\mathcal{L}_\perp^T\right)^{p_\perp} \mathcal{A}_{t,s}^{+,-}\,   d\xi dv,  \\
   \mathcal{A}_{t,s}^{+,-}  (z, \xi, \eta) &:=    {b}^-_{t,s}(z,\xi)\mathcal{F}_{x,v}G_s(\xi_x-\eta,\xi_v) \frac{1}{\eps} \sin\left(\frac{\eps \xi_v\cdot \eta} {2}\right)  W^+(\eps \eta).
\end{align*}

In the next lemma we establish a  result analogous to Lemma~\ref{lem-paH} in this new situation.
\begin{lem}
\label{lem-paH-highlow}
Let $\ell>d$,  $p_v, p_\perp \geq 2k_{d}$, if the assumption  ($\boldsymbol{A}_{2k_{d} +p_{\perp}+ p_{v}}$) holds, then  for every $0 \leq \vert \alpha\vert,\vert \beta\vert \leq k_{d}$,  $\partial^{\alpha}_{\eta}\partial^{\beta}_{x}H^{+,-}_{t,s}$ can be put under the  form 
\begin{equation}\label{Hder2-highlow}
\partial^{\alpha}_{\eta}\partial^{\beta}_{x}H^{+,-}_{t,s}(x,\eta)=  \int_{v} \int_{\xi}e^{-ix\cdot \eta}  e^{ i  \Phi_{t,s}(z,\xi)} d_{t,s}^{\alpha,\beta}(z,\xi,\eta) \, d\xi dv,
\end{equation}
where $d^{\alpha,\beta}$ satisfies
\begin{multline}
\label{EstD-highlow}
| d_{t,s}^{\alpha,\beta}(z,\xi,\eta)| \lesssim  
 \frac{\mathds{1}_{|v_\parallel|\leq \sqrt{2} \eps |\eta|}}{\eps} \frac{1}
{ \langle \xi_v\rangle^\ell \langle \xi_x-\eta\rangle^{\ell+1} \langle v_\perp\rangle^{p_\perp} \langle (t-s) \xi_x  \rangle^{p_{v}- 2 k_{d}}}  \\
\times \norm{b_{t,s}}_{L^\infty_{T}W^{p_\perp+2k_{d}+p_{v},\infty}_{z,\xi}} \norm{\langle \eps\nabla_v\rangle^{k_{d}}  \langle \eps \nabla_{x} \rangle^{k_{d}} G_{t,s}}_{T,  2\ell+p_v + k_d +1, p_\perp+k_{d}}.
\end{multline}
\end{lem}

\begin{proof}[Proof of Lemma~\ref{lem-paH-highlow}]
As previously, by using   $\mathcal{L}_v$ defined in~\eqref{eq:Lv} and the vector fields $X_{\eta}$, $X_{x}$ defined
in \eqref{defXetaj}, \eqref{defXxj}, we can  write 
$$
\partial^{\alpha}_{\eta}\partial^{\beta}_{x} H^-_{t,s}(x,\eta) = 
  \int_{v}\int_{\xi} e^{-ix\cdot \eta} e^{i    \Phi_{t,s} (z,\xi)} 
  (\mathcal{L}_v^T)^{p_{v}}  X_{x}^\beta X_{\eta}^\alpha  (\mathcal{L}_{\perp}^T)^{p_{\perp}} \mathcal{A}_{t,s}^{+, -}\, d \xi dv.  
$$
and we  set 
\begin{equation*}
d^{\alpha,\beta}_{t,s}(z,\xi,\eta) =   (\mathcal{L}_v^T)^{p_{v}}   X_{\eta}^\alpha X_{x}^\beta  (\mathcal{L}_{\perp}^T)^{p_{\perp}} \mathcal{A}_{t,s}^{+,-},
\end{equation*}
to get the form \eqref{Hder2-highlow}.
By using the expansion \eqref{expanL-perp} and again the expansions \eqref{XetaXxexpansion}, \eqref{expanL}
together with  \eqref{EstC-perp}, \eqref{estim-eij} and \eqref{EstC}, we get that
\begin{multline}\label{MajH1:plusmoins}
d_{t,s}^{\alpha,\beta}(z,\xi,\eta)=  \\ \sum_{\substack{0\leq\vert\gamma\vert\leq k_{d} \\ 0\leq \vert\sigma\vert+\vert\rho\vert + \vert \mu \vert \leq k_{d}}}
 \sum_{0 \leq | \alpha '| \leq p_v}
\sum_{ 0\leq |\beta'| \leq p_\perp } f_{t,s}^{\alpha,\beta, \gamma,\sigma,\rho,\mu, \alpha',\beta'}(z,\xi,\eta)  \partial^{\gamma}_x \,
\partial^{\sigma}_{\xi_x}\partial^{\rho}_\eta \left( (t-s) \partial_{\xi_{v}}\right)^\mu  \partial_{v}^{\alpha'}\partial_{{\xi_v}_\perp}^{\beta'} 
 \mathcal{A}^{+, -}_{t,s},
\end{multline}
where  the coefficients satisfy
\begin{equation}\label{MajH2:plusmoins}
\bigg\vert f_{t,s}^{\alpha,\beta,\gamma,\sigma,\rho,\mu,\alpha'\beta'}(z,\xi,\eta)  \bigg  \vert  \lesssim
 \frac{ \langle \xi_v\rangle^{p_{v}}
\langle \xi_x-\eta\rangle^{k_{d}}}{ \langle  v_{\perp}  \rangle^{p_\perp} \left( \langle \xi_{v} \rangle + \langle  (t-s) \xi_x  \rangle\right)^{p_{v}-k_{d}}} \lesssim  \frac{ \langle \xi_v\rangle^{p_{v}}
\langle \xi_x-\eta\rangle^{k_{d}}}{ \langle  v_{\perp}  \rangle^{p_\perp}  \langle  (t-s) \xi_x  \rangle^{p_{v}-k_{d}}} .
\end{equation}
By  using that 
$\partial_{{\xi_v}_\perp}  \sin\left(\frac{\eps \xi_v\cdot \eta} {2}\right)=0$, that
$$  \left| \partial_{\eta}^\rho \left(( t-s)  \partial_{\xi_{v}}\right)^\mu
 \left( \sin{ (\eps \xi_{v} \cdot \eta)} \right) \right|
  \lesssim  \left(  \langle \eps \xi_{v}\rangle^{|\rho|} \langle \eps (t-s) \eta \rangle^{| \mu|} \right)
   \lesssim   \left(  \langle \eps \xi_{v}\rangle^{|\rho|} \langle \eps (\xi_{x}- \eta) \rangle^{| \mu|} \langle (t-s) \xi_{x}\rangle^{|\mu|}  \right),$$
 and by  recalling that we are in the low velocity regime, we  have
\begin{multline}\label{MajH3:plusmoins}
\bigg\vert \partial^{\gamma}_x\,\partial^{\sigma}_{\xi_x}\,\partial^{\rho}_\eta \,\left( (t-s) \partial_{\xi_{v}}\right)^\mu\, \partial_v^{\alpha'}\,\partial_{{\xi_v}_\perp}^{\beta'}\,\mathcal{A}^{+,-}_{t,s}(z,\xi) \bigg \vert 
\lesssim   \frac{\mathds{1}_{|v_\parallel|\leq 2 \eps |\eta|}}{\eps} \langle (t-s) \xi_{x} \rangle^{k_{d}}
		  \\
		  \times
		  \sum_{0\leq |\alpha''| + |\beta''|\leq k_{d}+ p_\perp}  \langle \eps\xi_v\rangle^{k_{d}}  \langle \eps(\xi_{x}- \eta)\rangle^{k_{d}} |\pa^{\alpha''}_{\xi_x} \pa^{\beta''}_{\xi_v}  \mathcal{F}_{x,v}  G_{t,s}(\xi_x-\eta,\xi_v)|   
		  \norm{b}_{L^\infty_{T}W^{p_\perp+p_{v}+ 2k_{d},\infty}_{z,\xi}}. 
\end{multline}
Finally combining \eqref{MajH1:plusmoins}, \eqref{MajH2:plusmoins} and \eqref{MajH3:plusmoins}, we thus  obtain that
\begin{multline*}
| d_{t,s}^{\alpha,\beta}(z,\xi,\eta)| \lesssim  
 \frac{\mathds{1}_{|v_\parallel|\leq 2 \eps |\eta|}}{\eps} \frac{1}
{ \langle \xi_v\rangle^\ell \langle \xi_x-\eta\rangle^{\ell+1} \langle v_\perp\rangle^{p_\perp} \langle (t-s) \xi_x  \rangle^{p_{v}- 2 k_{d}}}  \\
\times  \norm{b_{t,s}}_{L^\infty_{T}W^{p_\perp+p_{v}+ 2k_{d},\infty}_{z,\xi}} \norm{\langle \eps\nabla_v\rangle^{k_{d}}  \langle \eps\nabla_x\rangle^{k_{d}} G}_{T,   2\ell + p_v +k_d  +1 ,p_\perp+k_{d}},
\end{multline*}
hence the result.
\end{proof}

\bigskip

To conclude the estimate for $H^{+, -}$, 
we  choose $p_\perp =d$, $p_{v}= 2 k_{d} + 2$  and use the previous Lemma. We observe that the integral in velocity contributes as
\begin{align*}
\int_v  \mathds{1}_{|v_\parallel|\leq 2 \eps |\eta|} \frac{1}{\langle v_\perp\rangle^{d} } \, dv 
 =   \int_{v_\parallel}  \mathds{1}_{|v_\parallel|\leq 2   \eps |\eta|}  \, dv_\parallel \int_{v_\perp} \frac{1}{\langle v_\perp\rangle^{d} }  \, dv_\perp \lesssim \eps |\eta| \leq \eps  \langle \xi_x \rangle \langle \eta-\xi_x\rangle.
\end{align*}
Therefore we can conclude in the high $\eta$, low $v$ regime. By Lemma~\ref{lem-paH-highlow}, estimate~\eqref{EstD-highlow} and the previous estimate, we get
\begin{multline*}
\vert \partial^{\alpha}_{\eta}\partial^{\beta}_{x}H_{t,s}(x,\eta)\vert
 \lesssim \int_{\xi_x} \frac{\langle  \xi_x\rangle }{  \langle \xi_x-\eta \rangle^{\ell}\langle (t-s)\xi_x\rangle^2} \, d\xi_x \\
 \times  \norm{b_{t,s}}_{L^\infty_{T}W^{4k_{d}+d+2,\infty}_{z,\xi}} \norm{\langle \eps\nabla_v\rangle^{k_{d}}  \langle \eps\nabla_x\rangle^{k_{d}} G}_{T,  2r+ 3k_{d} + 3,d+k_{d}}.
\end{multline*}

\bigskip

\noindent $\bullet$ {\bf Study of $H^{+,+}$: the high $\eta$, high $v$ regime.}
In the high $\eta$, high $v$ regime we shall also need to use   the operator $\mathcal{L}_\parallel$ which involves derivatives with respect to ${\xi_v}_\parallel$ in order to get integrability in $v_\parallel$ and to absorb the prefactor $\eps^{-1}$.

We shall use $\left(\mathcal{L}_\parallel^T\right)^{2}$, as previously, we can expand
\begin{equation}\label{expanL-para} 
\left(\mathcal{L}_\parallel^T\right)^{2} = \sum_{\vert \alpha\vert  \leq 2} c_{\parallel, t,s}^{\alpha}(z,\xi)\partial_{{\xi_v}_\parallel}^{\alpha} ,
\end{equation}
where 
thanks to  the lower bound~\eqref{eq:below-parallel} and the upper bound
\begin{equation}
\label{eq:upper-phi-parallel}
\sup_{0\leq |\alpha_z| +  |\alpha_\xi|\leq 2k_{d}+p_{v}} \sup_{0\leq |\alpha| \leq  2} \left|  \pa^{\alpha_z}_z \pa^{\alpha_\xi}_\xi  \pa^{\alpha}_{{\xi_v}} \nabla_{{\xi_v}_\parallel} \Phi_{t,s} \right| \lesssim \langle v_\parallel \rangle,
\end{equation}

the functions $c_{\parallel,t,s}^{\alpha}$  satisfy the estimate
\begin{equation}\label{EstC-para}
\sup_{0\leq |\alpha_z| +  |\alpha_\xi|\leq  2k_{d}+p_{v}+p_{\perp}  } \vert \pa^{\alpha_z}_z \pa^{\alpha_\xi}_\xi   c_{\parallel,t,s}^{\alpha}(z,\xi)\vert  \lesssim   \frac{1} {  \langle v_\parallel \rangle^{2}}
\end{equation}
when the assumption 
($\boldsymbol{A}_{2k_{d} + 2 + p_{v}+p_{\perp}}$) is matched.

The analogue of Lemmas~\ref{lem-paH} and \ref{lem-paH-highlow} reads in this case as follows.
\begin{lem}
\label{lem-paH-highhigh}
Let $\ell>d$, $p_{v},p_\perp \geq 2 k_{d}$, if the assumption  ($\boldsymbol{A}_{2k_{d} +p_{\perp}+ p_{v}+ 2}$) holds,  then for every $0 \leq \vert \alpha\vert,\vert \beta\vert \leq k_d$, $\partial^{\alpha}_{\eta}\partial^{\beta}_{x}H^{+,+}_{t,s}$ can be put under the  form 
\begin{equation}\label{Hder2-highhigh}
\partial^{\alpha}_{\eta}\partial^{\beta}_{x}H^{+,+}_{t,s}(x,\eta)=  \int_{v} \int_{\xi}e^{-ix\cdot \eta}  e^{i \Phi_{t,s}(z,\xi)}   d_{t,s}^{\alpha,\beta}(z,\xi,\eta) \, d\xi dv  
\end{equation}
where $d^{\alpha,\beta}$ satisfies
\begin{multline}
\label{EstD-highhigh}
| d_{t,s}^{\alpha,\beta}(z,\xi,\eta)| \lesssim  
 \mathds{1}_{|v_\parallel| \geq  \eps |\eta|} \frac{\langle \xi_{x} \rangle \langle \eps \eta \rangle}
{ \langle \xi_v\rangle^\ell \langle \xi_x-\eta\rangle^{\ell} \langle v_\perp\rangle^{p_\perp} \langle v_{\parallel} \rangle^2 \langle (t-s) \xi_x  \rangle^{p_{v}- 2 k_{d}}}  \\
\times  \norm{b_{t,s}}_{L^\infty_{T}W^{p_\perp+p_{v}+ 2k_{d} + 2 ,\infty}_{z,\xi}} \norm{\langle \eps\nabla_v\rangle^{k_{d}}  \langle \eps\nabla_x\rangle^{k_{d}} G}_{T,  2 \ell + p_v +k_d + 2  ,p_\perp+k_{d}+2}.
\end{multline}
\end{lem}

\begin{proof}
By using again    $\mathcal{L}_v$  as defined in~\eqref{eq:Lv} and the vector fields $X_{\eta}$, $X_{x}$ as defined
in \eqref{defXetaj}--\eqref{defXxj}, we can write 
$$
\partial^{\alpha}_{\eta}\partial^{\beta}_{x} H^-_{t,s}(x,\eta) = 
  \int_{v}\int_{\xi} e^{-ix\cdot \eta} e^{i    \Phi_{t,s} (z,\xi)} 
  (\mathcal{L}_v^T)^{p_{v}}  X_{x}^\beta X_{\eta}^\alpha(  \mathcal{L}_{\perp}^T)^{p_{\perp}} (\mathcal{L}_{\parallel}^T)^{2} \mathcal{A}_{t,s}^{+, +}\, d \xi dv.  
$$
and we  set 
\begin{equation*}
d^{\alpha,\beta}_{t,s}(z,\xi,\eta) =   (\mathcal{L}_v^T)^{p_{v}}  X_{\eta}^\alpha   X_{x}^\beta (\mathcal{L}_{\perp}^T)^{p_{\perp}}
 (\mathcal{L}_{\parallel}^T)^{2}
 \mathcal{A}_{t,s}^{+,+},
\end{equation*}
where
$$  \mathcal{A}_{t,s}^{+,+}  (z, \xi, \eta) =    {b}^+_{t,s}(z,\xi)\mathcal{F}_{x,v}G_s(\xi_x-\eta,\xi_v) \frac{1}{\eps} \sin\left(\frac{\eps \xi_v\cdot \eta} {2}\right)  W^+(\eps \eta),$$
to get the form \eqref{Hder2-highhigh}.

By using the expansion \eqref{expanL-para} and again the expansions  \eqref{expanL-perp} ,  \eqref{XetaXxexpansion}, \eqref{expanL}
together with  \eqref{EstC-para},   \eqref{EstC-perp}, \eqref{estim-eij} and \eqref{EstC}, we get that
\begin{multline*}d_{t,s}^{\alpha,\beta}(z,\xi,\eta)=  \\ \sum_{\substack{0\leq\vert\gamma\vert\leq k_{d} \\ 0\leq \vert\sigma\vert+\vert\rho\vert + \vert \mu \vert \leq k_{d}}}
 \sum_{0 \leq | \alpha '| \leq 2}
\sum_{ \substack{ 0\leq |\beta'| \leq p_\perp  \\ 0 \leq | \gamma '| \leq 2}} f_{t,s}^{\alpha,\beta, \gamma,\sigma,\rho,\mu,\alpha',\beta', \gamma'}(z,\xi,\eta)  \partial^{\gamma}_x \,
\partial^{\sigma}_{\xi_x}\partial^{\rho}_\eta \left( (t-s) \partial_{\xi_{v}}\right)^\mu  \partial_{v}^{\alpha'}\partial_{{\xi_v}_\perp}^{\beta'}
 \partial_{{\xi_{v}}_{\parallel}}^{\gamma'} 
 \mathcal{A}^{+, +}_{t,s},
\end{multline*}
where  the coefficients satisfy
\begin{equation}
\label{cellelaaussi}
\bigg\vert f_{t,s}^{\alpha,\beta,\gamma,\sigma,\rho,\mu,\alpha'\beta', \gamma'}(z,\xi,\eta)  \bigg  \vert  \lesssim
 \frac{ \langle \xi_v\rangle^{p_{v}}
\langle \xi_x-\eta\rangle^{k_{d}}}{ \langle  v_{\perp}  \rangle^{p_\perp} \langle v_{\parallel}\rangle^2 \left( \langle \xi_{v} \rangle + \langle  (t-s) \xi_x  \rangle\right)^{p_{v}-k_{d}}} \lesssim  \frac{ \langle \xi_v\rangle^{p_{v}}
\langle \xi_x-\eta\rangle^{k_{d}}}{ \langle  v_{\perp}  \rangle^{p_\perp} \langle v_{\parallel}\rangle^2  \langle  (t-s) \xi_x  \rangle^{p_{v}-k_{d}}}.
\end{equation}
Arguing as in the other cases, we can estimate
\begin{multline*}
\bigg\vert \partial^{\gamma}_x\,\partial^{\sigma}_{\xi_x}\,\partial^{\rho}_\eta \,\left( (t-s) \partial_{\xi_{v}}\right)^\mu\, \partial_v^{\alpha'}\,\partial_{{\xi_v}_\perp}^{\beta'}\, \mathcal{A}^{+,+}_{t,s}(z, \xi, \eta)\bigg \vert \\
\lesssim   \mathds{1}_{|v_\parallel|\geq  \eps |\eta|} \langle (t-s) \xi_{x} \rangle^{k_{d}} \langle \eta \rangle \langle \eps \eta \rangle \langle \xi_{v} \rangle   \norm{b}_{L^\infty_{T}W^{p_\perp+p_{v}+ 2k_{d}+ 2,\infty}_{z,\xi}} \cdot  \\
		 \sum_{0\leq |\alpha''| + |\beta''|\leq k_{d}+ p_\perp+2}  \langle \eps\xi_v\rangle^{k_{d}}  \langle \eps(\xi_{x}- \eta)\rangle^{k_{d}} |\pa^{\alpha''}_{\xi_x} \pa^{\beta''}_{\xi_v}  \mathcal{F}_{x,v}  G_{t,s}(\xi_x-\eta,\xi_v)|.  		  
\end{multline*}
We have used here that $\partial_{\xi_{\perp}} (\xi_{v} \cdot \eta)= 0$. Moreover, note that if no derivatives hit the $\sin$, 
 the inequality $| \sin x| \leq |x|$ allows to compensate  the prefactor $\eps^{-1}$.
 Otherwise, whenever a derivative hits the $\sin$, we directly gain a factor $\eps$.
 A crucial observation is that we have in the end at most one power of $|\eta|$ and one power of $\eps | \eta|$ because we use at most
  two derivatives $\partial_{\xi_{\parallel}}.$ By using $\langle \eta \rangle \leq \langle \xi_{x} \rangle \langle\xi_{x}- \eta\rangle$, 
  we end up with
  \begin{multline*}
| d_{t,s}^{\alpha,\beta}(z,\xi,\eta)| \lesssim  
 \mathds{1}_{|v_\parallel| \geq  \eps |\eta|} \frac{\langle \xi_{x} \rangle \langle \eps \eta \rangle}
{ \langle \xi_v\rangle^\ell \langle \xi_x-\eta\rangle^{\ell} \langle v_\perp\rangle^{p_\perp} \langle v_{\parallel} \rangle^2 \langle (t-s) \xi_x  \rangle^{p_{v}- 2 k_{d}}}  \\
\times  \norm{b}_{L^\infty_{T}W^{p_\perp+p_{v}+ 2k_{d} + 2 ,\infty}_{z,\xi}} \norm{\langle \eps\nabla_v\rangle^{k_{d}}  \langle \eps\nabla_x\rangle^{k_{d}} G}_{T,   2 \ell + p_v +k_d + 2  ,p_\perp+k_{d}+2}.
\end{multline*}

\end{proof}

\bigskip

To conclude the proof  for $H^{+, +}$, we choose $p_{v}= 2 k_{d} + 2$, $p_{\perp}=d$ in \eqref{EstD-highhigh} and we observe that
$$ \langle \eps \eta \rangle \int_{| v_{\parallel}| \geq \eps | \eta | }  { 1 \over \langle v_{\parallel} \rangle^2} \, d v_{\parallel}
 \lesssim 1.$$
This yields 
\begin{multline*}
\vert \partial^{\alpha}_{\eta}\partial^{\beta}_{x}H_{t,s}(x,\eta)\vert  \\
\lesssim \int_{\xi_x} \frac{\langle  \xi_x\rangle }{  \langle \xi_x-\eta \rangle^{\ell}\langle (t-s)\xi_x\rangle^2} \, d \xi \, 
 \norm{b}_{L^\infty_{T}W^{4 k_{d} + d + 4 ,\infty}_{z,\xi}} \norm{\langle \eps\nabla_v\rangle^{k_{d}}  \langle \eps\nabla_x\rangle^{k_{d}} G}_{T, 
  2\ell + 3 k_{d} + 4,  k_{d} + d + 2 }.
\end{multline*}
This finally ends the proof of Proposition \ref{prop-finH}.

\subsection {Improved variants of Theorem \ref{M-Pro} }\label{GenMoySub} 

We can first improve Theorem \ref{M-Pro} by allowing some polynomial growth of $b$ in $\xi$. Namely,  instead of the boundedness of
$ \|b\|_{L^\infty_{T} W_{z,\xi}^{d + 4k_{d} + 4, \infty} }$, we can require the boundedness  of  
$ \|b/\left( \langle \xi_{v} \rangle + \langle (t-s) \xi_{x} \rangle\right)^q\|_{L^\infty_{T} W_{z,\xi}^{p, \infty} }$
 for any $q \in \mathbb{N}$, when  $p$  is accordingly taken sufficiently large. 
 \begin{thm}\label{M-Pro2}
 Let  $q \in \mathbb{N}$. For every $T_{0}>0$, 
there exists $C_{0}>0$ such that for every $T \in [0, T_{0}]$,  if the assumption  ($\boldsymbol{A}_{4k_{d} + d + 4}$) holds, 
 we have  for every  $\eps \in (0,1)$ that
\begin{multline*}
\norm{\mathcal{U}_{[ \Phi, b, G]}}_{\mathscr{L}(L^2(0,T;L^2(\R^d)))}\leq  \\ C_{0}   \norm{{b \over \left( \langle \xi_{v} \rangle + \langle (t-s) \xi_{x} \rangle\right)^q}}_{L^\infty_{T} W^{q+ d+4k_{d}+4,\infty}_{z,\xi}}  
 \norm{\langle \eps\nabla_x \rangle^{k_{d}} \langle \eps\nabla_v \rangle^{k_{d}} G}_{T, q+ 3k_{d}+ 2 d + 6, k_{d}+ d+2}.
 \end{multline*}
 \end{thm}
For the proof, it suffices to notice that  thanks to the intermediate estimates in  \eqref{cellelaaussi}, \eqref{MajH2:plusmoins}
and \eqref{MajH2}, we can absorb  the additional powers of  $\langle \xi_{v} + \langle (t-s) \xi_{x} \rangle$
if we replace $p_{v}$ by $p_{v}+ q$. This directly yields the result.

\bigskip

Finally, in the more specific case when the phase $\Phi_{t,s}(z, \xi)$ is given by
$\Phi_{t,s}(z, \xi) = \Psi_{t,s}(z, \eps \xi)/\eps$, we can also extend the above  continuity on $L^2(0, T; L^2)$
 to a continuity result on  $L^2(0, T; H^0_{r})$.
 
 \begin{thm}\label{M-Pro3}
 Let  $q \in \mathbb{N}$, $r\in \mathbb{N}$ and assume that
$$\Phi_{t,s}(z, \xi) = { \Psi_{t,s}(z, \eps \xi) \over \eps }$$
 For every $T_{0}>0$, 
there exists $C_{0}>0$ such that for every $T \in [0, T_{0}]$,  if the assumption  ($\boldsymbol{A}_{4k_{d} + d + 4 +r}$) holds 
and if moreover
\begin{equation}
\label{lhypotheseenplus}
 \sup_{t, \, s \in [0, T]}
\underset{{\substack{0\leq\vert \alpha\vert + \vert \beta\vert \leq q+d+ 4k_{d}+ 3 + r
 }}}{\sup} \left|\partial^{\alpha}_z\partial^{\beta}_{\xi}  \left( \nabla_{x} \Psi_{t,s}(z,\xi) - \xi_{x}  \right)\right|\leq  1, 
\end{equation}
 then, 
 we have  for every  $\eps \in (0,1)$ that
\begin{multline*}
\norm{\mathcal{U}_{[ \Phi, b, G]}}_{\mathscr{L}(L^2(0,T;H^0_{r}))}\leq  \\ C_{0}   \norm{{ \langle \eps \nabla_x \rangle^r b \over \left( \langle \xi_{v} \rangle + \langle (t-s) \xi_{x} \rangle\right)^q}}_{L^\infty_{T} W^{q+ d+4k_{d}+4,\infty}_{z,\xi}}  
 \norm{\langle \eps\nabla_x \rangle^{k_{d}+r} \langle \eps\nabla_v \rangle^{k_{d}} G}_{T, q+ 3k_{d}+ 2 d + 6, k_{d}+ d+2}.
 \end{multline*}
 \end{thm}
 
 \begin{proof}
 We observe that
 \begin{multline*} \eps \partial_{x_{j}} \mathcal{U}_{[ \Phi, b, G]}(\rho)=\\ 
 \frac{1}{(2\pi )^{2d}} \int_{v}\int_{0}^t\int_{\xi} \int_{y} e^{i \Psi_{t,s} (z, \eps\xi) \over \eps } \left( \eps \partial_{x_j}b_{t,s}(z,\xi)
  +  i \partial_{x_{j}} \Psi_{t,s}(z, \eps \xi)  b_{t,s}(z, \xi) \right) \widehat{B[\varrho,G_{t,s}]}(\xi) d\xi ds dv
  \end{multline*}
 Next, by using \eqref{defBrecalled}, we also  have that
 $$ i \eps \xi_{x_{j}} \widehat{ B[\rho,G_{t, s}] }(\xi)=  \widehat{ B[ \eps \partial_{x_{j}}\rho,G_{t, s}] }(\xi)
  +  \widehat{ B[ \rho,\eps \partial_{x_{j}} G_{t, s}] }(\xi), $$
  therefore, we can write
  $$  \eps \partial_{x_{j}} \mathcal{U}_{[ \Phi, b, G]}(\rho)
  =  \mathcal{U}_{[ \Phi, b, G]}( \eps \partial_{x_{j}}\rho) +  \mathcal{U}_{[ \Phi,   b^{e_{j}}, G]}(\rho)
   +  \mathcal{U}_{[ \Phi, b,  \eps \partial_{x} G]}(\rho) $$
   where 
   $$ b^{e_{j}}_{t,s}(z, \xi)= \eps \partial_{x_{j}} b_{t,s}(z, \xi) +  i \left( \partial_{x_{j}} \Psi_{t,s} (z, \eps \xi) - \eps \xi_{x_{j}}\right) b_{t,s}(z, \xi).$$
   The result then follows by iterating this identity and by applying Theorem \ref{M-Pro2}. Note that  $b^{e_{j}}$ and its derivatives 
   can be controlled by using the assumption \eqref{lhypotheseenplus}.

 \end{proof}

\section{Higher estimates for the density}
\label{sec:penrose}

We move on to the last part of the proof. From now on, we always consider positive times $T \leq \min(T_\eps,T(M))$ so that Propositions~\ref{prop-HJ1}, \ref{EstB} and \ref{ConsFIO} apply. We start from the equation \eqref{LinT} for the solution $\mathrm{F}$ to the extended Wigner system~\eqref{IJ} and take the integral in $v$, by setting
\begin{equation}
\label{rhoFdef}
\rho_{\mathrm{F}}= \int_{v} \mathrm{F}\, dv, 
\end{equation}
we obtain
\begin{equation}\label{EqPrin}
\rho_{\mathrm{F}} = - \int_{v} \int_{0}^t U_{t,s} B[{\rho_{\mathrm{F}}}, f]\, ds dv + \int_v U_{t,0}  \mathrm{F}^0 \, dv  +  \int_v \int_0^t U_{t,s}  \mathcal{R}(s) dsdv.
\end{equation}

The philosophy will be to  simplify as much as possible~\eqref{EqPrin}, using the machinery developed in the previous sections. This will allow to reach a scalar semiclassical pseudodifferential equation, which we shall invert using the quantum Penrose stability condition.

Recalling that $\mathrm{F}$ is related to the solution  $f$ of the Wigner equation by the formula $\mathrm{F} = (\partial^{\alpha(i)} f)_{i \in \llbracket1, N_m\rrbracket}$, the outcome of this section will be
\begin{prop}
\label{prop:rho-final}
For all $T\in [0,\min (T(M),T_\eps)]$,
\begin{equation}
\label{estRho}
\norm{\rho}_{L^2(0,T; H^m_r)}\leq (T^{1/2} + \eps)  \Lambda( c_{0}^{-1},  \|f^0\|_{\H^m_{r}}, T,M).
\end{equation}
\end{prop}

\subsection{First reductions}

We observe that thanks to \eqref{embedH02} and \eqref{eq:estim-m-1} we have
$$\norm{\rho}_{L^2(0,T; H^{m-1}_r)} \lesssim  \| f \|_{L^2(0, T; \H^{m-1}_{r})}
 \leq T^{1 \over 2} \Lambda(T,M) ( \|f^0\|_{\H^{m-1}_{r}} + 1)$$
 so that we only have to estimate $\| \partial^\alpha_x \rho \|_{L^2(0,T; H^{0}_r)} $ for $| \alpha |=m$ or equivalently
 $ \| \rho_{\mathrm{F}} \|_{L^2(0,T; H^{0}_r)}$ thanks to the definition \eqref{rhoFdef} and thus indeed to study \eqref{EqPrin}.

We first estimate the  terms involving the initial condition and $\mathcal{R}$ in  \eqref{EqPrin}.
\begin{lem}
\label{lem:first}
The following estimate holds for all $T$, $T \leq T_\eps$, 
$$
 \norm{\int_v \int_0^t U_{t,s} \mathcal{R} \,  dsdv}_{L^2(0,T;H^0_{r})} + \norm{\int_v U_{t,0} \mathrm{F}^0  \,  dv}_{L^2(0,T;H^0_{r})}\leq T^{1/2} \Lambda (T,M) ( 1 + \|f^0\|_{\H^m_{r}}).
$$
\end{lem}
\begin{proof}
By using successively  \eqref{embedH02} and \eqref{MajOpU}, we have
\begin{align*}
 \norm{\int_0^t \int_v  U_{t,s} \mathcal{R} \, dv ds}_{L^2(0,T;H^0_{r})}
 &\leq \norm{\int_0^t \norm{ \int_v  U_{t,s} \mathcal{R} \, dv }_{H^0_{r}}  ds}_{L^2(0, T)} \\
  &\leq C  \norm{\int_0^t \|   U_{t,s} \mathcal{R}  \|_{\H^0_{r, 0}}  ds}_{L^2(0, T)}
  \\ 
  &\leq \Lambda(T,M)  \norm{\int_0^t  \| \mathcal{R}  \|_{\H^0_{r, 0}}}_{L^2(0, T)}
   \leq  T  \Lambda(T,M)  \norm{  \mathcal{R}  }_{L^2(0, T; \H^0_{r, 0})}.
  \end{align*}
  To conclude, we  use the estimate \eqref{EstV} for the remainder.
  
  In a similar way, by using again \eqref{embedH02} and \eqref{MajOpU}, we obtain that
  $$ \norm{\int_v U_{t,0} \mathrm{F}^0  \,  dv}_{L^2(0,T;H^0_{r})}
   \lesssim \norm{  \norm{U_{t,0} \mathrm{F}^0 }_{\H^0_{r, 0}} }_{L^2(0, T)}
    \leq T^{1 \over 2} \Lambda(T, M) \| \mathrm{F}^0 \|_{\H^0_{r, 0}} \leq  T^{1 \over 2} \Lambda(T, M)  \|f^0 \|_{\H^m_{r}}     $$
  where the final estimate just follows from the definition of $\mathrm{F}^0$.

\end{proof}
In this section, a \emph{remainder} will stand for a term, generically denoted by $R=R(t,x)$, satisfying an estimate of the form
\begin{equation}
\label{defremainder}
\| R\|_{L^2(0,T;H^0_{r})} \leq (T^{1/2} + \eps) \Lambda (T,M, \|f^0\|_{\H^m_{r}})
\end{equation}
for $T \leq \min(T_\eps,T(M))$. By using this notation,  owing to Lemma~\ref{lem:first}, we can recast \eqref{EqPrin} as 
\begin{equation}\label{EqPrin2}
\rho_{\mathrm{F}} = 
- { 1 \over \eps} \int_{v} \int_{0}^t U_{t,s} B[\rho_{\mathrm{F}}, f]\, ds dv+ R(t,x),
\end{equation}
where $R$ is a remainder. 

Next,  thanks to the results of Section \ref{sec:param}, namely Lemma~\ref{lem:reduction-propagator}and Proposition~\ref{ConsFIO}, we have obtained the approximation of the propagator of $\mathscr{T}+\mathcal{M}$ as
$$
U_{t,s}=U_{t,s}^{\mathrm{FIO}}+ \eps U^{\mathrm{rem}}_{t,s}.
$$
Let us show that the term involving $\eps U_{t,s}^{\mathrm{rem}}$ in the right-hand side of~\eqref{EqPrin2} can also be seen as  a remainder.
\begin{lem}
For every $T \leq (T_\eps,T(M))$, we have the estimate
$$
\norm{\int_v \int_0^t \eps U^{\mathrm{rem}}_{t,s} (B[\rho_{\mathrm{F}}, f])  dsdv}_{L^2(0,T;H^0_{r})}\leq T^2 \Lambda (T,M).
$$
\end{lem}
\begin{proof}
By using the same arguments as in the proof of the previous Lemma, and applying Proposition~\ref{ConsFIO} together with \eqref{MajUr} in  Lemma~\ref{lem:reduction-propagator}, we obtain  
$$  \norm{\int_v \int_0^t \eps U^{\mathrm{rem}}_{t,s}B[ \rho_{\mathrm{F}},  f]  dsdv}_{L^2(0,T;H^0_{r})}
 \leq T^2 \Lambda(T, M)  \varepsilon \norm{ B[ \rho_{\mathrm{F}},  f]  }_{L^2(0, T; \H^0_{r, 0})}.$$
 Since we have by definition of $\rho_{\mathrm{F}}$ that 
 $$    \norm{ B[ \rho_{\mathrm{F}},  f]  }_{L^2(0, T; \H^0_{r, 0})}
  \leq \sup_{| \alpha |= m}    \norm{ B[  \partial^\alpha \rho,  f]  }_{L^2(0, T; \H^0_{r, 0})}, $$
  we get from \eqref{estB1}  that
  $$  \eps \norm{ B[ \rho_{\mathrm{F}},  f]  }_{L^2(0, T; \H^0_{r, 0})} \leq  \norm{ \rho }_{L^2(0, T; H^m_{r})} \| f \|_{L^\infty(0, T; \H^m_{r})},
   \leq \Lambda(T, M).$$
  hence the lemma.

\end{proof}
As a consequence of this preliminary analysis, we have been able to reduce  \eqref{EqPrin2} 
to 
\begin{equation}\label{EqPrin3}
\rho_{\mathrm{F}} =  - \int_v \int_0^t U^{\mathrm{FIO}}_{t,s}B[\rho_{\mathrm{F}}, f]dsdv +R,
\end{equation}
where $R$ is a remainder.

\subsection{Further reductions using  a quantum averaging lemma}

By definition of  the Fourier Integral Operator $U_{t,s}^{\mathrm{FIO}}$, we have
\begin{multline}
\label{oldFIO++}
\int_v \int_0^t U_{t,s}^{\mathrm{FIO}} B[ \rho_{\mathrm{F}}(s), f(s)] dsdv\\ =\frac{1}{(2\pi)^{2d}} \int_{v}\int_{0}^t\int_{y}\int_{\xi} e^{\frac{i}{\varepsilon} \left( \varphi _{t,s}^{\eps}(z,\xi)-\langle y,\eps\xi\rangle\right)}\mathrm{B} _{t,s}^{\eps}(z,\xi) B[\rho_{\mathrm{F}}(s), f(s)] (y)d\xi dy ds dv.
\end{multline}
Let us introduce $\widetilde{U}^{\mathrm{FIO}}_{t,s}$ the Fourier integral operator associated with the phase $\varphi_{t,s}$ and the  amplitude $\mathrm{I}$, and we consider its action on the vector  $B[\rho_{\mathrm{F}}, f^0]$ where $f^0$ is the initial datum, which gives rise to the integral 
\begin{multline}
\label{newFIO}
\int_v \int_0^t \widetilde{U}^{\mathrm{FIO}}_{t,s}B[\rho_{\mathrm{F}}(s), f^0] dsdv \\ 
= \frac{1}{(2\pi)^{2d}} \int_{v} \int_{0}^t  \int_{\xi}  \int_{y}e^{\frac{i}{\varepsilon} \left( \varphi_{t,s}^{\eps} (z,\xi)-\langle y,\xi\rangle\right)} B[\rho_{\mathrm{F}}(s), f^0](y) dy d\xi    ds dv.
\end{multline}
The difference between the  terms~\eqref{oldFIO++} and~\eqref{newFIO} is  shown to be a remainder in the next lemma. To this end, we need to apply a quantum averaging lemma of Section \ref{sec:quantum}.

\begin{lem}
For $T \leq \min(T_\eps,T(M))$, we have the estimate
$$
\norm{\int_v \int_0^t \left(U_{t,s}^{\mathrm{FIO}}B[\rho_{\mathrm{F}}, f] - \widetilde{U}^{\mathrm{FIO}}_{t,s}B[\rho_{\mathrm{F}}, f^0]  \right)dsdv}_{L^2(0,T;H^0_{r})} \\ \leq T \Lambda (T,M).
$$
\end{lem}
\begin{proof}
By  the triangular inequality, we first write 
\begin{multline*}
\norm{\int_v \int_0^t \left(U_{t,s}^{\mathrm{FIO}}B[ \rho_{\mathrm{F}}, f] - \widetilde{U}^{\mathrm{FIO}}_{t,s} B[  \rho_{\mathrm{F}}, f^0]  \right) dsdv}_{L^2(0,T;H^0_{r})} \\
\leq \norm{\int_v \int_0^t U_{t,s}^{\mathrm{FIO}}B[ \rho_{\mathrm{F}}, f(s)- f^0] dsdv}_{L^2(0,T;H^0_{r})} \\
+\norm{\int_v \int_0^t \left(U_{t,s}^{\mathrm{FIO}}- \widetilde{U}^{\mathrm{FIO}}_{t,s}\right)B[ \rho_{\mathrm{F}}, f^0] dsdv}_{L^2(0,T;H^0_{r})}.
\end{multline*}
For the first term of the right-hand side, we apply the quantum averaging lemma adapted to the space $H^0_r$, namely Theorem \ref{M-Pro3}, with 
$$\Phi_{t,s}=\frac{1}{\eps}\varphi_{t,s}^{\eps}, \quad b_{t,s}=B_{t,s}^{\eps}, \quad G_{t, s}=f(s)- f^0,$$ 
and $q=0$.  Recall the notation $k_d =\lfloor d/2 \rfloor +2$. As already explained, the fact that the phase $\Phi_{t,s}$ satisfies the assumption $(A_{4k_d +d +4+r})$  comes from Proposition~\ref{prop-HJ1} and the fact that $m\geq 5k_d +d +4+r$. We obtain
\begin{multline*}
\norm{\int_v \int_0^t U_{t,s}^{\mathrm{FIO}}B[ \rho_{\mathrm{F}}, f(s)-f^0] dsdv}_{L^2(0,T;H^0_{r})}  \\
 \lesssim  \norm{\mathrm{B}^\eps_{t,s}}_{L^\infty_T W^{d+4k_d +4+r,\infty}_{z,\xi}} \norm{ \langle \eps \nabla_x \rangle^{k_d+r} \langle \eps \nabla_v \rangle^{k_d} (f(s)- f^0)}_{T, 3k_d + 2d+6, k_d+ d+2} \norm{\rho}_{L^2(0,T;H^m_r)}.
\end{multline*}
According to Proposition~\ref{EstB}, we have
$$
 \norm{\mathrm{B}^\eps_{t,s}}_{L^\infty_T W^{d+4k_d +4+r,\infty}_{z,\xi}}  \leq \Lambda(T,M),
$$
since $m \geq 5 k_d + d + 5 + r$.
 Furthermore, using Remark \ref{rem:norme-av}  and the fact that $f$ solves  the Wigner equation~\eqref{eq:wigner} with initial condition $f^0$, we  obtain that
\begin{align*}
\norm{ \langle \eps \nabla_x \rangle^{k_d+r} \langle \eps \nabla_v \rangle^{k_d} (f(s)- f^0)}_{T, 3k_d + 2d+6, k_d+ d+2}  &\leq \sup_s \| f(s) -f^0\|_{\H^{m-2}_{r-1}} \\
 &\leq \sup_s \int_0^s \norm{ \partial_\tau f  }_{\H^{m-2}_{r-1}}  \, d\tau  \leq T \Lambda(T,M),
\end{align*}
by the fact that  $m\geq 4k_d+2d+ 8 + r $, $r \geq 2k_d +2d +4$.

For the second term of the right-hand side, we apply again  Theorem \ref{M-Pro3}, still for $q=0$, with 
$$\Phi_{t,s}=\frac{1}{\eps}\varphi_{t,s}^{\eps}, \quad b_{t,s}={\mathrm{B}}_{t,s}^{\eps} - \mathrm{I}, \quad G=f^0.$$ 
This leads to 
\begin{align*}
\Big\|\int_v &\int_0^t \left(U_{t,s}^{\mathrm{FIO}} - \widetilde{U}^{\mathrm{FIO}}_{t,s}\right)B[\rho_{\mathrm{F}}, f^0]  dsdv\Big\|_{L^2(0,T;H^0_{r})}  \\
&\lesssim   \norm{ \mathrm{B}^\eps_{t,s}-\mathrm{I}}_{L^\infty_T W^{d+4k_d +4+r,\infty}_{z,\xi}} \norm{ \langle \eps \nabla_x \rangle^{k_d+r} \langle \eps \nabla_v \rangle^{k_d}f^0}_{T, 3k_d + 2d+6, k_d+ d+2} \norm{\rho}_{L^2(0,T;H^m_r)}.
\end{align*}
According to Proposition~\ref{EstB}, we have
$$
 \norm{\mathrm{B}^\eps_{t,s}- \mathrm{I}}_{L^\infty_T W^{d+4k_d +4+r,\infty}_{z,\xi}}  \leq T \Lambda(T,M),
$$
since $m \geq 5 k_d + d + 6+ r$. We also have thanks to Remark \ref{rem:norme-av} that
\begin{equation}
\label{eq:normavecf0}
 \norm{ \langle \eps \nabla_x \rangle^{k_d+r} \langle \eps \nabla_v \rangle^{k_d}f^0}_{T, 3k_d + 2d+6, k_d+ d+2}  \leq  \| f^0\|_{\H^{m}_{r}} \leq M_0,
\end{equation}
since $m\geq 4k_d+2d+ 6 + r $, $r \geq 2k_d +2d +3$. We thus get
$$
\Big\|\int_v \int_0^t \left(U_{t,s}^{\mathrm{FIO}} - \widetilde{U}^{\mathrm{FIO}}_{t,s}\right)B[\rho_{\mathrm{F}}, f^0]  dsdv\Big\|_{L^2(0,T;H^0_{r})}   \leq T \Lambda(T,M).
$$
Gathering these two estimates, we obtain the claimed result.

\end{proof}
At this point of the proof, we have therefore been able to recast \eqref{EqPrin} as
$$
\rho_{\mathrm{F}} = -  \int_v \int_0^t \widetilde{U}^{\mathrm{FIO}}_{t,s} B[ \rho_{\mathrm{F}}(s), f^0] dsdv + R,
$$
where $R$ is a remainder.  Observe that the matrix-valued FIO $\widetilde{U}^{\mathrm{FIO}}_{t,s}$ acts diagonally  and that by definition
 of $\mathrm{F}$, we have 
 $$ (\rho_{\mathrm{F}})_{k}=  \int_v\mathrm{F}_{k} dv= \partial^{\alpha(k)}_x \rho, \quad 1 \leq k \leq n_{m}.$$
 We can study  this diagonal system componentwise and thus  focus on the scalar equations
 $$ \partial^\alpha_x \rho (t)= - \int_v \int_0^t \widetilde{U}^{\mathrm{FIO}}_{t,s}B[\partial^{\alpha}_x  \rho(s), f^0] dsdv + R, \quad
  | \alpha |=m.$$
  Note that in the above expression we are abusing notation and still write  $\widetilde{U}^{\mathrm{FIO}}_{t,s}$ for the scalar FIO where the amplitude is now
   $1$ instead of $\mathrm{I}$ and that now $B$ acts on a scalar quantity and  is also scalar.

The next step is to relate the above integral  to 
\begin{multline*}
\int_v \int_0^t {U}^{\mathrm{free}}_{t,s}\left[ B[\partial^{\alpha}_x \rho(s), f^0]\right]dsdv\\ 
:=\frac{1}{(2\pi)^{2d}} \int_{v}\int_{0}^t\int_{\xi}\int_{y} e^{{i}\left(\left( x-(t-s)v\right)\cdot\xi_x+v\cdot\xi_v-y \cdot \xi \right)}  B[\partial^{\alpha}_x  \rho(s) , f^0]  dy d\xi  ds dv.
\end{multline*}
The operator ${U}^{\mathrm{free}}_{t,s}$ can be seen as a FIO, with the free  phase 
 $\varphi(z,\xi)=\left( x-(t-s)v\right)\cdot\xi_x+v\cdot\xi_v$, and amplitude $1$.
To compare  $\widetilde{U}^{\mathrm{FIO}}_{t,s}$ and  ${U}^{\mathrm{free}}_{t,s}$, 
we shall again use a quantum averaging lemma of Section~\ref{sec:quantum}.

\begin{lem}
For $T \leq \min(T_\eps,T(M))$, we have the estimate 
$$
\norm{ \int_v \int_0^t \widetilde{U}^{\mathrm{FIO}}_{t,s}B[\partial^{\alpha}_x  \rho(s), f^0]dsdv
- \int_v \int_0^t {U}^{\mathrm{free}}_{t,s} B[\partial^{\alpha}_x  \rho(s), f^0] dsdv }_{L^2(0,T;H^0_{r})} 
\leq T^{1/2} \Lambda (T,M).
$$
\end{lem}
\begin{proof}
Let us write
$$
\varphi_{t,s} (z,\xi)=\left(x-(t-s)v\right)\cdot\xi_x+v \cdot \xi_v + \widetilde{\varphi}_{t,s} (z,\xi).
$$
We aim at applying Theorem \ref{M-Pro}  with 
$$\Phi_{t,s}(z,\xi)=(x-(t-s)v)\cdot  \xi_x+ v\cdot \xi_v, \quad b_{ t,s}(z, \xi)= e^{ i { \widetilde{\varphi}_{t,s}^\eps(z,\xi) \over \eps} } - 1, \quad G_{t,s}= f^0
$$
and $q=1$. We obtain 
\begin{multline*}
\norm{ \int_v \int_0^t \widetilde{U}^{\mathrm{FIO}}_{t,s}B[\partial^{\alpha}_x  \rho(s), f^0]dsdv
- \int_v \int_0^t {U}^{\mathrm{free}}_{t,s} B[\partial^{\alpha}_x  \rho(s), f^0] dsdv }_{L^2(0,T;H^0_{r})} \\
 \lesssim   \norm{ \left(e^{ i { \widetilde{\varphi}_{t,s}^\eps \over \eps} } - 1\right)(\langle \xi_v \rangle + \langle (t-s) \xi_x\rangle)^{-1}}_{L^\infty_T W^{d+4k_d +5+r,\infty}_{z,\xi}} \norm{ \langle \eps \nabla_x \rangle^{k_d+r} \langle \eps \nabla_v \rangle^{k_d}f^0}_{T, 3k_d + 2d+7, k_d+ d+2} \\
 \times \norm{\rho}_{L^2(0,T;H^m_r)} 
\end{multline*}
We then use the sharp estimates of Lemma~\ref{lemestimtildephi}.
As $|b_{ t, s}(z, \xi)| \leq\frac{1}{\eps} | \widetilde{\varphi}_{t,s}^\eps |$, we can use
\eqref{tildephidegreun}
  to obtain
$$
|b_{ t, s}(z, \xi)|  \leq T^{1/2}\Lambda(T,M) (|\xi_v| + |t-s| |\xi_x|),
$$
since $m \geq k_d$ and $ \widetilde{\varphi}_{t,s}^\eps(z, \xi)= \widetilde{\varphi}_{t,s}(z, \eps \xi)$. Regarding $\partial^\alpha_z \partial^\beta_\xi b_{ t, s}$ for  $0<|\alpha|+|\beta|\leq d+4k_d +5+r$, since $m\geq d+5k_d +5+r$, 
\begin{itemize}
\item when $\beta\neq 0$,  we can use \eqref{tildephiunif} and the fact that $\widetilde\varphi$ is evaluated at $\eps \xi$  so that we  gain a factor $\eps$ when we take derivatives in $\xi$,
\item when $\beta=0$,  we can  use again  \eqref{tildephidegreun}. 
\end{itemize}
This yields in all the cases
$$
|\partial^\alpha_z \partial^\beta_\xi b_{t, s}(z, \xi)|  \leq  T^{1/2} \Lambda(T,M)  (\langle \xi_v\rangle + \langle(t-s)\xi_x\rangle),
$$
that is to say
$$
 \norm{ \left(e^{ i { \widetilde{\varphi}_{t,s}^\eps \over \eps} } - 1\right)(\langle \xi_v \rangle + \langle (t-s) \xi_x\rangle)^{-1}}_{L^\infty_T W^{d+4k_d +5+r,\infty}_{z,\xi}}  \leq T^{1/2} \Lambda(T,M).
$$
Finally using a variant of~\eqref{eq:normavecf0} to control the contribution of $f^0$, 
we obtain the claimed inequality.

\end{proof}
Thanks to the above Lemma, we have reached the point where we have been able to reduce  \eqref{EqPrin} to
\begin{equation}
\label{encoreunesimplif}
\partial^\alpha_x \rho  =-\int_v \int_0^t {U}^{\mathrm{free}}_{t,s}B[\partial^{\alpha}_x  \rho(s), f^0]dsdv+ R, \quad | \alpha |=m
\end{equation}
where $R$ is  a  remainder. 

\subsection{Final reduction}

It remains to further simplify the action of ${U}^{\mathrm{free}}_{t,s}$ on $B$, which is the object of the following lemma.
\begin{lem}
\label{dernierereduc}
 For $T \leq \min(T_\eps,T(M))$, we have
 \begin{align*}
 \Bigg\| \int_v \int_0^t &{U}^{\mathrm{free}}_{t,s}B[\partial^{\alpha}_x  \rho(s), f^0] \, dsdv \\
 &- 
   {2 \over (2\pi)^{d}}  \int_{\eta} \int_{0}^t  e^{ix\cdot \eta} {1 \over \eps} \sin \left( \eps(t-s) {|\eta|^2 \over 2} \right) \widehat{V}(\eps \eta) 
  \mathcal{F}_{v}f^0(x, (t-s) \eta) \widehat{ \partial^\alpha_x \rho}(s, \eta) ds d\eta
     \Bigg\|_{L^2(0, T;H^0_{r})} \\
      &\leq (T + \eps) \Lambda(T,M).
      \end{align*}
\end{lem}

\begin{proof}
We first observe that
$$\int_v \int_0^t {U}^{\mathrm{free}}_{t,s}B[\partial^{\alpha}_x  \rho(s), f^0] \, dsdv
=  \int_{v}   \int_{0}^t B[\partial^{\alpha}_x  \rho(s), f^0](x-(t-s)v, v) \, ds dv.$$
Next, by using \eqref{SinForm} together with the expression \eqref{eq:K3} of the symbol, 
 we obtain that
 \begin{multline*}
 \int_v \int_0^t {U}^{\mathrm{free}}_{t,s}B[\partial^{\alpha}_x  \rho(s), f^0] \, dsdv =  \\
 {1 \over (2\pi)^{d}} \int_{v} \int_{0}^t \int_{\xi_{x}} e^{i (x-(t-s)v ) \cdot \xi_{x}} \left( \int_{-{1 \over 2}}^{1 \over 2} \xi_{x} \cdot \nabla_{v} f^0(x-(t-s)v, v+ \lambda \eps  \xi_{x}) \, d\lambda\right) \widehat{ \partial^\alpha_x \mathrm{V}_{\rho} }(\xi_{x}) \, d\xi_{x} ds dv.
  \end{multline*}
We then use a Taylor expansion to write
\begin{multline*}
 \xi_{x} \cdot   \nabla_{v} f^0(x-(t-s)v, v+ \lambda \eps  \xi_{x})
=   \xi_{x} \cdot   \nabla_{v} f^0(x, v+ \lambda \eps  \xi_{x})
\\-  \int_{0}^1 D_xD_{v} f^0(x- \lambda'(t-s) v, 
 v+ \lambda \eps \xi_{x})\cdot[(t-s) v, \xi_{x}]\, d\lambda',
 \end{multline*}
 where we have denoted $D_xD_{v} f^0 = (\partial_{x_i} \partial_{v_j} f^0)_{i,j}$.
We thus get the expression
 \begin{multline*}
\int_v \int_0^t {U}^{\mathrm{free}}_{t,s}B[\partial^{\alpha}_x  \rho(s), f^0]dsdv \\ 
=  {1 \over (2\pi)^{d}} \int_{v} \int_{0}^t \int_{\xi_{x}} e^{i (x-(t-s)v ) \cdot \xi_{x}} \left( \int_{-{1 \over 2}}^{1 \over 2} \xi_{x} \cdot \nabla_{v} f^0(x, v+ \lambda \eps  \xi_{x}) \, d\lambda\right) \widehat{ \partial^\alpha_x \mathrm{V}_{\rho} }(s, \xi_{x}) \, d\xi_{x} ds dv - \mathcal{I} 
\end{multline*}
where
$$ \mathcal{I} =
{1 \over (2\pi)^{d}}  \int_{\xi_{x}} e^{i x \cdot \xi_{x}} \left( \int_{v} \int_{0}^t  \int_{-{1 \over 2}}^{1 \over 2} \int_{0}^1 e^{- i(t-s) v \cdot \xi_{x}} H_{\eps,  t-s}(x, v, \xi_{x}, \lambda, \lambda') d\lambda' d\lambda ds dv\right) \widehat{ \partial^\alpha_x \mathrm{V}_{\rho} }(s, \xi_{x}) \, d\xi_{x} 
 $$ 
and we have set
\begin{equation}
\label{Hdeffin}  H_{\eps, \tau}(x, v, \xi_{x}, \lambda, \lambda')=  D_x D_{v} f^0(x- \lambda'\tau v, 
 v+ \lambda \eps \xi_{x})\cdot[\tau v,\xi_{x}].
 \end{equation}
By using similar computations to those in the proof of Lemma \ref{lem:Keps}, we have that
\begin{multline}
\label{alternateL}
 {1 \over (2\pi)^{d}} \int_{v} \int_{0}^t \int_{\xi_{x}} e^{i (x-(t-s)v ) \cdot \xi_{x}} \left( \int_{-{1 \over 2}}^{1 \over 2} \xi_{x} \cdot \nabla_{v} f^0(x, v+ \lambda \eps  \xi_{x}) \, d\lambda\right) \widehat{ \partial^\alpha_x \mathrm{V}_{\rho} }(s, \xi_{x}) \, d\xi_{x} ds dv\\
=   {2 \over (2\pi)^{d}} \frac{1}{\eps} \int_{\xi_{x}} \int_{0}^t e^{i x \cdot \xi_{x}} \sin \left( { \eps(t-s) |\xi_{x}|^2  \over 2 } \right)
 \mathcal{F}_{v}f^0(x, (t-s) \xi_{x}) \widehat{ \partial^\alpha_x \mathrm{V}_{\rho}}(s, \xi_{x})\, d \xi_{x} ds,
 \end{multline}
 so that recalling the definition \eqref{def:Vrho} of $\mathrm{V}_{\rho}$,  to get Lemma \ref{dernierereduc}, it suffices  to prove that
$$\| \mathcal{I}\|_{L^2(0, T; H^0_{r})} \leq (T + \eps) \Lambda(T,M).$$
This estimate is reminiscent of the averaging Lemma proven in \cite{HKR} on the torus.
We shall follow here another approach based on the operator-valued  pseudodifferential calculus
developed in Appendix \ref{Pseudo} (the proof is thus close to that for the quantum averaging lemmas in Section~\ref{sec:quantum}).
We can write $\mathcal{I}$ under the form
$$ \mathcal{I}= \operatorname{Op}_{\mathrm{L}}(\partial^\alpha_x \mathrm{V}_{\rho})$$
where $\mathrm{L}(x, \eta)$ is  an operator-valued symbol acting on $L^2(0, T)$ and defined by  the convolution
$$\mathrm{L}(x, \eta)( \Upsilon) (t) = \int_{0}^t K_{\eps, t-s}(x, \eta) \Upsilon(s)\, ds,$$
where we have set
\begin{equation}
\label{Kdeffin} K_{\eps,  \tau}(x, \eta)=   \int_{v}   \int_{-{1 \over 2}}^{1 \over 2} \int_{0}^1 e^{- i t v \cdot \eta} H_{\eps, \tau}(x, v, \eta, \lambda, \lambda') d\lambda' d\lambda dv.
\end{equation}
By using the Calder\'on-Vaillancourt theorem of Appendix  \ref{Pseudo}, to obtain the estimate, we only have to show that
$$ \sup_{x, \eta} \| \partial_{x}^\alpha \partial_{\eta}^{\alpha'}  \mathrm{L}(x,\eta) \|_{\mathscr{L}(L^2(0, T))}
 \leq (T + \eps)\Lambda(T, M), \quad | \alpha | \leq k_d+r, \, | \alpha' | \leq k_d, \, k_d= 2 + \lfloor{ d \over 2} \rfloor.$$
 From the Young inequality for convolution in time, we have 
 $$ \| \partial_{x}^\alpha \partial_{\eta}^{\alpha'}  \mathrm{L}(x,\eta) \|_{\mathscr{L}(L^2(0, T))}
  \lesssim \sup_{x, \eta} \int_{0}^T  | \partial_{x}^\alpha \partial_{\eta}^{\alpha'} K_{\eps,  t}(x, \eta)| \, dt,$$
  so that the proof is reduced to showing that
  $$ \sup_{x, \eta} \int_{0}^T  | \partial_{x}^\alpha \partial_{\eta}^{\alpha'} K_{\eps,  t}(x, \eta)| \, dt
  \leq  (T + \eps)\Lambda(T, M).$$
  By using integration by parts in the $v$ integral, we get from the definition \eqref{Kdeffin}
of $K_{\eps, t}$ and  \eqref{Hdeffin}
 that for any $\alpha'' \in \mathbb{N}^d$, $| \alpha ''| \leq p$, and $t\leq T$,
\begin{multline*} |(t \eta)^{\alpha''}|\,  |  \partial_{x}^{\alpha} \partial_{\eta}^{\alpha'}K_{\eps, t}(x, \eta) |  \\
 \lesssim \sup_{|\beta| \leq p + 2k_d + r + 2  }
 \int_{-{1 \over 2}}^{1 \over 2} \int_{0}^1 \int_v |\partial^\beta_{x,v} f( x- \lambda't v, v + \lambda \eps \eta)  
 \langle t |v| \rangle^{k} (  1 + t |v| + t| \eta| + t | \eta|\, |v|) \, dv d\lambda d \lambda'.
\end{multline*}
Thanks to \eqref{embedH01}  and the  Sobolev embedding in $\R^{2d}$, we have the pointwise
estimate 
$$ \langle v \rangle^{r_{0}} | \partial^\beta_{x,v} f^0 (x,v)| \lesssim \|f^0 \|_{\H^{m-1}_{r_{0}}}$$
 if $m \geq  | \beta | +  2 + d$. 
 Therefore, we obtain for  $m \geq p+2k_d +r +d +4 $ that 
 $$
 |(t \eta)^{\alpha''}| \, | \partial_{x}^\alpha \partial_{\eta}^{\alpha'}K_{\eps, t}(x, \eta) | \lesssim
 \|f^0 \|_{\H^{m-1}_{r_{0}}} 
 \int_{v}   \int_{-{1 \over 2}}^{1 \over 2}   {  | \langle  t | v| \rangle^{k} \over  \langle v+ \lambda \eps | \eta|
  \rangle^{r_{0}} } (  1 + t |v| + t| \eta| + t | \eta|\, |v|) \, d\lambda dv.$$
  By using that $|v| \leq | v+ \lambda \eps  \eta |  + \lambda \eps | \eta |,$ we get  that  
 $$    | \partial_{x}^\alpha \partial_{\eta}^{\alpha'} K_{\eps, t}(x, \eta) | \leq \Lambda(T) \|f^0 \|_{\H^{m-1}_{r_{0}}}  \int_{v}   \int_{-{1 \over 2}}^{1 \over 2}  \left(  { 1 + \eps | \eta|  \over \langle t | \eta | \rangle^{p-1-k}
  \langle v + \lambda \eps \eta \rangle^{r_{0}-k-1}}
   \right) \, d\lambda dv .$$
   To conclude, we choose,  $p= 4 + k_d$,  $ r_{0} =   k_d + 2 + d$,
    which is justified since $r \geq k_d + 2  + d  $  and $m \geq 3k_d + d+ 8 + r$. 
 This finally yields
$$   | \partial_{x}^\alpha \partial_{\eta}^{\alpha'} K_{\eps, t}(x, \eta) | \leq \Lambda(T,M)     { 1 + t| \eta | + \eps | \eta|  \over \langle t | \eta | \rangle^{3}}, $$ 
and after integration in time
$$ \sup_{x, \eta} \int_{0}^T  | \partial_{x}^\alpha \partial_{\eta}^{\alpha'} K_{\eps,  t}(x, \eta)| \, dt
\lesssim   ( T + \eps) \Lambda(T,M),$$
concluding the proof.

 \end{proof}

By using Lemma \ref{dernierereduc}, we can thus further simplify  \eqref{encoreunesimplif}
into \begin{multline}\label{EqFFinal}
\partial^\alpha_x \rho(t,x) = -{2 \over (2\pi)^{d}} \int_{\eta}\int_{0}^t e^{{i} x\cdot \eta }\mathcal{F}_{v}  f^0(x,(t-s)\eta)\\ \left(  \frac{1}{\eps} \sin\left(\eps(t-s)\frac{\vert \eta \vert^2}{2}\right)  \widehat{V}(\eps \eta)\widehat {\partial^\alpha_x \rho} (s,\eta)\right)  ds    d\eta +R, \quad | \alpha |=m
\end{multline}
where $R$ is a remainder. 
We have therefore managed to turn the study of the initial identity~\eqref{EqPrin} to that of~\eqref{EqFFinal}.

\subsection{Quantum Penrose stability} 
To complete the proof of Proposition \ref{prop:rho-final}, we  need to prove a quantitative estimate for  a solution $\mathbf{h} \in L^2(0, T;L^2(\mathbb{R}^d))$  
 to the  scalar equation 
\begin{multline}\label{Eqh}
{\mathbf{h}} (t,x) = -{2 \over (2\pi)^{d}} \int_{\eta}\int_{0}^t e^{{i} x\cdot \eta }\mathcal{F}_{v} f^0(x,(t-s)\eta))\\ \left(  \frac{1}{\eps} \sin\left(\eps(t-s)\frac{\vert \eta \vert^2}{2}\right)  \widehat{V}(\eps \eta)\widehat{ {\mathbf{h}}} (s,\eta)\right)  ds   d\eta +\mathbf{R}(t,x),
\end{multline}
where $\mathbf{R}$ is a given source term and $\widehat{ {\mathbf{h}}}$ stands for the Fourier transform of $\mathbf{h}$ with respect to $x$.

\begin{definition}Let us define the operator acting on $ \mathbf{h} \in L^2( \mathbb{R}; L^2(\mathbb{R}^d))$ by
\begin{multline}
\label{defLf0} \mathcal{L}_{\eps, f^{0}}{\mathbf{h}} (t,x)= \\
-{2 \over (2\pi)^{d}} \int_{\eta}\int_{0}^t e^{{i} x\cdot \eta }\mathcal{F}_{v} f^0(x,(t-s)\eta)) \left(  \frac{1}{\eps} \sin\left(\eps(t-s)\frac{\vert \eta \vert^2}{2}\right)  \mathds{1}_{t-s \geq 0} \widehat{V}(\eps \eta)\widehat{ {\mathbf{h}}} (s,\eta)\right)  ds   d\eta.
\end{multline}
\end{definition}
We shall  first relate $\mathcal{L}_{\eps, f^0}$ to a space-time pseudodifferential operator with parameter.

\begin{lem}\label{Eqpseudo}
For all $h\in \mathscr{S}(\R \times \R^d)$ satisfying ${h}|_{t<0}=0$, and every $\gamma \geq 0$,  
we have
\begin{align}
\label{conjugform}
 e^{-\gamma t }\mathcal{L}_{\eps, f^0}( e^{\gamma t } h)=
\mathbf{Op}_{\mathcal{P}_{\mathrm{quant}}}^{\eps, \gamma}(h),
\end{align}
where $\mathbf{Op}_{\mathcal{P}_{\mathrm{quant}}}^{\eps, \gamma}$ is the pseudodifferential operator in time and space associated with the symbol
\begin{equation}\label{DefASymb}
  \mathcal{P}_{\mathrm{quant}}(x,  \gamma, \tau,  \eta)=  -2 \widehat V( \eta)  \int_0^{+\infty}  e^{-(\gamma +i\tau) s}\mathcal{F}_{v}f^0(x,s\eta)    \sin\left( s\frac{\vert \eta \vert^2}{2}\right)  ds,
\end{equation}
which is  the quantum Penrose function introduced in~\eqref{Penrose-quantique}.
\end{lem}
\begin{proof}
Since ${h}|_{t<0}=0$, we first note that 
$$
\mathcal{L}_{\eps, f^{0}}{{h}}={1 \over (2\pi)^{d}} \int_{\eta}\int_{-\infty}^t e^{{i} x\cdot \eta }  e^{-\gamma (t-s)}\mathcal{F}_{v}f^0(x,(t-s)\eta) \left(  \frac{-2}{\eps} \sin\left(\eps(t-s)\frac{\vert \eta \vert^2}{2}\right)  \widehat{V}(\eps \eta)\widehat h (s,\eta)\right)dsd\eta.
$$
Taking the inverse Fourier transform in time, we can write
$$
\widehat h (s,\eta)= \frac{1}{2\pi}\int_{\tau}e^{i\tau s} \mathcal{F}_{t,x}h(\tau,\eta) d\tau
$$
Plugging in this identity, we reach the formula
\begin{multline*}
\mathcal{L}_{\eps, f^{0}}{{h}}
= \frac{1}{(2\pi)^{d+1}}\int_{\tau}\int_{\eta} e^{{i}\left( x\cdot \eta +\tau t \right)}  \int_{-\infty}^t  e^{-(\gamma +i\tau) (t-s)}\mathcal{F}_{v}f^0(x,(t-s)\eta) \\ \left(  \frac{-2}{\eps} \sin\left(\eps(t-s)\frac{\vert \eta \vert^2}{2}\right)  \widehat{V}(\eps \eta) \mathcal{F}_{t,x}h(\tau,\eta)\right)ds d\eta d\tau 
\end{multline*}
Changing variable in the integral in $s$, we eventually obtain 
$$
\mathcal{L}_{\eps, f^{0}}{{h}}= \frac{1}{(2\pi)^{d+1}} \int_{\tau}\int_{\eta} e^{{i}\left( x\cdot \eta +\tau t \right)}  \mathcal{P}_{\mathrm{quant}}(x,\eps\gamma,\eps\tau,\eps \eta)  \mathcal{F}_{t,x}h(\tau,\eta) d\eta d\tau = \mathbf{Op}_{ \mathcal{P}_{\mathrm{quant}}}^{\eps,\gamma}(h),
$$
recalling the quantization \eqref{PseudoParamDef}.
\end{proof}
To save space, we will denote $\zeta=(\gamma,\tau,\eta)$; also, since no confusion is possible, we denote  from now on $\mathcal{P}$ instead of $\mathcal{P}_{\mathrm{quant}}$ for the quantum Penrose function.

Recall  the notation $k_d = \lfloor d/2 \rfloor +2$. We provide in  Appendix~\ref{sec:PseuDiff} the  required pseudodifferential calculus associated with this quantization. Namely, we shall rely on 
\begin{prop}\label{prop:pseudoparam}
 There exists $C>0$ such that for every $\eps \in (0,1]$ and every $\gamma>0$, we have 
\begin{itemize}
\item  for every symbol $a$ such that $ |a|_{k_d,0}<+\infty$
$$
\norm{\mathbf{Op}_a^{\eps, \gamma}}_{\mathscr{L}(L^2(\mathbb{R} \times \mathbb{R}^d))}\leq C |a|_{k_d,0}, 
$$
\item  for every symbol $a,b$ such that $ |a|_{k_d,1}<+\infty,  |b|_{k_d+1,0}<+\infty$
$$
\norm{\mathbf{Op}_a^{\eps,\gamma}\mathbf{Op}_b^{\eps,\gamma}-\mathbf{Op}_{ab}^{\eps,\gamma}}_{\mathscr{L}(L^2(\mathbb{R} \times \mathbb{R}^d))}\leq \frac{C}{\gamma} |a|_{k_d,1}|b|_{k_d+1,0} .
$$
\end{itemize}
The seminorms $|\cdot|_{k,0}$ and $|\cdot|_{k,1}$  are defined for any $k \in \mathbb{N}$ as
$$
\begin{aligned}
\vert c \vert _{k,0} &= \sup_{\vert \alpha\vert \leq k } \norm{\mathcal{F}_x(\partial^{\alpha}_x c)}_{L^1(\mathbb{R}^d;L^{\infty}_\zeta)}, \\
\vert c \vert _{k,1} &= \sup_{\vert \alpha\vert \leq k } \norm{\gamma \mathcal{F}_x(\partial^{\alpha}_x \nabla_\xi c)}_{L^1(\mathbb{R}^d;L^{\infty}_\zeta)},
\end{aligned}
$$
where $\xi=( \tau,\kappa)$.
\end{prop}

The symbol $\mathcal{P}$ is a good symbol for this calculus, as checked in the next lemma.
\begin{lem}
\label{lem:a-symbol}
For the quantum Penrose function $\mathcal{P}$, we have  for every $k\in \mathbb{N}$ such that $m\geq k+ 6$  the estimates
\begin{align*}
|\mathcal{P}|_{k,0} &\leq C \norm{f^0}_{\mathcal{H}^{k+4}_{k_d-1}}, \\
|\mathcal{P}|_{k,1} &\leq C \norm{f^0}_{\mathcal{H}^{k+6}_{k_d}},
\end{align*}
\end{lem}
Lemma~\ref{lem:a-symbol} will be specifically used for $k= k_d$ or $k_d+1$; we therefore use that $m\geq k_d +7$ and $r \geq k_d+1$.
To ease readability, the proof of Lemma~\ref{lem:a-symbol} is postponed to the end of the section.
This lemma implies, thanks to the first item of~Proposition \ref{prop:pseudoparam}
that $\mathbf{Op}_{\mathcal{P}}^{\eps,\gamma} \in \mathscr{L}(L^2(\R \times \R^d))$ with norm uniform in $\eps$.

In order to study \eqref{Eqh} on $[0, T]$, we shall first study
the global (that is for all $t \in \mathbb{R}$) pseudodifferential equation
\begin{equation}\label{EqPseudo}
h=\mathbf{Op}_{\mathcal{P}}^{\eps,\gamma}(h)+\mathcal{F},
\end{equation}
for a given source term  $\mathcal{F}$.
We have reached the point of the proof where the quantum Penrose stability condition~\eqref{QCrit} plays a crucial role.

\begin{prop}\label{MajhFin}
Under the $c_0$-quantum Penrose stability condition~\eqref{QCrit},
 we have the following properties:
 \begin{itemize}
 \item[i)] 
  there exists $\gamma_0\geq 1$ depending only on $\norm{f^0}_{\mathcal{H}^m_{r}}$ and $c_0$ such that, for $\gamma\geq \gamma_0$, the operator $\mathrm{I}-\mathbf{Op}_{\mathcal{P}}^{\eps,\gamma}$ is invertible on $L^2(\R \times \R^d)$:  there exists  $\Lambda(c_0^{-1}, \norm{f^0}_{\mathcal{H}^m_{r}})$ such that
  for every $\mathcal{F} \in L^2(\mathbb{R} \times \mathbb{R}^d)$,  $\gamma \geq \gamma_{0}$,  
  $ \eps \in (0, 1)$, there exists
a  unique  solution $h_{\gamma, \eps}$ to \eqref{EqPseudo}, and  we have the estimate
$$
\norm{h_{\gamma, \eps}}_{L^2(\mathbb{R}\times \mathbb{R}^d)}\leq \Lambda(c_0^{-1}, \norm{f^0}_{\mathcal{H}^m_{r}}) \norm{{\mathcal{F}}}_{L^2(\mathbb{R}\times \mathbb{R}^d)}.
$$
\item[ii)] Consider  $ \mathbb{F} \in L^2( \mathbb{R}; L^2(\mathbb{R}^d))$ such that
 $F_{|t<0}=0$.
Then, the fonction
$$ \mathbf{h}= e^{\gamma t } ( \mathrm{I}-\mathbf{Op}_{\mathcal{P}}^{\eps,\gamma})^{-1}
(  e^{-\gamma t} \mathbb{F})$$
vanishes for $t<0$ and 
 does not depend on $\gamma$ for $\gamma \geq \gamma_{0}$.
 \item
 [iii)] Consider  $ \mathbb{F} \in L^2( \mathbb{R}; L^2(\mathbb{R}^d))$ such that
 $\mathbb{F} _{|t \leq T}=0$ for some $T>0$.
Then,  for $\gamma \geq \gamma_{0}$,  $ \mathbf{h}= e^{\gamma t } ( \mathrm{I}-\mathbf{Op}_{\mathcal{P}}^{\eps,\gamma})^{-1}
(  e^{-\gamma t} \mathbb{F})$
 vanishes for $t \leq T$.

\end{itemize}

\end{prop}

\begin{proof}

For i), we consider the symbol $c=\frac{\mathcal{P}}{1-\mathcal{P}}$. 
This is a good symbol for our pseudodifferential calculus with parameters, as by the  quantum Penrose stability condition~\eqref{QCrit} and thanks to the same arguments as in the proof  of Lemma~\ref{lem:a-symbol} for $\mathcal{P}$, we have
 $$\vert c\vert_{k_d+1,0} \leq \Lambda( c_0^{-1}, \norm{f^0}_{\mathcal{H}^m_{r}}), \qquad \vert c\vert_{k_d,1} \leq \Lambda( c_0^{-1}, \norm{f^0}_{\mathcal{H}^m_{r}}).
$$ 
Therefore, owing to Proposition~\ref{prop:pseudoparam}, 
\begin{equation}
\label{eq:pseu1}
\|\mathbf{Op}_{\frac{\mathcal{P}}{1-\mathcal{P}}}^{\eps,\gamma}\|_{\mathscr{L}(L^2(\R \times \R^d))} \leq \Lambda( c_0^{-1}, \norm{f^0}_{\mathcal{H}^m_{r}}).
\end{equation}
Let us  consider 
$$
\left(\mathrm{I}+\mathbf{Op}_{\frac{\mathcal{P}}{1-\mathcal{P}}}^{\eps,\gamma}\right)(\mathrm{I}-\mathbf{Op}_{\mathcal{P}}^{\eps,\gamma}) = \left[ \mathrm{I} - \left(\mathbf{Op}_{\frac{\mathcal{P}}{1-\mathcal{P}}}^{\eps,\gamma}\mathbf{Op}_{\mathcal{P}}^{\eps,\gamma}- \mathbf{Op}_{\frac{\mathcal{P}^2}{1-\mathcal{P}}}^{\eps,\gamma}\right)\right]. 
$$
Again by Proposition~\ref{prop:pseudoparam}, it holds
\begin{equation}
\label{eq:pseu2}
\left\|\mathbf{Op}_{\frac{\mathcal{P}}{1-\mathcal{P}}}^{\eps,\gamma}\mathbf{Op}_{\mathcal{P}}^{\eps,\gamma}- \mathbf{Op}_{\frac{\mathcal{P}^2}{1-\mathcal{P}}}^{\eps,\gamma}\right\|_{\mathscr{L}(L^2(\R \times \R^d))} \leq \frac{1}{\gamma} \Lambda(c_0^{-1},\norm{f^0}_{\mathcal{H}^m_{r}}).
\end{equation}
We deduce that  there exists $\gamma_0>0$ depending only on $\norm{f^0}_{\mathcal{H}^m_{r}}$ and $c_0$ such that, for $\gamma\geq \gamma_0$, the operator $\left[ \mathrm{I} - \left(\mathbf{Op}_{\frac{\mathcal{P}}{1-\mathcal{P}}}^{\eps,\gamma}\mathbf{Op}_{\mathcal{P}}^{\eps,\gamma}- \mathbf{Op}_{\frac{\mathcal{P}^2}{1-\mathcal{P}}}^{\eps,\gamma}\right)\right]$ is invertible on $L^2(\R \times \R^d)$, and so  $(\mathrm{I}-\mathbf{Op}_{\mathcal{P}}^{\eps,\gamma})$ is left-invertible. Similarly it is also right-invertible and hence it is invertible. The claimed estimate follows from~\eqref{eq:pseu1}--\eqref{eq:pseu2}.

For $ii)$, we shall crucially use the following Lemma which relies on the Paley-Wiener Theorem.
\begin{lem}
\label{lemPaleyWiener}
Consider $a(x, \zeta)$ a symbol such that $|a|_{k_{d},0}<+\infty$ is  finite, assume in addition that $a(x, \zeta)= \mathfrak{a}(x, \xi, \tau-i \gamma)$
 where  $\mathfrak{a}(x, \xi, z)$ is holomorphic in $\mbox{Im } z <0$,  continuous
  on $ \mbox{Im } z \leq 0$.
  Then, for   every 
   $ F \in L^2(\mathbb{R}\times \mathbb{R}^d)$ such that $F_{|t<0}=0$, 
   we have  for every $\eps \in (0, 1]$ that 
  $u=  e^{  \gamma t }\mathbf{Op}_{a}^{\eps,\gamma} (e^{-\gamma t } F)$
  is independent of $\gamma \geq 0 $. Moreover, we have $u \in  L^2(\mathbb{R} \times \mathbb{R}^d)$
   and $u_{|t<0}=0$.
\end{lem}
Note that $u$ as defined in the Lemma depends on $\eps$ but since $\eps$ plays only the role of a parameter, we do not stress this dependence.
Let us postpone the proof of this Lemma and first finish the proof of ii).

From the proof of i), we have for $\gamma \geq\gamma_0$,
$$h_\gamma:=( \mathrm{I}-\mathbf{Op}_{\mathcal{P}}^{\eps,\gamma})^{-1} ( e^{-\gamma t} \mathbb{F})
=  \left( \mathrm{I} -  \mathbf{R}_{\gamma}\right)^{-1}  \left(\mathrm{I}+\mathbf{Op}_{\frac{\mathcal{P}}{1-\mathcal{P}}}^{\eps,\gamma}\right) ( e^{-\gamma t} \mathbb{F}), $$
where we have set
\begin{equation}
\label{Rgammadef} \mathbf{R}_{\gamma}=  \left(\mathbf{Op}_{\frac{\mathcal{P}}{1-\mathcal{P}}}^{\eps,\gamma}\mathbf{Op}_{\mathcal{P}}^{\eps,\gamma}- \mathbf{Op}_{\frac{\mathcal{P}^2}{1-\mathcal{P}}}^{\eps,\gamma}\right).
\end{equation}
By definition, $h_\gamma$ is the unique $L^2$ solution to
\begin{equation}
\label{equationPgamma} ( \mathrm{I}-\mathbf{Op}_{\mathcal{P}}^{\eps,\gamma}) h_{\gamma}= e^{-\gamma t} \mathbb{F}.
\end{equation}
Thanks to the expression \eqref{DefASymb}, we observe that we can write
$\mathcal{P}(x, \gamma, \tau, \xi )= \mathfrak{P}(x, \xi, \tau-i \gamma)$ where
 $ \mathfrak{P}(x, \xi,\cdot)$ is holomorphic in $\mbox{Im } z <0$.
From the  Penrose condition and Lemma  \ref{lem:a-symbol},  we can use 
Lemma \ref{lemPaleyWiener} with  the symbol ${ \mathcal{P} \over  1 - \mathcal{P} }$, 
 to  first get that
 $G_{\gamma}=  \left(\mathrm{I}+\mathbf{Op}_{\frac{\mathcal{P}}{1-\mathcal{P}}}^{\eps,\gamma}\right) ( e^{-\gamma t} \mathbb{F})$ is such that  $G_{\gamma}= e^{-\gamma t} G$
 with $G \in L^2(\mathbb{R} \times \mathbb{R}^d)$ and $G_{|t<0}= 0$.
 
   Then since the operator norm of $\mathbf{R}_{\gamma}$ is small enough, we can write
   $$ (\mathrm{I}-  \mathbf{R}_{\gamma})^{-1} G_{\gamma}
   =  \sum_{n \geq 0} (\mathbf{R}_{\gamma})^n (e^{-\gamma t}  G).$$
   By using  \eqref{Rgammadef} and Lemma \ref{lemPaleyWiener} repeatedly with the symbols
   $\mathcal{P}/(1- \mathcal{P})$, $\mathcal{P}$ and $\mathcal{P}^2/(1- \mathcal{P})$, 
   we get that
   $$  (\mathrm{I}-  \mathbf{R}_{\gamma})^{-1} G_{\gamma}
   =  \sum_{n \geq 0} u^{(n)}_{\gamma},$$
   where  $u^{(n)}_{\gamma}= e^{- \gamma t  } u^{(n)}$ with 
    $u^{(n)} \in L^2(\mathbb{R}\times \mathbb{R}^d)$  and $u^{(n)}_{|t<0}=0$. 
    Since the series converges in $L^2$  for $\gamma \geq \gamma_{0}$, this yields in particular that
     $h_\gamma=(\mathrm{I}-  \mathbf{R}_{\gamma})^{-1} G_{\gamma}$ vanishes for negative times.

Since $h_{\gamma_{0}}$ vanishes for negative times,
we have  for $\gamma \geq \gamma_{0}$ that $ e^{-(\gamma- \gamma_{0})}h_{\gamma_{0}}
 \in L^2(\mathbb{R} \times \mathbb{R}^d)$ and 
 we can use the conjugation formula \eqref{conjugform}
 to get that  for $\gamma \geq \gamma_{0}$, 
 $$ ( \mathrm{I}-\mathbf{Op}_{\mathcal{P}}^{\eps,\gamma}) (e^{-(\gamma- \gamma_{0}) t} h_{\gamma_{0} } )
 = e^{-(\gamma- \gamma_{0}) t }  ( \mathrm{I}-\mathbf{Op}_{\mathcal{P}}^{\eps,\gamma_{0} }) h_{\gamma_{0}}= e^{-\gamma t} \mathbb{F}.$$
 By  uniqueness of the $L^2$ solution of \eqref{equationPgamma}, 
 we thus deduce that 
  $ h_{\gamma}=  e^{-(\gamma- \gamma_{0}) t}h_{\gamma_{0}}$.
   This ends the proof of ii).
 
 Let us prove iii).   From i) and ii), we first get that $\mathbf{h}=e^{\gamma t} h_\gamma$ vanishes for negative times,  is independent of $\gamma$ for $\gamma \geq \gamma_{0}$,  and such that
 $$ \| e^{-\gamma t } \mathbf{h} \|_{L^2(\mathbb{R} \times \mathbb{R}^d)} \leq C 
  \| e^{-\gamma t } \mathbb{F}\|_{L^2(\mathbb{R} \times \mathbb{R}^d)}
   \leq C   \| e^{-\gamma t } \mathbb{F}\|_{L^2([T, + \infty)  \times \mathbb{R}^d)},$$
 since $\mathbb{F}$ vanishes for $t  \leq T$,  with  $C$ independent of $\gamma \geq \gamma_{0}.$
  This yields 
   $$  \|  \mathbf{h} \|_{L^2((0, T)  \times \mathbb{R}^d)}
    \leq C   \| e^{-\gamma (T- t )} \mathbb{F}\|_{L^2([T, + \infty)  \times \mathbb{R}^d)}.$$
    By letting $\gamma $ go to infinity, the right-hand side tends to zero by dominated convergence
     and consequently,  $h= 0$ also on $(0, T)$.

\bigskip   

It only remains to prove Lemma \ref{lemPaleyWiener}.
\begin{proof}[Proof of Lemma \ref{lemPaleyWiener}]
We first consider $x$ and $\xi$ as parameters. For almost every $x$, we have that  $F(\cdot, x) \in L^2(\mathbb{R})$ and that it vanishes for negatives times therefore its Fourier transform in time $\widehat{F}(\tau, x)$ extends
into an holomorphic function  on $\mbox{Im z}<0$ such that
$\sup_{\gamma>0} \| \widehat{F}(\cdot- i \gamma, x) \|_{L^2(\mathbb{R})} <+\infty$.
By the boundedness assumption on the symbol, we also have that
$$ \sup_{\gamma>0} \| \mathfrak{a}( x,  \eps \xi, \eps( \cdot- i  \gamma)) \widehat{F}(\cdot- i \gamma, x) \|_{L^2(\mathbb{R})} <+\infty.$$ 
By the Paley-Wiener Theorem,  there therefore exists a function
$ H_{x, \xi} \in L^2(\mathbb{R})$ which vanishes for negative times such that
$\mathfrak{a}(  x, \eps \xi,\eps( \cdot- i  \gamma)) \widehat{F}(\cdot- i \gamma, x) = \widehat{H}_{x, \xi}(\tau- i\gamma).$  We deduce from the definition of the pseudodifferential operator that 
$$ \mathbf{Op}_{a}^{\eps,\gamma} (e^{-\gamma t } F)
= (2\pi)^{-d}\int_{\xi} e^{ix \cdot \xi}  e^{-\gamma t} H_{x, \xi}(t)\, d\xi$$
and thus that  $  \mathbf{Op}_{a}^{\eps,\gamma} (e^{-\gamma t } F_{\gamma})$ 
vanishes for negative times.
Moreover, we also get  from the last expression   that
\begin{align*} \mathbf{Op}_{a}^{\eps,\gamma} (e^{-\gamma t } F)
&= (2\pi)^{-d- 1}\int_{\xi} \int_\tau e^{ix \cdot \xi}  e^{i \tau t}  e^{-\gamma t}\widehat{H}_{x, \xi} (\tau)
\, d\tau d \xi \\
&= (2\pi)^{-d- 1} \int_{\xi} \int_\tau e^{ix \cdot \xi}  e^{i \tau t}  e^{-\gamma t}
a( x,0, \eps \tau, \eps \xi) \widehat{F}(\cdot, x)\, d\tau d \xi
=  e^{-\gamma t }\mathbf{Op}_{a}^{\eps, 0}  F.
\end{align*}
 This yields that  $u= e^{\gamma t} \mathbf{Op}_{a}^{\eps,\gamma} (e^{-\gamma t } F_{\gamma})$
 is independent of $\gamma$  and such that  $u \in L^2(\mathbb{R} \times \mathbb{R}^d)$
  since $F \in L^2(\mathbb{R} \times \mathbb{R}^d)$ and $\mathbf{Op}_{a}^{\eps, 0} $
  is continuous on $L^2(\mathbb{R} \times \mathbb{R}^d)$.

\end{proof}

    We are now in position to prove Proposition \ref{prop:rho-final}.

 \subsection{Proof of Proposition  \ref{prop:rho-final} }

We have to study the equation \eqref{EqFFinal} which reads by using the definition \eqref{defLf0}, 
\begin{equation}
\label{cellededepart} \partial^\alpha_x \rho(t) = \mathcal{L}_{\eps, f^0} \partial^\alpha_x \rho+ R, \quad | \alpha | =m
\end{equation}
where $R$ is a remainder and thus enjoys the estimate  \eqref{defremainder}.

\noindent $\bullet$ {\bf Step 1.}
We shall first prove the estimate \eqref{estRho} for $r= 0$, that is to say
\begin{equation}
\label{pourr=0} \| \partial^\alpha_x \rho \|_{L^2((0, T) \times \mathbb{R}^d)}
 \lesssim  \Lambda (c_{0}^{-1},  \norm{f^0}_{\mathcal{H}^m_{r}}, T)
  \| R \|_{L^2((0, T) \times \mathbb{R}^d)}.
  \end{equation}
Let us define $\mathbf{h}_{1}$  as $\mathbf{h}_1 = \partial_{x}^{\alpha}\rho$ on $[0,T]$ and $\mathbf{h}_1 =0$ on $(-\infty,0)\cup (T,+\infty)$ so that
 $\mathbf{h}_1 \in L^2(\R\times \R^d)$.
 Then $\mathbf{h}_1 $ solves for $t \in \mathbb{R}$  the equation
 $$ \mathbf{h}_{1} = \mathcal{L}_{\eps, f^0} \mathbf{h}_{1}+ \mathbf{R}_{1},$$
 which can be seen as the  definition of the  source term $\mathbf{R}_{1}$.
  Since $\mathbf{h}_{1}$ vanishes for negative times  and is in $L^2(\mathbb{R} \times \mathbb{R}^d)$,
  we also have that 
  $\mathbf{R}_{1} \in L^2(\mathbb{R} \times \mathbb{R}^d)$. Indeed,  by using
   Lemma \ref{Eqpseudo}, we have  $ \mathcal{L}_{\eps, f^0} \mathbf{h}_{1}
   = \mathbf{Op}^{\eps, 0}_{\mathcal{P}} \mathbf{h}_{1}$ and $  \mathbf{Op}^{\eps, 0}_{\mathcal{P}}$
   is continuous on $L^2(\mathbb{R} \times \mathbb{R}^d)$.
 Moreover we have that $\mathbf{R}_{1}$ coincides with $R$ on $[0, T]$  and vanishes for negative times.
 By setting $h_{1}=  e^{-\gamma t} \mathbf{h}_{1}$
 and by using again Lemma \ref{Eqpseudo}, we get that $h_{1}$ is a $L^2$ 
   solution of 
 \begin{equation}
 \label{defh1}
  h_{1}=   \mathbf{Op}^{\eps, \gamma}_{\mathcal{P}} h_{1}  + e^{-\gamma t } \mathbf{R}_{1}
  \end{equation}
 which vanishes for negative times.
 
 We can also define a source term $\mathbf{R_{2}}$ by setting $\mathbf{R_{2}}= R$ on $[0, T]$ where
  $R$ is the original source term  in \eqref{cellededepart} and $\mathbf{R_{2}}= 0$
   for $t \leq 0$ and $t \geq T$.  Thanks to Proposition \ref{MajhFin},  i),  for $\gamma \geq \gamma_{0}$
 we can set 
  \begin{equation}
  \label{defh2}
  h_{2}= ( \mathrm{I}-\mathbf{Op}_{\mathcal{P}}^{\eps,\gamma})^{-1}
(  e^{-\gamma t} \mathbf{R_{2}})
\end{equation}
  and get
 \begin{equation}
 \label{labonneestim} 
 \norm{h_{2}}_{L^2(\mathbb{R}\times \mathbb{R}^d)}\leq \Lambda(c_0^{-1}, \norm{f^0}_{\mathcal{H}^m_{r}}, T) \norm{e^{-\gamma t } \mathbf{R_{2}}}_{L^2(\mathbb{R}\times \mathbb{R}^d)}
  =  \Lambda(c_0^{-1}, \norm{f^0}_{\mathcal{H}^m_{r}}, T) \norm{e^{-\gamma t } \mathbf{R_{2}}}_{L^2((0, T)\times \mathbb{R}^d)}.
\end{equation}
We also know from  Proposition \ref{MajhFin},  ii) that $h_{2}$ vanishes for negative times.

Thanks to \eqref{defh1} and \eqref{defh2}, we obtain that $h= h_{1}-h_{2} \in L^2(\mathbb{R} \times \mathbb{R}^d)$ vanishes for negative times and solves
$$ h =   \mathbf{Op}^{\eps, \gamma}_{\mathcal{P}} h +  e^{-\gamma t }(  \mathbf{R}_{1} - \mathbf{R}_{2})
$$
with  $\mathbf{R}_{1} - \mathbf{R}_{2} \in   L^2(\mathbb{R} \times \mathbb{R}^d)$
 and $\mathbf{R}_{1} - \mathbf{R}_{2}= 0$ for $t \leq T$.
  Thanks to Proposition \ref{MajhFin} iii), we get that $h_{1}= h_{2}$ on $[0, T]$; 
  this yields that $e^{-\gamma t}  \partial^\alpha_x \rho$ also enjoys the estimate
   \eqref{labonneestim}, hence we get \eqref{pourr=0}.
 
 \bigskip 
   
\noindent $\bullet$ {\bf Step 2.} We will finally get by induction  that
 $$   \| \partial^\alpha_x \rho \|_{L^2(0, T; H_{r}^0)}
 \lesssim  \Lambda (c_{0}^{-1},  \norm{f^0}_{\mathcal{H}^m_{r}},T)
  \| R \|_{L^2((0, T) \times \mathbb{R}^d)}.$$
  Indeed,  from \eqref{cellededepart}, we have that for every $j \in \llbracket 1, d\rrbracket$, 
  $$ \eps \partial_{x_{j}} \partial^\alpha_x \rho=\mathcal{L}_{\eps, f^0}  (\eps \partial_{x_{j}} \partial^\alpha_x \rho)
  + \eps  \mathcal{L}_{\eps,\partial_{x_{j}} f^0} \partial^\alpha_x \rho+ \eps \partial_{x} R.$$
  From Lemma \ref{Eqpseudo}, we obtain that 
  $$ 
  \mathcal{L}_{\eps,\partial_{x_{j}} f^0} \partial^\alpha_x \rho
  = \mathbf{Op}^{\eps, 0}_{\partial_{x_{j}}\mathcal{P}} h,
  $$
  where we have set $h= \partial^\alpha_x \rho$ on $[0, T]$ and $0$ elsewhere.
  From Proposition~\ref{prop:pseudoparam} and Lemma \ref{lem:a-symbol}, we know that    $\mathbf{Op}^{\eps, 0}_{\partial_{x_{j}}\mathcal{P}}$
   is continuous on $L^2$ and thus we get from \eqref{pourr=0} that
   $$ \|\mathcal{L}_{\eps,\partial_{x_{j}} f^0} \partial^\alpha_x \rho\|_{L^2(0, T; L^2(\mathbb{R}^d))}
    \leq   \Lambda (c_{0}^{-1},  \norm{f^0}_{\mathcal{H}^m_{r}}) \| R \|_{L^2((0, T) \times \mathbb{R}^d)}.$$
  We consequently obtain that $\eps \partial_{x_{j}} \partial^\alpha_x \rho$ solves 
  $$   \eps \partial_{x_{j}} \partial^\alpha_x \rho=\mathcal{L}_{\eps, f^0}  (\eps \partial_{x_{j}} \partial^\alpha_x \rho)
   + R_{1}$$
   where the source term $R_{1}$ enjoys the estimate
   $$  \|R_{1}\|_{L^2(0, T; L^2(\mathbb{R}^d))}
    \leq   \Lambda (c_{0}^{-1},  \norm{f^0}_{\mathcal{H}^m_{r}},T)  \|R\|_{L^2(0, T; H^0_{1})}.$$
    Hence, from {\bf Step 1},  we deduce that
    $$   \| \partial^\alpha_x \rho \|_{L^2((0, T; H^0_{1})}
 \lesssim   \Lambda (c_{0}^{-1},  \norm{f^0}_{\mathcal{H}^m_{r}},T)
  \| R \|_{L^2(0, T; H^0_{1})}.$$
  The general case
  $$ \| \partial^\alpha_x \rho \|_{L^2((0, T; H^0_{r})}
 \lesssim   \Lambda (c_{0}^{-1},  \norm{f^0}_{\mathcal{H}^m_{r}},T)
  \| R \|_{L^2(0, T; H^0_{r})}.$$
   follows similarly by induction,  since $m >5 + r + {d \over 2}$.
   Since $R$ is a remainder, we finally get  \eqref{estRho} by recalling \eqref{defremainder}. 
   The proof of Proposition \ref{prop:rho-final} is finally complete.

   \end{proof}

To conclude this section, it only remains to prove  Lemma~\ref{lem:a-symbol}.

\subsection{Proof of Lemma~\ref{lem:a-symbol}}

Let us first treat the first estimate, which is fairly straightforward. 
For all $\vert \alpha\vert \leq k$, using the inequality $|\sin x|\leq |x|$, we get 
\begin{align*}
&\vert\left(\mathcal{F}_x\partial^{\alpha}_x \mathcal{P} \right) (\kappa,\zeta)\vert = \bigg\vert   \int_0^{+\infty}e^{-(\gamma +i\tau) s}  2 \sin\left( s\frac{\vert \eta \vert^2}{2}\right)  \widehat{V}( \eta) \cdot   \mathcal{F}_{x,v}(\partial^\alpha_x f^0)(\kappa,s\eta) ds\bigg\vert
	\\&\qquad \lesssim \left( \int_v  (1+\vert v\vert )^{2 (k_{d}-1)}\left[ \vert( \mathcal{F}_x  \partial^{\alpha}_x (I-\Delta_v)^2 f^0)(\kappa,v)\vert \vert  \widehat{V}( \eta) \vert   \right]^2 dv \right)^{1/2} \int_0^{+\infty} \frac{s\vert \eta\vert^2}{(1+ \vert s  \eta\vert^2)^2}ds
	\\&\qquad \lesssim \left( \int_v  (1+\vert v\vert )^{2 (k_{d}-1)}\left[ \vert ( \mathcal{F}_x  \partial^{\alpha}_x (I-\Delta_v)^2 f^0)(\kappa,v)\vert \vert  \widehat{V}( \eta) \vert \right]^2 dv \right)^{1/2} \int_0^{+\infty} \frac{1}{(1+s^2)}ds,
\end{align*}
where we recall $k_d = \lfloor d/2\rfloor +2$.
Consequently, by the Bessel-Parseval identity and the fact that $\widehat{V}$ is bounded,
$$
\norm{\mathcal{F}\partial^{\alpha}_x \mathcal{P} }_{L^2(\mathbb{R}^d; L^\infty_\zeta)}\leq C \norm{
\langle v \rangle^{k_{d}-1} f^0}_{{H}^{k+4}_{x,v}}
 \leq C \norm{ f^0}_{\H^{k+4}_{k_{d}-1}},
$$
where the last inequality comes from \eqref{embedH01}.
 
Let us focus on the second item. We want to estimate of $\norm{\gamma \mathcal{F}_x(\partial^{\alpha}_x \nabla_\xi c)}_{L^2(\mathbb{R}^d;L^{\infty}_\zeta)}$ for all $|\alpha|\leq k $. Denote $\xi=(\tau,\eta)$. 
We have  
\begin{multline*}
\frac{1}{2}\nabla_\xi \mathcal{P}=   i\left(\int_0^{+\infty}  e^{-(\gamma +i\tau) s} \mathcal{F}_{v}f^0\left(x,s\eta\right)  s  \sin\left( s \frac{\vert \eta \vert^2}{2}\right)  \widehat{V}( \eta)    ds\right) e_0 
	\\ + \sum_{j=1}^{d} \bigg(-\int_0^{+\infty}  e^{-(\gamma +i\tau) s}  \partial_{\eta_j} \mathcal{F}_{v} f^0 \left(x,s\eta\right)  s  \sin\left(s \frac{\vert \eta \vert^2}{2}\right)  \widehat{V}( \eta)   ds 
 	\\  - \int_0^{+\infty}  e^{- (\gamma +i\tau) s } \mathcal{F}_{v}f^0\left(x,s\eta\right)  \ s  \eta_j \cos\left( s \frac{\vert \eta \vert^2}{2}\right)  \widehat{V}( \eta)   ds 
	\\  - \int_0^{+\infty}  e^{- (\gamma +i\tau) s}   \mathcal{F}_{v}f^0\left(x,s\eta\right)   \sin\left(s \frac{\vert \eta \vert^2}{2}\right)   \partial_{\eta_j}  \widehat  V( \eta)   ds\bigg) e_{j},
\end{multline*}
where $(e_j)_{j\in\llbracket0,d\rrbracket}$ represents the canonical basis of $\mathbb{R}^{d+1}$. 
We can then make the change of variable $s'=s\langle \zeta\rangle$, where $\langle \zeta\rangle=\left( \gamma^2+\tau^2+\vert \eta\vert ^2 \right)^{1/2}$,  in this formula.
Consequently, we are left to consider four types of symbols that we denote by 
\begin{align*}
I^{\alpha}_1&=  \frac{\gamma}{\langle \zeta\rangle}\int_0^{+\infty}  e^{-\frac{(\gamma +i\tau)}{\langle \zeta\rangle} s} \partial^{\alpha}_x\mathcal{F}_{v}f^0\left(x,s\frac{\eta}{\langle \zeta\rangle}\right)  \frac{s}{ \langle \zeta\rangle}  \sin\left( \frac{s}{\langle \zeta\rangle} \frac{\vert \eta \vert^2}{2}\right)  \widehat{V}( \eta) ds, \\
I^{\alpha}_2&=  \frac{\gamma}{\langle \zeta\rangle}\int_0^{+\infty}  e^{-\frac{(\gamma +i\tau)}{\langle \zeta\rangle} s}  \partial^{\alpha}_x\partial_{\eta_{j}}\mathcal{F}_{v}f^0\left(x,s\frac{\eta}{\langle \zeta\rangle}\right)  \frac{s}{ \langle \zeta\rangle}  \sin\left( \frac{s}{\langle \zeta\rangle} \frac{\vert \eta \vert^2}{2}\right)  \widehat{V}( \eta)   ds, \\
I^{\alpha}_3&=  \frac{\gamma}{\langle \zeta\rangle} \int_0^{+\infty}  e^{-\frac{(\gamma +i\tau)}{\langle \zeta\rangle} s} \partial^{\alpha}_x \mathcal{F}_{v}f^0\left(x,s\frac{\eta}{\langle \zeta\rangle}\right)  \frac{s \eta_j}{ \langle \zeta\rangle}  \cos\left( \frac{s}{\langle \zeta\rangle} \frac{\vert \eta \vert^2}{2}\right)  \widehat{V}(\eta) ds, \\
I^{\alpha}_4&=  \frac{\gamma}{\langle \zeta\rangle} \int_0^{+\infty}  e^{-\frac{(\gamma +i\tau)}{\langle \zeta\rangle} s}\partial^{\alpha}_x  \mathcal{F}_{v}f^0\left(x,s\frac{\eta}{\langle \zeta\rangle}\right)   \sin\left( \frac{s}{\langle \zeta\rangle} \frac{\vert \eta \vert^2}{2}\right)  \partial_\eta \widehat{V}( \eta)  ds.
\end{align*}

 \noindent$\bullet$ {\bf Estimate for $I^{\alpha}_1$.}
Let us rewrite $I^{\alpha}_1$ with the new variables $(\widetilde{\gamma},\widetilde{\tau},\widetilde{\eta})=(\gamma,\tau,\eta)/\langle \zeta\rangle$  on the unit sphere:
$$
I^{\alpha}_1(x,\widetilde{\gamma},\widetilde{\tau},\widetilde{\eta},\langle \zeta\rangle)= \widetilde{\gamma} \int_0^{+\infty}  e^{-(\widetilde{\gamma} +i\widetilde{\tau}) s} \partial^{\alpha}_x\mathcal{F}_{v}f^0\left(x,s\widetilde{\eta}\right)  \frac{s}{ \langle \zeta\rangle}  \sin\left( \langle \zeta\rangle s \frac{\vert \widetilde{\eta} \vert^2}{2}\right)  \widehat{V}({\langle \zeta\rangle}\widetilde{ \eta}) ds.
$$
We first  consider the case  $\vert \widetilde{\eta}\vert\geq 1/2$, for which we can follow
 the same lines as in the proof of the first item:
\begin{align*}
\vert \mathcal{F}_x I^{\alpha}_1 &(\kappa,\widetilde{\gamma},\widetilde{\tau},\widetilde{\eta},\langle \zeta\rangle) \vert \\
&= \bigg\vert \left(  \int_0^{+\infty} e^{-(\widetilde{\gamma} +i\widetilde{\tau}) s}   (\mathcal{F}_x \partial^\alpha_x)\left(\mathcal{F}_v f^0\right)  \left(\kappa,s\widetilde{\eta}\right)  \frac{s}{ \langle \zeta\rangle}  \sin\left( \langle \zeta\rangle s \frac{\vert \widetilde{\eta} \vert^2}{2}\right)  \widehat{V}({\langle \zeta\rangle}\widetilde{ \eta}) ds \right)\bigg\vert
	\\& \lesssim \left( \int_v    (1+\vert v\vert )^{2 k_d}  \vert \mathcal{F}_x  \partial^{\alpha}_x (I-\Delta_v)^2 f^0(\kappa,v)\vert^2  dv \right)^{1/2} \int_0^{+\infty} \frac{s^2\vert \widetilde{\eta}\vert^2}{(1+\vert s \widetilde{\eta}\vert^2)^2}ds
	\\& \lesssim \left( \int_v   (1+\vert v\vert )^{2 k_d} \vert  \mathcal{F}_x  \partial^{\alpha}_x (I-\Delta_v)^2 f^0(\kappa,v)\vert^2 dv \right)^{1/2} \int_0^{+\infty} \frac{|\widetilde{\eta}|}{(1+s^2\vert \widetilde{\eta}\vert ^2)}ds.
\end{align*}
On the other hand, if $\vert \widetilde{\eta}\vert <1/2$, we must have $\vert \widetilde{\gamma}\vert^2+\vert \widetilde{\tau}\vert ^2\geq 3/4$. Writing $ e^{-(\widetilde{\gamma} +i\widetilde{\tau}) s}= \frac{-1}{\widetilde{\gamma} +i\widetilde{\tau}}\partial_s ( e^{-(\widetilde{\gamma} +i\widetilde{\tau}) s})$, this allows to perform an integration by parts in $s$ to obtain
\begin{align*}
\vert \mathcal{F}_x I^{\alpha}_1 (\kappa,\widetilde{\gamma},\widetilde{\tau},\widetilde{\eta},\langle \zeta\rangle) \vert 	&\lesssim  \int_0^{+\infty}\bigg\vert  (\mathcal{F}_x \partial^\alpha_x)\left(\partial_\eta\mathcal{F}_v f^0\right) \left(\kappa,s\widetilde{\eta}\right)  \frac{s\vert \widetilde{\eta}\vert}{ \langle \zeta\rangle}  \sin\left( \langle \zeta\rangle s \frac{\vert \widetilde{\eta} \vert^2}{2}\right)  \widehat{V}({\langle \zeta\rangle}\widetilde{ \eta}) \bigg\vert ds
	\\&+   \int_0^{+\infty}\bigg\vert  (\mathcal{F}_x \partial^\alpha_x)\left(\mathcal{F}_v f^0\right) \left(\kappa,s\widetilde{\eta}\right)  \frac{1}{\langle \zeta\rangle}  \sin\left( \langle \zeta\rangle s \frac{\vert \widetilde{\eta} \vert^2}{2}\right)  \widehat{V}({\langle \zeta\rangle}\widetilde{ \eta}) \bigg\vert ds
	\\&+   \int_0^{+\infty}\bigg\vert  (\mathcal{F}_x \partial^\alpha_x)\left(\mathcal{F}_v f^0\right) \left(\kappa,s\widetilde{\eta}\right)  s\vert \eta\vert ^2  \cos\left( \langle \zeta\rangle s \frac{\vert \widetilde{\eta} \vert^2}{2}\right)  \widehat{V}({\langle \zeta\rangle}\widetilde{ \eta}) \bigg\vert ds.
\end{align*}
For the first term, we have 
\begin{multline*}
  \int_0^{+\infty}\bigg\vert (\mathcal{F}_x \partial^\alpha_x)\left(\partial_\eta\mathcal{F}_v f^0\right) \left(\kappa,s\widetilde{\eta}\right)  \frac{2s\vert \widetilde{\eta}\vert}{ \langle \zeta\rangle}  \sin\left( \langle \zeta\rangle s  \frac{\vert \widetilde{\eta} \vert^2}{2}\right)  \widehat{V}({\langle \zeta\rangle}\widetilde{ \eta}) \bigg\vert  ds
	\\ \lesssim \left( \int_v (1+\vert v\vert )^{2 k_d  }  \left[ \vert( \mathcal{F}_x  \partial^{\alpha}_x (I-\Delta_v)^3 f^0)(\kappa,v)\vert  \right]^2 dv \right)^{1/2} \int_0^{+\infty} \frac{s^2\vert \eta\vert^3}{(1+\vert s \eta\vert^2)^3}ds
\end{multline*}
and similarly for the others, it holds 
\begin{align*}
   \int_0^{+\infty}&\bigg\vert (\mathcal{F}_x \partial^\alpha_x) \left(\mathcal{F}_v f^0\right) \left(\kappa,s\widetilde{\eta}\right)   \frac{2}{ \langle \zeta\rangle}  \sin\left( \langle \zeta\rangle s \frac{\vert \widetilde{\eta} \vert^2}{2}\right)  \widehat{V}({\langle \zeta\rangle}\widetilde{ \eta}) \bigg\vert ds
	\\ &+   \int_0^{+\infty}\bigg\vert  (\mathcal{F}_x \partial^\alpha_x)\left(\mathcal{F}_v f^0\right) \left(\kappa,s\widetilde{\eta}\right)  s\vert \eta\vert ^2  \cos\left( \langle \zeta\rangle s \frac{\vert \widetilde{\eta} \vert^2}{2}\right)  \widehat{V}({\langle \zeta\rangle}\widetilde{ \eta}) \bigg\vert  ds
\\ &\lesssim \left( \int_v (1+\vert v\vert )^{2 (k_d-1)}  \left[ \vert ( \mathcal{F}_x  \partial^{\alpha}_x (I-\Delta_v)^2 f^0)(\kappa,v)\vert \right]^2 dv \right)^{1/2} \int_0^{+\infty} \frac{1}{(1+s^2)}ds.
\end{align*}

 \noindent$\bullet$  {\bf Estimate of $I^{\alpha}_2$.}
The term  $I^{\alpha}_2$ is similar  to $I^{\alpha}_1$, 
we just change $f^0$ into $v_{j} f^0$ and thus 
  we obtain the same type of estimate where we only change the weight in v of order $k_d-1$ into a weight of order $k_d$.

 \noindent$\bullet$ {\bf Estimate of $I^{\alpha}_3$.} We once again write  that
$$
I^{\alpha}_3(x,\widetilde{\gamma},\widetilde{\tau},\widetilde{\eta},\langle \zeta\rangle)= - \int_0^{+\infty} e^{-(\widetilde{\gamma} +i\widetilde{\tau}) s} \partial^{\alpha}_x \mathcal{F}_{v}f^0\left(x,s\widetilde{\eta}\right)  s\widetilde{\eta} \cos\left( s{\langle \zeta\rangle}  \frac{\vert \widetilde{\eta} \vert^2}{2}\right)   \widehat{V}({\langle \zeta\rangle}\widetilde{ \eta}) ds.
$$
Here we must be more careful about the precise structure of the integrand and use that the $\cos$ term is  oscillatory.
Since $$\widetilde{\gamma} e^{-\left[\widetilde{\gamma} +i\left(\widetilde{\tau} \pm \langle \zeta\rangle  \frac{\vert \widetilde{\eta} \vert^2}{2}\right)\right] s}= -  \frac{\widetilde{\gamma} }{\widetilde{\gamma} +i(\widetilde{\tau} \pm \langle \zeta\rangle  \frac{\vert \widetilde{\eta} \vert^2}{2})} \partial_s \left( e^{-\left[\widetilde{\gamma} +i\left(\widetilde{\tau} \pm \langle \zeta\rangle  \frac{\vert \widetilde{\eta} \vert^2}{2}\right)\right] s} \right),$$
and $\left| \frac{\widetilde{\gamma} }{\widetilde{\gamma} +i(\widetilde{\tau} \pm \langle \zeta\rangle  \frac{\vert \widetilde{\eta} \vert^2}{2})}\right| \leq 1$, we have by integration by parts
\begin{align*}
\vert \mathcal{F}_x I^{\alpha}_3 (\kappa,\widetilde{\gamma},\widetilde{\tau},\widetilde{\eta},\langle \zeta\rangle) \vert 	
	&\lesssim  \int_0^{+\infty}\bigg\vert (\mathcal{F}_x \partial^\alpha_x)\left(\partial_\eta \mathcal{F}_v f^0\right) \left(\kappa,s\widetilde{\eta}\right)   s\vert \widetilde{\eta}\vert^2   \bigg\vert  ds
	\\&+   \int_0^{+\infty}\bigg\vert  (\mathcal{F}_x \partial^\alpha_x)\left(\mathcal{F}_v f^0\right) \left(\kappa,s\widetilde{\eta}\right) \vert \widetilde \eta \vert     \bigg\vert ds,
\end{align*}
which can be estimated as above.

 \noindent$\bullet$ {\bf Estimate for $I^{\alpha}_4$.}  We estimate this integral as in the previous item, we split
 the sin term, regroup the exponentials and  integrate by parts in $s$.
 
\bigskip

Summing up the four estimates,  taking the $L^2$ norm in $\kappa$
and using again \eqref{embedH01}
we  obtain that 
$$
|\mathcal{P}|_{k,1} \leq C \norm{f^0}_{\mathcal{H}^{k+6}_{k_d}}.
$$
This ends the proof of  the lemma.

\section{End of the proof}
\label{sec:end}

\subsection{Proof of  Theorem~\ref{thm:main2}}We are in position to close the bootstrap argument initiated in Section \ref{bootstrap}. We start by fixing  $T(M)$ small enough such that all the results from the previous sections hold for $T \in (0,\min(T_\eps,T(M))]$. 
By Lemma \ref{PE} (for what concerns $f$) and Proposition~\ref{prop:rho-final} (for what concerns $\rho$), for all all $\eps \in (0,1)$ and  $T \in (0,\min(T_\eps,T(M)))$, it holds 
$$
\mathcal{N}_{m,r}(T,f) \leq C M_0 + (T^{1/2} + \eps)  \Lambda( c_{0}^{-1},  \|f^0\|_{\H^m_{r}}, T,M). 
$$
Let us fix $M$ by setting $M= 2C M_0 +1$, so that
$$
\frac{1}{2} M > C M_0.
$$ 
Then by continuity, we can find  $T^{\#}  \in (0,T(M)]$ independent of $\eps$ and $\eps_0 \in (0,1)$ small enough such that for all $T\in[0,T^{\#}]$ and all $\eps \in (0,\eps_0)$,
$$
(T^{1/2} + \eps)  \Lambda( c_{0}^{-1},  \|f^0\|_{\H^m_{r}}, T,M)<  \frac{1}{2}M.
$$
This means that 
for all  $\eps \in (0,\eps_0)$, for all $T\in[0,\min(T_\eps,T^{\sharp}))$, $\mathcal{N}_{m,r}(T,f)<M$ and  
therefore, we must have $T_{\eps}> T^{\#}$ (otherwise this would contradict the definition of $T_\eps$, as we are in the case when $T_\eps<T^\ast$, the maximal time of existence, recall Section~\ref{bootstrap}): Theorem~\ref{thm:main2} is thus proved.

\subsection{Proof of Theorem~\ref{thm:main}}
\label{sec:end-conv}


Let us finally prove the second part of Theorem~\ref{thm:main} as a consequence of Theorem~\ref{thm:main2}. To enhance readability, let us assume that either $\widehat{V}(0)=1$ or $\widehat{V}(0)=-1$. We focus on the case  $\widehat{V}(0)=1$, we will discuss the other case $\widehat{V}(0)=-1$ in the end.
Let us also put back the subscripts $\eps$ in the unknowns of the Wigner equation.
Applying Theorem~\ref{thm:main2} , we  fix $\eps_0>0$,  $M>0$ and $T>0$ such that $\sup_{\eps\in(0,\eps_0]} \mathcal{N}_{m,r}(T,f_\eps)\leq M$.

Recall Definition~\ref{def:Hmr-droit} for the weighted Sobolev space $\mathrm{H}^m_r$. Thanks to~\eqref{embedH01}, we have for all $m,r \in \mathbb{N}$ that  $\| \cdot \|_{\mathrm{H}^m_r} \lesssim \| \cdot \|_{\H^m_r}$.
The family $(f_\eps)_{\eps \in (0,\eps_0)}$ is therefore uniformly bounded in $L^\infty(0,T;\mathrm{H}^{m-1}_r)$ and up to taking a subsequence (that we do not explicitly write for readability), there exists $f \in L^\infty(0,T;\mathrm{H}^{m-1}_r)$ such that $f_\eps$ weakly-$\ast$ converges to $f$ in $L^\infty(0,T;{H}^{m-1}_r)$.  Furthermore, still by weak compactness, we have that $\rho_f \in L^2(0,T;H^m)$.

By (a slight variant of) the estimate~\eqref{estBreste} of Lemma~\ref{lembilinB} (since $m>5+ d/2$) and thanks to~\eqref{embedH02},  we obtain for all $t \in [0,T]$, (using $\rho_{\eps}$ instead of $\rho_{f_{\eps}}$ for the sake
of readability)
$$
\| B_\eps[\rho_\eps, f_\eps] \|_{\mathcal{H}^{m-2}_{r}} \lesssim  \|  \rho_\eps\|_{H^{m-1}_r} \| f_\eps\|_{\mathcal{H}^{m-1}_r}  \lesssim  \| f_\eps\|_{\mathcal{H}^{m-1}_{r}}^2.
$$
Therefore, since $f_\eps$ satisfies the Wigner equation~\eqref{eq:wigner}, we infer that $(\partial_t f_\eps)_{\eps \in (0,\eps_0)}$ is uniformly bounded in $L^\infty(0,T;\mathrm{H}^{m-2}_{r-1})$. By the Ascoli theorem, we first deduce that $f_\eps$ actually converges strongly to $f$ in $L^\infty(0,T;L^2)$, and thus by interpolation that $f_\eps$  converges strongly to $f$ in $L^\infty(0,T;\mathrm{H}^{m-1-\delta}_{r-\delta})$ for all $\delta>0$. 
Moreover,  thanks to \eqref{perturbpenrose}, we  also have that 
$$\sup_{(\gamma,\tau, \eta) \in (0,+\infty)\times \R \times \R^d} |\mathcal{P}_{\mathrm{quant}}(\gamma, \tau, \eta, f_{\eps}(t)) - \mathcal{P}_{\mathrm{quant}}(\gamma, \tau, \eta,  f_{\eps}^0)| \lesssim T \Lambda(T, M), $$
therefore, by taking $T$ smaller if necessary, we  can get that  
 $( f_\eps)_{\eps \in (0,\eps_0)}$ satisfies the $c_{0}/2$ Penrose stability condition uniformly for all $t \in [0,T]$, 
and by passing to the limit, that $f$ also satifies the $ c_{0}/2$ Penrose stability condition uniformly for all $t \in [0,T]$.

Let us now show that $f$ satisfies the Vlasov-Benney equation~\eqref{eq:vlasov-benney-intro} by passing to the limit in the Wigner equation~\eqref{eq:wigner}. The only term that deserves a proper study is $B_\eps[ \rho_{\eps}, f_\eps]$.
We write the decomposition
\begin{align*}
B_\eps[ \rho_\eps, f_\eps] + \nabla_x \rho_f  \cdot \nabla_v f = 
B_\eps[\rho_\eps-\rho_f, f_\eps]  + B_\eps[\rho_f, f_\eps - f] 
+ B_\eps[\rho_f, f] + \nabla_x \rho_f  \cdot \nabla_v f.
\end{align*}
The first two terms are estimated as in the proof of Lemma~\ref{lembilinB}, using~\eqref{usefulK}:
\begin{align*}
\| B_\eps[\rho_\eps-\rho_f, f_\eps] \|_{L^2} 
&\lesssim  \left\|  \int_\eta \frac{1}{\eps} \sin\left( \frac{\eps (\xi_x - \eta)\cdot \xi_v}{2}\right) (\widehat{\rho_{f_\eps}}- \widehat{\rho_f})(\xi_x-\eta)  \widehat{f_\eps}(\eta, \xi_v) \, d\eta \right\|_{L^2_{\xi}} \\
&\lesssim \| \rho_\eps-\rho_f\|_{H^{1}} \| f_\eps\|_{{H}^{m-1}} \lesssim \| f_\eps-f\|_{\mathrm{H}^{1}_{r-1}}  \| f_\eps\|_{{H}^{m-1}},
\end{align*}
where we have used $m-1>1+d/2$ and $r-1>d/2$.
Similarly, we obtain
$$
\| B_\eps[\rho_f, f_\eps - f] \|_{L^2} \lesssim   \| f\|_{\mathrm{H}^{1}_r}  \| f_\eps-f\|_{{H}^{m-1}}.
$$
For the remaining term,  we write 
\begin{multline*}
\| B_\eps[\rho_f, f] +\nabla_x \rho_f  \cdot \nabla_v f \|_{L^2_x} \\
 \lesssim
  \left\|  \int_\eta \left[(\widehat{V}(\eps (\xi_x-\eta))-1) (\xi_x - \eta)\cdot \xi_v\right] \widehat{\rho_f}(\xi_x-\eta)  \widehat{f}(\eta, \xi_v) \, d\eta \right\|_{L^2_{\xi}} \\
 + \left\|  \int_\eta \left[\frac{2}{\eps} \sin\left( \frac{\eps (\xi_x - \eta)\cdot \xi_v}{2}\right) - (\xi_x - \eta)\cdot \xi_v\right] \widehat{\rho_f}(\xi_x-\eta)  \widehat{f}(\eta, \xi_v) \, d\eta \right\|_{L^2_{\xi}}.
\end{multline*}
For the  first term in the right-hand side, we use $|\widehat{V}( x)-\widehat{V}(0)| \lesssim |x|$ and therefore obtain a control by
$
\eps \| f\|_{\mathrm{H}^{2}_{r-1}}  \| f\|_{\mathrm{H}^{m-1}}.
$
For the second one, by the elementary inequality $|\sin x - x| \lesssim  |x|^3$ which holds for all $x \in \R$, we are left to estimate 
$$ 
\left\|\eps^2 \int_{\eta} |\xi_x -\eta|^3 |\widehat{\rho_f}(\xi_x-\eta)|  |\xi_v |^3 |\widehat{f}(\eta, \xi_v)| \, d \eta\right\|_{L^2_{\xi}} \lesssim \eps^2 \| f\|_{\mathrm{H}^{3}_r} \| f\|_{{H}^{m-1}},
$$
where we have  used $m>4+ d/2$. Gathering all pieces together, we conclude that $B_\eps[ \rho_{\eps}, f_\eps] $ converges strongly to $- \nabla_x \rho_f  \cdot \nabla_v f$ in $L^2(0,T; L^2)$, and consequently $f$ satisfies the Vlasov-Benney equation. 

Eventually, by weak compactness, $f \in L^\infty(0,T;\mathrm{H}^{m-1}_r)\cap \mathscr{C}_w([0,T];\mathrm{H}^{m-1}_r)$ and we already know that $\rho_f \in L^2(0,T;H^m)$: since $f$ satisfies the Vlasov-Benney equation, by a standard argument based on an energy estimate, we get that  $f \in \mathscr{C}([0,T]; \mathrm{H}^{m-1}_{r})$.

The following holds.

\begin{lem}
If $f$ satisfies the $c_{0}/2$ quantum Penrose condition on $[0, T]$, we also have that:
\begin{equation}\label{lem:Penrose}
\underset{t \in [0,T], \, x \in \mathbb{R}^d}{\inf}\underset{(\gamma,\tau,\eta)\in(0,+\infty)\times \mathbb{R}\times \mathbb{R}^d}{\inf} \vert1- \widehat{V}(0) \mathcal{P}_{\mathrm{VB}}(\gamma,\tau,\eta,{f}(t,x,\cdot))\vert \geq c_0/2.
\end{equation}

\end{lem}

\begin{proof}
We use polar coordinates and write $(\gamma,\tau,\eta) = (r \widetilde\gamma,r \widetilde\tau, r\widetilde\eta)$, with $r = (|\gamma|^2 + |\tau|^2 + |\eta|^2)^{1/2} >0$ and 
$$(\widetilde\gamma, \widetilde\tau, \widetilde\eta) \in S_+ := \left\{ \widetilde\gamma>0, \widetilde\tau \in \R, \widetilde\eta \in \R^d, \,|\widetilde\gamma|^2 + |\widetilde\tau|^2 + |\widetilde\eta|^2=1 \right\}.$$
Introducing 
$$
\widetilde{\mathcal{P}_{\mathrm{quant}}}(r,\widetilde\gamma,\widetilde\tau,\widetilde\eta,f)=- 2 \widehat{V}( r \widetilde\eta)  \int_{0}^{+\infty}  e^{-(\widetilde\gamma+i\widetilde\tau)s} \frac{1}{r}{\sin\left(\frac{r s  |\widetilde\eta|^2}{2}\right)} (\mathcal{F}_v  f) (t,x, s\widetilde\eta ) ds,$$
the $c_{0}/2$ quantum Penrose condition  implies that for all $t \in [0,T]$, $x \in \R^d$, $r>0$ and $(\widetilde\gamma, \widetilde\tau, \widetilde\eta) \in S_+$, $|1-\widetilde{\mathcal{P}_{\mathrm{quant}}}| \geq c_0/2$. But $\widetilde{\mathcal{P}_{\mathrm{quant}}}$ extends as a continuous function on $[0,+\infty) \times S_+$ with 
$\widetilde{\mathcal{P}_{\mathrm{quant}}}(0,\widetilde\gamma,\widetilde\tau,\widetilde\eta,f) = \widehat{V}(0) \mathcal{P}_{\mathrm{VB}}(\widetilde\gamma,\widetilde\tau,\widetilde\eta,f)$ and $\mathcal{P}_{\mathrm{VB}}$ is homogeneous of order $0$ with respect to $(\gamma,\tau,\eta)$, so we deduce the lemma.
\end{proof}

Consequently, by uniqueness of the solution to Vlasov-Benney  in $\mathscr{C}([0,T]; \mathrm{H}^{m-1}_{r})$ that satisfies the Penrose stability condition
\begin{equation}\label{Pen-fin1}
\underset{t \in [0,T], \, x \in \mathbb{R}^d}{\inf}\underset{(\gamma,\tau,\eta)\in(0,+\infty)\times \mathbb{R}\times \mathbb{R}^d}{\inf} \vert1- \mathcal{P}_{\mathrm{VB}}(\gamma,\tau,\eta,{f}(t,x,\cdot))\vert \geq c_0/2,
\end{equation}
as obtained in \cite[Theorem 1.3]{HKR}\footnote{As a matter of fact, this result is proved for the equation set on $\mathbb{T}^d \times \R^d$, but extends straightforwardly to $\R^d \times \R^d$.}, we finally conclude that no subsequence is actually required and the whole family $(f_\eps)_{\eps \in (0,\eps_0)}$ converges to $f$. This concludes the proof of Theorem~\ref{thm:main} in the defocusing case.

\bigskip

The proof is  similar in the case $\widehat{V}(0)<0$, except that the formal limit is the singular Vlasov equation
\begin{equation}\label{vlasov-singular-focusing}
\partial_t f + v\cdot \nabla_x f  + \nabla_x\rho_f \cdot \nabla_v f=0, 
\end{equation}
which has not (as far as we know) been studied \emph{per se} in the mathematical literature, except in \cite{CHKR}. However, the estimates  of \cite{HKR} devised for Vlasov-Benney transpose perfectly, as soon as the right Penrose condition is considered, namely
\begin{equation}\label{Pen-fin-2}
\underset{t \in [0,T], \, x \in \mathbb{R}^d}{\inf}\underset{(\gamma,\tau,\eta)\in(0,+\infty)\times \mathbb{R}\times \mathbb{R}^d}{\inf} \vert1 + \mathcal{P}_{\mathrm{VB}}(\gamma,\tau,\eta,{f}(t,x,\cdot))\vert \geq c_0/2.
\end{equation}
By Lemma~\ref{lem:Penrose}, the quantum Penrose condition  implies \eqref{Pen-fin-2} when $\widehat{V}(0)=-1$.
It can be readily checked that under~\eqref{Pen-fin-2}, the uniqueness result of \cite[Theorem 1.3]{HKR}  holds as well for~\eqref{vlasov-singular-focusing}, hence allowing to conclude the proof as in the case $\widehat{V}(0)=1$.

\appendix

\section{Pseudodifferential and Fourier integral operators}\label{sec:pseudo}

The goal  of this section is to gather the  various results on  pseudodiffential 
and Fourier integral operators that are needed in the proof of the main result.

\subsection{$L^2$ continuity of pseudodifferential operators for operator-valued symbols. }\label{Pseudo}

Let $n \in \mathbb{N}$. Let $\mathbf{H}$ be a separable Hilbert space. Consider a symbol 
$$
\mathrm{L}(y,\eta): \, \R^n \times\R^n \to \mathscr{L}(\mathbf{H})
$$
where  $\mathscr{L}(\mathbf{H})$ stands for the set of linear bounded operators on $\mathbf{H}$, 
the pseudodifferential operator associated with the symbol $\mathrm{L}$ is defined as
$$
\operatorname{Op}_{\mathrm{L}} u := (2 \pi)^{-n} \int_\eta e^{i y\cdot \eta} \mathrm{L}(y,\eta) \mathcal{F}u(\eta) d\eta,
$$
for all  smooth functions $u$ from $\R^n$ to $\mathbf{H}$. For $\mathbf{H} =\R$, we recover the standard pseudodifferential calculus. In this work we will specifically consider the case $\mathbf{H}= L^2(0,T)$.
The Calder\'on-Vaillancourt theorem reads for such operators as:
\begin{prop}\label{SymTim}
Let $k_n=\lfloor n/2\rfloor+2$.  Assume that 
$$
\sup_{\vert\alpha\vert, \,  \vert \beta\vert\leq k_n} \sup_{y,\eta \in \R^n} \norm{\partial^{\alpha}_{y}\partial^{\beta}_\eta  \mathrm{L}}_{\mathscr{L}(\mathbf{H})} < +\infty.
$$
Then the operator $\operatorname{Op}_{\mathrm{L}}  $ is bounded on $L^2(\R^n;\mathbf{H})$ and there exists  $C>0$ such that 
$$
\norm{\operatorname{Op}_{\mathrm{L}}}_{\mathscr{L}(L^2(\R^n;\mathbf{H}))}\leq C   \sup_{\vert\alpha\vert, \,  \vert \beta\vert\leq k} \sup_{y,\eta \in \R^n} \norm{\partial^{\alpha}_{y}\partial^{\beta}_\eta \mathrm{L}}_{\mathscr{L}(\mathbf{H})}.
$$
\end{prop}

\begin{Rem}
As readily seen from the upcoming proof, in dimension $n=4k+j$, $j=2,3$, Proposition~\ref{SymTim} holds when replacing $k_n$ by $\lfloor n/2 \rfloor +1$.
\end{Rem}

\begin{proof} 
We prove this proposition by a duality argument, closely following the approach of  
\cite[Proof of Theorem 1.1.4]{Lerner}. Since $\mathscr{S}(\mathbb{R}^{n}; \mathbf{H})$ is dense in $L^2(\mathbb{R}^{n}; \mathbf{H})$, it is enough to prove that for all $F,G\in \mathscr{S} (\mathbb{R}^{n}; \mathbf{H})$,
$$
\vert \langle \operatorname{Op}_{\mathrm{L}}  F,G \rangle_{L^2(\mathbb{R}^n; \mathbf{H})} \vert \leq C \norm{F}_{L^2(\mathbb{R}^{n};\mathbf{H})}\norm{G}_{L^2(\mathbb{R}^{n};\mathbf{H})}.
$$
For $n= 4p + j, \, j= 0, 1$, we set $k=\lfloor n/2\rfloor+2$ while for $n= 4p + j, \, j= 2, 3$, we set $k=\lfloor n/2\rfloor+1$. Note that $k$ is always an even integer.
Following \cite{Lerner}, let us  introduce the polynomial function $P_{k}(x)$ of degree $k$ defined by  
$$P_k(x) = (1+\vert x\vert^2)^{k/2}.$$
We shall consider for any  $F \in  \mathscr{S} (\mathbb{R}^{n}; \mathbf{H})$, 
   the function
$$
Z_F(x,\eta)=\int_{\mathbb{R}^n} F(y)P_k(x-y)^{-1} e^{-iy\cdot \eta} dy.
$$
Notice that $Z_F$ can be seen (up to a multiplicative factor depending only on dimension) as the partial Fourier transform of $(x,y)\mapsto F(y) P_k(x-y)^{-1}$. With the choice of $P_k$, since $k>n/2$, we infer that  $1/P_k \in L^2(\R^n)$ and
\begin{equation}
\label{wignerappendix1}
\norm{Z_F}_{L^2(\mathbb{R}^{2n}; \mathbf{H})}=c_k\norm{F}_{L^2(\mathbb{R}^{n}; \mathbf{H})},
\end{equation} 
and that $Z_F$  is  $\mathscr{C}^\infty(\mathbb{R}^{2n}; \mathbf{H})$ and has localization properties that are suitable to justify the following computations (we refer to \cite[p.4]{Lerner},
see also below for a quantitative estimate which is needed).
Starting from  the above scalar  product, we write 
\begin{align*}
\langle  \operatorname{Op}_{\mathrm{L}}  F,G\rangle_{L^2(\mathbb{R}^n; \mathbf{H})} &=\int_{x} \int_{\eta}  \langle e^{ix\cdot\eta}  \mathrm{L} (x,\eta) \widehat F(\eta), G(x) \rangle_{\mathbf{H}} d\eta dx 
\\& =   \int_{x} \int_{\eta}  \langle\mathrm{L} (x,\eta) P_{k}(D_{\eta}) \int_{y} e^{i(x-y) \cdot \eta} P_{k}(x-y)^{-1} F(y)\,dy, G(x) \rangle_{\mathbf{H}} d\eta dx  \\
 &  =   \int_{x} \int_{\eta}  \langle \mathrm{L} (x,\eta) P_{k}(D_{\eta})\left( e^{ix \cdot \eta} Z_{F}(x, \eta)\right), G(x) \rangle_{\mathbf{H}} d\eta dx  \\
 & = \int_{x} \int_{\eta}  \langle  P_{k}(D_{\eta})\left( e^{ix \cdot \eta} Z_{F}(x, \eta)\right),  \mathrm{L} (x,\eta)^*G(x) \rangle_{\mathbf{H}} d\eta dx,
\end{align*}
where $\mathrm{L}(x,\eta)^*$ stands for the adjoint operator of $\mathrm{L}(x,\eta)$ in $\mathbf{H}$.
Thanks to the regularity and the decay of $Z_{F}$, we can integrate by parts to get
\begin{align*}
\langle  \operatorname{Op}_{\mathrm{L}}  F,G\rangle_{L^2(\mathbb{R}^n; \mathbf{H})} &= (-1)^{k}
\int_{x} \int_{\eta}  \langle   e^{ix \cdot \eta} Z_{F}(x, \eta), P_{k}(D_{\eta}) \mathrm{L} (x,\eta)^*G(x) \rangle_{\mathbf{H}} d\eta dx  \\
& = (-1)^{k}  \int_{x} \int_{\eta}  \langle   \left(P_{k}(D_{\eta}) \mathrm{L} (x,\eta)\right) e^{ix \cdot \eta} Z_{F}(x, \eta) ,  G(x) \rangle_{\mathbf{H}} d\eta dx.
\end{align*}
Next, we write
\begin{align*}
&\langle \operatorname{Op}_{\mathrm{L}}  F,G \rangle_{L^2(\mathbb{R}^n; \mathbf{H})} \\
&=(-1)^{k}   \int_{x}\int_{\eta}  \bigg \langle  (P_k(D_\eta)\mathrm{L}(x,\eta) )  Z_F(x,\eta)  ,   P_k(D_x)\left(\int_{\xi}e^{ix\cdot(\xi-\eta)}P_k(\xi-\eta)^{-1} {\mathcal{F}_x{G}}(\xi)d\xi\right)\bigg\rangle_{\mathbf{H}} d\eta dx
	\\&= c_n (-1)^{k} \int_{x}  \int_{\eta} \bigg \langle    (P_k(D_\eta)\mathrm{L}(x,\eta))  Z_F(x,\eta)   ,  P_k(D_x) \left( e^{-ix\cdot\eta} Z_{{\mathcal{F}_x^{-1}{G}}}(-\eta,x)\right)\bigg\rangle_{\mathbf{H}}   d\eta dx.
	\end{align*}
By integrating by parts, this yields
\begin{multline*}
	\langle \operatorname{Op}_{\mathrm{L}}  F,G \rangle_{L^2(\mathbb{R}^n; \mathbf{H})} 
=  \\
c_n (-1)^{k} \int_{x}  \int_{\eta} \bigg \langle   P_{k}(D_{x}) \left[ (P_k(D_\eta)\mathrm{L}(x,\eta))  Z_F(x,\eta) \right] , e^{-ix\cdot\eta} Z_{{\mathcal{F}_x^{-1}{G}}}(-\eta,x) \bigg\rangle_{\mathbf{H}}   d\eta dx
\end{multline*}
and hence, by expanding the polynomials into monomials and by using the Leibniz formula, we obtain
$$ \langle \operatorname{Op}_{\mathrm{L}}  F,G \rangle_{L^2(\mathbb{R}^n; \mathbf{H})} 
	=\sum_{\underset{\vert \beta\vert+\vert \gamma\vert\leq k}{\vert\alpha\vert\leq k} }c_{\alpha,\beta,\gamma} \int_{x} \int_{\eta} \bigg\langle    \partial^{\alpha}_{\eta}\partial^{\beta}_x \mathrm{L} (x,\eta) \partial^\gamma_x Z_F(x,\eta), e^{-ix\cdot\eta} Z_{ {\mathcal{F}_x^{-1}{G}}}(-\eta,x)  \bigg\rangle_{\mathbf{H}}  d\eta dx.
$$
Using the Cauchy-Schwarz inequality, we obtain that 
\begin{align*}
\vert \langle \operatorname{Op}_{\mathrm{L}}F,G\rangle  \vert &\leq \sum_{\underset{\vert \beta\vert+\vert \gamma\vert\leq k}{\vert\alpha\vert\leq k} }c_{\alpha,\beta,\gamma}\norm{  Z_{ {\mathcal{F}_x^{-1}{G}}}}_{L^2(\mathbb{R}^{2n}; \mathbf{H})} \norm{   \partial^{\alpha}_{\eta}\partial^{\beta}_x \mathrm{L} (x,\eta) \partial^\gamma_x Z_F(x,\eta) }_{L^2(\mathbb{R}^{2n}; \mathbf{H})}
	\\&\lesssim \norm{G}_{L^2(\mathbb{R}^{n}; \mathbf{H}) } \sup_{\vert \gamma\vert\leq k}\norm{\partial^{\gamma}_x Z_F}_{L^2(\mathbb{R}^{2n}; \mathbf{H})} \sup_{\underset{\vert \beta\vert\leq k}{\vert\alpha\vert\leq k} }\norm{\partial^{\alpha}_{\eta}\partial^{\beta}_x \mathrm{L}}_{\mathscr{L}(\mathbf{H})},
\end{align*}
where we have also used \eqref{wignerappendix1} and Bessel-Parseval to get the last estimate.
To conclude the proof we are left to estimate $\norm{\partial^{\gamma}_x Z_F}_{L^2(\mathbb{R}^{2n}; \mathbf{H})}$. But since $k>n/2$, we still have that $\partial^{\gamma}(1/P_k)\in L^2(\mathbb{R}^{n})$ and hence as for \eqref{wignerappendix1}, we get
$$
\norm{\partial^{\gamma}_x Z_F}_{L^2(\mathbb{R}^{2n}; \mathbf{H})}=\norm{\int_y F(y) \partial^{\gamma}(1/P_k)(x-y)e^{-iy\cdot\eta}dy}_{L^2(\mathbb{R}^{2n}; \mathbf{H})} \lesssim \norm{F}_{L^2(\mathbb{R}^n;\mathbf{H})}.
$$
This allows to conclude the proof.
\end{proof}

\subsection{Weighted $L^2$ continuity of Fourier Integral Operators}\label{SFIO}

\begin{definition}
Let $n \in \mathbb{N}$. Given an amplitude function $b_{t,s}(z,\xi)$ and a real phase function $\varphi_{t,s}(z,\xi)$, we define
the semiclassical Fourier Integral Operators $U^{\mathrm{FIO}}_{t,s}$ acting on  a function $u\in \mathscr{S}(\mathbb{R}^{n})$ as 
$$
U_{t,s}^{\mathrm{FIO}}u(z)=\frac{1}{(2\pi)^{n  }} \int_{\xi} e^{\frac{i}{\varepsilon}  \varphi_{t,s}^{\eps} (z,\xi) }b_{t,s}^{\eps}(z,\xi)   \widehat u(\xi) d\xi.
$$
We recall the notation
$$ \varphi_{t,s}^{\eps} (z,\xi)= \varphi_{t, s}(z, \eps \xi), \quad b_{t,s}^{\eps}(z,\xi)= b_{t,s}(z, \eps \xi).$$
\end{definition}

We shall first obtain the following general $L^2$ continuity result.

\begin{prop}\label{MajFIO0}
Let $k =\lfloor n/2 \rfloor +1$. 
Let $b_{t,s}(z,\xi)$ and $\varphi_{t,s}(z,\xi)$  be an amplitude and a real phase and assume that there exist $T>0$ and $C>0$ such that the following estimate hold:
\begin{align}\label{Hypb}
\sup_{t,s \in [0,T]}\norm{\partial^{\alpha}_z\partial^{\beta}_{\xi} b_{t,s}(z,\xi)}_{L^{\infty}_{z ,\xi}} &\leq C, \quad | \alpha | \leq k, \, | \beta | \leq k, \\
\label{Hyppsi2}
 \sup_{t,s \in [0,T]} \norm{\partial^{\alpha}_z\partial^{\beta}_{\xi}  \varphi_{t,s}(z,\xi)}_{L^\infty_{z,\xi}}&\leq C, \quad | \alpha | \leq k+2, \, |\beta| \leq k+ 2,   \, | \alpha |+ |\beta | \geq 2.
\end{align}
Assume moreover that
\begin{equation}\label{Hyppsi}
\sup_{t,s \in [0,T]} \norm{ \left(\partial_z\partial_{\xi} \varphi_{t,s} - \mathrm{I} \right) \left(z,\xi\right)}_{L^{\infty}_{z,\xi}}\leq \frac{1}{2}.
\end{equation}
Then the operator $U_{t,s}^{\mathrm{FIO}}  $ is bounded on $L^2(\mathbb{R}^n)$:  there exists  $C_0>0$ such that for every $\eps \in (0, 1]$, 
\begin{equation}\label{eq:continuite-FIO-0}
\sup_{t,s \in [0,T]} \| U_{t,s}^{\mathrm{FIO}}\|_{\mathscr{L}(L^2(\mathbb{R}^n))} \leq C_0.
\end{equation}

\end{prop}

\begin{Rem}
Note that this result applies  as well for standard pseudodifferential operators, as one can choose the phase $\varphi_{t,s}(z,\xi) =  z\cdot \xi$.  Note also that the regularity assumption for the symbol in Proposition~\ref{MajFIO0} is (slightly) better than the one of Proposition~\ref{SymTim}. However, the proof of Proposition~\ref{MajFIO0} involves the  use of  properties of the  Fourier transform of  the  symbol  
which do not extend  to  operator-valued symbols. This is why we needed to resort to a more robust proof for Proposition~\ref{SymTim}, which is unfortunately less  sharp when it comes to regularity assumptions. 
\end{Rem}

By using this general result, we will be able to obtain a more specific form which is tailored for our needs (see Section~\ref{sec:param}).
We focus on the case $n= 2d$, so that we use as in the rest of the paper the notation $z= (x,v)$, $\xi= (\xi_{x}, \xi_{v})$. We namely obtain a sharp continuity result in the weighted space $\mathcal{H}^0_{r,0}$ (recall the definition in~\eqref{defH0rbisbis}), for phases and amplitudes of limited regularity. 

\begin{prop}
\label{MajFIO}
For $r \in \mathbb{N}^*$, assume that \eqref{Hyppsi} holds, that we have
\begin{equation}\label{hypFIO-weighted1}
\sup_{t,s \in [0,T]}\norm{\langle \eps \nabla_x \rangle^r  \langle \eps \nabla_{\xi_v} \rangle^r  \partial^{\alpha}_z\partial^{\beta}_{\xi} b_{t,s}(z,\xi)}_{L^{\infty}_{z ,\xi}} \leq C, \quad | \alpha | + | \beta| \leq 2 (1 + d),
\end{equation}
and assume in addition that
\begin{equation}\label{hypFIO-weighted2}
 \| \partial_{z}^\alpha \partial_{\xi}^\beta (\nabla_{\xi_{v}} \varphi_{t,s}(z, \xi) - v )\|_{L^\infty_{z,\xi}} +  \| \partial_{z}^\alpha \partial_{\xi}^\beta (\nabla_{x} \varphi_{t,s}(z, \xi) - \xi_{x} )\|_{L^\infty_{z,\xi}} \leq C,  \quad |\alpha | +  | \beta |  \leq 2 (1 + d) + 2r-1.
\end{equation}
Then, the operator  $U_{t,s}^{\mathrm{FIO}}  $ is bounded on $\H^0_{r, 0}(\mathbb{R}^{2d})$: there exists  $C_0>0$ such that for every
 $\eps \in (0, 1]$, 
\begin{equation}\label{eq:continuite-FIO}
\sup_{t,s \in [0,T]} \| U_{t,s}^{\mathrm{FIO}}\|_{\mathscr{L}(\H^0_{r, 0})} \leq C_0.
\end{equation}
\end{prop}

\subsubsection*{Proof of Proposition \ref{MajFIO0}}
We shall omit the dependence in $t,s$, in the proof, all the estimates will be uniform on $[0, T]$. 
We notice that we can write
$$ U_{t,s}^{\mathrm{FIO}}= S_{\eps}^{-1} A_{\eps} S_{\eps}$$
where $S_{\eps}$ is the scaling operator
$$ S_{\eps} f (z) = \eps^{n \over 4} f( \sqrt{\eps} z) $$
which is in a isometry on $L^2(\mathbb{R}^n)$ and $A_{\eps}$ is defined by 
$$ A_{\eps} u(z)= \frac{1}{(2\pi)^{n} } \int_{\xi} e^{ i  \varphi_{\eps} (z,\xi) } b_\eps (z,\xi)  \widehat u(\xi) d\xi,$$
where $\varphi_{\eps}$ and $b_{\eps}$ are defined by 
\begin{equation}
\label{defnouveaux} \varphi_{\eps} (z,\xi)= {1 \over \eps} \varphi_{t,s}( \eps^{1 \over 2} z, \eps^{1 \over 2} \xi), \quad b_{\eps}(z, \xi)= b_{t,s} ( \eps^{1 \over 2} z, \eps^{1 \over 2} \xi).
\end{equation}
We deduce from \eqref{defnouveaux} and the assumptions \eqref{Hypb}, \eqref{Hyppsi2}, \eqref{Hyppsi},  that 
\begin{equation}
\label{hypmeilleure1}
\norm{\partial^{\alpha}_z\partial^{\beta}_{\xi} b_{\eps}(z,\xi)}_{L^{\infty}_{z ,\xi}} \leq C, \, | \alpha |, \, |\beta| \leq k, \quad 
\   \norm{\partial^{\alpha}_z\partial^{\beta}_{\xi}  \varphi_{\eps}(z,\xi)}_{L^\infty_{z,\xi}}\leq C, \, | \alpha |,  \, |\beta| \leq k+2, \, | \alpha | + | \beta| \geq 2,
\end{equation}
and that
\begin{equation}
\label{hypmeilleure2}
 \quad  \norm{ \left(\partial_z\partial_{\xi} \varphi_{\eps} - \mathrm{I} \right) \left(z,\xi\right)}_{L^{\infty}_{z,\xi}}\leq \frac{1}{2}.
\end{equation}
To prove the $L^2$ continuity of $U^{\mathrm{FIO}}_{t,s}$ it is now  equivalent to prove the $L^2$ continuity of $A_{\eps}$. With the properties
 \eqref{hypmeilleure1}, \eqref{hypmeilleure2} we can rely on the approach of \cite{Bou} to get uniform estimates  in $\eps$.
There is another classical proof relying on a $TT^*$ argument (see for example \cite{Robert}) but which is much more demanding
in terms of regularity.  

For any $v \in L^2(\mathbb{R}^n)$, we shall estimate:
$$ I := \int_{z}   \int_{\xi} e^{ i  \varphi_{\eps} (z,\xi) } b_\eps (z,\xi)  \widehat u(\xi)  v(z) \, d\xi dz.$$
Let us take $\chi$ a smooth compactly supported function such that $\int_{\mathbb{R}^n} \chi (x) \, dx = 1$. We write
\begin{align*}
I &=  \int_{z}   \int_{\xi} \int_{m} \int_{l} e^{ i  \varphi_{\eps} (z,\xi) } b_\eps (z,\xi)  \widehat u(\xi)  v(z) \chi(z-m)  \chi(\xi-l) \, dl dm  d\xi dz
\\&=  \int_{z}   \int_{\xi} \int_{m} \int_{l}  e^{ i  \varphi_{\eps} (z + m,\xi + l) } b_{\eps}(z+m, \xi+l)  \chi(z)  \chi(\xi) v(z+m)\widehat{u}(\xi+l) \, dl dm d\xi dz
\end{align*}
and we finally obtain
$$ I =  \int_{z}   \int_{\xi} \int_{m} \int_{l}  e^{ i  \varphi_{\eps} (z + m,\xi + l) } b_{m,l}(z , \xi)  v_{m}(z) {u}_{l}(\xi)  dldm d\xi dz $$
with
\begin{equation}
\label{deflocales}  b_{m,l}(z , \xi) =  b_{\eps}(z+m, \xi+l)  \chi(z)  \chi(\xi), \quad v_{m}(z)= v(z+m)\widetilde{\chi}(z) , \quad u_{l}(\xi)= \widehat{u}(\xi+l) \widetilde{\chi}(\xi),
\end{equation}
where $\widetilde{\chi}$ is a smooth compactly supported function which is equal to one on the support of $\chi$.

We shall now use a Taylor expansion of the phase, by writing
$$  \varphi_{\eps} (z + m,\xi + l) = \varphi_{\eps}(m,l) + \nabla_{z} \varphi_{\eps}(m,l)\cdot z + \nabla_{\xi} \varphi_{\eps}(m,l)\cdot \xi + R_{m,l}(z,\xi),$$
where
\begin{equation}
\label{defRkl} R_{m,l}(z , \xi)=  \int_{0}^1 (1-t) D^2\varphi_{\eps}(m+tz, l+t \xi) \cdot(z, \xi)^2 \, dt.
\end{equation}
Let us then define
$$ a_{m,l}(z, \xi)= e^{i R_{m,l}(z, \xi)} b_{m,l}(z,\xi).$$
Thanks to the definition  \eqref{deflocales}, we observe that $a_{m,l}$ is compactly supported in $z,$ $\xi$ and consequently, 
we can deduce from  \eqref{hypmeilleure1} and \eqref{hypmeilleure2} that
\begin{equation}
\label{amplitudefinale}
  \sup_{m,l, z, \xi} | \partial_{z}^\alpha \partial_{\xi}^\beta a_{m,l}(z, \xi)| \leq C, \quad |\alpha|\leq k, \, |\beta | \leq k.
\end{equation}
We have
\begin{align*} I &=  \int_{z}   \int_{\xi} \int_{m} \int_{l}  e^{ i  \varphi_{\eps}(m,l) +i  \nabla_{z}\varphi_{\eps}(m,l) \cdot z + i \nabla_{\xi} \varphi_{\eps}(m,l)\cdot \xi} a_{m,l}(z , \xi)  v_{m}(z) {u}_{l}(\xi) \, dl dm d\xi dz 
\\ &= c_n \int_{m} \int_{l}  \int_{\eta} \int_{y}  e^{i  \varphi_{\eps}(m,l)}  \widehat{a_{m,l}}(\eta, y) \widehat{v_{m} } \left(- \eta - \nabla_{z} \varphi_{\eps}(m,l)\right) 
 \widehat{u_{l}}\left(-y- \nabla_{z} \varphi_{\eps}(m,l)\right) \, dy d\eta  dl dm,
\end{align*}
where $\widehat{a_{m,l}}$ stands for the Fourier transform with respect to both sets of variables, and $c_n$ is a normalizing constant.
By using Cauchy-Schwarz, we get that
$$ | I | \lesssim \int_{m ,l}   \| a_{m,l} \|_{H^{k,k}} \left( \int_{\eta, y} {| \widehat{v_{m}} (- \eta - \nabla_{z} \varphi_{\eps}(m,l))|^2 \over \langle
 \eta \rangle^{2k}} {| \widehat{u_{l}} (- y - \nabla_{\xi} \varphi_{\eps}(m,l))|^2 \over \langle
 y \rangle^{2k}} \, d\eta dy \right)^{1 \over 2} \, dm dl,$$
 where  the Sobolev norm  $ \| a_{m,l} \|_{H^{k,k}}$ is defined by
 $$  \| a_{m,l} \|_{H^{k,k}}^2 = \int_{ \eta, y }  \langle \eta \rangle^{2k} \langle y \rangle^{2k}  |\widehat{a_{m,l}}(\eta, y) |^2\, d\eta dy.$$
 Note that from \eqref{amplitudefinale} and the fact that $a_{m,l}$ is compactly supported, we get that 
 $$ \sup_{m,l}  \| a_{m,l} \|_{H^{k,k}} \lesssim 1.$$
  This yields by using again Cauchy-Schwarz
  \begin{align*}  | I | &\lesssim \left(\int_{m,l, \eta}  {| \widehat{v_{m}} (- \eta - \nabla_{z} \varphi_{\eps}(m,l))|^2 \over \langle
 \eta \rangle^{2k}} \, d \eta dm dl  \int_{m,l,y}   {| \widehat{u_{l}} (- y - \nabla_{\xi} \varphi_{\eps}(m,l))|^2 \over \langle
 y \rangle^{2k}} \, dy dm dl \right)^{1 \over 2}
 \\
&\lesssim   \left(\int_{m,l, \eta}  {| \widehat{v_{m}} (\eta)|^2 \over \langle
 \eta + \nabla_{z} \varphi_{\eps}(m,l) \rangle^{2k}} \, d \eta dm dl  \int_{m,l,y}   {| \widehat{u_{l}} (y)|^2 \over \langle
 y + \nabla_{\xi} \varphi_{\eps}(m,l) \rangle^{2k}} \, dy dm dl \right)^{1 \over 2}.
 \end{align*}
Thanks to  \eqref{hypmeilleure2}, we know that 
$ l\mapsto \nabla_{z} \varphi_{\eps}(m,l)$ and  $m \mapsto  \nabla_{\xi} \varphi_{\eps}(m,l)$ are diffeomorphisms with controlled Jacobians.
We can thus use them to  change variables to get that
$$ | I | \lesssim   \left(\int_{m,l', \eta}  {| \widehat{v_{m}} (\eta)|^2 \over \langle
 \eta + l' \rangle^{2k}} \, d \eta dm dl'  \int_{m',l,y}   {| \widehat{u_{l}} (y)|^2 \over \langle
 y + m' \rangle^{2k}} \, dy dm' dl \right)^{1 \over 2}
 $$
 and as a result, we finally obtain by using Bessel-Parseval and $k>n/2$ that
 $$ 
 | I | \lesssim   \left(\int_{m,z}  | v_{m} (z)|^2  \, dz  dm  \int_{l,\xi}   |u_{l} (\xi)|^2 d\xi dl \right)^{1 \over 2} \lesssim \|u\|_{L^2} \|v\|_{L^2},
 $$
 where the final estimate comes from the definition  \eqref{deflocales} and the fact that $\widetilde \chi$ is compactly supported. This ends the proof 
  of Proposition \ref{MajFIO0}.
  
  \subsubsection*{Proof of Proposition \ref{MajFIO}}
  Let us set
  $$ r_{1}(t,s, z, \xi)= 2\left( \nabla_{\xi_{v}} \varphi_{t,s}(z,\xi) - v\right),  \quad r_{2}(t,s, z, \xi)= \nabla_{x} \varphi_{t,s}(z, \xi) - \xi_{x}.$$
  We observe that we can write 
  \begin{multline*} V_{\pm} (U_{t,s}^{\mathrm{FIO}}u)
  \\
  ={ 1 \over (2\pi)^d}  \int_{\mathbb{R}^{2d}} e^{ \frac{i}{\eps} \varphi_{t,s}^\eps } \left[ i \left(   \eps \xi_{x}
   \pm \frac{2}{\eps} \nabla_{\xi_{v}} \varphi_{t,s}^\eps \pm r_{2}^\eps (t,s,z,  \xi)  \mp r_{1}^\eps(t,s, z, \xi)  \right) b^\eps + \eps  \nabla_{x}b^\eps\right]
    \widehat{u}(\xi) \, d\xi.
    \end{multline*}
    Next,   by integrating by parts we have
    $$ { 1 \over (2\pi)^d}  \int_{\mathbb{R}^{2d}}  \frac{i}{\eps} \nabla_{\xi_{v}} \varphi_{t,s}^\eps  e^{ \frac{i}{\eps} \varphi_{t,s}^\eps }
     b^\eps \widehat u \, d\xi
      = -  { 1 \over (2\pi)^d}   \int_{\mathbb{R}^{2d}}  e^{ \frac{i}{\eps} \varphi_{t,s}^\eps } \left(\nabla_{\xi_{v}}b^\eps  \widehat u + b^\eps \nabla_{\xi_{v}} \widehat u
      \right) \, d\xi, $$
      and therefore, we finally get the identity
  $$  V_{\pm} (U_{t,s}^{\mathrm{FIO}}u)
  =  { 1 \over (2\pi)^d}  \int_{\mathbb{R}^{2d}} e^{ \frac{i}{\eps} \varphi_{t,s}^\eps } b^\eps \widehat{ V_{\pm} u}\, d\xi  
   +   { 1 \over (2\pi)^d}  \int_{\mathbb{R}^{2d}} e^{ \frac{i}{\eps} \varphi_{t,s}^\eps } \left(  \pm  r_{2}^\eps \mp r_{1}^\eps  + \eps  \left[(\nabla_{x} \mp 2\nabla_{\xi_{v}} )
    b \right]^\eps \right)  \widehat u  \, d\xi.$$
  The result then follows by iterating this identity and by  applying Proposition \ref{MajFIO0}.

\subsection{Pseudifferential calculus with parameter }\label{sec:PseuDiff}
 In this section, we present some useful  results for pseudodifferential calculus with parameter $\gamma>0$, following \cite{HKR} {(see also \cite{MetZum})}.  Here we do not need only $L^2$ continuity results but also  calculus results for the composition of operators, for this reason, 
 we shall use different norms of symbols compared to  Section \ref{Pseudo}, the main interest is that  they are less demanding in terms of regularity when dealing with composition formulas when we apply them to our specific setting.
 
  We consider symbols $a(x,\gamma,\tau,\kappa)=a(x,\zeta)$ on $\mathbb{R}^d \times ]0,+\infty [\times \mathbb{R} \times  \mathbb{R}^d\backslash \{0\}, \gamma >0$ is a parameter. We introduce the following seminorms, for $k \in \mathbb{N}$,
\begin{equation}
\label{def:semi-norms-gamma}
\begin{aligned}
\vert a \vert _{0} &= \sup_{\vert \alpha\vert \leq k } \norm{\mathcal{F}_x(\partial^{\alpha}_xa)}_{L^2(\mathbb{R}^d;L^{\infty}_\zeta)}, \\
\vert a \vert _{k,1} &= \sup_{\vert \alpha\vert \leq k } \norm{\gamma \mathcal{F}_x(\partial^{\alpha}_x \nabla_\xi a)}_{L^2(\mathbb{R}^d;L^{\infty}_\zeta)},
\end{aligned}
\end{equation}
where $\xi=( \tau,\kappa)$.

\begin{Rem}
Note that, we are considering pseudodifferential operators  acting on functions defined
 on $\mathbb{R} \times \mathbb{R}^d$ and that denoting by $t$ the first variable of $\mathbb{R}\times \mathbb{R}^d$,  the symbols that  we consider here do not depend on $t$
  so that they act as Fourier multipliers on this component. This class is the one actually needed
  for the analysis  in the paper and this simplification allows to slightly lower the level of regularity
  needed on the symbols in order to have a good calculus.
We also point out that the semi-norm $\vert \cdot \vert _{k_d,1}$ is slightly different from the one used in \cite{HKR}, as the weight here is $\gamma$ whereas it was
$
\langle \zeta \rangle= (\gamma^2+\tau^2+\vert \kappa\vert^2)^{1/2}
$
in  \cite{HKR}. This is because when $\widehat{V}$ is not decaying, the symbols that we consider in this work only have finite semi-norm for the one defined here.
\end{Rem}
The  continuity results that  we will need in this work are given below.
\begin{prop}\label{CalcPseudoParam}
Let $k_d := \lfloor d/2 \rfloor +2$. There exists $C>0$ such that for every $\gamma>0$, we have 
\begin{itemize}
\item for every symbol $a$ such that $ |a|_{k_d,0}<+\infty$,
$$
\norm{\mathbf{Op}_a^{\gamma}}_{\mathscr{L}(L^2(\mathbb{R} \times \mathbb{R}^d))}\leq C |a|_{k_d,0}, 
$$
\item for every symbol $a,b$ such that $ |a|_{k_d,1}<+\infty,  |b|_{k_d+1,0}<+\infty$,
$$
\norm{\mathbf{Op}_a^{\gamma}\mathbf{Op}_b^{\gamma}-\mathbf{Op}_{ab}^{\gamma}}_{\mathscr{L}(L^2(\mathbb{R} \times \mathbb{R}^d))}\leq \frac{C}{\gamma} |a|_{k_d,
1}|b|_{k_d+1,0} .
$$
\end{itemize}
\end{prop}

\begin{Rem}
Exactly as for Proposition~\ref{SymTim},  in dimension $d=4k+j$, $j=2,3$,  Proposition~\ref{CalcPseudoParam} holds when replacing $k_d$ by $\lfloor d/2 \rfloor +1$.
\end{Rem}

Note that the first item above is the same as in  \cite{HKR}, we shall reproduce the proof
for the sake of completeness. For the second item, there is a slight difference due 
to the different definition of the seminorm $|\cdot|_{k,1}$ compared to  \cite{HKR}.

%
\begin{proof}
We expand the operator $\mathbf{Op}_a^{\gamma}$ as
\begin{align*}
\mathbf{Op}_a^{\gamma} u&=(2\pi)^{-d-1} \int_\eta\int_\tau e^{i(\tau t + x\cdot \eta)}a(x,\zeta)\widehat u(\tau,\eta)d\tau d\eta 
		       \\&=(2\pi)^{-2d-1}\int_\kappa e^{ix\cdot \kappa}\int_\eta\int_\tau e^{i(\tau t + x\cdot \eta)}\mathcal{F}_x a(\kappa,\zeta)\widehat u(\tau,\eta)d\tau d\eta d\kappa
			 \\&=(2\pi)^{-2d-1}  \int_\kappa\int_\tau e^{i(\tau t + x\cdot \kappa)}\left(\int_\eta \mathcal{F}_x a(-\eta+ \kappa,\zeta)\widehat u(\tau,\eta) d\eta \right) d\tau d\kappa.
\end{align*}
Using the Bessel-Parseval identity, this yields 
$$
\norm{\mathbf{Op}_a^{\gamma}u}_{L^2(\mathbb{R} \times \mathbb{R}^d)}\lesssim \norm{\norm{\int \mathcal{F}_x a(-\eta+\kappa,\gamma,\tau,\eta)\widehat u(\tau,\eta) d\eta}_{L^2_\kappa(\mathbb{R}^d)}}_{L^2_\tau(\mathbb{R})}.
$$
Then, by Cauchy-Schwarz and Fubini, we obtain 
\begin{align*}
&\norm{\int \mathcal{F}_x a(\eta-\kappa,\gamma,\tau,\eta)\widehat u(\tau,\eta) d\eta}_{L^2_\kappa(\mathbb{R}^d)}^2 
\\& \lesssim \norm{\sup_\kappa \vert \mathcal{F}_x a(\cdot, \gamma,\tau,\kappa)\vert}_{L^1(\mathbb{R}^d)} \int_{\eta} \int_{\kappa} \vert \mathcal{F}_x a(-\eta+\kappa, \gamma,\tau,\eta)\vert \vert \widehat u (\tau,\eta)\vert^2 d\eta d\kappa 
\\& \lesssim \norm{\sup_\kappa \vert \mathcal{F}_x a(\cdot, \gamma,\tau,\kappa)\vert}_{L^1(\mathbb{R}^d)} \norm{\sup_\kappa \vert \mathcal{F}_x a(\cdot, \gamma,\tau,\kappa)\vert}_{L^1(\mathbb{R}^d)}\norm{\widehat u(\tau,\cdot)}^2_{L^2(\mathbb{R}^d)}
\\& \lesssim \norm{\sup_\kappa \vert \mathcal{F}_x a(\cdot, \gamma,\tau,\kappa)\vert}_{L^1(\mathbb{R}^d)} ^2 \norm{\widehat u(\tau,\cdot)}^2_{L^2(\mathbb{R}^d)}.
\end{align*}
We finally take the integral in $\tau$ to obtain 
$$
\norm{\mathbf{Op}_a^{\gamma}u}_{L^2(\mathbb{R} \times \mathbb{R}^d)}\lesssim\norm{ \mathcal{F}_x a}_{L^1(\mathbb{R}^d;L^\infty_\zeta)} \norm{\widehat u(\tau,\cdot)}_{L^2(\mathbb{R}\times \mathbb{R}^d)},
$$
 As in the proof of Proposition~\ref{SymTim},  if $d= 4p + j, \, j= 0, 1$, we set $k=\lfloor d/2\rfloor+2$ while for $d= 4p + j, \, j= 2, 3$, we set $k=\lfloor d/2\rfloor+1$. 
By the Cauchy-Schwarz inequality, we have $\norm{ \mathcal{F}_x a}_{L^1(\mathbb{R}^d;L^\infty_\zeta)} \lesssim |a|_{k,0}$ and we obtain the first item.

\bigskip

We then study the second estimate. Provided that $ab$ belongs to a suitable class of symbols, we can use the composition formula for pseudodiffential operators 
$$
\mathbf{Op}_a^{\gamma}\mathbf{Op}_b^{\gamma}=\mathbf{Op}_c^{\gamma},
$$
where 
\begin{align*}
c(x,\zeta)&=\int_{\kappa'} e^{i\kappa'\cdot x} a(x,\gamma,\tau,\kappa+\kappa')\mathcal{F}_xb(\kappa',\zeta)d\kappa'
	   \\&= \int_{\kappa'} e^{i\kappa'\cdot x} a(x,\gamma,\tau,\kappa)\mathcal{F}_xb(\kappa',\zeta)d\kappa'+\int_{\kappa'} e^{i\kappa'\cdot x} \int_0^1 \nabla_\kappa a(x,\gamma,\tau,\kappa+r\kappa')dr \cdot \kappa' \mathcal{F}_xb(\kappa',\zeta)d\kappa'
	  \\&=  a(x,\zeta)b(x,\zeta)+\int_{\kappa'} e^{i\kappa'\cdot x} \int_0^1 \nabla_\kappa a(x,\gamma,\tau,\kappa+r\kappa')dr \cdot \kappa' \mathcal{F}_xb(\kappa',\zeta)d\kappa'
	   \\&=  a(x,\zeta)b(x,\zeta)+\frac{1}{\gamma}d(x,\zeta),
\end{align*}
defining $d(x,\zeta)$ by 
$$
d(x,\zeta)=\gamma \int_{\kappa'} \int_0^1  e^{i\kappa'\cdot x} \nabla_\kappa a(x,\gamma,\tau,\kappa+r\kappa')dr \cdot \kappa' \mathcal{F}_xb(\kappa',\zeta)d\kappa'.
$$
Let us now estimate $|d|_{k_d,0}$. We have 
$$
(\mathcal{F}_x \partial^{\alpha}_x d) (\eta,\gamma,\tau,\kappa)=\gamma\int_0^1  \int_{\kappa'}  (i\eta)^{\alpha} (\mathcal{F}_x \nabla_\kappa a)(\eta-\kappa',\gamma,\tau,\kappa+r\kappa')dr \cdot \kappa' \mathcal{F}_xb(\kappa',\zeta)d\kappa',
$$
and taking the $L^{\infty}_\zeta$ norm, it holds
\begin{align*}
\norm{(\mathcal{F}_x \partial^{\alpha}_x d) (\eta,\cdot)}_{L^{\infty}_\zeta}&\lesssim \bigg( \int_{\kappa'}  \vert \eta-\kappa' \vert^{\vert \alpha\vert} \norm{\gamma (\mathcal{F}_x \nabla_\kappa a)(\eta-\kappa',\gamma,\tau,\cdot)}_{L^{\infty}_\zeta} \vert \kappa'\vert \norm{ \mathcal{F}_xb(\kappa',\cdot)}_{L^{\infty}_\zeta}d\kappa'
\\&+  \int_{\kappa'}   \norm{\gamma (\mathcal{F}_x \nabla_\kappa a)(\eta-\kappa',\gamma,\tau,\cdot)}_{L^{\infty}_\zeta} \vert \kappa'\vert^{\vert \alpha\vert +1} \norm{ \mathcal{F}_xb(\kappa',\cdot)}_{L^{\infty}_\zeta}d\kappa'\bigg).
\end{align*}
Using convolution estimates  we deduce 
$$
|d|_{k_d,0}\lesssim\left(  |a|_{k_d,1}\norm{\mathcal{F}_x \nabla_x b }_{L^1(\mathbb{R}^d;L^{\infty}_\zeta) } +  |b|_{k_d+1,0}\norm{\gamma \mathcal{F}_x \nabla_\kappa a }_{L^1(\mathbb{R}^d;L^{\infty}_\zeta) }  \right),
$$
which precisely means that 
$$
|d|_{k_d,0}\lesssim|a|_{k_d,1}|b|_{k_d+1,0}.
$$
The continuity result of the first item hence shows that 
$$
\norm{\mathbf{Op}_a^{\gamma}\mathbf{Op}_b^{\gamma}u-\mathbf{Op}_{ab}^{\gamma}u}_{L^2(\mathbb{R} \times \mathbb{R}^d)}=\norm{\frac{1}{\gamma}\mathbf{Op}_d^{\gamma}u}_{L^2(\mathbb{R} \times \mathbb{R}^d)} \lesssim \frac{1}{\gamma} |a|_{k_d,1}|b|_{k_d+1,0} \norm{u}_{L^2(\mathbb{R} \times \mathbb{R}^d)},
$$
which concludes the proof of the proposition.
\end{proof}

 We finally deal with the  semiclassical version of the above calculus.
     For any symbol $a(x, \zeta)$ as above,  we set for $\eps \in (0, 1]$,  $a^\eps(x, \zeta)= a(x, \eps \zeta)= a(x, \eps \gamma, \eps \tau, \eps \kappa)$ and we define for $\gamma \geq 1$, 
 \begin{equation}
 (\mathbf{Op}^{\eps, \gamma}_{a} u)(t,x)=
  (\mathbf{Op}^{\gamma}_{a^\eps} u)(t,x).
 \end{equation}   
 For this calculus, we have the following result:
 \begin{prop}
 \label{propsemifini}
  There exists $C>0$ such that for every $ \eps \in (0, 1]$ and for every $\gamma \geq 1$, we have 
  \begin{itemize}
  \item 
    for every symbol $a$ such that $ |a|_{k_d,0}<+\infty$,
  $$ \|\mathbf{Op}_{a}^{\eps, \gamma} \|_{\mathscr{L}(L^2(\mathbb{R} \times \mathbb{R}^d))}  \leq C |a|_{k_d,0}, $$
   \item  for every symbol $a,b$ such that $ |a|_{k_d,1}<+\infty,  |b|_{k_d+1,0}<+\infty$,
  $$  \|\mathbf{Op}_{a}^{\eps, \gamma} \mathbf{Op}_{b}^{\eps, \gamma}  - \mathbf{Op}_{ab}^{\eps, \gamma} \|_{\mathscr{L}(L^2(\mathbb{R} \times \mathbb{R}^d))} \leq {C \over \gamma} |a|_{k_d, 1} |b|_{k_d+1, 0}.$$
  
  \end{itemize}
 
 \end{prop}
 
 \begin{proof}[Proof of Proposition~\ref{propsemifini}]
  The proof is a direct consequence of Proposition  \ref{CalcPseudoParam}
   since for any symbol $a$, we have by definition of $a^\eps$  that  for all $k\in \mathbb{N}$,
   $$ |a^\eps|_{k, 0}= |a|_{k,0}, \quad  |a^\eps |_{k, 1}=   |a|_{k,1}.$$
 
 \end{proof}

\bibliographystyle{abbrv}
\bibliography{biblio.bib}

\end{document}